\documentclass[11pt]{report}

\usepackage{multind}
\makeindex{A}
\makeindex{B}

\usepackage{type1cm}
\usepackage{amssymb}
\usepackage{amsmath}
\usepackage{amscd}
\usepackage{latexsym}
\usepackage{amsthm}
\usepackage{eucal}
\usepackage{enumerate}
\usepackage{epic}
\usepackage{eepic}
\usepackage{graphicx}
\usepackage{color}
\usepackage{float}
\makeatletter
\renewcommand{\@makecaption}[2]{%
  \vskip\abovecaptionskip
  \sbox\@tempboxa{#1 #2}%
  \ifdim \wd\@tempboxa >\hsize
    #1: #2\par
  \else
    \global \@minipagefalse
    \hb@xt@\hsize{\hfil\box\@tempboxa\hfil}%
  \fi
  \vskip\belowcaptionskip}
\makeatother



\numberwithin{equation}{section}
\numberwithin{figure}{section}

\newtheorem{thmSec}{Theorem}[section]
\newtheorem{theorem}[thmSec]{Theorem}
\newtheorem{proposition}[thmSec]{Proposition}

\newtheorem{claim}[thmSec]{Claim}
\newtheorem{fact}[thmSec]{Fact}
\newtheorem{remNonumber}[thmSec]{Remark}
\newtheorem{corollary}[thmSec]{Corollary}
\newtheorem{lemSec}[thmSec]{Lemma}
\newtheorem{example}[thmSec]{Example}
\newtheorem{problem}[thmSec]{Problem}
\newtheorem{definition}[thmSec]{Definition}
\newtheorem{program}[thmSec]{Program}
\newenvironment{texteqn}
{\begin{equation} \begin{minipage}[b]{0.8\linewidth}}
{\end{minipage} \end{equation} \ignorespacesafterend}
\newtheorem*{main theorem}{Main Theorem}
\newcommand{\rarrowsim}{\smash{\mathop{\,\rightarrow\,}\limits
  ^{\lower1.5pt\hbox{$\scriptstyle\sim$}}}}
\newcommand{\larrowsim}{\smash{\mathop{\,\leftarrow\,}\limits
  ^{\lower1.5pt\hbox{$\scriptstyle\sim$}}}}
\newcommand{\Har}[2]{\mathcal{H}^{#1}(\mathbb{R}^{#2})}
\newcommand{\set}[2]{ \{ {#1}: {#2} \}  }

\newcommand{\sgn}{\operatorname{sgn}}
\newcommand{\trF}{{ {}_2 \kern-.05em F_1}}

\newcommand{\trans}{{}^t \!}

\newcommand{\mmax}{{\mathfrak m}^{\max}}

\newcommand{\nmax}{{\mathfrak n}^{\max}}
\newcommand{\pmax}{{\mathfrak p}^{\max}}

\newcommand{\tilLap}[1]{{\widetilde {\Delta}}_{#1}}

\newcommand{\Mmax}{M^{\max}}

\newcommand{\Nmax}{N^{\max}}
\newcommand{\Pmax}{P^{\max}}
\newcommand{\Vpq}{V^{p,q}}

\newcommand{\n}{p+q-2}

\newcommand{\nbar}[1]{\overline{n}_{#1}}
 
\newcommand{\wL}{L^2(\mathbb R_+, r^{p+q-5}dr)}

\newcommand{\B}{\mathcal T}
\newcommand{\Mp}{Mp}
\newcommand{\tilC}{\widetilde{C}}
\newcommand{\hatotimes}{\mathbin{\widehat{\otimes}}}
\newcommand{\xzeta}{x}
\usepackage[dvips,hyperindex]{hyperref}\hypersetup{linkbordercolor=0.2 1 0.5}

\begin{document}

\def\today{July 18, 2008}

\title{%
The Schr\"{o}dinger model\\
for the minimal representation of\\
the indefinite orthogonal group $O(p, q)$}

\author{Toshiyuki \textsc{Kobayashi}$^*$\thanks{$^*$ Partially supported by
        Grant-in-Aid for Scientific Research (B) (18340037), Japan
        Society for the Promotion of Science.}
\ and
Gen \textsc{Mano}%
\footnote{2000 Mathematics Subject Classification. Primary 22E30;
Secondary 22E46, 43A80\newline
\indent \indent Key words and phrases:
minimal representation, Schr\"odinger model, 
generalization of the Fourier transform,
Weil representation, indefinite orthogonal group, 
unitary representation, isotropic cone,
Bessel functions, 
Meijer's $G$-functions
}
}

{
\renewcommand{\thefootnote}{}

\maketitle

}

\pagenumbering{roman}
\setcounter{page}{2}
\tableofcontents

\begin{abstract}
\thispagestyle{plain}
\pagenumbering{roman}
\setcounter{page}{4}
We introduce a generalization of the Fourier transform, 
denoted by $\mathcal{F}_C$, on the isotropic cone $C$ 
associated to an indefinite quadratic form of signature $(n_1,n_2)$ on
$\mathbb{R}^n$ ($n=n_1+n_2$: even).
This transform is in some sense the
unique and natural unitary operator on $L^2(C)$,
as is the case with the Euclidean Fourier transform
$\mathcal{F}_{\mathbb{R}^n}$ on $L^2(\mathbb{R}^n)$.
Inspired by recent developments of algebraic representation theory of
reductive groups,
we shed new light on classical analysis 
on the one hand,
and give the global formulas for the $L^2$-model 
of the minimal representation of the
simple Lie group $G=O(n_1+1,n_2+1)$ on the other hand.

The transform $\mathcal{F}_C$ 
 expands functions on $C$ into joint
eigenfunctions of \textit{fundamental differential operators}
which are mutually commuting, self-adjoint, and of
second order.
We decompose $\mathcal{F}_C$ into the singular Radon transform and the
 Mellin--Barnes integral, 
find its distribution kernel,
and
 establish the inversion and the Plancherel formula.
The transform
$\mathcal{F}_C$ reduces to the Hankel transform if
 $G$ is  $O(n,2)$ or
$O(3,3) \approx SL(4,\mathbb{R})$.

The unitary operator $\mathcal{F}_C$ together with multiplications and
translations coming from 
the conformal transformation group 
$CO(n_1,n_2)\ltimes\mathbb{R}^{n_1+n_2}$
 generates the minimal representation of the indefinite orthogonal
group $G$.
Various different models of the same representation 
have been constructed by
 Kazhdan, Kostant, Binegar--Zierau, Gross--Wallach, Zhu--Huang, Torasso,
Brylinski, and Kobayashi--{\O}rsted, and others.
Among them, our model gives the global formula of 
the whole group action on the simple Hilbert space $L^2(C)$,
and generalizes the classic
Schr\"odinger
 model $L^2(\mathbb R^n)$
of the Weil representation. 
Here, $\mathcal{F}_C$ plays a similar role to $\mathcal{F}_{\mathbb{R}^n}$.

Yet another motif is special functions.
Large group symmetries in the minimal representation 
yield functional equations of various special functions.
We find explicit $K$-finite vectors
on $L^2(C)$, 
and give a new
proof of the Plancherel formula
for Meijer's $G$-transforms.
\end{abstract}

\pagenumbering{arabic}

\chapter{Introduction}
\label{sec:intro}

\setcounter{page}{1}

This book is a continuation of a series of 
our research projects 
\cite{xkcheck,xkmano1, xkmano2, 
xkors1,xkors2,xkors3}.
Our motif is to open up and develop geometric analysis
 of a single infinite dimensional representation, namely, 
the minimal
representation  $\pi$ of the indefinite (even) orthogonal group.

This representation is surprisingly rich in
its different models,
through which we have cross-fertilization and interactions
 with various areas of mathematics such
as conformal geometry and the Yamabe operator,
Fourier analysis, 
ultra-hyperbolic equations and their conserved quantities,
 the Kepler problem,
holomorphic semigroups, and analysis on isotropic cones. 
Among them,
this book is devoted to the $L^2$-model
(\textit{Schr\"{o}dinger model}), 
for which the local formula was established in a previous paper
 \cite{xkors3} with B. {\O}rsted.
The global formula of the whole group action 
is the subject of this book.

We have limited ourselves to the very representation $\pi$,
although some of our results could be 
generalized to other settings by the ideas developed here.
This is primarily because we believe that
 geometric analysis of
this specific minimal representation
is of interest in its own right, 
and might open up an unexpected direction of research bridging different fields
 of mathematics,
 as in the case of the Weil representation 
 (e.g.\ \cite{xfolland,xHo,xHoTan,xKaVe,xRalSch}).

Bearing this in mind, we will not only
\begin{itemize}
\item
formalize our main results by means of representation theory, 
\end{itemize}
but  also
\begin{itemize}
\item
formalize our main results without group theory.
\end{itemize}

We have made effort to expound the theory in a self-contained fashion 
as much as possible.

\medskip

For $n=n_1+n_2$,
we denote by $\mathbb{R}^{n_1,n_2}$ the Euclidean space $\mathbb{R}^n$
endowed with the flat pseudo-Riemannian structure
\begin{equation*}
  ds^2=dx_1^2+\dots+dx_{n_1}^2 - dx_{n_1+1}^2 - \dots - dx_n^2,
\end{equation*}
and define the isotropic cone $C$ by
\begin{equation*}
   C:= \{ x \in \mathbb{R}^n \setminus \{0\}:
          x_1^2+\dots+x_{n_1}^2 - x_{n_1+1}^2 - \dots - x_n^2 = 0 \}.
\end{equation*}

In this book, we will introduce the `Fourier transform'
$\mathcal{F}_C$ on 
 the isotropic cone $C$ for $n$ even. 
This transform $\mathcal{F}_C$ is in some sense the
unique and natural unitary operator on $L^2(C)$,
as it is the case with the Euclidean Fourier transform
$\mathcal{F}_{\mathbb{R}^n}$ on $L^2(\mathbb{R}^n)$.

Here is a brief guide to the three motivations of this book,
with emphasis on the role of the unitary operator $\mathcal{F}_C$.

The first motivation comes from analysis on the isotropic cone $C$ itself.
Different from non-isotropic hypersurfaces
(e.g.\ hyperboloids) in $\mathbb{R}^{n_1,n_2}$, 
the restriction of $ds^2$ to $C$ is degenerate,
and we do not have a natural pseudo-Riemannian structure on $C$.
Consequently, there is no natural single operator on $C$ 
such as the 
\index{B}{Laplace--Beltrami operator}%
Laplace--Beltrami operator.
However, it turns out that 
there are commuting, self-adjoint, second order differential
operators $P_1,\dots,P_n$
 that we call
\index{B}{fundamental differential operator}%
\textit{fundamental differential operators} 
 on $C$
 satisfying the algebraic relation
$P_1^2+\dots+P_{n_1}^2-P_{n_1+1}^2-\dots-P_n^2=0$.
Then, what we want is to understand
 how an arbitrary function on $C$
(of appropriate class) is expanded 
into joint eigenfunctions of $P_1,\dots,P_n$.

We will find explicit joint eigendistributions
 for  $P_1,\dots,P_n$,
 and construct a (well-defined) transform, 
to be denoted by $\mathcal F_C$,
by
 means of these eigenfunction. 
The transform 
$\mathcal{F}_C$ intertwines the multiplication by coordinate functions
 with the differential operators $P_j$.
Moreover, we prove that we can normalize $\mathcal{F}_C$ such that it
 is involutive, i.e. $\mathcal F_C^2=
  \operatorname{id}$ and unitary.
Thus, we establish its inversion formula and the Plancherel type theorem.
It is noteworthy that 
the kernel function $K(x,x')$ of $\mathcal{F}_C$ involves singular
distributions 
(e.g.\ normal derivatives of Dirac's delta function with respect to a hypersurface)
but yet that the operator $\mathcal{F}_C$ is unitary in the general
case where $n_1,n_2>1$ and $n_1+n_2>4$.
In the case $n_1=1$, $n_2=1$ or $(n_1,n_2)=(2,2)$,
$\mathcal{F}_C$ reduces to the 
\index{B}{Hankel transform}%
Hankel transform
 composed by a
(singular) Radon transform.

The second motivation comes from representation theory of real reductive groups, in particular, 
from minimal representations.

\index{B}{minimal representation}%
Minimal representations 
 are infinite 
dimensional unitary representations that are the `closest' to the trivial one-dimensional
 representation.
The Weil representation of the metaplectic group
\index{A}{MpnR@$\Mp(n, \mathbb R)$}%
$\Mp(n, \mathbb R)$, which plays a prominent
role in the construction of theta series,
is a classic example.
Most minimal representations
 are isolated among the set of
irreducible unitary representations,
and cannot be built up from the exsisting
induction techniques of representation theory.

A multitude of different models of minimal representations have been
investigated recently by many people (see Sections \ref{subsec:1.4}
and \ref{subsec:1.4b}).
Each model known so far has its own advantages indeed but also has some disadvantages. 
For instance, the inner product of the Hilbert space is not explicit in some models,
whereas the whole group action 
is not clear in some other models.

A challenge to surmount that `disadvantage' 
 may turn up
 as a natural problem in other areas of mathematics.
In order to give its flavor,
let us consider two geometric models
 of minimal representations
of the indefinite orthogonal group
$G = O(n_1+1,n_2+1)$:
One is in the solution space to the Yamabe equation
 (\textit{conformal model}%
\index{B}{conformal model}%
),
 and the other is in 
 $L^2(C)$ (\textit{Schr\"{o}dinger model}%
\index{B}{Schr\"{o}dinger model}%
).

In the conformal model,
the whole group action is very clear,
whereas the inner product is not. 
The problem of finding the explicit inner product
was solved in the previous paper \cite{xkors1} as the theory of conserved
quantities for ultra-hyperbolic
 equations, such as the energy for
 the wave equation.

In the Schr\"{o}dinger model $L^2(C)$, 
the unitary structure is clear, 
whereas the
whole group action is not.
The understanding of the whole group action 
was a missing piece of \cite{xkors1}.
This problem is reduced to finding
 the generalization of the Fourier--Hankel transform on the
isotropic cone $C$, namely,
the above mentioned operator $\mathcal{F}_C$. 
By finding an explicit formula of $\mathcal{F}_C$, 
 we shall settle this problem. 
The role of $\mathcal{F}_C$ in our minimal representation is in
parallel to that of the Euclidean Fourier transform
$\mathcal{F}_{\mathbb{R}^n}$ in the Weil
representation, 
summarized as below:
\begin{center}
\begin{tabular}{l|cc}
  simple group
   & $Mp(n,\mathbb{R})$
   & $O(n_1+1,n_2+1)$
\\
   & (type $C$) 
   & (type $D$)
\\
\hline
\\
  Minimal representation
   & 
\index{B}{Weil representation}%
     Weil representation
   & $\pi$
\\
   & &
\\
\hline
\\
  $L^2$-model
   & $L^2(\mathbb{R}^n)$
   & $L^2(C)$
\\
  (Schr\"{o}dinger model)
  &&
\\
\hline
\\[-\medskipamount]
  unitary inversion
   & $e^{\frac{\sqrt{-1}n\pi}{4}} \mathcal{F}_{\mathbb{R}^n}$
   & $\mathcal{F}_C$
\\[\medskipamount]
\end{tabular}
\end{center}

The third motivation comes from special functions.
We note that the isotropic cone $C$ is so small that the group 
$G = O(n_1+1,n_2+1)$ cannot act on $C$ continuously and non-trivially.
This feature is reflected by the fact that the
 Gelfand--Kirillov dimension of the representation of $G$ on $L^2(C)$ 
attains its minimum amongst all infinite dimensional representations of
$G$.
Thus,
the representation space $L^2(C)$ is
extremely `small' with
respect to the group $G$.
In turn,
we could expect a very concrete theory of global
analysis on $C$ by using abundant symmetries of 
the group $G$ or its Lie algebra.

It turns out that
special functions in the Schr\"{o}dinger model $L^2(C)$ arise in a somewhat
different way from the well-known cases such as analysis on symmetric
spaces (e.g.\ \cite{xHe}) or its variants.
For instance, the Casimir operator of $K$ acts on $L^2(C)$ as a \textit{fourth}
differential operator. 

In this book,
we encounter many classically known special
functions (e.g.\ Bessel functions, Appell's hypergeometric
functions, Meijer's $G$-functions, etc.).
Special functions are a part of our method 
for the analysis of the minimal representation,
and conversely, by decomposing the operator $\mathcal{F}_C$ we provide
 a representation theoretic proof of 
[inversion, Plancherel, \dots] formulas of special functions
including Meijer's $G$-functions.

\medskip

Encouraged by a suggestion of R. Stanton,
we have decided to write a considerably
 long introduction.
What follows is divided into three parts
according to the aforementioned three motivations and new perspectives.
In Sections \ref{subsec:1.1}--\ref{subsec:FC},
we state key properties of the involutive unitary operator 
$\mathcal{F}_C$ on $L^2(C)$ from analytic perspectives,
in comparison with the well-known case of the Euclidean Fourier
transform $\mathcal{F}_{\mathbb{R}^n}$ on $L^2(\mathbb{R}^n)$.
Sections \ref{subsec:1.4}--\ref{subsec:uncertainty} 
give representation theoretic perspectives,
and we explain the role of $\mathcal{F}_C$ in the Schr\"{o}dinger
model of the minimal representation of the indefinite orthogonal group
in comparison with the role of $\mathcal{F}_{\mathbb{R}^n}$ for the
Weil representation.
Thus, 
we compare $\mathcal{F}_C$ again with $\mathcal{F}_{\mathbb{R}^n}$,
and correspondingly, the simple Lie algebra 
$\mathfrak{o}(n_1+1,n_2+1)$ with $\mathfrak{sp}(n,\mathbb{R})$.
In Section \ref{subsec:1.7},
we give a flavor of the interactions of the analysis on the minimal
representations with special functions.

\section{Differential operators on the isotropic cone}
\label{subsec:1.1}

Consider an indefinite quadratic form on
$\mathbb{R}^n=\mathbb{R}^{n_1+n_2}$: 
\begin{equation}\label{def:Q}
   Q(x)
   := x_1^2+\dots+x_{n_1}^2-x_{n_1+1}^2-\dots-x_{n_1+n_2}^2.
\end{equation}
Throughout the Introduction,
we assume $n_1,n_2>1$ and $n=n_1+n_2$ is an even integer greater than
two. 
(From Chapter \ref{sec:rev}, we will use the following 
notation: $p=n_1-1$, $q = n_2-1$.)

Associated to the quadratic form $Q$, we define the 
isotropic cone%
\index{B}{isotropic cone}
\begin{equation*}
   C := \{ x\in\mathbb{R}^n \setminus \{0\}: Q(x) = 0 \},
\index{A}{C@$C$}
\end{equation*}
and endow $C$ with the volume form $d\mu$ characterized by 
\begin{equation*}
  dQ \wedge d\mu = dx_1 \cdots dx_n.
\end{equation*}

Our object of study is the Hilbert space 
\index{A}{L2C@$L^2(C)$|main}%
$L^2(C)\equiv L^2(C,d\mu)$
consisting of square integrable functions on $C$.

A differential operator $P$ on $\mathbb{R}^n$ is said to be
\textit{tangential to} the submanifold $C$ if $P$ satisfies
\begin{equation} \label{eqn:tg}
   \psi_1|_C  = \psi_2|_C
   \Rightarrow
   (P\psi_1)|_C = (P\psi_2)|_C
\end{equation}
for any smooth functions $\psi_1, \psi_2$ defined in neighborhoods of
$C$ in $\mathbb{R}^n$. 
Then, we can `restrict' $P$ to $C$,
and get a differential operator $P|_C$ on $C$.

For instance, 
the following vector fields are tangential to $C$:
\begin{equation*}
\index{A}{E@$E$|main}%
   E := \sum_{i=1}^n x_i \frac{\partial}{\partial x_i}
   \quad\text{(the Euler operator)},%
\index{B}{Euler operator}
\end{equation*}
\begin{equation*}
   X_{ij}
   := \epsilon_i\epsilon_j x_i \frac{\partial}{\partial x_j}
      - x_j \frac{\partial}{\partial x_i}
  \qquad (1 \le i < j \le n),
\end{equation*}
where we set $\epsilon_j=1$ or $-1$ according as $1\le j\le n_1$ or 
$n_1+1\le j \le n$.
This is because the vector fields $E$ and $X_{ij}$ $(1\le i<j\le n)$ are
obtained as the differential of the conformal linear transformation
group 
\begin{equation*}
   CO(Q) := \{ g\in GL(n,\mathbb{R}):
               \text{$Q(gx)=cQ(x)$
                     $({}^\forall x\in \mathbb{R}^n)$ for some $c>0$} \},
\index{A}{CO(Q)@$CO(Q)$}
\end{equation*}
which preserves the isotropic cone $C$.

Let 
\index{A}{Rx@$\mathbb{R}[x,\frac{\partial}{\partial x}]$}%
$\mathbb{R}[x,\frac{\partial}{\partial x}]$ be the
$\mathbb{R}$-algebra of differential operators with polynomial
coefficients (the \textit{Weyl algebra}%
\index{B}{Weyl algebra}%
), namely, the non-commutative
ring generated by the multiplication by $x_1,\dots,x_n$ and the vector
fields 
$\frac{\partial}{\partial x_1},\dots,\frac{\partial}{\partial x_n}$.

We denote by
\index{A}{RxC@$\mathbb{R}[x,\frac{\partial}{\partial x}]^C$}%
$\mathbb{R}[x,\frac{\partial}{\partial x}]^C$ 
the subalgebra consisting of operators that are tangential to $C$. 
The multiplication by coordinate functions $x_k$ clearly
satisfies the condition \eqref{eqn:tg}.
Thus, we have seen
\begin{equation*}
   x_k, E, X_{ij} \in \mathbb{R} \Bigl[ x,\frac{\partial}{\partial x} \Bigr]^C
   \quad (1 \le k \le n, \  1 \le i < j \le n).
\end{equation*}
However, there exist yet other operators which are tangential
to $C$,
but are not generated by $x_k,E, X_{ij}$  in the Weyl
algebra (see Remark \ref{rem:Pjout}). 

Among them are the 
\index{B}{fundamental differential operator}%
\textit{fundamental differential operators} 
of second order,
to be denoted by $P_1,\dots,P_n$, which are defined by
\begin{equation} \label{eqn:Pjdef}
\index{A}{P_j@$P_j$}%
   P_j
   := \epsilon_j x_j \square - (2E+n-2) \frac{\partial}{\partial x_j}.
\end{equation}
Here, 
\index{A}{2square@\par\indexspace$\square$}%
$\square$ is the
Laplace--Beltrami operator associated to $Q$, namely,
\begin{equation*}
   \square := \sum_{j=1}^n \epsilon_j
                 \frac{\partial^2}{\partial x_j^2}
          \equiv \frac{\partial^2}{\partial x_1^2}+\dots+
                 \frac{\partial^2}{\partial x_{n_1}^2}-
                 \frac{\partial^2}{\partial x_{n_1+1}^2} -\dots-
                 \frac{\partial^2}{\partial x_n^2}.
\end{equation*}
In the degenerate case $n_1=n_2=1$,
our operators $P_1$ and $P_2$ take the following form:
we set $y_1:=x_1+x_2$, \thinspace  $y_2:=x_1-x_2$, 
\begin{equation*}
   P_1 + P_2 = -4y_1\frac{\partial^2}{\partial y_1^2},
   \quad 
   P_1 - P_2 = -4y_2\frac{\partial^2}{\partial y_2^2},
\end{equation*}
see Remark \ref{rem:Pn2}.
In general, 
these operators $P_1,\dots,P_n$ satisfy the following properties 
(see Theorem \ref{thm:Pj}):
\begin{enumerate}[{\bf {P}1}]
\index{A}{1zpropertiesP1-P5@\textbf{P1}--\textbf{P5}|(}%
\item \label{item:Pj1}
  $P_iP_j = P_jP_i$ for any $1 \le i,j \le n$.
\item \label{item:Pj2}
   $P_j \in \mathbb{R} [x,\frac{\partial}{\partial x}]^C$ for any
   $1 \le j \le n$.
\item \label{item:Pj3}
   The induced differential operators $P_j|_C$ on $C_0^\infty(C)$
   extend to self-adjoint
   operators on the Hilbert space $L^2(C)$.
\item \label{item:Pj4}
   $(P_1^2+\dots+P_{n_1}^2-P_{n_1+1}^2-\dots-P_n^2)|_C = 0$.
\item \label{item:Pj5}
The Lie algebra generated by $x_i$,
$P_i$ $(1 \le i \le n)$ contains the vector fields 
$E, X_{ij}$ $(1\le i<j\le n)$.
\index{A}{1zpropertiesP1-P5@\textbf{P1}--\textbf{P5}|)}%
\end{enumerate}

{}From now on, we simply write $P_j$ for $P_j|_C$.
Thus, we have commuting self-adjoint, 
second-order differential operators $P_1,\dots,P_n$ on $L^2(C)$.

We are brought naturally to the following:
\begin{problem} \label{prob:Pj}
{\upshape 1)}\enspace
Find joint eigenfunctions of the differential operators
$P_1,\dots,P_n$ on the isotropic cone $C$.

{\upshape 2)}\enspace
Given a function $f$ on $C$,
find an explicit expansion formula of $f$ into joint
eigenfunctions of $P_1,\dots,P_n$.
\end{problem}

\section{`Fourier transform' $\mathcal{F}_C$ on the isotropic cone
$C$}
\label{subsec:1.2}

In this book,
we shall give a solution to Problem \ref{prob:Pj}
by introducing a unitary operator $\mathcal{F}_C$ on $L^2(C)$.

To elucidate the operator $\mathcal{F}_C$,
let us consider first much simpler operators
\begin{equation*}
   p_j := -\sqrt{-1} \frac{\partial}{\partial x_j}
   \quad (1 \le j \le n)
\end{equation*}
in place of $P_j$.
Then, $p_1,\dots,p_n$ form a commuting family of differential
operators which extend to self-adjoint operators on
$L^2(\mathbb{R}^n)$.
Analogously to Problem \ref{prob:Pj},
consider the question of finding the explicit eigenfunction
 expansion for the operators
$p_1,\dots,p_n$. 
Then, as is well-known, 
this is done by using
the (Euclidean) Fourier transform 
\index{A}{FRn@$\mathcal{F}_{\mathbb{R}^n}$|main}%
$\mathcal{F}\equiv\mathcal{F}_{\mathbb{R}^n}$ on
$\mathbb{R}^n$. 
In what follows, we normalize $\mathcal{F}_{\mathbb{R}^n}$ as
\begin{equation}\label{eqn:Fourier}
\index{A}{FRn@$\mathcal{F}_{\mathbb{R}^n}$|main}%
  \mathcal{F}_{\mathbb{R}^n} u(\xi) :=
     \frac{1}{(2\pi)^{\frac{n}{2}}}
     \int_{\mathbb{R}^n} u(x) e^{\sqrt{-1}\langle x,\xi\rangle} dx,
\end{equation}
where
$\langle x,\xi\rangle = \sum_{i=1}^n x_i \xi_i$
and $dx = dx_1\cdots dx_n$.
We note that the signature of the power here is opposite from the
usual convention. 
Obviously, the kernel 
\begin{equation*}
   k(x,\xi) := \frac{1}{(2\pi)^{\frac{n}{2}}}
               \, e^{\sqrt{-1}\langle x,\xi \rangle}
\end{equation*}
of the Fourier transform $\mathcal{F}_{\mathbb{R}^n}$ is real analytic
on the direct product space $\mathbb{R}^n \times \mathbb{R}^n$.

We recall the following key properties of the 
Euclidean Fourier transform%
\index{B}{Euclidean Fourier transform}%
: 
\begin{enumerate}[{\bf {F}1}]
\index{A}{1zpropertiesF1-F6@\textbf{F1}--\textbf{F8}|(}%
\item \label{item:F1}
$p_j \, k(x,\xi) = \xi_j \, k(x,\xi)$.
\item \label{item:F2}
$k(x,\xi) = k(\xi,x)$.
\item \label{item:F3}
$\mathcal{F}_{\mathbb{R}^n} (C_0^\infty(\mathbb{R}^n))
 \subset C^\infty(\mathbb{R}^n) \cap L^2(\mathbb{R}^n)$.
\item \label{item:F4}
$\mathcal{F}_{\mathbb{R}^n}$ extends to a unitary operator on
$L^2(\mathbb{R}^n)$.
\item \label{item:F5}
$\mathcal{F}_{\mathbb{R}^n} \circ x_j
 = p_j
   \circ \mathcal{F}_{\mathbb{R}^n},
$

$\mathcal{F}_{\mathbb{R}^n} \circ 
   p_j
 = -\xi_j \circ \mathcal{F}_{\mathbb{R}^n}.
$
\item \label{item:F6}
$(\mathcal{F}_{\mathbb{R}^n}^2 u)(x) = u(-x)$,
$\mathcal{F}_{\mathbb{R}^n}^4 = \operatorname{id}$.
\item \label{item:F7}
$\mathcal{F}_{\mathbb{R}^n} u(\xi)
 = \displaystyle\frac{1}{(2\pi)^{\frac{n}{2}}}
   \int_{\mathbb{R}} e^{\sqrt{-1}t} (Ru)(\xi,t)dt
$.
\item \label{item:F8}
$\mathcal{F}_{\mathbb{R}^n}$ gives an automorphism of each of the following
topological vector spaces:
\begin{center}
$\mathcal{S}(\mathbb{R}^n) \subset L^2(\mathbb{R}^n) \subset
\mathcal{S}'(\mathbb{R}^n).$
\end{center}
\index{A}{1zpropertiesF1-F6@\textbf{F1}--\textbf{F8}|)}%
\end{enumerate}
Here, \textbf{F\ref{item:F7}} gives the plane wave expansion of the Fourier
transform by means of the Radon transform $R$ defined by
\begin{equation*}
   Ru(\xi,t) := \int_{\mathbb{R}^n} u(x) \delta
   (\langle x,\xi \rangle -t) dx.
\end{equation*}
In \textbf{F\ref{item:F8}},
we denote by
\index{A}{S1@$\mathcal{S}(\mathbb{R}^n)$}%
$\mathcal{S}(\mathbb{R}^n)$
the space of rapidly decreasing $C^\infty$-functions on $\mathbb{R}^n$
(the Schwartz space%
\index{B}{Schwartz space} 
 endowed with the Fr\'{e}chet topology),
and by 
\index{A}{S1'@$\mathcal{S}'(\mathbb{R}^n)$}%
$\mathcal{S}'(\mathbb{R}^n)$
the dual space consisting of tempered distributions.

\textbf{F\ref{item:F4}} is the Plancherel theorem,
and \textbf{F\ref{item:F8}} gives the Paley--Winer theorem for the
Schwartz space $\mathcal{S}(\mathbb{R}^n)$
(and its dual $\mathcal{S}'(\mathbb{R}^n)$).
By \textbf{F\ref{item:F1}} and \textbf{F\ref{item:F2}},
the inversion formula
\begin{align*}
   f(x)
   & = \mathcal{F}_{\mathbb{R}^n} \circ
       \mathcal{F}_{\mathbb{R}^n}^{-1} f(x)
\\
   & = \int_{\mathbb{R}^n} (\mathcal{F}_{\mathbb{R}^n}^{-1} f)
       (\xi) k(x,\xi) d\xi
\end{align*}
gives an expansion of a function $f$ into joint eigenfunctions
$k(x,\xi)$ of the commuting self-adjoint operators $p_j$ $(1 \le j \le n)$.

Moreover, the property \textbf{F\ref{item:F5}} characterizes 
the operator $\mathcal{F}_{\mathbb{R}^n}$ up
to scalar.
We pin down this algebraic statement in two ways as follows:
\begin{proposition} \label{prop:Apxj}
Let $A$ be a continuous operator on $L^2(\mathbb{R}^n)$ satisfying the
following identities:
\begin{equation} \label{eqn:Apxj}
   A\circ x_j = p_j \circ A, \ 
   A \circ p_j = -x_j \circ A
   \quad (1 \le j \le n).
\end{equation}
Then, $A$ is a scalar multiple of $\mathcal{F}_{\mathbb{R}^n}$.
\end{proposition}

\begin{proposition} \label{prop:Apxij}
Let $A$ be a continuous operator on $L^2(\mathbb{R}^n)$ satisfying
\begin{equation} \label{eqn:Apxij}
   A \circ x_i x_j = p_i p_j \circ A, \ 
   A \circ p_i p_j = x_i x_j \circ A
   \quad (1 \le i,j \le n).
\end{equation}
Then, $A$ is of the form 
$A = a \mathcal{F}_{\mathbb{R}^n} + b \mathcal{F}_{\mathbb{R}^n}^{-1}$
for some $a,b \in \mathbb{C}$.
\end{proposition}
Here, 
\eqref{eqn:Apxij} is obviously a weaker condition than
\eqref{eqn:Apxj}.

We did not go into details about the domain of definition for 
\eqref{eqn:Apxj} and \eqref{eqn:Apxij} in the above propositions.
The domain could be 
$D:=\{ f\in L^2(\mathbb{R}^n): x_j f, p_j f \in L^2(\mathbb{R}^n) 
  \ (1 \le {}^\forall j \le n) \}$,
on which we regard the identities
$A x_j f = p_j A f$  those of distributions
in the case of 
Proposition \ref{prop:Apxj}.
Likewise for Proposition \ref{prop:Apxij}.

\medskip
\noindent
\textbf{Intertwining characterization of $\mathcal{F}_C$}
\nopagebreak

Back to the setting in Section \ref{subsec:1.1},
we consider the differential operator $P_j$
(of second order!).
Then, 
it turns out that
the intertwining relation between $P_j$
and the multiplication by the coordinate function
$x_j$ again characterizes our operator $\mathcal{F}_C$ up to scalar:
\begin{thmSec}[see Theorem \ref{thm:FCunique}]
 \label{thm:APxj}
\begin{enumerate}[\upshape1)]
\item 
There exists a unitary operator $\mathcal{F}_C$ on $L^2(C)$ satisfying
the following relation:
\begin{equation} \label{eqn:APxj}
   A \circ P_j = 4x_j \circ A, \ 
   A \circ x_j = 4P_j \circ A
   \quad (1 \le j \le n).
\end{equation}
\item 
Continuous operators $A$ on $L^2(C)$ satisfying \eqref{eqn:APxj} 
are unique up to scalar.
In particular, 
any such operator $A$ is a scalar multiple of the unitary operator 
$\mathcal{F}_C$,
and $A^2$ is a scalar multiple of the identity operator.
\end{enumerate}
\end{thmSec}

\medskip
\noindent
\textbf{Joint eigendistributions on $C$}
\nopagebreak

Next, we consider the system of differential
equations on $C$: 
\begin{equation} \label{eqn:Pjsys}
   P_j \psi = 4\xi_j \psi
   \qquad (1 \le j \le n).
\end{equation}
The coefficient $4$ in the right-hand side is just for 
simplifying later notation. 

We shall deal with solutions in an appropriate class of distributions on $C$
(the dual space 
\index{A}{L2Cminusinfty@$L^2(C)^{-\infty}$}%
$L^2(C)^{-\infty}$ 
of smooth vectors%
\index{B}{smooth vector} 
 $L^2(C)^\infty$,
see \eqref{eqn:5C})
and prove the following:
\begin{thmSec}[see Theorem \ref{thm:Pjunique}]
 \label{thm:SolPx}
Fix $\xi = (\xi_1,\dots,\xi_n) \in \mathbb{R}^n \setminus \{0\}$. 
\begin{enumerate}[\upshape1)]
\item 
If $Q(\xi)\ne0$, then any distribution $\psi$ on $C$ satisfying
\eqref{eqn:Pjsys} is zero.
\item 
If $\xi\in C$, then the solution space of \eqref{eqn:Pjsys} in
$L^2(C)^{-\infty}$ is one-dimensional.
\end{enumerate}
\end{thmSec}
The first statement is an immediate consequence of 
{\bf P\ref{item:Pj4}}.
By the explicit formula given in Theorem \ref{thm:Kformula},
we shall see that the unique solution in Theorem \ref{thm:SolPx} (2) is
not real analytic if $n_1, n_2 > 1$.

\medskip
\noindent
\textbf{Abstract properties of $\mathcal{F}_C$}
\nopagebreak

We will prove in this book that the distribution solution $\psi(x)$ 
in Theorem \ref{thm:SolPx} (2) can be normalized depending on $\xi\in C$,
which we denote by $K(x,\xi)$ for now,
in such a way that the following key properties 
 are fulfilled:
\begin{enumerate}[{\bf {K}1}]
\index{A}{1zpropertiesK1-K6@\textbf{K1}--\textbf{K8}|(}%
\item \label{item:K1}
For each fixed $\xi\in C$,
$K(\cdot,\xi)$ is a distribution solution on $C$ to
\eqref{eqn:Pjsys}. 
\item \label{item:K2}
$K(x,\xi) = K(\xi,x)$ as a distribution on $C \times C$.
\item \label{item:K3}
We define
\begin{equation} \label{eqn:FCdef}
(\mathcal{F}_Cf)(\xi) := \int_C K(x,\xi)f(x)d\mu(x)
\qquad\text{for $f\in C_0^\infty(C)$}.
\end{equation}
Then, \eqref{eqn:FCdef} is well-defined, and we have
a linear map
$\mathcal{F}_C: C_0^\infty(C) \to C^\infty(C) \cap L^2(C)$.
\item \label{item:K4}
\index{A}{FC@$\mathcal{F}_C$}%
$\mathcal{F}_C$
 extends to a unitary operator on $L^2(C)$.
\item \label{item:K5}
$
   \mathcal{F}_C \circ 4x_j = P_j \circ \mathcal{F}_C
$,

$   \mathcal{F}_C \circ P_j  = 4\xi_j \circ \mathcal{F}_C
$.
\item \label{item:K6}
$\mathcal{F}_C^2 = \operatorname{id}$.
\item \label{item:K7}
$\mathcal{F}_C u(\xi)
 = \int_{\mathbb{R}} \Psi(t) \mathcal{R}f(\xi,t) dt
$
\item \label{item:K8}
$\mathcal{F}_C$ gives the automorphism of each of the following topological
vector spaces:
\begin{center}
$   L^2(C)^\infty \subset L^2(C) \subset L^2(C)^{-\infty}.$
\end{center}
\index{A}{1zpropertiesK1-K6@\textbf{K1}--\textbf{K8}|)}%
\end{enumerate}
These properties \textbf{K\ref{item:K1}}--\textbf{K\ref{item:K8}} are
stated in parallel to the Euclidean case 
\textbf{F\ref{item:F1}}--\textbf{F\ref{item:F8}}. 
In \textbf{K\ref{item:K7}},
$\mathcal{R}$ is the (singular) Radon transform on the isotropic cone $C$ which
will be defined in \eqref{eqn:RadonC},
and $\Psi(t)$ is a distribution on $\mathbb{R}$ which will be defined
in Theorem \ref{thm:Kformula}.
We note that the transform in \textbf{K\ref{item:K7}} by $\Psi(t)$ 
collapses to the Hankel if $n_2 = 1$.

In \textbf{K\ref{item:K8}},
we have the following inclusive relation
\begin{equation*}
   C_0^\infty(C) \subset L^2(C)^\infty \subset L^2(C)
   \subset L^2(C)^{-\infty} \subset \mathcal{D}'(C)
\end{equation*}
as in the Euclidean case
(see \textbf{F\ref{item:F8}}): 
\begin{equation*}
   C_0^\infty(\mathbb{R}^n) \subset \mathcal{S}(\mathbb{R}^n)
   \subset L^2(\mathbb{R}^n) \subset \mathcal{S}'(\mathbb{R}^n)
   \subset \mathcal{D}'(\mathbb{R}^n).
\end{equation*}

In summary, \textbf{K\ref{item:K4}} is a Plancherel type theorem of
$\mathcal{F}_C$ on $L^2(C)$,
\textbf{K\ref{item:K6}} gives its inversion formula,
\textbf{K\ref{item:K7}} expresses $\mathcal{F}_C$ by `plane wave'
decomposition, 
and \textbf{K\ref{item:K8}} gives a Paley--Winer type theorem for the
`Schwartz space' $L^2(C)^\infty$.

The above formulation brings us naturally to the following program:
\begin{program}\label{prog:FCanalysis}
Develop a theory of `Fourier analysis' on the isotropic cone $C$ by
means of $\mathcal{F}_C$.
\end{program}
We expect that this program could be enhanced by a solid foundation
and concrete formulas of the transform $\mathcal{F}_C$.

For this, the first step is to find 
explicit formulas for the (normalized) joint eigenfunctions
$K(x,\xi)$.
We prove that they are given by means of \textit{Bessel distributions} 
(see Theorem \ref{thm:Kformula}).
In particular, \textbf{K\ref{item:K2}} follows readily from the formulas.
The properties \textbf{K\ref{item:K1}}, \textbf{K\ref{item:K4}},
\textbf{K\ref{item:K5}}, and \textbf{K\ref{item:K6}} will be proved in Theorem
\ref{thm:FC} based on a representation theoretic
interpretation that $\mathcal{F}_C$ is the 
\index{B}{unitary inversion operator}%
`unitary inversion operator' 
on $L^2(C)$ for the minimal representation of the
indefinite orthogonal group $O(n_1+1,n_2+1)$.
By \textbf{K\ref{item:K6}}, we get the inversion formula just as 
$\mathcal{F}_C^{-1}=\mathcal{F}_C$, which gives an explicit solution to 
the problem of joint eigenfunction expansions
(see Problem \ref{prob:Pj}):
\begin{equation*}
   f(x) = \int_C (\mathcal{F}_C f)(\xi) K(x,\xi) d\mu(\xi).
\end{equation*}

The kernel $K(x,\xi)$ is not locally integrable but is a distribution in
general.
To see the convergence of the right-hand side \eqref{eqn:FCdef},
we note that $K(x,\xi)$ depends only on
$\langle x,\xi \rangle = \sum_{i=1}^n x_i \xi_i$
(see Theorem \ref{thm:Kformula}).
This fact leads us to the factorization \textbf{K\ref{item:K7}}
 through the (singular) Radon
 transform 
\index{A}{R@$\mathcal{R}$}%
$\mathcal{R}$ 
on the isotropic cone $C$, which is defined by the integration over the
intersection of $C$ with the hyperplane
\[
  \{ x \in \mathbb{R}^n : \langle x,\xi \rangle = t \}.
\]

For a quick summary of the transform $\mathcal{R}$ 
(see Chapter \ref{sec:A} for details),
we identify a compactly supported smooth function $f$ on $C$
 with a measure $f d\mu$.
It is a tempered distribution 
 on $\mathbb{R}^n$ $(n > 2)$.
Then, the 
Radon transform%
\index{B}{Radon transform} 
$\mathcal{R}$ of $f d\mu$ is defined by
\begin{equation}\label{eqn:RadonC}
   \mathcal{R} f(\xi,t)
   := \int_C f(x) \delta(\langle x,\xi \rangle -t) d\mu (x)
\end{equation}
for $(\xi,t)\in(\mathbb{R}^n\setminus\{0\})\times(\mathbb{R}\setminus\{0\})$.
The point here is that the integration is taken over the isotropic
cone.
In other words,
$\mathcal{R}f(\xi,t)$ is obtained by the integration over submanifolds
which are generically of codimension two in $\mathbb{R}^n$. 
Consequently,
$\mathcal{R} f(\xi,t)$ satisfies the 
\index{B}{ultra-hyperbolic equation}%
ultra-hyperbolic differential equation 
 of the $\xi$-variable:
\begin{equation}\label{eqn:RfPDE}
   \Bigl( \sum_{j=1}^{n_1} \frac{\partial^2}{\partial\xi_j^2}
     \, - \sum_{j=n_1+1}^{n} \frac{\partial^2}{\partial\xi_j^2}
   \Bigr)
   \mathcal{R} f(\xi,t)
   = 0.
\end{equation}

Next, in order to see the regularity of $\mathcal{R}f(\xi,t)$ at
$t=0$, we fix $\xi$. Then,
the intersection of the isotropic cone $C$ with the hyperplane
$\{ x\in \mathbb{R}^n: \langle x,\xi \rangle = t \}$
forms a one parameter
 family of submanifolds of codimension two for $t \ne 0$,
which have singularities at $t=0$.
Accordingly, the Radon transform
$\mathcal{R} f(\xi,t)$ is not of $C^\infty$ class at 
$t = 0$ even for 
\index{A}{C0infty@$C_0^\infty(C)$}%
$f \in C_0^\infty(C)$.
The regularity of $\mathcal{R}f(\xi,t)$ 
at $t = 0$ is the principal object of the paper \cite{xkmano5},
where it is proved that $\mathcal{R} f(\xi,t)$ is
$[\frac{n-5}{2}]$ times continuously differentiable at $t = 0$.
Here, $[x]$ denotes the greatest integer that does not exceed $x$.
This regularity is exactly sufficient for what we need to prove
 that the singular integral
\eqref{eqn:FCdef} makes sense for $f \in C_0^\infty(C)$.
See Section \ref{subsec:Radon} for details. 

The reverse direction,
namely, the application of our results on $\mathcal{F}_C$ to the
results on the singular Radon transform $\mathcal{R}$ includes:
\begin{corollary} \label{cor:Rfinverse}
Any compactly supported smooth function
$f \in C_0^\infty(C)$ can be recovered only from the restriction of
the Radon transform $\mathcal{R} f(\xi,t)$ to $C \times \mathbb{R}$.
\end{corollary}

\medskip
\noindent
\textbf{Underlying algebraic structures}
\nopagebreak

The underlying algebraic structure of Propositions \ref{prop:Apxj} and
\ref{prop:Apxij} and Theorem \ref{thm:APxj} will be revealed by the Lie
algebras generated by the differential operators in each setting as follows:

In Proposition \ref{prop:Apxj},
the Lie algebra generated by $p_i, x_i$ $(1 \le i \le n)$ is the
Heisenberg Lie algebra%
\index{B}{Heisenberg Lie algebra}%
.

In Proposition \ref{prop:Apxij},
the Lie algebra generated by $p_i p_j, x_i x_j$ $(1 \le i, j \le n)$
is the 
symplectic Lie algebra%
\index{B}{symplectic Lie algebra} 
\index{A}{spnR@$\mathfrak{sp}(n,\mathbb{R})$}%
 $\mathfrak{sp}(n,\mathbb{R})$.

In Theorem \ref{thm:APxj},
the Lie algebra generated by $P_i, x_i$ $(1 \le i \le n_1+n_2)$ is the
indefinite orthogonal Lie algebra%
\index{B}{indefinite orthogonal Lie algebra} 
 $\mathfrak{o}(n_1+1,n_2+1)$.

These actions of the Lie algebras lift to 
unitary representations of the corresponding Lie groups:
On $L^2(\mathbb{R}^n)$,
the Schr\"{o}dinger representation of the Heisenberg group,
and the 
\index{B}{Weil representation}%
Weil representation of the metaplectic group
$Mp(n,\mathbb{R})$,
namely, a double cover of the symplectic group 
$Sp(n,\mathbb{R})$%
\index{A}{SpnR@$Sp(n, \mathbb R)$} 
(see Section \ref{subsec:1.4});
on $L^2(C)$,
the minimal representation of the indefinite orthogonal group
$O(n_1+1,n_2+1)$
(see Section \ref{subsec:1.4b}).
In Sections \ref{subsec:1.4}--\ref{subsec:uncertainty}, 
we shall discuss some perspectives from representation theory.
Before entering representation theory, 
we continue an account from the viewpoints of analysis
in the next section.

\section{Kernel of $\mathcal{F}_C$ and Bessel distributions}
\label{subsec:FC}

In this section,
we give an explicit formula of the kernel distribution $K(x,\xi)$ of
the transform $\mathcal{F}_C$ on $L^2(C)$,
and observe the similarities to and differences
from the kernel 
$k(x,\xi) =
(2\pi)^{-\frac{n}{2}} \, e^{\sqrt{-1}\langle x,\xi \rangle}$
of the Euclidean Fourier transform $\mathcal{F}_{\mathbb{R}^n}$ on
$L^2(\mathbb{R}^n)$.

Here is the highlight of this book.
\begin{theorem}[see Theorem \ref{thm:A}]
 \label{thm:Kformula}
Suppose $n=n_1+n_2$ is even, $>2$.
The unique unitary operator 
$\mathcal{F}_C$%
\index{A}{FC@$\mathcal{F}_C$} 
 on $L^2(C)$ in Theorem
\ref{thm:APxj} is given by the distribution kernel
\index{A}{Kxx'@$K(x,x')$}%
$K(x,\xi) := \Psi(\langle x,\xi \rangle)$, where
\begin{align*}
\Psi(t)
={}& 2(-1)^{\frac{n_1(n_1-1)}{2}} \pi^{-\frac{n-2}{2}}
\\
&\times
 \begin{cases}
    \Phi_{\frac{n-4}{2}}^+ (t)
    &\text{if\/ $\min(n_1,n_2)=1$},
    \\[1ex]
    \Psi_{\frac{n-4}{2}}^+ (t)
    &\text{if\/ $n_1,n_2>1$ are both odd},
    \\[1ex]
    \Psi_{\frac{n-4}{2}} (t)
    &\text{if\/ $n_1,n_2>1$ are both even}.
 \end{cases}
\end{align*}
\end{theorem}  
\noindent
As for the normalization of a scalar constant in the above theorem,
we note that the intertwining property \eqref{eqn:APxj} determines
$K(x,\xi)$ up to scalar,
and moreover \textbf{K\ref{item:K6}} determines $K(x,\xi)$ up to
signature. 
The signature is taken
to be compatible with the action on the Schr\"{o}dinger model of 
the minimal representation of
$O(n_1+1,n_2+1)$ which will be discussed in Section
\ref{subsec:uncertainty}. 

In Theorem \ref{thm:Kformula}, 
$\langle \ , \ \rangle$ denotes the standard (positive definite)
inner product on $\mathbb{R}^n$.
$\Phi_m^+(t)$, 
$\Psi_m^+(t)$ and $\Psi_m(t)$ are tempered distributions on
$\mathbb{R}$, 
defined below in \eqref{def:Psi0}, \eqref{def:Psi+}, and
\eqref{def:Psi}, respectively.

\medskip
\noindent
\textbf{Bessel distributions}%
\index{B}{Bessel distribution|main}
\nopagebreak

Let $J_\nu(x)$, $Y_\nu(x)$ and $K_\nu(z)$ be 
the (modified) Bessel functions (see Appendix \ref{subsec:B}).
We use the following notational convention:
$$
  f(t_+) := \begin{cases}  f(t) &(t>0)  \\
                              0 &(t\le0),
            \end{cases}
  \quad 
  f(t_-) := \begin{cases}  0  &(t\ge0)  \\
                        f(|t|)&(t<0),
            \end{cases}
$$
for a function (or a `generalized function') $f(t)$ on $\mathbb{R}$.
Then, 
$\Phi_m^+$, $\Psi_m^+$, and $\Psi_m$ 
in Theorem \ref{thm:Kformula} are the distributions 
on $\mathbb{R}$ given by 
\begin{alignat}{2} \label{def:Psi0}
\index{A}{1Phi@$\Phi_m^+(t)$}%
\Phi_m^+(t)&:= (2t)_+^{-\frac{m}2} J_m(2\sqrt{2t_+}),\\
\label{def:Psi+}
\index{A}{1Psi@$\Psi_m^+(t)$}%
\Psi_m^+(t)&:=  (2t)_+^{-\frac{m}2} J_m(2\sqrt{2t_+})
   -\sum_{k=1}^{m} \frac{(-1)^{k-1}}{2^k(m-k)!}\, \delta^{(k-1)}(t),  \\
\label{def:Psi}
\index{A}{1Psi@$\Psi_m(t)$}%
\Psi_m(t)&:=   (2t)_+^{-\frac{m}2} Y_m(2\sqrt{2t_+})
       +\frac{2(-1)^{m+1}}\pi (2t)_-^{-\frac{m}2} K_m(2\sqrt{2t_-}).
\end{alignat}
Here,
$(2t_+)^{-\frac{m}{2}} J_m (2\sqrt{2t_+})$
makes sense as a locally integrable function on $\mathbb{R}$.
On the other hand,
$(2t)_+^{-\frac{m}{2}} Y_m(2\sqrt{2t_+})$
and
$(2t)_-^{-\frac{m}{2}} K_m(2\sqrt{2t_-})$
are defined as regularized distributions,
and $\Psi_m(t)$ have the following singularity:
$$
\Psi_m(t) = (\text{locally integrable function})
   + \frac{-1}{\pi}
   \sum_{k=1}^m \frac{t^{-k}(k-1)!}{2^k(m-k)!}.
$$
We shall say that $\Phi_m^+$, $\Psi_m^+$, and $\Psi_m$ are
\textit{Bessel distributions}. 

\medskip
\noindent
\textbf{Integral expressions and differential equations}
\nopagebreak

The Bessel distributions
$\Phi_m^+(t)$, $\Psi_m^+(t)$, and $\Psi_m^+(t)$ are real analytic on 
$\mathbb{R} \setminus \{0\}$,
and satisfy the following differential equation:
\begin{equation*}
   t \frac{d^2\Psi}{dt^2} + (m+1)\frac{d\Psi}{dt} + 2\Psi = 0,
\end{equation*}
or equivalently,
\begin{equation*}
   (\theta^2+m\theta+2t)\Psi = 0,
\end{equation*}
where $\theta := t\frac{d}{dt}$.
Furthermore, all of the three solutions satisfy the following 
\index{B}{asymptotic behavior!Bessel distribution@---, Bessel distribution}
asymptotic behavior
\begin{equation*}
   \Psi(t) = O(t^{-\frac{2m+1}{4}})
   \quad \text{as $t \to +\infty$}.
\end{equation*}
 
In this book, we adopt 
an alternative definition of $\Phi_m^+$, $\Psi_m^+$ and $\Psi_m$ in 
Section \ref{subsec:intPsi} by means of the Mellin--Barnes type integral for
 distributions.
The expressions \eqref{def:Psi0}--\eqref{def:Psi} will be explained there.
Another (slightly different) expression of $\Phi_m^+$, $\Psi_m^+$, and
$\Psi_m$ by means of `normalized' Bessel functions $\widetilde{J}_m$,
$\widetilde{K}_m$, and $\widetilde{Y}_m$ is given in Remark \ref{rem:PhiPsi}.

Chapter \ref{sec:Diffeq} is devoted to these
Bessel distributions.
We shall discuss their integral formulas and differential equations.

\medskip
\noindent
\textbf{Support of the kernel $K(x,\xi)$}
\nopagebreak

Unlike the kernel 
$k(x,\xi) = (2\pi)^{-\frac{n}{2}} e^{\sqrt{-1}\langle x,\xi\rangle}$
for the Euclidean Fourier transform $\mathcal{F}_{\mathbb{R}^n}$,
our formula in Theorem \ref{thm:Kformula} shows that 
the support of the kernel $K(x,\xi)$ for $\mathcal{F}_C$ 
differs according to the signature
($n_1, n_2$). 

To see this, we set the `half' space of the direct product manifold 
$C \times C$ by
\begin{equation*}
   (C\times C)_+ := \{ (x,\xi) \in C \times C:
                       \langle x,\xi\rangle \ge 0 \}.
\end{equation*}
Then, we have the following mysterious
phenomenon: 
\begin{corollary} \label{cor:suppC}
The kernel 
$K(x,\xi)$ of the unitary operator $\mathcal{F}_C$
on $L^2(C)$ satisfies
\begin{equation*}
   \operatorname{supp} K(x,\xi)
   = \begin{cases}
        (C \times C)_+
        & \text{if $n_1,n_2$ both odd},
     \\
        C \times C
        & \text{if $n_1,n_2$ both even}.
     \end{cases}
\end{equation*}
\end{corollary}
In particular,
\begin{equation*}
   \operatorname{supp} K(x,\xi) \subsetneqq C \times C
   \quad\text{if $n_1,n_2$ both odd}.
\end{equation*}

If $n_1=1$ or $n_2=1$,
the 
isotropic cone%
\index{B}{isotropic cone} 
 $C$ is the 
light cone%
\index{B}{light cone} 
 in the relativistic cosmology, 
 which splits into two connected components,
namely, the 
forward light cone%
\index{B}{forward light cone} 
 $C_+$ and 
backward light cone%
\index{B}{backward light cone} 
 $C_-$,
and correspondingly,
we have a direct sum decomposition as Hilbert spaces:
\begin{equation}\label{eqn:Ltwosplit}
   L^2(C) = L^2(C_+) \oplus L^2(C_-).
\end{equation}
Then, $\mathcal{F}_C$ leaves $L^2(C_+)$ and $L^2(C_-)$ invariant,
respectively (see \cite{xkmano1,xkmano2}).
This gives a geometric explanation of Corollary \ref{cor:suppC}
because
\begin{equation}\label{eqn:CCdeco}
   (C \times C)_+ = (C_+ \times C_+) \cup (C_-\times C_-)
\end{equation}
in this case.
On the other hand, if $n_1,n_2 > 1$,
then $C$ is connected and we do not have a natural decomposition of
$(C \times C)_+$ like \eqref{eqn:CCdeco}.
Moreover, the representation $\pi$ 
of the indefinite orthogonal group
$O(n_1+1,n_2+1)$ on $L^2(C)$ (discussed later) 
stays irreducible
when restricted to the identity component $SO_0(n_1+1,n_2+1)$.
Nevertheless, the support of $K(x,\xi)$ is half the space of 
 $C \times C$ when both $n_1$ and $n_2$ are
 odd integers greater than one.
We do not see its intrinsic reason. 

\medskip
\noindent
\textbf{Singularities of the kernel $K(x,\xi)$}
\nopagebreak

Another distinguishing feature of our kernel $K(x,\xi)$ for
$\mathcal{F}_C$ is that it is not real analytic on $C \times C$,
whereas the kernel  
$k(x,\xi) = (2\pi)^{-\frac{n}{2}} e^{\sqrt{-1}\langle x,\xi\rangle}$ 
for the Euclidean Fourier transform 
$\mathcal{F}_{\mathbb{R}^n}$ is obviously real analytic on
$\mathbb{R}^n \times \mathbb{R}^n$.

Among the three Bessel distributions introduced in
\eqref{def:Psi0}--\eqref{def:Psi}, 
$\Phi_m^+(t)$ is a locally integrable function on $\mathbb{R}$,
whereas $\Psi_m^+(t)$ and $\Psi_m(t)$ $(m\ge1)$ are not.
The singular part of the distribution $\Psi_m^+(t)$ is given as a
linear combination of the Dirac delta function $\delta(t)$ and its
$l$th derivative $\delta^{(l)}(t)$ $(l=1,2,\dots,m-1)$.
The singular part of the 
distribution $\Psi_m(t)$ is given by a linear combination of the
distribution $t^{-k}$ $(k=1,2,\dots,m)$
(see Theorem \ref{prop:Psiint}).

Some readers might wonder why the kernel function of a unitary
operator involves such singularities.
So, let us examine to which extent the regularity of the kernel $K(x,\xi)$
is required from 
 the general theory of functional analysis.

By the Schwartz kernel theorem,
any continuous operator
$T: L^2(C) \to L^2(C)$
is expressed as
\begin{equation*}
   (Tf)(\xi) = \int_C K(x,\xi) f(x) d\mu(x)
\end{equation*}
by some distribution kernel
$K(x,\xi) \in \mathcal{D}'(C \times C)$.
Here, we regard $K(x,\xi)$ as a generalized function by using the
measure $d\mu$ on $C$
(see \cite{xGeGrVi}).

If $T$ is a Hilbert--Schmidt operator,
then $K\in L^2(C \times C)$.
If $T$ is the identity operator,
then $K$ is Dirac's delta function $\delta(x-\xi)$.
In general, the continuity of $T$ 
forces 
 any such $K$ to be 
at most $(\dim C+2)$ times derivatives 
 of a locally integrable function on 
$C \times C$ 
(see \cite[pp.~296--299]{xTaylor} for the argument 
 using the Sobolev space theory 
in the compact
torus case).

\medskip
\noindent
\textbf{`Laurent series expansions' of Bessel distributions}
\nopagebreak

We end this section with an interesting observation on 
\index{B}{Bessel distribution!Laurent series expansion@---, Laurent series expansion}%
`Laurent series expansions' 
of Bessel distributions $\Phi_m^+(t)$ and $\Psi_m^+(t)$:
\begin{align}
   & 
\index{A}{1Phi@$\Phi_m^+(t)$}%
\Phi_m^+(t) = \sum_{j=0}^\infty
                   \frac{(-1)^j 2^j t_+^j}{\Gamma(m+j+1)\Gamma(j+1)},
\label{eqn:Psipower}
\\
   & 
\index{A}{1Psi@$\Psi_m^+(t)$}%
\Psi_m^+(t) = (\sum_{j=0}^\infty+\sum_{j=-\infty}^{-1})
                   \frac{(-1)^j 2^j t_+^j}{\Gamma(m+j+1)\Gamma(j+1)}.
\label{eqn:Psi+power}
\end{align}
We note that $\Phi_m^+$ arises as the kernel for $\mathcal{F}_C$ when
$(n_1,n_2)=(2m+3,1)$,
and that $\Psi_m^+$ arises when
$n_1,n_2,\ge3$
are both odd and 
$n_1+n_2=2m+4$.
The first formula \eqref{eqn:Psipower}
 is a usual Taylor expansion.
But the second formula \eqref{eqn:Psi+power} involves negative terms,
for which we need a justification.
For this,
we think of 
\index{A}{xlambda@$x_+^\lambda$}%
$t_+^\lambda$ 
as a distribution meromorphically dependent on $\lambda$, 
and then we get
\begin{equation*}
   \frac{t_+^\lambda}{\Gamma(1+\lambda)} \,\Big|_{\lambda=-k}
   = \delta^{(k-1)}(t)
\end{equation*}
by \eqref{eqn:resxlmd}.
Therefore, we have 
\begin{equation*}
   \frac{(-1)^\lambda 2^\lambda t_+^\lambda}
        {\Gamma(m+\lambda+1)\Gamma(1+\lambda)} \,\Big|_{\lambda=-k}
   = \begin{cases}
        \dfrac{(-1)^k}{2^k(m-k)!} \delta^{(k-1)}(t)
        & (1 \le k \le m),
     \\[2ex]
        0 & (m+1 \le k).
     \end{cases}
\end{equation*}

In this sense, 
the series \eqref{eqn:Psi+power} contains only finitely many negative
terms,
and is equal to \eqref{def:Psi+}.

\section{Perspectives from representation theory --
finding smallest objects}
\label{subsec:1.4}

The philosophy of analysis and synthesis asks for the understanding of
the smallest objects and of how things are built from these objects.

For a Hilbert space $\mathcal{H}$ over $\mathbb{C}$,
we denote by $U(\mathcal{H})$ the group consisting of unitary
operators on $\mathcal{H}$.
By a unitary representation on a Hilbert space $\mathcal{H}$ of a
topological group $G$,
we mean a group homomorphism
\begin{equation*}
   \pi: G \to U(\mathcal{H})
\end{equation*}
such that $G \times \mathcal{H} \to \mathcal{H}$,
$(g,v) \mapsto \pi(g)v$ is continuous.

The `smallest objects' of unitary representations are irreducible
unitary representations.
By a theorem of Mautner and Teleman,
any unitary representation of a locally compact group $G$
(e.g.\ a Lie group) can be decomposed into the direct integral of
irreducible unitary representations of $G$, see \cite{xwal}.
The classification of irreducible unitary representations of Lie
groups has been a long standing unsolved problem since 1940s,
originally arising from quantum mechanics.

The `smallest objects' of Lie groups consist of simple Lie groups
such as $SL(n,\mathbb{R})$, $O(p,q)$, and $Sp(n,\mathbb{R})$,
and one-dimensional abelian Lie groups such as $\mathbb{R}$ and $S^1$. 
Loosely speaking,
a theorem of Duflo \cite{xduflo} asserts that all irreducible unitary
representations of general (real algebraic) Lie groups
are built up from those of simple Lie groups.

However, irreducible unitary representations of simple Lie groups
are not fully understood despite huge efforts for many decades and
also significant results.
Among them, 
powerful algebraic machinery including the theory of 
cohomological induction%
\index{B}{cohomological induction} 
 has been largely developed in 1980s by Zuckerman, Vogan,
Wallach and others \cite{xvuni,xwal}.
As a result,
the problem of classifying irreducible unitary representations
has been focused on those representations
of simple Lie groups that cannot be `induced up' from other
representations.
Such representations may be regarded as `atoms' of unitary
representations of Lie groups,
and they are still mysterious creatures.
See Vogan \cite{xvuni}, for example,
for a discussion on how to understand them as the theory of
`unipotent representations' of reductive Lie groups.

`Minimal representations' are the simplest,
infinite dimensional `unipotent representations'.
There has been active study on minimal representations,
mostly by algebraic methods since 1990s
(see Gan and Savin \cite{xgansavin} and J.-S. Li \cite{xLi} for
surveys both in the real and in the $p$-adic fields).
In contrast to these existing algebraic approaches,
new geometric analysis of the minimal representation of the simple Lie
group $O(p,q)$ is a motif of this book.
We will discuss minimal representations in the next section 
followed by Sections \ref{subsec:1.6a}--\ref{subsec:uncertainty}
 in more details.

\section{Minimal representations of simple Lie groups}
\label{subsec:1.4b}

To formulate the `smallness' of an irreducible representation $\pi$ of
a simple Lie group $G$,
algebraic representation theory usually appeals to the `largeness' of
the annihilator 
\index{A}{Annpi@$\operatorname{Ann}(\pi)$}%
$\operatorname{Ann}(\pi)$ 
in $U(\mathfrak{g}_{\mathbb{C}})$ 
of the differential
representation $d\pi$.
Here, $U(\mathfrak{g}_{\mathbb{C}})$ is the universal enveloping algebra
of the complexified Lie algebra
$\mathfrak{g}_{\mathbb{C}} = \mathfrak{g} \otimes_{\mathbb{R}} \mathbb{C}$.
We recall:
\begin{definition}\label{def:minrep}
An irreducible unitary representation $\pi$ of a simple Lie group $G$ is a
\textit{minimal representation}%
\index{B}{minimal representation|main} 
 if the annihilator
$\operatorname{Ann}(\pi)$ is
equal to the Joseph ideal \cite{xgansavinJ,xjoseph} of 
$U(\mathfrak{g}_{\mathbb{C}})$. 
\end{definition}

The Joseph ideal $\operatorname{Ann}(\pi)$ is a completely prime
ideal whose 
associated variety%
\index{B}{associated variety} 
 $\mathcal{V}(\operatorname{Ann}\pi)$ 
is the closure of
$\mathcal{O}_{\min}^{\mathbb{C}}$.
Here, 
$\mathcal{O}_{\min}^{\mathbb{C}}$%
\index{A}{OminC@$\mathcal{O}_{\min}^{\mathbb{C}}$} 
 is the complex minimal nilpotent
orbit in $\mathfrak{g}_{\mathbb{C}}^*$ for the coadjoint
representation. 
Therefore, 
the 
Gelfand--Kirillov dimension%
\index{B}{Gelfand--Kirillov dimension} 
 of $\pi$, to be denoted by 
$\operatorname{DIM}\pi$%
\index{A}{DIMpi@$\operatorname{DIM}\pi$}%
, satisfies
\begin{equation*}
   \operatorname{DIM}\pi = \frac{1}{2} \dim_{\mathbb{C}}
   \mathcal{O}_{\min}^{\mathbb{C}},
\end{equation*}
and in particular, 
$\pi$ has the smallest possible Gelfand--Kirillov dimension.

The Weil representation $\varpi$ of the
metaplectic group $\Mp(n,\mathbb{R})$
is a classic example of minimal representations 
(to be precise, the Weil representation decomposes into a direct sum of
two irreducible representations, both of which are minimal
representations),
and $\operatorname{DIM}\varpi = n$. 

The study of minimal representations of other reductive groups is relatively
new, and it is only in the last two decades that various models of minimal
representations have been proposed and studied extensively
(see  \cite{xBiZi,xbk,xDvSa,xfks,xgansavin,xGrWal,xKa,xkazsav,xkmano2,xkors1,xkors2,xkors3,
xKo,xLi,xSa,xTo,xHuZu} for instance).

We note that not every reductive Lie group admits minimal
representations. 
For instance,
the minimal nilpotent orbit $\mathcal{O}_{\min}^{\mathbb{C}}$ 
of $SO(m,\mathbb{C})$ has the
dimension $2(m-3)$,
and it was found by Howe and Vogan (see \cite{xV})
 that there exists no
representation of any covering group of $SO_0(p,q)$
(or its finite extension) 
whose Gelfand--Kirillov dimension equals
$p+q-3$ 
if $p+q$ is odd and $p,q>3$. 
This is the primary reason that we dealt with the case
$n=n_1+n_2$ is even in Sections \ref{subsec:1.1}--\ref{subsec:FC}.
Here, $p=n_1+1$ and $q=n_2+1$.
(It is another story that there exists an `infinitesimally unitary'
and `minimal' representation of the
Lie algebra $\mathfrak{so}(p,q)$ for $p+q$ odd, see \cite{xkors3}.)

In accordance with the philosophy of the orbit method advocated by
A. Kirillov, B. Kostant, and M. Duflo,
irreducible unitary representations of Lie groups $G$ 
are supposed to be attached to coadjoint orbits in $\mathfrak{g}^*$.
This works perfectly for simply-connected nilpotent Lie groups.
However, for simple Lie groups,
the orbit method does not work very well,
though the orbit method still gives
 an approximation of the unitary dual.
In particular,
it is not known how to attach unitary representations to
nilpotent orbits of simple Lie groups.

The reverse direction is easier.
As we have seen at the beginning of this section,
we can attach nilpotent orbits to admissible representations $\pi$
(in particular, to irreducible unitary representations):
\begin{equation}\label{eqn:Asspi}
   \pi \rightsquigarrow \operatorname{Ann}(\pi) \rightsquigarrow
   \mathcal{V}(\operatorname{Ann}(\pi)).
\end{equation}
There are also several ways to attach (a union of) real nilpotent
orbits in $\mathfrak{g}^*$ to admissible representations $\pi$ such as
\begin{alignat*}{2}
    &\mathfrak{g}^* \supset \operatorname{WF}_{\operatorname{N}}(\pi): {}
    &&\text{wave front set%
\index{B}{wave front set} 
 by Howe, \cite{xhoweWF}}
   \\
    &\mathfrak{g}^* \supset \operatorname{WF}(\pi):
    &&\text{asymptotic support%
\index{B}{asymptotic support} 
 of the character Trace $\pi$,
     \cite{xbavo}}
   \\
    &\mathfrak{p}_{\mathbb{C}}^* \supset \operatorname{Ass}(\pi):
    &&\text{associated variety%
\index{B}{associated variety} 
 of the underlying 
     $(\mathfrak{g}_{\mathbb{C}},K)$-module of $\pi$, \cite{xvoass}}
\end{alignat*}
Here, $\mathfrak{g} = \mathfrak{k}+\mathfrak{p}$ is a Cartan
decomposition
and
$\mathfrak{g}_{\mathbb{C}}=\mathfrak{k}_{\mathbb{C}}+\mathfrak{p}_{\mathbb{C}}$
is its complexification. 
$\operatorname{WF}_{\operatorname{N}}(\pi)$ and 
$\operatorname{WF}(\pi)$ are a union of nilpotent orbits of $G$ in
$\mathfrak{g}^*$, 
and $\operatorname{Ass}(\pi)$ is a union of nilpotent orbits of $K_{\mathbb{C}}$ on
$\mathfrak{p}_{\mathbb{C}}^*$.
It was conjectured by Barbash and Vogan \cite{xbavo},
and proved by Schmid and Vilonen \cite{xschvil} that
$\operatorname{WF}(\pi)$ corresponds to 
$\operatorname{Ass}(\pi)$ via the 
\index{B}{Kostant--Sekiguchi correspondence}%
Kostant--Sekiguchi correspondence. 
For a minimal representation $\pi$,
these invariants are contained in the closure of 
the intersection $\mathfrak{g}^* \cap
\mathcal{O}_{\min}^{\mathbb{C}}$,
$\mathfrak{p}_{\mathbb{C}}^* \cap \mathcal{O}_{\min}^{\mathbb{C}}$,
respectively.
Here, 
$\mathcal{O}_{\min}^{\mathbb{C}}$ is 
the complex minimal
nilpotent orbit in
$\mathfrak{g}_{\mathbb{C}}^*$, 
and we identify $\mathfrak{g}_{\mathbb{C}}$ with
$\mathfrak{g}_{\mathbb{C}}^*$ by the Killing form.

Conversely, it is much more difficult to construct irreducible unitary
representations from nilpotent orbits in general.
If the complex minimal nilpotent orbit
$\mathcal{O}_{\min}^{\mathbb{C}}$ 
has a non-empty intersection with $\mathfrak{g}^*$,
then $\mathcal{O}_{\min}^{\mathbb{C}} \cap \mathfrak{g}^*$ consists of
equi-dimensional $\operatorname{Ad}^*(G)$-orbits, namely,
minimal nilpotent orbits in $\mathfrak{g}^*$.
An optimistic picture is that minimal representations are attached to
such orbits,
however, this is false in general
(see \cite{xtorassosl} for the $\mathfrak{sl}(3,\mathbb{R})$ case).
In this direction Brylinski and Kostant \cite{xbk} constructed
minimal representations from $\operatorname{Ass}(\pi)$ on the
$(\mathfrak{g}_{\mathbb{C}},K)$-module level,
but both the Hilbert structure and the whole group action are not
given globally,
 depending on each $K$-type.

In this book,
we deal with 
$G = O(n_1+1,n_2+1)$
($n_1,n_2 > 1, n_1+n_2$ even $> 4$), 
where 
$\mathcal{O}_{\min}%
\index{A}{Omin@$\mathcal{O}_{\min}$} 
 := \mathcal{O}_{\min}^{\mathbb{C}} \cap
\mathfrak{g}^*$ 
is a single $G$-orbit, namely, 
a minimal nilpotent orbit. 
Our minimal representation is realized on the concrete Hilbert space $L^2(C)$.
On the other hand,
 the isotropic cone $C$ becomes a Lagrangian submanifold of
$\mathcal{O}_{\min}$ which is endowed with the natural symplectic structure. 
In this sense, 
our Schr\"{o}dinger model $L^2(C)$ 
may be viewed as a 
\textit{geometric quantization}%
\index{B}{geometric quantization} 
 of the real 
minimal coadjoint orbit $\mathcal{O}_{\min}$.
Our main theorem enables us to give the
 whole group $G$-action on $L^2(C)$ explicitly and 
independently of $K$-types.

\section{Schr\"{o}dinger model for the Weil representation}
\protect\index{B}{Schr\"{o}dinger model}%
\label{subsec:1.6a}

In Section \ref{subsec:1.6},
we shall discuss the minimal representation of $G=O(n_1+1,n_2+1)$,
which is a simple group of type $D$ if $n_1+n_2$ is even, $>2$.
In this section,
we recall the best understood minimal representation of a
simple group of type $C$,
that is, the 
Segal--Shale--Weil representation%
\index{B}{Segal--Shale--Weil representation} 
 $\varpi$,
simply, the 
Weil representation%
\index{B}{Weil representation}%
, 
or sometimes referred to as the 
oscillator representation%
\index{B}{oscillator representation}%
, 
or 
harmonic representation%
\index{B}{harmonic representation}%
, 
of the 
\index{B}{metaplectic group}%
metaplectic group 
\index{A}{MpnR@$\Mp(n, \mathbb R)$}%
$\Mp(n, \mathbb R)$.
Here, $Mp(n,\mathbb{R})$ is
the twofold cover of the real 
\index{B}{symplectic group}%
symplectic group 
\index{A}{SpnR@$Sp(n, \mathbb R)$}%
$Sp(n, \mathbb R)$.
Let $\xi_0$ denote the (unique) non-trivial element in the kernel of
the homomorphism $Mp(n,\mathbb{R}) \to Sp(n,\mathbb{R})$.
That is, we have an exact sequence of Lie groups:
\begin{equation*}
   1 \to \{e,\xi_0\} \to Mp(n,\mathbb{R}) \to
   Sp(n,\mathbb{R}) \to 1.
\end{equation*}
We will set $G' := Mp(n,\mathbb{R})$.

Among various realizations of the Weil representation
(see \cite[\S3]{xLi} for a brief survey),
the Schr\"{o}dinger model gives
 a realization of the Weil representation $\varpi$ on the Hilbert space
$L^2(\mathbb{R}^n)$.
Since our model $(\pi,L^2(C))$ of the minimal representation of 
$G$ has a strong resemblance to $(\varpi,L^2(\mathbb{R}^n))$ of 
$G'$,
we list
 some important aspects of the Schr\"odinger model of $\varpi$
(see e.g.\ \cite{xfolland,xGel,xHo,xKaVe}):

\begin{enumerate}[{\bf {C}1}]
\index{A}{1zpropertiesC1-C4@\textbf{C1}--\textbf{C4}|(}%
\item
\label{C1}
The representation is realized on a very concrete
Hilbert space, that is, $L^2(\mathbb R^n)$.

\item
\label{C2}
The restriction of $\varpi$ to the 
\index{B}{Siegel parabolic subgroup}%
Siegel parabolic subgroup 
$P_\text{Siegel}$ is still irreducible. 
 The restriction $\varpi_{P_\text{Siegel}}$ has a relatively simple form
 (translations and multiplications by unitary characters).

\item
\label{C3}
The infinitesimal action $d \varpi$ of the Lie algebra 
 $\mathfrak{sp}(n, \mathbb R)$ is given by
   differential operators of at most second order. 
     
\item \label{C4}
There is a distinguished element $w'_0$ 
 of $G'$ that sends $P_\text{Siegel}$ to the opposite parabolic
subgroup. 
The corresponding unitary operator $\varpi(w'_0)$ on $L^2(\mathbb R^n)$
 is equal to $e^{\frac{\sqrt{-1}n\pi}{4}}
\index{A}{FRn@$\mathcal{F}_{\mathbb{R}^n}$}%
\mathcal F_{\mathbb{R}^n}$.
Correspondingly to the fact that $(w'_0)^4=\xi_0$ and
 $\varpi(\xi_0)=-\operatorname{id}$, 
the Fourier transform $\mathcal{F}_{\mathbb{R}^n}$ is of order four.
\index{A}{1zpropertiesC1-C4@\textbf{C1}--\textbf{C4}|)}%
\end{enumerate}

We write $Sp(n,\mathbb{R})$ in the matrix form as 
\begin{equation*}
  Sp(n,\mathbb{R})
  = \{ g \in GL(2n,\mathbb{R}):
       {}^t \! g J_n g = J_n \},
\end{equation*}
where
$J_n = \begin{pmatrix}0 &-I_n \\ I_n &0 \end{pmatrix}$.
Then, down to
$
Sp(n,\mathbb{R}) \simeq Mp(n,\mathbb{R})
\mod \{ e,\xi_0 \}
$,
we may take $P_\text{Siegel}$ and $w'_0$ as
\begin{align*}
\index{A}{Psiegel@$P_\text{Siegel}$}%
   P_\text{Siegel} \mod \{ e,\xi_0 \}
   & \simeq
     \Bigl\{ \begin{pmatrix} A &B \\ 0 &{}^t\! A^{-1} \end{pmatrix} :
             A \in GL(n,\mathbb{R}), \, 
             A \, {}^t\! B = B \, {}^t\! A \Bigr\},
\\
\index{A}{w0'@$w'_0$}%
   w'_0 \mod \{ e,\xi_0 \}
   & = J_n.
\end{align*}

Since $G'$ is generated by $P_\text{Siegel}$ and $w'_0$,
\textbf{C\ref{C2}} and \textbf{C\ref{C4}} determine the action of $G'$ on 
$L^2(\mathbb R^n)$
(see \cite{xPe} for an explicit formula of the whole group $G'$-action
on $L^{2}(\mathbb{R}^{n})$). 
As for \textbf{C\ref{C3}},
if $X \notin \mathfrak{p}_\text{Siegel}$
(the Lie algebra of $P_\text{Siegel}$), 
then $d \varpi(X)$ contains a differential operator of second order,
and is not given by a
vector field. This reflects the fact that $G'$ acts only on
$L^2(\mathbb{R}^n)$, and not on
$\mathbb{R}^n$.

In various places of this book (usually, as Remarks), 
we will compare 
our results with the corresponding results for the Weil representation.
We will see that the case $n_1=1$ or $n_2=1$ 
 in our setting $G = O(n_1+1,n_2+1)$ has similar features to the case of
 the Weil representation in both analytic and representation theoretic aspects,
 and that the general case $n_1, n_2>1$ often
  provides new analytic features.

\section{Schr\"{o}dinger model for the minimal representation of
$O(p,q)$}
\protect\index{B}{Schr\"{o}dinger model}%
\label{subsec:1.6}

In this section,
we consider 
the 
\index{B}{indefinite orthogonal group|main}%
indefinite orthogonal group
\begin{equation} \label{eqn:Opq}
  G:= 
\index{A}{O(p,q)@$O(p,q)$|main}%
O(p,q)= \set{ g \in GL(p+q, \mathbb R)}{ {}^t \! g I_{p,q} g=I_{p,q}},
\end{equation}
where 
$
I_{p,q}:= \begin{pmatrix}  I_p  &  0  \\  0    & -I_q      \end{pmatrix}
$,
\thinspace
$p,q \ge 2$ and $p+q$ is an even integer greater than four. 
Then, 
there exists a distinguished irreducible unitary representation $\pi$
of $G$ with the following properties:
\begin{enumerate}[\bfseries{M}1]
\index{A}{1zpropertiesM1-M5@\textbf{M1}--\textbf{M5}|(}%
\item  \label{item:M1}
$\pi$ is still irreducible when restricted to the identity component
$SO_0(p,q)$ of $G$ if and only if $p,q > 2$.
\item  \label{item:M2}
If $p=2$ or $q=2$, then $\pi$ is a direct sum of an irreducible
\index{B}{highest weight module}%
unitary highest weight representation 
and a lowest weight
representation.
\item  \label{item:M3}
$\pi$ is 
\index{B}{spherical representation|main}%
spherical 
(i.e.\ has a non-trivial vector fixed by 
$K = O(p) \times O(q)$)
if and only if $p = q$.
\item   \label{item:M4}
$\pi$ is a 
\index{B}{minimal representation}%
minimal representation 
in the sense of Definition
\ref{def:minrep} if $p+q>6$ (see \cite{xBiZi,xKo}).
\item   \label{item:M5}
$\operatorname{WF}(\pi) = \mathcal{O}_{\min} =
 (\mathcal{O}_{\min}^{\mathbb{C}} \cap \mathfrak{g}^*)$,
$\operatorname{Ass}(\pi) = \mathcal{O}_{\min}^{\mathbb{C}} \cap
 \mathfrak{p}_{\mathbb{C}}^*$, 
and $\operatorname{DIM}(\pi) = p+q-3$,
(see \cite{xkors1}).
\index{A}{1zpropertiesM1-M5@\textbf{M1}--\textbf{M5}|)}%
\end{enumerate}
The above properties \textbf{M2} and \textbf{M3} show that our minimal
representation $\pi$ is neither a spherical representation nor a
highest weight representation in the generic case where $p \ne q$,
\thinspace 
$p,q > 2$, and $p+q > 6$.

Various realizations of the minimal representation $\pi$ have been
proposed so far by a number of people. 
For example, 
Kazhdan in \cite{xKa}, and Kostant in \cite{xKo} for $p=q=4$, 
and Binegar and Zierau \cite{xBiZi}
for general $p, q \ge 2$, 
constructed $\pi$ 
as a subrepresentation of a maximally degenerate principal series
representation 
(see also Howe and Tan \cite{xHoTan} for a full discussion on its
composition series); 
Zhu and Huang \cite{xHuZu} constructed $\pi$ 
as the 
\index{B}{theta correspondence}%
theta correspondence of the trivial one-dimensional
representation of $SL(2,\mathbb{R})$ for the reductive dual pair
$O(p,q)\cdot SL(2,\mathbb{R}) \subset Sp(p+q,\mathbb{R})$, 
see also \cite{xRalSch};
and Kobayashi and \O rsted
\cite{xkors1} constructed $\pi$ as the solution to the Yamabe equation
(\textit{conformal model});
in \cite{xkors3} the \textit{Schr\"{o}dinger model} (an $L^2$-model)
of $\pi$.
Yet another construction has been proposed in Brylinski and Kostant
\cite{xbk}, Gross and Wallach \cite{xGrWal}, and Torasso \cite{xTo}.

Among various realizations of the minimal representations,
our concern is with the Schr\"{o}dinger model realized on $L^2(C)$. 
Here, $C$ is the isotropic cone in 
$\mathbb{R}^{n_1+n_2} = \mathbb{R}^n$ 
discussed in Sections \ref{subsec:1.1}--\ref{subsec:FC},
and the relation between the above parameters $p$ and $q$ is given by
\begin{equation*}
   p=n_1+1, \quad q=n_2+1, \quad n=p+q-2.
\end{equation*}

The conformal linear transformation group $CO(Q)$
$(\subset GL(n_1+n_2,\mathbb{R}))$
acts on the isotropic cone $C$ $(\subset \mathbb{R}^{n_1+n_2})$,
and then induces a unitary representation on $L^2(C)$
(see \eqref{eqn:rM}--\eqref{eqn:rA} for a concrete formula).
Much more than that,
this action on $L^2(C)$ can be extended to a unitary representation $\pi$ of
the indefinite orthogonal group $G=O(n_1+1,n_2+1)$. 
To explain its idea,
we also recall another geometric model, namely,
the conformal model.

The conformal model of the minimal representation $\pi$ is realized, for example,
in the solution space to the 
ultra-hyperbolic equation%
\index{B}{ultra-hyperbolic equation}
\begin{equation*}
   \Bigl( \frac{\partial^2}{\partial x_1^2} + \dots +
          \frac{\partial^2}{\partial x_{n_1}^2} -
          \frac{\partial^2}{\partial x_{n_1+1}^2} - \dots -
          \frac{\partial^2}{\partial x_{n}^2} \Bigr) u = 0
\end{equation*}
in $\mathcal{S}'(\mathbb{R}^n)$ 
\index{B}{tempered distribution}%
(tempered distributions),
and also in the solution space
\begin{equation*}
   \widetilde{\Delta}_{S^{n_1}\times S^{n_2}} \, v = 0
   \quad \text{on $S^{n_1}\times S^{n_2}$}
\end{equation*}
where 
$\widetilde{\Delta}_{S^{n_1}\times S^{n_2}} 
 = \Delta_{S^{n_1}} - \Delta_{S^{n_2}} - (\frac{n_1-1}{2})^2
    + (\frac{n_2-1}{2})^2$
is the 
Yamabe operator%
\index{B}{Yamabe operator} 
 on 
the direct product manifold 
$S^{n_1} \times S^{n_2}$ equipped with the
pseudo-Riemannian structure of signature $(n_1,n_2)$.
These two models are isomorphic to each other by the general theory of
conformal geometry,
and the intertwining operator is given by the `twisted pull-back' 
$\widetilde{\Psi}^*$ of 
the conformal map
$\Psi: \mathbb{R}^{n_1+n_2} \to S^{n_1} \times S^{n_2}$,
the inverse of the stereographic projection 
(see \cite{xkcheck} for an elementary account).
Then, taking the Fourier transform $\mathcal{F}_{\mathbb{R}^n}$
of the conformal model, 
we get the Schr\"{o}dinger model $L^2(C)$.

The intertwining operator
$\mathcal{T}: L^2(C) \to
\mathcal{S}ol(\widetilde{\Delta}_{S^{n_1}\times S^{n_2}})$
is defined in Section \ref{subsec:Sch} so that the following diagram
commutes: 
\begin{figure}[H]
\begin{equation*}
\begin{matrix}
  & & \index{B}{conformal model}%
\text{Conformal model}
  \\ 
  && 
\index{B}{K-picture@$K$-picture}%
\text{($K$-picture)}
  \\[1ex]
  && \mathcal{S}ol(\widetilde{\Delta}_{S^{n_1}\times S^{n_2}})
\index{A}{Soldelta@$\mathcal{S}ol(\widetilde{\Delta}_{S^{n_1}\times S^{n_2}})$}
  \\
  & \hidewidth\rlap{\raisebox{1ex}
    {\kern.3em$\scriptstyle \widetilde{\Psi}^*$}}
    \rotatebox[origin=c]{-135}{$\xrightarrow{\kern2em}$} \hidewidth
  && \hidewidth\rotatebox[origin=c]{-45}{$\xleftarrow{\kern2em}$} 
     \rlap{\raisebox{1ex}
     {$\kern-.8em\scriptstyle \mathcal{T}$}}\hidewidth 
  \\
  \mathcal{S}ol(\square_{\mathbb{R}^{n_1,n_2}}) 
\index{A}{Solsquare@$\mathcal{S}ol(\square_{\mathbb{R}^{n_1,n_2}})$}
  && \xrightarrow[\mathcal{F}_{\mathbb{R}^n}]{\kern4em}
  && L^2(C)
  \\[1ex]
  \hidewidth\text{Conformal model}\hidewidth
  &&&& \hidewidth\text{Schr\"{o}dinger model}\hidewidth
\index{B}{Schr\"{o}dinger model}
  \\
\index{B}{N-picture@$N$-picture}%
  \hidewidth\text{($N$-picture)}\hidewidth
\end{matrix}
\end{equation*}
\renewcommand{\figurename}{Diagram}
\caption{}
\label{diag:1.7.1}
\end{figure}
We remark that the isotropic cone 
$C$ is defined as a hypersurface in $\mathbb{R}^{n} = \mathbb{R}^{n_1+n_2}$,
and the group $G=O(n_1+1,n_2+1)$ cannot act (non-trivially) on $C$.
(In fact, any (non-trivial) $G$-space is of dimension at least 
$n_1+n_2=\dim C+1$.)

In the $L^2$-model $L^2(C)$ of the indefinite
 orthogonal group $G$, 
the action of a maximal parabolic subgroup
\begin{equation*}
   \overline{\Pmax}
\index{A}{Pmax=@$\overline \Pmax:= \Mmax A \overline \Nmax$}%
 \simeq \mathbb{Z}_2 \cdot CO(Q) \ltimes
   \mathbb{R}^{n_1+n_2}
\end{equation*}
on $L^2(C)$ is of
 a simple form 
 (see \eqref{eqn:rM}--\eqref{eqn:rN} for definition).
This parabolic subgroup $\overline{\Pmax}$ 
plays a similar role of the Siegel parabolic subgroup $P_\text{Siegel}$,
and analogous results to the properties
\textbf{C\ref{C1}}, \textbf{C\ref{C2}} and \textbf{C\ref{C3}} hold
(see Section \ref{subsec:act}). 
If we set
\begin{equation*}
\index{A}{w0@$w_0$}%
  w_0 := I_{n_1+1,n_2+1},
\end{equation*}
then $w_0$ sends $\overline{\Pmax}$ to the opposite parabolic
subgroup $\Pmax$,
and $G$ is generated by $w_0$ and $\overline{\Pmax}$.

In light of the Bruhat decomposition
$G=\overline{\Pmax} \amalg \overline{\Pmax}w_0
\overline{\Pmax}$,
we can get directly the concrete form of the action of the whole group
$G$ once we know
$\pi(w_0)$ explicitly.

In the degenerate case
$(n_1,n_2)=(2,0)$, $G = O(3,1)$
is locally isomorphic to $SL(2,\mathbb{C})$ acting on
$\mathbb{R}^{2+0}$ $(\simeq \mathbb{C})$
as 
\index{B}{Mobius transform@M\"{o}bius transform}%
M\"{o}bius transforms 
 (linear fractional transforms),
and $w_0$ acts on $\mathbb{C}$ as the 
conformal inversion%
\index{B}{conformal inversion}
\begin{equation*}
   \mathbb{C} \to \mathbb{C}, 
   \quad
   z \mapsto -\frac{1}{z}.
\end{equation*}

Although $(n_1,n_2)=(2,0)$ is beyond the parameter in our
consideration of the representation,
this feature of the conformal inversion $w_0$ is valid for any
$(n_1,n_2)$ 
(see Section \ref{subsec:w} for a list of key properties of this
 element $w_0$). 

In this book, we establish an analogous result
 to \textbf{C\ref{C4}} for $G=O(n_1+1,n_2+1)$,
that is, we find the unitary operator $\pi(w_0)$ 
on $L^2(C)$ for the conformal inversion $w_0$.

Then, here is our main result in this context.
\begin{main theorem}[see Theorem \ref{thm:A}]
\label{mainthm}
Let $n_1, n_2\ge 1$ and $n = n_1+n_2 \ge 4$ is even.
Then the unitary operator
$\pi(w_0):L^2(C) \to L^2(C)$ takes the form:
\begin{equation*}
   \pi(w_0) = \mathcal{F}_C,
\end{equation*}
where $\mathcal{F}_C$ is the involutive unitary operator given in
Section \ref{subsec:1.2}.
That is, 
\begin{equation}\label{eqn:mainthm}
\pi(w_0)u(\xzeta)=\int_C K(\xzeta, \xzeta')u(\xzeta')d\mu(\xzeta'),
\qquad  u \in L^2(C),
\end{equation}
where the distribution kernel $K(x,x')$ is given in Theorem
\ref{thm:Kformula}. 
\end{main theorem}

As obvious corollaries of the representation theoretic
 interpretation of $\mathcal F_C$ as above,
we have:
\begin{corollary}
\label{cor:FCinvunitary}
\index{B}{Plancherel formula!FC@---, $\mathcal {F}_C$}%
\upshape{(Plancherel and inversion formulas, see Corollaries \ref{cor:A1} and \ref{cor:A2})}
\begin{alignat*}{2}
 \| \mathcal F_C u \|_{L^2(C)}&=\| u\|_{L^2(C)} 
 \quad&&\text{for}~ u \in L^2(C),\\
 \mathcal F_C^{-1}&=\mathcal{F}_C \qquad &&\text{on}~L^2(C).
\end{alignat*}
\end{corollary}
The inversion formula $\mathcal{F}_C^{-1} = \mathcal{F}_C$ implies
that the following relation
\begin{align*}
   (\mathcal{F}_C u)(x)
   &=
   \int_C K(x,x') u(x') d\mu(x'), 
\\
   u(x)
   &=
   \int_C K(x,x') (\mathcal{F}_C u) (x') d\mu(x')
\end{align*}
\index{B}{reciprocal formula!FC@---, $\mathcal {F}_C$}%
is reciprocal.
Such an inversion formula
 is sometimes referred to
as a 
\textit{reciprocal formula}%
\index{B}{reciprocal formula} 
 (see Titchmarsh \cite{xTi}, for this terminology in a general setting).

It is noteworthy that Corollary \ref{cor:FCinvunitary} is stated
without any language of group theory. 
It would be an interesting problem to give a straightforward proof of
 Corollary \ref{cor:FCinvunitary} from the definition of Bessel
 distributions \eqref{def:Psi0}--\eqref{def:Psi} without group theory.

By \textbf{M\ref{item:M2}} and \textbf{M\ref{item:M4}}, 
$\pi$ is a non-highest weight,
 minimal representation if and only if 
$n_1,n_2 > 1$ and $n > 4$.
Therefore, we have discovered another mysterious phenomenon:
\begin{corollary} \label{cor:singK}
The kernel $K(x,x')$ of the unitary inversion $\pi(w_0)$ is a
not locally integrable function if and only if $\pi$ is a non-highest
weight,
minimal representation.
\end{corollary}

\medskip
\noindent
\textbf{The 
highest weight module%
\index{B}{highest weight module} 
 case}
($n_1 = 1$ or $n_2 = 1$)
\nopagebreak

In the case $n_2=1$ (likewise $n_1=1$), 
$\pi$ splits into the direct sum of a highest weight module $\pi_+$ 
and a lowest 
weight module $\pi_-$ 
when restricted to the identity component $G_0=SO_0(n_1+1,2)$ of 
$G = O(n_1+1,2)$
according to the decomposition \eqref{eqn:Ltwosplit}.
Both $\pi_+$ and $\pi_-$ are minimal representations of $G_0$.

We note that $G$ is 
the conformal group $O(n_1+1,2)$ of the 
Minkowski space%
\index{B}{Minkowski space} 
 $\mathbb{R}^{n_1,1}$, namely, the Euclidean space $\mathbb{R}^{n_1+1}$
equipped with the flat Lorentz metric of signature $(n_1,1)$.
In this case 
our representation $\pi$ has been studied also in physics.
The minimal representation $\pi_+$
may be interpreted as the symmetry of 
the solution space to the mass-zero spin-zero wave equation.
The representation $\pi_+$ arises also on the Hilbert space of
 bound states of the Hydrogen atom.

It is known that highest weight representations can be extended to holomorphic
semigroups of a complexified Lie group $G_{\mathbb{C}}$.
This theory has been initiated by Olshanski \cite{xOls} and
Stanton \cite{xStanton}, among others,
in connection with the 
Gelfand--Gindikin program%
\index{B}{Gelfand--Gindikin program} 
 to realize a family
of representations in a geometrically unified manner.
In this context the unitary operator $\pi(w_0)$ may be 
 regarded as the 
boundary value%
\index{B}{boundary value} 
 of a 
holomorphic semigroup%
\index{B}{holomorphic semigroup}%
. 
We then ask an explicit form of the holomorphic semigroup. 
This idea was first exploited by Howe \cite{xHo} for the Weil
representation, 
where he showed that $\varpi(w_0)$ 
$(= e^{\frac{\sqrt{-1}n\pi}{4}}\mathcal{F}_{\mathbb{R}^n})$ 
is given as the boundary value of the 
Hermite semigroup%
\index{B}{Hermite semigroup}%
, 
i.e.,
the holomorphic semigroup with the 
Mehler kernel%
\index{B}{Mehler kernel} 
 (see also Folland \cite{xfolland} for an exposition).
The same idea also works in our setting of 
$O(n_1+1,n_2+1)$ with $n_2=1$,
and the explicit formula of a certain holomorphic semigroup 
(the `Laguerre semigroup'%
\index{B}{Laguerre semigroup}%
) yields the formula of the
 unitary inversion operator  $\pi(w_0)$ by taking its 
 boundary value
 in \cite{xkmano1,xkmano2}.
In this book,
as a special case (i.e.\ $n_2=1$) of Theorem \ref{thm:A},
we give a new proof of the formula of $\pi(w_0)$.

\section{Uncertainty relation -- inner products and $G$-actions}
\label{subsec:uncertainty}

In this section,
we consider the models of representations in the previous section in
a more general setting,
and formalize two representation theoretic questions
(see Problem \ref{prob:Uncertainty}).

Let $P = LN = MAN$
be a parabolic subgroup of a real reductive Lie group $G$,
$\overline{P} = L\overline{N}$ its opposite parabolic subgroup, 
and
$\mathfrak{g} = \mathfrak{n} + \mathfrak{l} + \overline{\mathfrak{n}}$
the corresponding Gelfand--Naimark decomposition of the Lie algebra
$\mathfrak{g}$.
Assume that the nilradical $\mathfrak{n}$ is abelian,
and in particular $P$ is a maximal parabolic subgroup.

Take a (non-unitary) one-dimensional representation
$\chi: L \to \mathbb{C}^{\times}$,
and consider the induced representation
$W := \operatorname{Ind}_P^G(\chi \otimes \mathbb{C})$.
Then the space $W^\infty$ consisting of its smooth vectors can be 
regarded as a subspace of
$(C^\infty \cap \mathcal{S}')(\overline{\mathfrak{n}})$
by the restriction to
$\overline{\mathfrak{n}} \simeq \overline{N} \subset G/P$. 
Here, $\mathcal{S}'(\overline{\mathfrak{n}})$ denotes the space of 
\index{B}{tempered distribution}%
tempered distributions 
on $\overline{\mathfrak{n}}$ (regarded as the Euclidean space).
Then, taking the 
\index{B}{Euclidean Fourier transform}%
Euclidean Fourier transform,
we have
\begin{equation*}
   W^\infty \subset \mathcal{S}'(\overline{\mathfrak{n}}) 
   \xrightarrow[\mathcal{F}]{\sim} \mathcal{S}'(\mathfrak{n}),
\end{equation*}
where $\mathfrak{n}$ is identified with the dual space of
$\overline{\mathfrak{n}}$. 
We let $G$ act on $\mathcal{F}(W^\infty)$ through
$\mathcal{F}$. 
This $G$-action cannot be extended to $\mathcal{S}'(\mathfrak{n})$, 
but its restriction to the parabolic subgroup
 $\overline{P} = L\overline{N}$ can be extended to 
$\mathcal{S}'(\mathfrak{n})$ because $\overline{P}$ acts 
on $\mathcal{S}'(\mathfrak{n})$ just 
by translations and
multiplications of unitary characters.
Let $w_0$ be an element of $K$ such that 
$w_0 L w_0^{-1} = L$
and 
$w_0 N w_0^{-1} = \overline{N}$ (a conformal inversion).

Now we consider the following setting:
\begin{enumerate}[1)]
\item  
Let $(\pi,\mathcal{H})$ be an irreducible unitary representation of $G$,
such that the underlying $(\mathfrak{g}_{\mathbb{C}}, K)$-module
$\mathcal{H}_K$ is a subrepresentation of $W$,
and the Hilbert space $\mathcal{H}$ is realized in 
$\mathcal{S}'(\overline{\mathfrak{n}})$.
\item  
Let $C$ be an $L$-orbit in $\mathfrak{n}$,
having an $M$-invariant measure $d\mu$ such that
$L^2(C,d\mu) \subset \mathcal{S}'(\mathfrak{n})$.
Further,
$\overline{P} = L\overline{N}$ leaves $L^2(C,d\mu)$ invariant,
and acts as a unitary representation.
\item  
$\mathcal{F}(\mathcal{H}) = L^2(C,d\mu)$.
\end{enumerate}

We note that the condition (2) determines the absolute value
$|\chi(a)|$ for $a \in A$.
As we mentioned,
the $\overline{P}$-action on $L^2(C)$ is given just by translations
and multiplications of unitary characters.
Since $G$ is generated by $\overline{P}$ and $w_0$,
the action $\pi$ of 
the whole group $G$ on $L^2(C,d\mu)$ 
is determined by finding the formula of $\pi(w_0)$. 

The (Euclidean) Fourier transform $\mathcal{F}$ transfers the
 defining ideal for the affine variety $\overline{C}$ in $\mathfrak{n}$
 to the system,
to be denoted by $\mathcal{M}$,
of differential equations on $\overline{\mathfrak{n}}$
such that the space $\mathcal{H}$ is contained in the solution space: 
\begin{equation*}
\mathcal{S}ol(\mathcal{M}) := \{ f \in
 \mathcal{S}'(\overline{\mathfrak{n}}): Pf = 0
 \ \, \text{for any $P\in\mathcal{M}$} \}.
\end{equation*}
In the previous example,
$\mathfrak{n} \simeq \mathbb{R}^n$,
$C$ is the isotropic cone, 
and the system $\mathcal{M}$ is generated by
$\square_{\mathbb{R}^{n_1,n_2}}
 = \frac{\partial^2}{\partial x_1^2} +\dots+
   \frac{\partial^2}{\partial x_{n_1}^2} -
   \frac{\partial^2}{\partial x_{n_1+1}^2} -\dots-
   \frac{\partial^2}{\partial x_n^2}
$.

Thus, we have two models of the irreducible unitary representation
with the following nature:

\medskip 
\noindent
\fbox{\begin{minipage}{0.97\linewidth}
\vspace*{.7ex}
\textbf{Solution model}%
\index{B}{solution model} 
 on $\mathcal{S}ol(\mathcal{M})$ in
$\mathcal{S}'(\overline{\mathfrak{n}})$. 
  \begin{itemize}
  \item
  The `intrinsic inner product' on the Hilbert space
  $\mathcal{H}$ $(\subset \mathcal{S}ol(\mathcal{M}))$
  is not clear.
  \item
  The $G$-action on $\mathcal{H}$,
  to be denoted by $\varpi(g)$, is simple 
  (essentially, the translations by the M\"{o}bius transform of $G$ on
   $G/P \fallingdotseq \overline{\mathfrak{n}}$).
  \end{itemize}
\vspace*{.1ex}
\end{minipage}
}

\medskip 
\noindent
\fbox{\begin{minipage}{0.97\linewidth}
\vspace*{.7ex}
\textbf{$L^2$-model}%
\index{B}{L2model@$L^2$-model} 
 on $L^2(C,d\mu) \subset \mathcal{S}'(\mathfrak{n})$.
  \begin{itemize}
  \item
  The inner product on the Hilbert space $L^2(C,d\mu)$ is very clear.
  \item
  The $G$-action on $L^2(C,d\mu)$ is not simple except for the 
  $\overline{P}$-action.
  \end{itemize}
\vspace*{.1ex}
\end{minipage}
}

\medskip
\noindent
Now, we have a kind of 
`uncertainty relation'%
\index{B}{uncertainty relation} 
 in the sense that it is hard to find a
single model having 
explicit descriptions of both $G$-actions and inner products.
This feature in the above two models is symbolically
 summarized
as follows:
\begin{center}
\begin{tabular}{l|cc}
   & $\mathcal{S}ol(\mathcal{M})$
   & $L^2(C,d\mu)$ 
\\
\hline
\\
   inner product
   & ?
   & simple
\\[\medskipamount]
   $G$-action
   & simple
   & ?
\\
\end{tabular}
\end{center}

\medskip

Then, we ask:
\begin{problem} \label{prob:Uncertainty}
\quad
\begin{enumerate}[\upshape 1)]
\item  
\index{B}{conserved quantity}%
\textup{(conserved quantity for solution model)}%
\enspace
Find an 
intrinsic inner product%
\index{B}{intrinsic inner product} 
 on the solution space
$\mathcal{S}ol(\mathcal{M})$ that is invariant by $G$.
\item  
\textup{(generalization of the 
\index{B}{Fourier--Hankel transform}%
Fourier--Hankel transform for $L^2$-model)}\enspace
Find an explicit formula for the 
\index{B}{unitary inversion operator}%
unitary inversion operator
 $\pi(w_0)$.
\end{enumerate}
\end{problem}

1)
\textbf{Solution model}. 
In order to clarify the meaning of `intrinsic inner product',
we list three approaches to describe the inner product on
$\mathcal{H}$ in $\mathcal{S}ol(\mathcal{M})$.
\begin{enumerate}[{1}-a)]
\item  
\index{B}{Parseval type formula}%
(Parseval type formula)\enspace
Describe the inner product on $\mathcal{H}$ according to the
$K$-type decomposition.
\item  
\index{B}{Green function}%
(Green function)\enspace
Give an integral expression of the solutions to $\mathcal{M}$,
and describe the inner product on $\mathcal{H}$ by means of the
integral expression.
\item  
(Conserved quantities)\enspace
Find an inner product formula in terms of only solutions.
\end{enumerate}
The approaches (1-a) and (1-b) give explicit inner products in a
sense and are usually sufficient for representation theoretic purposes
(e.g.\ showing the unitarizability), 
however, do not give an intrinsic formula in the sense that the
formula depends on the $K$-type decomposition or on the integral
expression of solutions.
The approach (1-c) seeks for an intrinsic formula based purely on
solutions.
Here are 
a few remarks on (1-a)--(1-c) in order.

The approach (1-a) is algebraic.
We note that $(G,L)$ forms a reductive symmetric pair under the
assumption that the nilpotent radical $\mathfrak{n}$ is abelian.
Consequently,
$(K,L\cap K)$ is a compact symmetric pair,
and therefore 
$\operatorname{Ind}_P^G(\chi \otimes \mathbb{C})$ 
is 
\index{B}{multiplicity free}%
$K$-multiplicity free.
Then, 
the unitary inner product on $\mathcal{H}$ 
is a scalar multiple of the $L^2$-inner product on 
$L^2(K/L\cap K)$ on each $K$-type by Schur's lemma.
Thus, the unitary inner product on $\mathcal{H}$ is expressed by the 
`weight function'
\begin{equation*}
   m: \widehat{K}_{L\cap K} \to \mathbb{R}_+,
\end{equation*}
where 
$\widehat{K}_{L\cap K} = \{ \tau \in \widehat{K}:
 \operatorname{Hom}_K (\tau, L^2(K/L\cap K)) \ne 0 \}$
is explicitly known by the 
\index{B}{Cartan--Helgason theorem}%
Cartan--Helgason theorem.
The weight function $m$ for unitarizable subquotients have been
computed for a 
 number
of degenerate principal series representations, 
especially since the
influential paper \cite{xHoTan} by Howe and Tan.
More generally, 
\index{B}{discretely decomposable branching law}%
discretely decomposable branching laws 
for non-compact subgroups $H$ give an 
extension of the approach (1-a) (see \cite{xkors2}).

The approaches (1-b) and (1-c) are analytic.
The integral expression of the solutions in (1-b) corresponds to the
\index{B}{Knapp--Stein intertwining operator}%
Knapp--Stein intertwining operator 
$A$ in representation theory.
The $G$-invariance of the resulting bilinear form 
$f_1,f_2 \mapsto (f_1,Af_2)$ 
is clear,
and the positivity of the bilinear form implies the unitarizability of
$\operatorname{Image}A$. 
This formula of the inner product is explicit, 
however,
the formula is not written directly in terms of solutions.
In fact, 
 it is non-trivial to find a preimage $f$
such that $u=Af$ when a solution $u$ is given.

What we seek for in (1-c) is to describe \textit{directly} the inner
product on solutions.
The energy for the wave equation is invariant under the time
translation, 
and is a classic example of conserved quantities.
The unitarizability of a solution space to $\mathcal{M}$ predicts the
existence of positive definite conserved quantities.
If the solution space is an irreducible $G$-module,
then such conserved quantities must be unique up to scalar.
This is what we call the `intrinsic inner product' on the solution
space.
The uniqueness and the existence is predicted by representation
theory. 
Finding its explicit form would be a challenging problem in analysis,
arising naturally from representation theory.

2)
\textbf{$L^2$-model}.
The existence of $L^2$-models 
of small representations has been found 
 for some other reductive groups. 
See \cite{xclerc,xhno,xrove} for unitary highest weight
representations, 
\cite{xDvSa} for spherical cases
 by using Jordan algebras;
\cite{xSab} for $SO(4,3)$, \cite{xTo} for minimal representations 
of general 
reductive groups by using amalgamation of maximal parabolic subgroups.

In the setting that we discussed as $L^2$-model,
there is a simple action of a parabolic subgroup
$\overline{P}$ on $L^2(C)$,
and we know the existence of the $G$-action on $L^2(C)$ by
some other reasons. 
On the other hand, 
it is often the case that the missing piece is the \textit{explicit}
global formula on how to extend the action on $L^2(C)$ 
 from $\overline{P}$ to the whole group
$G$. 

Our subject stated in Problem \ref{prob:Uncertainty} (2) 
is to fill this missing piece
 by finding the formula of the `unitary inversion'
$\pi(w_0)$ for the conformal inversion element $w_0\in G$.
We remark that the operator $\pi(w_0)$ in the $L^2$-model can be
 written as
\begin{equation}\label{eqn:piFpi}
   \pi(w_0) = \mathcal{F} \circ \varpi(w_0) \circ
              \mathcal{F}^{-1},
\end{equation}
where $\varpi(w_0)$ is a simple action on the solution model
$\mathcal{S}ol(\mathcal{M})$.
However,
the right-hand side of the formula \eqref{eqn:piFpi} does not give 
a solution to 
Problem \ref{prob:Uncertainty} (2) readily, because 
it is not easy to carry out the computation of the composition of
integral forms in general.

On the other hand,
finding the formula of the unitary inversion $\pi(w_0)$ in the $L^2$
model $L^2(C)$ has a significant meaning.
To see this, we recall that there is a Bruhat decomposition
$G = \overline{P} \cup \overline{P} w_0 \overline{P}$.
Therefore, once we get a formula of $\pi(w_0)$,
then the whole group action of $G$ can be written by using $\pi(w_0)$ 
\textit{at most once} 
(without any further composition of integral operators).
Thus, the formula of $\pi(w_0)$ is critical for finding the global
formula of the whole action of $G$.

We have seen in Sections \ref{subsec:1.6a} and \ref{subsec:1.6} that
Problem \ref{prob:Uncertainty} (2) is settled for the Weil
representation for the metaplectic group $Mp(n,\mathbb{R})$ 
and the minimal representation of the indefinite
orthogonal group $O(n_1+1,n_2+1)$, respectively. 
More generally, 
in the case that $(\pi,\mathcal{H})$ is a minimal representation of a
reductive group $G$,
we expect that the operator could be described
 by means of some `special function'
  of one variable.

\section{Special functions and minimal representations}
\label{subsec:1.7}

Yet another theme is special functions.

In this book, we shall see special functions arise from the minimal
representation. 
For example, 
\index{B}{K-Bessel function@$K$-Bessel function}%
$K$-Bessel functions appear as the radial part of $K$-finite
vectors in the Schr\"{o}dinger model $L^2(C)$.
\index{B}{Meijer's $G$-function}%
Meijer's $G$-functions appear as the radial part of the integral
kernel of $\mathcal{F}_C$. 
\index{B}{Appell's hypergeometric function}%
Appell's hypergeometric functions bridge two models of the minimal
representation, namely, the Schr\"{o}dinger model and the conformal
model. 
All together, we develop a new line of investigation on various
special functions in connection with the minimal
representation.

\medskip
\noindent
\textbf{The `radial part' of the unitary inversion $\mathcal{F}_C$}
\nopagebreak

We begin with the 
\index{B}{Euclidean Fourier transform}%
(Euclidean) Fourier transform  
\index{A}{FRn@$\mathcal{F}_{\mathbb{R}^n}$}%
$\mathcal{F}_{\mathbb{R}^n}$ 
 as an illustrative example.
As we already discussed, 
this corresponds to the unitary inversion operator for the 
(original)
 Schr\"{o}dinger model of the Weil representation of
$Mp(n,\mathbb{R})$.

According to the polar coordinate
\begin{equation*}
   \mathbb{R}_+ \times S^{n-1} \to \mathbb{R}^n,
   \quad 
   (r,\omega) \mapsto r\omega,
\end{equation*}
we have a unitary equivalence:
\begin{align} \label{eqn:RnRs}
   L^2(\mathbb{R}^n)
   & \simeq L^2(\mathbb{R}_+, r^{n-1}dr) \hatotimes
            L^2(S^{n-1})
\nonumber
\\
   & \simeq \sideset{}{^\oplus}\sum_{l=0}^\infty  
            L^2(\mathbb{R}_+, r^{n-1}dr) \otimes
            \mathcal{H}^l(\mathbb{R}^n),
\end{align}
where 
\index{A}{HjRm@$\mathcal{H}^j(\mathbb{R}^m)$}%
$\mathcal{H}^l(\mathbb{R}^n)$ is the space of 
\index{B}{spherical harmonics}%
spherical harmonics of
degree $l$,
$ \{\varphi \in C^\infty(S^{n-1}):
 \Delta_{S^{n-1}} \varphi
 = -l(l+n-2)\varphi \}$
(see Appendix \ref{subsec:H}).
Here, 
\index{A}{3w@$\widehat{\otimes}$}%
$\widehat{\otimes}$
stands for the Hilbert completion of the tensor product space,
and
\index{A}{3sideset@$\sideset{}{^\oplus}\sum$}%
$\sideset{}{^\oplus}\sum$ 
stands for the Hilbert completion of an algebraic direct sum.

Correspondingly to the direct sum decomposition \eqref{eqn:RnRs},
the Fourier transform $\mathcal{F}_{\mathbb{R}^n}$ is decomposed as 
\begin{equation*}
\index{A}{FRn@$\mathcal{F}_{\mathbb{R}^n}$}%
   \mathcal{F}_{\mathbb{R}^n} 
   = \sideset{}{^\oplus}\sum_{l=0}^\infty
     T_l \otimes \operatorname{id}.
\end{equation*}
Here, $T_l$ is the 
\index{B}{Hankel transform}%
Hankel transform 
 of the following form
(see Remark \ref{rem:4.1.3}):
\begin{equation*}
   (T_l f)(r)
   = \frac{c}{r^n} \int_0^\infty x^{\frac{n}{2}}
     J_{\frac{n-2+2l}{2}} (x) f \Bigl(\frac{x}{r}\Bigr) dx.
\end{equation*}
Next, we consider our minimal representation
of $O(n_1+1,n_2+1)$ realized on $L^2(C)$.
Then, the 
\index{B}{bipolar coordinate}%
bipolar coordinate on the isotropic cone $C$,
\begin{equation*}
   \mathbb{R}_+ \times S^{n_1-1} \times S^{n_2-1} \to C
\end{equation*}
induces a unitary equivalence
\begin{align*}
\index{A}{L2C@$L^2(C)$}%
   L^2(C)
   &\simeq L^2(\mathbb{R}_+, \frac{1}{2} r^{n-3}dr) \hatotimes
           L^2(S^{n_1-1}) \otimes L^2(S^{n_2-1})
\\
   &\simeq \sideset{}{^\oplus}\sum_{l,k=0}^\infty
           L^2(\mathbb{R}_+, \frac{1}{2} r^{n-3} dr) \otimes
           \mathcal{H}^l(\mathbb{R}^{n_1}) \otimes
           \mathcal{H}^k(\mathbb{R}^{n_2}).
\end{align*}
Then, the unitary inversion operator $\mathcal{F}_C$ is decomposed as
\begin{equation*}
   \mathcal{F}_C
   = \sideset{}{^\oplus}\sum_{l,k=0}^\infty 
\index{A}{Tlk@$T_{l,k}$}%
T_{l,k} \otimes
     \operatorname{id} \otimes \operatorname{id},
\end{equation*}
where $T_{l,k}$ is a unitary operator on 
$L^2(\mathbb{R}_+, \frac{1}{2} r^{n-3} dr)$.
The unitary operators $T_{l,k}$ may be regarded
 as a generalization of Hankel transforms.

It turns out that the kernel function of $T_{l,k}$ is real analytic 
(see Theorem \ref{thm:TlkG} below). 
This is a good
 contrast to the fact that the unitary operator $\mathcal{F}_C$ 
on $L^2(C)$ is given by a
distribution kernel for general $n_1,n_2 > 1$
(see Theorem \ref{thm:Kformula}).

\begin{thmSec}[see Theorem \ref{thm:C}]
 \label{thm:TlkG}
Let $G_{04}^{20}$ be 
\index{B}{Meijer's $G$-function}%
Meijer's $G$-function 
 (see Appendix \ref{subsec:G} for definition), 
and we define a real analytic function $K_{l,k}$ by
\begin{equation*}
\index{A}{Klkt@$K_{l,k}(t)$}%
   K_{l,k}(t)
   := 4(-1)^{a+\frac{n_1+k+l}{2}} G_{04}^{20}
      (t^2 \mid \frac{l+k}{2}, a+2, \frac{-n_1-n_2+4-l-k}{2}, b+2),
\end{equation*}
where
\begin{equation*}
    a := \max \Bigl(\tfrac{-n_1-l+k}{2}, \tfrac{-n_2+l-k}{2}\Bigr),
 \
    b := \min \Bigl(\tfrac{-n_1-l+k}{2}, \tfrac{-n_2+l-k}{2}\Bigr).
\end{equation*}
Then, we have
\begin{equation*}
   (T_{l,k} f)(r)
   = \frac{1}{2} \int_0^\infty K_{l,k} (rr') f(r') r'{}^{n-3} dr'.
\end{equation*}
\end{thmSec}

It is noteworthy that Meijer's $G$-functions
arise in the representation theory of reductive Lie groups.
We observe the Casimir operator of a maximal compact
subgroup $K$ acts on $L^2(C)$ as a
fourth order differential operator. 
Correspondingly,
 Meijer's $G$-functions $G_{04}^{20} (x|b_1,b_2,b_3,b_4)$ solve 
ordinary differential equations of order four
(see \eqref{eqn:diffeqG24}):
\begin{equation*}
   \prod_{j=1}^4 (x \frac{d}{dx} -b_j) u = 0.
\end{equation*}
In the case $n_1=1$ (or $n_2=1$),
our minimal representation $\pi$ is a direct sum of a highest weight
representation and a lowest weight representation.
In this case,
the kernels $K_{l,k}$ collapse to Bessel functions,
and the unitary operators
$T_{l,k}$ are reduced to 
\index{B}{Hankel transform}%
Hankel transforms.

The group law $w_0^2=1$ in $G$
implies $\pi(w_0)^2= \operatorname{id}$,
and consequently,
 $T_{l,k}^2= \operatorname{id}$ for every
$l,k\in\mathbb{N}$.
Hence, 
Theorem \ref{thm:TlkG} gives
a group theoretic proof for 
the 
\index{B}{Plancherel formula!0for Meijer's $G$-transform@---, for Meijer's $G$-transform}%
Plancherel 
 and 
\index{B}{reciprocal formula!0for Meijer's $G$-transform@---, for Meijer's $G$-transform}%
reciprocal formulas 
on 
\index{B}{Meijer's $G$-transform}%
Meijer's $G$-transforms 
which were first proved by C. Fox \cite{xFo}
by a completely different method.

\begin{corollary}[see Corollary \ref{cor:C1}]
\label{cor:GPlancherel}
Let $b_1, b_2, \gamma$ be half-integers such that $b_1 \ge 0$, 
$\gamma \ge 1$, $\frac{1-\gamma}2 \le b_2 \le \frac{1}2+b_1$.
Then, the integral transform 
\begin{equation*}
  S_{b_1,b_2,\gamma}: f(x) \mapsto 
\frac{1}\gamma 
\int_0^\infty \!\!\! G_{04}^{20} 
((xy)^{\frac{1}\gamma} \mid b_1, b_2, 1-\gamma-b_1, 1-\gamma-b_2)f(y)dy
\end{equation*}
is a unitary operator on $L^2(\mathbb R_+)$.
\end{corollary}

\begin{corollary}[see Corollary \ref{cor:C2}]
\label{cor:Greciprocal}
The unitary operator $S_{b_1,b_2,\gamma}$ is of order two in 
$L^2(\mathbb R_+)$, that is, 
$(S_{b_1,b_2,\gamma})^{-1} = S_{b_1,b_2,\gamma}$.
\end{corollary}

A special case of the above corollaries
(i.e.\ $n_2=1$ case)
yields the classic formulas of the Hankel
transform (see Remark \ref{rem:4.1.6}).

\medskip
\noindent
\textbf{$K$-finite vectors in $L^2(C)$}
\nopagebreak

By the general theory due to Vogan
\cite{xV}, 
$K$-types of minimal representations $\pi$ 
are indexed by a natural number $a=0,1,2,\dotsc$
(by this property, $\pi$ is an example of the so-called 
\index{B}{ladder representation}%
\textit{ladder representation}). 
In contrast to the conformal model on
$\mathcal{S}ol(\widetilde{\Delta}_{S^{_p-1}\times S^{q-1}})$
where explicit $K$-finite vectors are given readily by spherical
harmonics, 
it is not clear a priori what $K$-finite vectors look like in the
Schr\"{o}dinger model $L^2(C)$ because the whole group $K$ cannot act
on the isotropic cone $C$.

Our idea is to compute explicitly the intertwining integral operator
$\mathcal{T}^{-1}$ between these two models
 in Diagram \ref{diag:1.7.1}. 
Then, by using a reduction formula of 
\index{B}{Appell's hypergeometric function}%
Appell's hypergeometric functions, 
we have (loosely):
\begin{equation*}
 \mathcal{T}^{-1} \ 
\index{B}{Gegenbauer polynomial}%
\text{(Gegenbauer polynomials)}
 =
\index{B}{K-Bessel function@$K$-Bessel function}%
\text{$K$-Bessel functions}
\end{equation*}
and prove the following result:
\begin{thmSec}[see Corollary \ref{cor:Kfinite}]
 \label{thm:KfiniteC}
Let\/ $n_1 \ge n_2$.
For $a=0,1,2,\dotsc$, 
\begin{equation*}
   r^{a-\frac{n_2-1}{2}}
   K_{\frac{n_2-1}{2}} (2r) \phi(w)
   \quad (\phi \in \mathcal{H}^a(\mathbb{R}^{n_1}))
\end{equation*}
is a $K$-finite vector in $L^2(C)$.
In the $K$-type formula
(see \eqref{eqn:Ktype}),
this vector belongs to the $K$-type
$\mathcal{H}^a(\mathbb{R}^{n_1+1}) \otimes
 \mathcal{H}^{a+\frac{n_1-n_2}{2}}(\mathbb{R}^{n_2+1})$.
\end{thmSec}
In Theorem \ref{thm:KfiniteC},
the $a=0$ case corresponds to the minimal $K$-type,
and was previously proved in \cite[Theorem 5.8]{xkors3}.
We note that $\pi$ is spherical if $n_1 = n_2$.
Even in the case of spherical representations,
finding explicit forms of the $K$-fixed vectors in $L^2$-model is
non-trivial.
See \cite{xDvSa} for
similar formulas of the $K$-fixed vectors
 in $L^2$-models for some other groups by means of $K$-Bessel
functions.

\section{Organization of this book}

This book is organized as follows. 
We review quickly the $L^2$-realization
(a generalization of the classic Schr\"odinger model) of the minimal representation of 
$O(p,q)$ in the first half of 
Chapter \ref{sec:rev}, 
and develop a basic theory of fundamental differential operators on
the isotropic cone $C$ in the latter half.
Then we find some $K$-finite vectors on $L^2(C)$
explicitly by means of the $K$-Bessel function $K_\nu(z)$ in Chapter  
\ref{sec:K}. 
Chapter \ref{sec:C} is devoted entirely
to the integral formula of the unitary operator 
$T_{l,k}$ on $\wL$ for double spherical harmonics expansions 
(see Theorem \ref{thm:C}). 
In Chapter \ref{sec:A}, 
building on the results of Chapter \ref{sec:C}, 
 we complete the proof of our main theorem
(see Theorem \ref{thm:A}). 
In order to make the proof readable
 as much as possible, we collect in Appendix the
formulas and the properties of various special functions used in this book.

\section{Acknowledgements}

A large part of the results here was obtained while both  authors were
at the Research Institute for Mathematical Sciences, Kyoto. 
The materials for Sections \ref{subsec:1.1}--\ref{subsec:FC} 
and \ref{subsec:Pjb}--\ref{subsec:w} 
were
developed while the first author was participating in the program
\lq\lq Representation Theory, Complex Analysis and Integral
Geometry\rq\rq \ 
 organized by B. Kr\"otz and S. Gindikin
 at Max-Planck-Institut f\"ur Mathematik, Bonn in 2007. 
We are very grateful to the colleagues and the staff for the wonderful
atmosphere of research at these institutes.

The authors have benefited from various discussions with
J. Faraut, S. Gindikin, S. Hansen, M. Kashiwara, T. Koorwinder, B. \O
rsted, H. Sabourin, S. Sahi, P. Torasso, and D. Vogan.

Some of the results of this book were presented
at the workshop on ``Harmonic Analysis and Homogeneous Spaces"
in honor of Professor G. van Dijk at Lorentz Center in Leiden in 2004,
at the Faraut Seminar in Paris in 2005,
at the workshop on ``Representation Theory and Prehomogeneous Vector Spaces"
at IRMA, Strasbourg
 in 2006, and
 at the workshop in Poitiers and
at the Analysis Seminar of Aarhus University in 2007.
We express our deep gratitude to the organizers of these workshops and seminars, 
and to the
participants for helpful and stimulating comments on various occasions.

A part of the results here was announced in \cite{xkmanoprojp} with a
sketch of proof.

\medskip
The first author would like to thank Ms. Suenaga
for indispensable help in preparing \LaTeX{} manuscript.
\medskip

Notation: $\mathbb R_+:=\set{x \in \mathbb R}{x>0}, \ 
\mathbb N:=\{0,1,2,\cdots\}$.

\chapter{Two models of the minimal representation of $O(p,q)$}
\label{sec:rev}

This chapter gives an account of the connection of
 the following two topics:
\newline\indent 
 1) Analysis on the isotropic cone $C$ for
   commuting differential operators $P_j$
 ($1 \le j \le n$) associated to the quadratic
 form of signature $(n_1, n_2)$
 (see Introduction \ref{subsec:1.1}--\ref{subsec:FC}).
\newline\indent
2) Minimal representation of the indefinite
 orthogonal group $O(p,q)$.

Throughout this chapter, we shall use the following
 notation:
$$
   p=n_1+1, \quad q=n_2+1, \quad n = n_1 + n_2 = p+q-2.
$$
 
The first half of this chapter is
 a review from \cite{xkors1,xkors3}
about two concrete models of the minimal
representation of the group $G=O(p,q)$, namely,
the conformal model $(\varpi^{p,q}, \overline{V^{p,q}})$
 using the Yamabe operator \eqref{eqn:Yamabe} in Section
\ref{subsec:Kpic} 
and 
the $L^2$-model 
(the Schr\"odinger model) $(\pi, L^2(C))$
in Section \ref{subsec:Sch}.
In the terminology of representation theory of reductive Lie groups
(e.g.\ \cite{xKnapp,xwal}), 
 the former realization gives a subrepresentation (the 
\index{B}{K-picture@$K$-picture}%
$K$-picture,
the $N$-picture, etc.)\ of a degenerate principal
series representation, 
whereas the latter corresponds to the dual of
the 
\index{B}{N-picture@$N$-picture}%
$N$-picture via the (Euclidean) Fourier transform $\mathcal{F}$.  

The intertwining operator $\mathcal{T}$ between these two
models will be given in \eqref{def:B},
which is summarized as the following diagram:
\[
\begin{matrix}
   ~~L^2(C)    &  \overset{T}\hookrightarrow 
                   & \mathcal S'(\mathbb R^{p+q-2})   \\
  \mathcal{T} 
  \begin{picture}(10,10)
  \put(5,0){\dashbox(0,8.5)}
  \put(5,-1){\vector(0,-1){1}}
  \end{picture}
  &       &   \uparrow \mathcal{F}  \\
  ~~~K\text{-picture}   &  \underset{\widetilde\Psi^*}\to   
                           &  N\text{-picture} .
\end{matrix}
\]
Here, $T$ is the identification map between functions on $C$
and distributions supported on $C$
  by the canonical measure on $C$,
$\widetilde\Psi^*$ is the $G$-intertwining operator between the
$K$-picture and the $N$-picture,
and is interpreted as the twisted pull-back for the conformal
map $\Psi : \mathbb{R}^{p+q-2} \to S^{p-1} \times S^{q-1}$.

The latter half of this chapter is new.
In Section \ref{subsec:Pjb},
we analyze a commuting family of differential operators 
$P_j(b)$ $(1 \le j \le n)$ of second order with parameter $b$ in
$\mathbb{R}^n$, 
and prove that they are tangential to the isotropic cone $C$ if
$b=1$. 
The resulting differential operators $P_j := P_j(1)|_C$ 
on $C_0^\infty(C)$ 
extend to
self-adjoint operators on $L^2(C)$
(see Theorem \ref{thm:Pj}).
Thus,
we get a family of commuting differential operators
$P_j$ $(1 \le j \le n)$ of second order,
which we call \textit{fundamental differential operators}
on $C$.
In Section \ref{subsec:w},
we see that the unitary inversion operator $\mathcal{F}_C = \pi(w_0)$
diagonalizes $P_j$ $(1 \le j \le n)$
and that the intertwining relation of $P_j$ and the multiplication by
coordinate function $x_j$ characterizes $\mathcal{F}_C$ up to scalar.
Thus, 
we develop an abstract theory of the unitary operator $\mathcal{F}_C$
on $L^2(C)$ by taking the (Euclidean) Fourier transform
$\mathcal{F}_{\mathbb{R}^n}$ on $L^2(\mathbb{R}^n)$ as a prototype
(see Theorems \ref{thm:FC}, \ref{thm:FCunique}, \ref{thm:Pjunique}).

\section{Conformal model}
\label{subsec:Kpic}

This section summarizes the conformal model of
the minimal representation of the indefinite orthogonal group
$G=O(p,q)$ ($p+q: $ even).
The advantage of the 
\index{B}{conformal model}%
conformal model is that the group action on the
representation is simple and that its geometric idea is clear.
Since this conformal model corresponds to a subrepresentation of the most degenerate
principal series representations (with a very special parameter),
the same representation can be studied also by the purely algebraic
method of $(\mathfrak{g},K)$-modules.
See \cite{xBiZi,xHoTan,xKo} in this direction.
The same subrepresentation can also be captured by the theta
correspondence arising from the dual pair
$O(p,q) \cdot SL(2,\mathbb{R}) \subset Sp(p+q,\mathbb{R})$
(see \cite{xHuZu}).
Our approach in this section is geometric, and the basic reference
here is \cite{xkors1}.
See also \cite{xkcheck} for an elementary exposition from viewpoints
of conformal transformation groups.

The general geometric idea here is summarized as follows.
Let $X$ be an $n$-dimensional manifold equipped with a Riemannian (or
more generally, pseudo-Riemannian) structure $g$.
Then, associated to $g$, we define:
\begin{align}
  & \kappa: \text{the scalar curvature on $X$},
\nonumber
\\
  & \Delta_X: 
\index{B}{Laplace--Beltrami operator}%
\text{the Laplace--Beltrami operator on $X$},
\nonumber
\\
  & 
\index{A}{1Delta_X@$\protect\widetilde{\Delta}_X$}%
\widetilde{\Delta}_X := \Delta_X + \frac{n-2}{4(n-1)}\kappa
\index{B}{Yamabe operator}%
\quad\text{(Yamabe operator)}.
\label{eqn:Yamabe}
\end{align}
Then, although these objects depend on the (pseudo-)Riemannian
structure $g$,
the solution space
\begin{equation*}
\index{A}{SolDeltaX@$\mathcal{S}ol(\widetilde{\Delta}_X)$}%
   \mathcal{S}ol(\widetilde{\Delta}_X)
   := \{ f\in C^\infty(X): \widetilde{\Delta}_Xf = 0 \}
\end{equation*}
is conformally invariant, namely, if
$\varphi:X\to X$ is a conformal diffeomorphism with a conformal factor
$c_\varphi \in C^\infty(X)$ satisfying
\begin{equation*}
   \varphi^* g = c_\varphi^2 g
\end{equation*}
then
\begin{equation} \label{eqn:conf}
   C^\infty(X) \to C^\infty(X),
   \quad
   f \mapsto c_\varphi^{-\frac{n-2}{2}} f \circ \varphi
\end{equation}
leaves $\mathcal{S}ol(\widetilde{\Delta}_X)$ invariant,
and hence we get a representation of the conformal transformation
group 
$\operatorname{Conf}(X,g)$ on $\mathcal{S}ol(\widetilde{\Delta}_X)$
(see \cite[Theorem A]{xkors1}).
The point here is that the above construction is functional under
conformal maps,
and in particular,
if two pseudo-Riemannian manifolds are conformally equivalent
(not necessarily isometric),
then the resulting two representations are isomorphic.

A special case applied to pseudo-Riemannian manifolds which are
conformally equivalent to flat pseudo-Riemannian space forms
 gives rise to the minimal representations of the
indefinite orthogonal groups.
Let us explain this specific case in more details.

We denote by 
\index{A}{Rpq@$\mathbb R^{p,q}$}%
$\mathbb R^{p,q}$ 
the Euclidean space $\mathbb R^{p+q}$ equipped with 
the pseudo-Riemannian structure $g_{\mathbb R^{p,q}}$ of signature $(p,q)$:
$$
ds^2=dx_1^2+\cdots+dx_p^2-
dy_1^2-\cdots-dy_q^2.
$$ 
Then, the restriction of $ds^2$ to the submanifold 
\begin{align}\label{def:M}
\index{A}{MS@$M\simeq S^{p-1}\times S^{q-1}$}%
  M  & := \set{(x,y) \in \mathbb{R}^{p+q}}{|x|=|y|=1, \ x \in \mathbb{R}^p
           , \ y \in \mathbb{R}^q} \\
        & \simeq S^{p-1}\times S^{q-1}  \nonumber
\end{align}
is non-degenerate,
and defines a pseudo-Riemannian structure on $M$ of signature
$(p-1,q-1)$. 
Here, $|\cdot | $ stands for the usual Euclidean norm.
The resulting pseudo-Riemannian structure $g_M$ on $M$ is nothing but the
direct product of the standard unit sphere $S^{p-1}$ (positive
definite metric) and the unit sphere $S^{q-1}$ equipped with the negative definite metric 
($(-1) \times $
the standard metric).

Then,
the Yamabe operator $\tilLap{M}$ of $M$ takes the following form
(see \cite[(3,4,1)]{xkors3}):
\begin{equation}
\label{eqn:Ya}
      \tilLap{M}=\Delta_{S^{p-1}}-\Delta_{S^{q-1}}
        -\left(\frac{p-2}{2}\right)^2+\left(\frac{q-2}{2}\right)^2, 
\end{equation}
where $\Delta_{S^{p-1}}$ and $\Delta_{S^{q-1}}$ are the Laplace--Beltrami
operators on $S^{p-1}$ and $S^{q-1}$ respectively.

The indefinite orthogonal group 
\index{A}{O(p,q)@$O(p,q)$}%
$G=O(p,q)$ 
acts naturally on $\mathbb R^{p,q}$
 as isometries.  
This action preserves the cone
$$
     \Xi :=\{(x,y)\in \mathbb R^{p,q}: |x|=|y| \ne 0\}
$$
 but does not preserve $M$.  
In order to let $G$ act on $M$, 
we set a function $\nu$ on
$\mathbb{R}^{p,q}$ by
$$   \nu : \mathbb{R}^{p,q} \to \mathbb{R}, \quad  (x,y)\mapsto |x|.  
$$  
If $v \in M(\subset \Xi)$ and $h \in G$, 
then $h \cdot v \in \Xi$, 
and consequently
 $\frac{h \cdot v}{\nu(h \cdot v)} \in M$.  
Thus, 
we can define the action of $G$ on $M$: 
$$   
     L_h: M \to M,  \quad v \mapsto \frac{h\cdot v}{\nu(h\cdot v)} \quad (h \in G).
$$
Then, 
 we have $L_h^{\ast} g_M=\frac{1}{\nu(h \cdot v)^2} g_M$ at $T_vM$
 and thus the diffeomorphism $L_h$ is conformal with respect to the
pseudo-Riemannian metric on $M$.
Conversely, any conformal diffeomorphism of $M$ is of the form $L_h$
for some $h\in G$
 (see \cite[Chapter IV]{xskob}).

By the general theory (see \eqref{eqn:conf}) of conformal geometry,
we can construct a representation, denoted by
$\varpi^{p,q}$, of $G$
on the solution space to $\tilLap{M}$ in $C^\infty (M)$:
\begin{align}
 &
\index{A}{Vpq@$\Vpq$|main}%
\Vpq 
:= \mathcal{S}ol \tilLap{M}=
\set{f \in C^\infty(M)}{\tilLap{M} f =0},  \notag
\\
\intertext{where we set}
\label{def:vpi}
    &(
\index{A}{1pipq@\protect$\varpi^{p,q}$|main}%
\varpi^{p,q}
(h^{-1})f) (v) := \nu(h \cdot v)^{-\frac {p+q-4} 2}
                                     f(L_h v),  
\end{align}
for $h \in  G, \ v \in M$, and $f \in V^{p,q}$.
The following theorem was proved in \cite{xkors3}
in this geometric framework.
There are also algebraic proofs 
(see Remark \ref{rem:subrep}).  
\begin{fact}[{see \cite[Theorem~3.6.1]{xkors3}}]
\label{fact:A}
Let $p, q \ge 2$ and $p+q \ge 6$ be even.

{\rm 1)}(irreducibility)
$(\varpi^{p,q},V^{p,q})$ is an irreducible unitary representation of $G$.

{\rm 2)}(unitarizability)
There exists a $G$-invariant inner product
$(\ ,\ )_M$ on $V^{p,q}$.
Such a $G$-invariant inner product is unique up to a scalar multiple,
and we shall normalize it in \eqref{eqn:innerM}.

We write $\overline{V^{p,q}}$ for the Hilbert completion of $V^{p,q}$,
and use the same letter $\varpi^{p,q}$ to denote the resulting
irreducible unitary representation.

{\rm 3)}($K$-type formula) 
Let $K \simeq O(p)\times O(q)$ be a maximal compact subgroup of $G$.
Then, the restriction
of $(\varpi^{p,q}, \overline \Vpq)$ to $K$ decomposes into irreducible
representations of $K$ as follows:
\begin{equation}
\label{eqn:Ktype}
\overline \Vpq \simeq 
\sideset{}{^\oplus}\sum_{\substack{a+\frac{p-q}2=b,\\
                                   a,b \in \mathbb N}}^\infty 
\Har{a}{p} \otimes \Har{b}{q}.
\end{equation}
Here, 
\index{A}{HjRm@$\mathcal{H}^j(\mathbb{R}^m)$}%
$\mathcal{H}^a(\mathbb{R}^p)$ denotes the irreducible
representation of $O(p)$ on the space of spherical harmonics of degree
$a$
(see Section \ref{subsec:H}).

{\rm 4)} (Parseval--Plancherel formula) 
On each $K$-type $\Har{a}{p} \otimes \Har{b}{q}$
for $(a,b) \in \mathbb{N}^2$ such that
$a + \frac{p-q}{2} = b$, or equivalently,
$a + \frac{p-2}{2} = b + \frac{q-2}{2}$,
the unitary inner product $(\cdot, \cdot)_M$ is of the form:
\begin{equation}\label{eqn:innerM}
  (F, F)_M = \Bigl(a+\frac{p-2}2\Bigr) \|F\|_{L^2(M)}^2.
\end{equation}
\end{fact}

Next, we consider the following 
 injective map (see \cite[(2.8.2)]{xkors3}) by
$$
\Psi: \mathbb{R}^{p+q-2} \to M, \quad
z \mapsto \tau(z)^{-1} \iota(z),
$$
where for $z = (z',z'') \in \mathbb{R}^{p-1} \oplus \mathbb{R}^{q-1}$ 
we set
\begin{align*}
&\tau(z) := \biggl(1+\Bigl(\frac{|z'|+|z''|}{z}\Bigr)^2\biggr)^{\frac{1}{2}}
           \biggl(1+\Bigl(\frac{|z'|-|z''|}{2}\Bigr)^2\biggr)^{\frac{1}{2}},
\\
&\iota: \mathbb{R}^{p+q-2} \to \mathbb{R}^{p+q}, \ 
 (z',z'') \mapsto \biggl(1-\frac{|z'|^2-|z''|^2}{4}, z', z'',
                   1+\frac{|z'|^2-|z''|^2}{4}\biggr).
\end{align*}
Then, 
$\Psi$ is a conformal map
 such that $\Psi^{\ast}g_M=\tau(z)^{-2} g_{\mathbb R^{p-1,q-1}}$.  
According to \cite[Definition 3.4]{xkors1},
the twisted pull-back $\widetilde{\Psi}^*$ of the conformal map $\Psi$
is a linear map
\begin{equation*}
\index{A}{1Psi@$\widetilde{\Psi}^*$}%
   \widetilde{\Psi}^*
 : C^\infty(M) \to C^\infty(\mathbb{R}^{p+q-2})
\end{equation*}
given by
\begin{equation} \label{eqn:Psitilde}
   (\widetilde{\Psi}^* f)(z)
   := \tau(z)^{-\frac{p+q-4}{2}} (f\circ \Psi)(z).
\end{equation}

The image $M_+$ of $\Psi$ is roughly the half of $M$:
$$
\index{A}{Mz@$M_+$}%
M_+ := \{ u = (u_0,u',u'',u_{p+q-1}) \in M:
u_0+u_{p+q-1} > 0 \}.  
$$
We note that $\Psi$ induces a 
\index{B}{conformal compactification}%
conformal compactification of the flat space $\mathbb R^{p-1,q-1}$:
$$
     \mathbb R^{p-1,q-1} \hookrightarrow (S^{p-1} \times S^{q-1})/ \sim {\mathbb Z}_2, 
$$
 where $\sim {\mathbb Z}_2$ denotes the equivalence relation in $M=S^{p-1} \times S^{q-1}$
 defined by $u \sim -u$.  

The inverse of $\Psi: \mathbb{R}^{p+q-2} \overset{\sim}{\to} M_+$ 
is given by
$$
\Psi^{-1} (u_0,u',u'',u_{p+q-1}) =
\Bigl(\frac{u_0+u_{p+q-1}}{2}\Bigr)^{-1}
(u',u'').
$$
We note that $\Psi^{-1}$ is the ordinary
stereographic projection of the sphere $S^{p-1}$
if $q=1$.  

We write $(\Tilde{\Psi}^\ast)^{-1}=\widetilde{(\Psi^{-1})}^*$ for
the twisted pull-back 
(in the sense of \cite[Definition 2.3]{xkors1})
of the conformal map $\Psi^{-1}: M_+ \to \mathbb{R}^{p+q-2}$,
that is, 
$$
(\Tilde{\Psi}^\ast)^{-1}: C^\infty(\mathbb R^{\n}) \to
C^\infty(M_+)
$$ 
is given by
\begin{equation}
\label{def:Psiinv}
  (\Tilde{\Psi}^\ast)^{-1}(F)(v) :=
  (\frac{v_0+v_{p+q-1}}2)^{-\frac{p+q-4}2} 
    F(\frac{2}{v_0+v_{p+q-1}} \begin{pmatrix} v' \\ v'' \end{pmatrix} ),
\end{equation}
where $v={}^t \! (v_0, v', v'', v_{p+q-1}) \in M_+$, 
$v_0, v_{p+q-1} \in \mathbb R$,
$v' \in \mathbb R^{p-1}$, $v'' \in \mathbb R^{q-1}$.
In the group language (e.g.\ \cite{xKnapp}), 
$(\widetilde{\Psi}^*)^{-1}$ is the standard intertwining operator from the $N$-picture to the $K$-picture.  
The map $(\widetilde{\Psi}^*)^{-1}$ will be applied also to other
classes of functions.

\begin{remNonumber}
\label{rem:subrep}
Our manifold $M$ is a double cover of the generalized flag variety
 $G/\overline{\Pmax}$ by a maximal parabolic subgroup
$\overline{\Pmax}$ (see \eqref{eqn:Pmax}).  
Then, 
 $(\varpi^{p,q}, V^{p,q})$ is identified with a subrepresentation
 of the degenerate principal series representation induced from a certain
 one-dimensional representation of $\overline{\Pmax}$.  
In this framework, 
 Fact \ref{fact:A} was proved by Kostant \cite{xKo} for $p=q=4$
 and by Binegar and Zierau \cite{xBiZi}, for general
$p,q$ satisfying the condition that $p, q \ge 2$ and $p+q>4$.
Zhu and Huang \cite{xHuZu} identified
 this subrepresentation with the local 
\index{B}{theta correspondence}%
theta correspondence 
associated to the dual pair
 $O(p,q) \times SL(2,\mathbb R)$ in $Sp(p+q, \mathbb R)$
 (to be more precise, 
its metaplectic cover)
 and the trivial one-dimensional representation of $SL(2,\mathbb R)$.  
\end{remNonumber}
  
\begin{remNonumber}
If $p+q \ge 8$, then $\varpi^{p,q}$ becomes a 
\index{B}{minimal representation}%
minimal representation in the sense 
of Definition \ref{def:minrep} 
 (see \cite{xBiZi}).
\end{remNonumber}

\section{$L^2$-model (the Schr\"odinger model)}
\label{subsec:Sch}

In Sections \ref{subsec:Sch} and \ref{subsec:act}, 
we summarize the known results on 
the Schr\"odinger model of the 
minimal representation $(\pi, L^2(C))$ of $G=O(p,q)$.
The basic reference is \cite{xkors3}.

A naive idea here is the following.
Since $\Psi: \mathbb{R}^{p-1,q-1}\to M$ is a conformal map between two
pseudo-Riemannian manifolds, 
we have
\begin{equation*}
\index{A}{Soldelta@$\mathcal{S}ol(\widetilde{\Delta}_{S^{n_1}\times S^{n_2}})$}
   \widetilde{\Psi}^*(\mathcal{S}ol(\widetilde{\Delta}_M))
   \subset 
\index{A}{Solsquare@$\mathcal{S}ol(\square_{\mathbb{R}^{n_1,n_2}})$}
\mathcal{S}ol(\square_{\mathbb{R}^{p-1,q-1}}).
\end{equation*}
Since $M$ is compact,
it follows from the formula \eqref{eqn:Psitilde} of
$\widetilde{\Psi}^*$ that
\begin{equation*}
   \widetilde{\Psi}^*(C^\infty(M))
   \subset \mathcal{S}' (\mathbb{R}^{p+q-2}),
\end{equation*}
where $\mathcal{S}'(\mathbb{R}^{p+q-2})$ denotes the space of tempered distributions.
By taking the Euclidean Fourier transform
$\mathcal{F}_{\mathbb{R}^{p+q-2}}$,
we get
\begin{equation*}
   (\mathcal{F}_{\mathbb{R}^{p+q-2}} \circ \widetilde{\Psi}^*)
   (\mathcal{S}ol(\widetilde{\Delta}_M))
   \subset \{ u \in \mathcal{S}'(\mathbb{R}^{p+q-2}):
              \operatorname{Supp} u \subset \overline{C} \}.
\end{equation*}
It was proved in \cite{xkors3} that the right-hand side is contained
in $L^2(C)$ (and in fact a dense
subspace of $L^2(C)$).
Let us fix some notation to formalize this fact.

We define a quadratic form by 
\begin{equation}
\index{A}{Qx@$Q(\xzeta)$|main}%
 Q(\xzeta) := \xzeta_1^2 +\cdots +\xzeta_{p-1}^2
-\xzeta_p^2-\cdots-\xzeta_{p+q-2}^2,
\end{equation}  
which is the defining polynomial of the isotropic cone $C$ in $\mathbb R^{p+q-2}$.

The substitution $\delta(Q)$ of $Q$ into the 
\index{B}{Dirac delta function}%
Dirac delta function $\delta$
of one variable
defines a distribution on $\mathbb R^{\n} \setminus \{0\}$,
which is represented as a measure,
denoted by $d\mu$ on $C$. 
Alternatively, this measure is obtained as the volume form $\alpha|_C$
where
$\alpha$ is any $(p+q-3)$ form such that 
$$
d Q \wedge \alpha = d \xzeta_1 \wedge d \xzeta_2 \wedge \cdots \wedge d \xzeta_{p+q-2}
$$
(see \cite[Chapter III, Chapter 2]{xGeSh}). 
An example of such $\alpha$ is given by
\begin{equation*}
   \sum_{j=1}^n
   \frac{(-1)^{j-1} \epsilon_j x_j dx_1\wedge\cdots\wedge \widehat{dx}_j
         \wedge\cdots\wedge dx_n}
        {2 \|x\|^2} .
\end{equation*}

In the bipolar coordinate:
\begin{equation}
\label{def:pol}
  \mathbb R_+ \times S^{p-2} \times S^{q-2} \overset{\sim}{\to} C, \quad 
  (r, \omega, \eta) \mapsto 
 \begin{pmatrix} r\omega \\ r\eta \end{pmatrix},
\end{equation}
the distribution 
\index{A}{1delta@$\delta (Q)$|main}%
$\delta (Q)$
 is given by
\begin{equation}
\label{eqn:Q}
  \langle \delta(Q), \varphi \rangle 
  = \frac{1}2  \int_0^\infty \int_{S^{p-2}} \int_{S^{q-2}}
     \varphi(\begin{pmatrix} r\omega  \\ r\eta\end{pmatrix}) 
     r^{p+q-5} dr d\omega d\eta
\end{equation}
for a test function $\varphi$ on $\mathbb{R}^{p+q-2}$.
Here,
$d\omega$ and $d\eta$ denote the standard measures on 
$S^{p-2}$ and  $S^{q-2}$, respectively.
By this formula,
we see that 
$\delta(Q)$ extends to a Schwartz 
distribution on $\mathbb R^{\n}$ of measure class
if $p+q >4$ because $r^{p+q-5}dr$ is locally integrable,  
then.
Obviously, we have 
\begin{equation*}
\operatorname{supp}~\delta(Q) = C \cup \{0\}.
\end{equation*}
Let 
\index{A}{L2C@$L^2(C)$}%
$L^2(C) \equiv L^2(C,d\mu)$ be the Hilbert space consisting of
square integrable functions on $C$.
Thus, for a function $\varphi$ on $C$,
the $L^2$-norm of $\varphi$ is given by
\begin{equation}\label{eqn:Cnorm}
\| \varphi \|_{L^2 (C)}
= \frac{1}{2} \int_0^\infty \int_{S^{p-2}} \int_{S^{q-2}}
  |\varphi(r\omega, r\eta)|^2
  r^{p+q-5} dr d\omega d\eta .
\end{equation}
Correspondingly to the coordinates, 
we have an isomorphism of Hilbert spaces:
\begin{equation}\label{eqn:LtwoC}
L^2(\mathbb{R}_+, \frac{1}{2} r^{p+q-5} dr) \hatotimes
L^2(S^{p-2}) \hatotimes L^2(S^{q-2})
\simeq L^2(C).
\end{equation}
Here,
we employ the usual notation
 $\hatotimes$ for the Hilbert completion of the tensor
product. 

If $p+q > 4$, 
then $u\mapsto u\delta(Q)$ defines a continuous,
 injective map from the Hilbert space
$L^2(C)$ into
the space 
\index{A}{S1'@$\mathcal{S}'(\mathbb{R}^n)$}%
$\mathcal S'(\mathbb R^{p+q-2})$ 
of tempered distributions on $\mathbb R^{p+q-2}$:
\begin{equation}\label{def:T}
\index{A}{TL2CS@$T: L^2(C) \to \mathcal{S}'(\mathbb{R}^{n})$}%
  T: L^2(C) \to \mathcal S' (\mathbb R^{\n}), 
  \quad u \mapsto u \delta (Q).
\end{equation}
See \cite[\S 3.4]{xkors3}. 

Now, 
 we are ready to introduce a key map 
 which will give an intertwining operator between the conformal model and the $L^2$-model.  
\begin{equation}
\label{def:B}
\index{A}{TL2C@$\B: L^2(C) \to \overline{V^{p,q}}$|main}%
\B :=(\mathcal{F}_{\mathbb{R}^{p+q-2}}\circ\Tilde{\Psi}^*)^{-1}\circ T = 
(\Tilde{\Psi}^\ast)^{-1}\circ\mathcal{F}_{\mathbb{R}^{p+q-2}}^{-1} \circ T.
\end{equation}

For $u \in C_0^\infty(C)$, $\B u \in C^\infty(M_+)$.
We extend $\B u$ to a function on 
\index{A}{Mz@$M_+$}%
$M_+ \cup (-M_+)$ by
$$
(\B u)(-v) =
(-1)^{\frac{p-q}{2}} (\B u)(v)
\quad (v \in M_+).
$$

We recall from Fact \ref{fact:A} that the inner product on 
 $\overline{V^{p,q}}$ is given by the formula \eqref{eqn:innerM}.  
Then, 
the main ingredient of \cite[Theorem 4.9]{xkors3} can be restated as:
\begin{fact}
\label{fact:B}
$\B$ extends to a unitary operator (up to scalar) from  
$L^2(C)$ onto $\overline{V^{p,q}}$.
\end{fact}

\begin{remNonumber} \label{rem:mathcalT}
The definition \eqref{eqn:Fourier} of the (Euclidean) Fourier transform 
\index{A}{FRn@$\mathcal{F}_{\mathbb{R}^n}$}%
$\mathcal{F}_{\mathbb{R}^n}$
adopted
in this book involves the scalar multiplication by
$(2\pi)^{-\frac{n}{2}}$.
Accordingly, the normalization of\/ $\mathcal{T}$ is different from that
of \cite{xkors3} by a scalar multiplication.
In our normalization, we have
\begin{equation*}
  \|\mathcal{T} u \|^2
  = \frac{1}{2} \|u\|_{L^2(C)}^2
 \quad (u \in L^2(C))
\end{equation*}
as we shall see in \eqref{eqn:Tnorm}.
\end{remNonumber}

Through the unitary operator $\B$,
we can transfer the unitary representation 
 $(\varpi^{p,q}, \overline{V^{p,q}})$ of $G=O(p,q)$ to a unitary 
representation, to be denoted by $\pi^{p,q}$, 
 on the Hilbert space $L^2(C)$ by 
\begin{equation}\label{eqn:piTpi}
 \pi^{p,q}(g):=\B^{-1}\circ \varpi^{p,q}(g)\circ \B, \quad g\in G.
\end{equation}
Hereafter, we shall write 
\index{A}{1pi@$\pi$}%
$\pi$ for $\pi^{p,q}$ for simplicity.
Obviously, $\pi$ is irreducible because so is 
$\varpi^{p,q}$ (see Section \ref{subsec:Kpic}).
We note that the unitary inner product of $\pi$
is nothing but the $L^2$-inner product of $L^2(C)$.
Naming after the classic
Schr\"{o}dinger model $L^2(\mathbb R^n)$ for the 
\index{B}{Weil representation}%
Weil representation
 of the metaplectic group (e.g.\ \cite{xfolland}),
we shall say the resulting irreducible unitary representation 
$(\pi, L^2(C))$ is
the 
\index{B}{Schr\"{o}dinger model|main}%
{\it Schr\"odinger model} 
for the minimal representation of $G=O(p,q)$.

In the philosophy of the orbit method due to Kirillov and Kostant,
the Schr\"{o}dinger model may be regarded as a geometric quantization
of the minimal nilpotent coadjoint orbit 
\index{A}{Omin@$\mathcal{O}_{\min}$}%
 $\mathcal{O}_{\min}$ in
$\mathfrak{g}^*$, the dual of $\mathfrak{g} = \mathfrak{o}(p,q)$.
We note that the isotropic cone $C$ is a Lagrangean
variety of the symplectic manifold 
\index{A}{Omin@$\mathcal{O}_{\min}$}%
$\mathcal{O}_{\min}$.

So far, we have introduced two models of 
$(\varpi^{p,q}, \overline \Vpq)$ and $(\pi, L^2(C))$
for the minimal representation of $G$. 
In the realization of $\overline \Vpq$, 
the $K$-structure is very clear to see, while on 
$L^2(C)$, its $K$-structure is not clear {\it a priori}. 
Generalizing the idea in \cite{xkors3} where we found explicitly an
$L^2$-function (essentially, a $K$-Bessel function) on $C$ belonging
to the minimal $K$-type,
we shall explicitly find 
a formula of $K$-finite vectors of $L^2(C)$ in Chapter \ref{sec:K}.
This computation is carried out by reducing the intertwining operator $\mathcal{T}$
to that of the 
Hankel transform of the $K$-Bessel functions. 

\section{Lie algebra action on $L^2(C)$}
\label{subsec:act}

We continue the review of an easy part of the Schr\"{o}dinger model
$L^2(C)$ of $G=O(p,q)$ from Section \ref{subsec:Sch}.
We shall explain how the Lie algebra $\mathfrak{g}=\mathfrak{o}(p,q)$
acts on smooth vectors of $L^2(C)$.
The action of a maximal parabolic subgroup $\overline{\Pmax}$ on
$L^2(C)$ will be also described.
The main reference of this section is \cite[\S3]{xkors3}.

Our notation here has the following relation with the notation in
Introduction (Sections \ref{subsec:1.1}--\ref{subsec:FC}): 
\begin{align*}
   & p = n_1+1, \ q = n_2+1,
\\
   & n = n_1+n_2 = p+q-2.
\end{align*}

Let $e_0,\dots,e_{p+q-1}$ be the standard basis of $\mathbb{R}^{p+q}$,
on which $G=O(p,q)$ (see \eqref{eqn:Opq}) acts naturally.
First, we define subgroups $\Mmax_+, \Mmax$, 
$K$, and $K'$ of $G$ as follows:
\begin{alignat*}{2}
\index{A}{m0@$m_0$|main}
m_0 &:= -I_{p+q}, 
\\
\index{A}{Mmax+@$\Mmax_+ \simeq O(p-1,q-1)$|main}%
\Mmax_+ &:= \set{g \in G}{g \cdot e_0 = e_0, \ g \cdot
e_{p+q-1} = e_{p+q-1}} &\ \simeq \  & O(p-1,q-1),
\\ 
\index{A}{Mmax@$\Mmax\simeq O(p-1,q-1)\times\protect\mathbb{Z}_2$|main}%
\Mmax &:= \Mmax_+ \cup m_0 \Mmax_+
&\ \simeq \  & O(p-1,q-1) \times \mathbb Z_2, \\
\index{A}{K1Opq@$K \simeq O(p)\times O(q)$|main}%
K &:= G \cap O(p+q) &\  \simeq \  & O(p) \times O(q), \\
\index{A}{K1Opqdash@$K' \simeq O(p-1)\times O(q-1)$|main}%
K'&:= K\cap \Mmax_+ &\ \simeq \ & O(p-1)\times O(q-1).
\end{alignat*}

Then, $K$ is a maximal compact subgroup of $G$ as we already used
in Section \ref{subsec:Kpic},
and $K'$ is a maximal compact subgroup of $\Mmax_+$.
Corresponding to the maximal compact subgroup $K$,
the 
\index{B}{Cartan involution}%
Cartan involution 
\index{A}{1theta@$\theta$}%
$\theta$ of $G$ is given by
$\theta(g) = {}^t \! g^{-1}$,
and its differential (by the same notation) is given by
$\theta(X) = -{}^t\! X$ in the matrix form.

We note that the group $\Mmax_+$ acts on 
the isotropic cone $C$ in $\mathbb{R}^{p+q-2}$ transitively, and leaves 
the measure $d\mu$ (see Section \ref{subsec:Sch}) invariant.

Next we set
\begin{equation}\label{eqn:ej1n}
 \epsilon_j := \begin{cases}
                       1 & (1 \le j \le p-1), 
                 \\
                      -1 & (p \le j \le p+q-2) .
                \end{cases}  
\end{equation}
Let $N_j, \overline N_j~(1 \le j \le p+q-2)$ and $H$ 
be elements of the Lie algebra
$\mathfrak{g} = \mathfrak{o}(p,q)$ given by
\begin{alignat}{2} \label{def:Nj} 
\index{A}{Nj@$N_j$}%
  N_j&:=E_{j, 0} - E_{j, p+q-1} - \epsilon_j E_{0, j} 
            - \epsilon_{j} E_{p+q-1, j}, \\
 \label{def:barNj}
\index{A}{Njbar@$\overline{N}_j$}%
  \overline{N}_j &:= E_{j, 0} + E_{j, p+q-1} - \epsilon_j E_{0, j} 
            + \epsilon_{j} E_{p+q-1, j},  \\ 
 \label{eqn:Edef}
\index{A}{H@$H$}%
 H&:= E_{0, p+q-1} + E_{p+q-1,0}.
\end{alignat}
We note that for $1 \le j \le p+q-2$,
\begin{align} \label{eqn:NNbar}
   \theta(N_j) &= \epsilon_j \overline{N}_j,
\nonumber
\\
   [N_j,\overline{N}_j] &= -2\epsilon_j H.
\end{align}
Then, we define abelian Lie algebras of $\mathfrak{g}$ by 
\begin{alignat*}{2}
   & 
\index{A}{nsmallmax@$\protect\mathfrak{n}^{\protect\max}$|main}%
\nmax &&:= \sum_{j=1}^{p+q-2} \mathbb{R} N_j,
\\
   & 
\index{A}{nsmallmaxbar@$\protect\overline{\protect\mathfrak{n}^{\protect\max}}$|main}%
\overline{\nmax} &&:= \sum_{j=1}^{p+q-2} \mathbb{R} \overline{N}_j,
\\
\index{A}{a-frak@$\protect\mathfrak{a}$|main}%
   & \mathfrak{a} &&:= \mathbb{R} H,
\end{alignat*}
and give coordinates for the corresponding
abelian Lie subgroups
\index{A}{Nmax@$\protect\Nmax$}%
$\Nmax$, 
\index{A}{Nmaxbar@$\overline \Nmax$}%
$\overline \Nmax$ and 
\index{A}{A@$A$}%
$A$ by
\begin{alignat}{3}
  \mathbb R^{p+q-2} \simeq {}&\Nmax,
  \quad &&a=(a_1,a_2, \cdots, a_{p+q-2})
  &&\mapsto 
\index{A}{na@$n_a$}%
n_a:= \exp (\sum_{j=1}^{p+q-2}a_j N_j),  \notag\\
\label{eqn:nbar}
  \mathbb R^{p+q-2} \simeq {}&\overline{\Nmax},  
  \quad &&a=(a_1, a_2, \cdots, a_{p+q-2}) 
  &&\mapsto 
\index{A}{nabar@$\overline{n}_a$}%
\overline{n}_a:=\exp (\sum_{j=1}^{p+q-2} a_j \overline N_j ), \\ \notag
  \mathbb{R} \simeq {}&A, &&t 
  &&\mapsto \exp (t H).
\end{alignat}
Since $M_+^{\max}$ normalizes $\Nmax$,
we have a semidirect product group $M_+^{\max} \Nmax$,
which has the following matrix form:
\begin{equation}
  \label{eqn:MN}
  \Mmax_+ \Nmax = \{g \in G : g(e_0 +e_{p+q-1})
                = e_0 + e_{p+q-1}\}.
\end{equation}
On the other hand,
the natural action of $G$ on $\mathbb{R}^{p+q}$ induces a transitive
action on 
\begin{equation}
  \label{eqn:GMN}
\index{A}{C0@$\widetilde{C}$|main}%
  \widetilde{C} := \{(\xzeta_0,\cdots,\xzeta_{p+q-1}) \in 
  \mathbb{R}^{p+q} \setminus \{0\}:
  \sum^{p-1}_{j=0} \xzeta^2_j - \sum^{p+q-1}_{j=p} \xzeta^2_j = 0 \} ,
\end{equation}
with the isotropy subgroup $M_+^{\max} \Nmax$ at
$e_0 + e_{p+q-1}$.
Thus, we get
a diffeomorphism
\begin{equation*}
  G/\Mmax_+ \Nmax \simeq \widetilde{C}.
\end{equation*}

We define a maximal parabolic subgroup by
\begin{equation}
\label{eqn:Pmax}
\index{A}{Pmax=@$\overline \Pmax:= \Mmax A \overline \Nmax$|main}%
   \overline \Pmax
:= \Mmax A \overline \Nmax.
\end{equation}
Then the action of $\overline{\Pmax}$
on $L^2(C)$ is described explicitly as follows (\cite[\S 3.3]{xkors3}):
\begin{alignat}{2}
\label{eqn:rM}
   (\pi(m) \psi)(\xzeta)
   &= \psi(\trans{m} \xzeta)
   && \qquad (m \in \Mmax_+),
\\
\label{eqn:rm}
 (\pi(m_0) \psi)(\xzeta)
   &= (-1)^{\frac{p-q}2}\psi(\xzeta),
   &&
\\
\label{eqn:rA}
   (\pi(e^{t H}) \psi)(\xzeta)
   &= e^{-\frac{p+q-4}2 t} \psi(e^{-t} \xzeta)
   && \qquad (t \in \mathbb R),
\\
\label{eqn:rN}
   (\pi(\nbar{a}) \psi)(\xzeta)
   &= e^{2\sqrt{-1} 
  (a_1 \xzeta_1 + \dots + a_{p+q-2} \xzeta_{p+q-2})} \psi(\xzeta)
   && \qquad (a \in \mathbb R^{p+q-2}).
\end{alignat}

For $g\notin\overline{\Pmax}$,
the action $\pi(g)$ on $L^2(C)$ was not given in \cite{xkors3}.
(In fact, the formula $\pi(g)$ for general $g\in G$ is the main issue
of the book.)
Instead, we obtaind a formula for the differential action of the Lie
algebra $\mathfrak{g}$ in \cite{xkors3}.
To state the formula,
we recall some general terminology for infinite dimensional
representations adapted to our special setting. 
\begin{definition} \label{def:smoothvec}
Let $\psi \in L^2(C)$.
We say $\psi$
is a 
\index{B}{differentiable vector|main}%
\emph{differentiable vector} 
if
\begin{equation*}
   d\pi(X)\psi := \lim_{t\to 0}
                  \frac{\pi(e^{tX})\psi-\psi}{t}
\end{equation*}
exists for any $X\in\mathfrak{g}$.
Iterating this process,
we say $\psi$ is a 
\index{B}{smooth vector|main}%
\emph{smooth vector} 
if
$d\pi(X_1)\cdots d\pi(X_k)\psi$ is a differentiable vector for any 
$k\ge 1$ and $X_1,\dots,X_k \in \mathfrak{g}$.
(The notion of smooth vectors is defined for continuous
representations on complete, locally convex topological vector spaces.)

We say $\psi$ is $K$-\emph{finite} if
\begin{equation*}
   \dim_{\mathbb{C}} \text{$\mathbb{C}$-\textup{span}}
   \{ \pi(k)\psi: k \in K \} < \infty.
\end{equation*}
\end{definition}
Let
\index{A}{L2Cinfty@$L^2(C)^\infty$|main}%
 $L^2(C)^\infty$ 
be the space of smooth vectors of the unitary
representation $(\pi,L^2(C))$ of $G$,
and 
\index{A}{L2CK@$L^2(C)_K$}%
$L^2(C)_K$ the space of $K$-finite vectors of $L^2(C)$.
Applying the general theory of representations of real reductive Lie groups (see
\cite{xKnapp,xwal}, for example) to our irreducible unitary
representation $(\pi,L^2(C))$,
we have
\begin{enumerate}
    \renewcommand{\labelenumi}{\upshape \theenumi)}
\item  
$L^2(C)^\infty$ has naturally a Fr\'{e}chet topology,
on which both the group $G$ and the Lie algebra $\mathfrak{g}$ act
continuously.
\item  
$L^2(C)_K \subset L^2(C)^\infty \subset L^2(C)$.
Moreover, $L^2(C)_K$ is dense in $L^2(C)^\infty$,
and $L^2(C)^\infty$ is dense in $L^2(C)$, in each topology.
\end{enumerate}

{}From now, we shall use the notation which is compatible with the
notation of Chapter \ref{sec:intro}:
\begin{equation*}
  n = n_1+n_2 = (p-1)+(q-1).
\end{equation*}
It follows from \eqref{eqn:rA} and \eqref{eqn:rN} that 
the differential action of $A$ and $\overline{\Nmax}$ 
on $L^2(C)^\infty$ is given as 
\begin{align}
& d\pi(H) = -(E+\frac{n-2}{2}),
\label{eqn:dpiH}
\\
& d\pi(\overline{N}_j)=2\sqrt{-1}\xzeta_j 
 \quad (1\le j \le n).
\label{eqn:barNj}
\end{align}
Here, 
\index{A}{E@$E$}%
$E := \sum_{i=1}^n x_i \frac{\partial}{\partial x_i}$
is the 
\index{B}{Euler operator}%
Euler operator,
and $x_j$ in \eqref{eqn:barNj} is the multiplication operator by the
coordinate function $x_j$.

On the other hand, 
the differential action of $\Nmax$ on $L^2(C)^\infty$ 
is more complicated.
In \cite[Lemma 3.2]{xkors3}, we gave its formula
 by means of second order differential
operators in the ambient space $\mathbb{R}^n=\mathbb{R}^{p+q-2}$ via 
the inclusion map 
$T: L^2(C) \hookrightarrow \mathcal S'(\mathbb R^{n})$,
\thinspace
$u \mapsto u\delta(Q)$ (see \eqref{def:T}) as follows.
Let $D_j$ be a differential operator on 
$\mathbb R^n$ 
(in the notation loc.\ cit.,
$D_j = d\hat \varpi_{\lambda, \epsilon}(N_j)$  with 
$\lambda =\frac{p+q-4}2$) given by:
\begin{equation}\label{eqn:Dj}
  D_j = \sqrt{-1} \biggl(-\frac{n-2}2 
          \epsilon_j \frac{\partial}{\partial \xzeta_j}
       - \bigl(\sum_{k=1}^{n}
          \xzeta_k \frac{\partial}{\partial \xzeta_k}\bigr) 
         \epsilon_j \frac{\partial }{\partial \xzeta_j}+\frac{1}2 \xzeta_j 
         \bigl(\sum_{k=1}^{n}\epsilon_k 
         \frac{\partial^2}{\partial \xzeta_k^2}\bigr) \biggr).  
\end{equation}
Then the differential action $d\pi(N_j)~(1 \le j \le n)$ 
is characterized
by the commutative diagram:
\begin{equation}
\label{diag:D}
\begin{matrix}
   ~~~~~~L^2(C)^\infty   &  \stackrel{T}\longrightarrow & \mathcal S'(\mathbb R^{n})   \\
  d\pi(N_j) \Big\downarrow  &       &   \Big\downarrow D_j  \\
  ~~~~~~L^2(C)^\infty   &  \stackrel{T}\longrightarrow   &  \mathcal S'(\mathbb R^{n}).
\end{matrix}
\end{equation}

In Section \ref{subsec:Pjb}, 
we shall treat these differential operators more systematically by
introducing the following
differential operators $P_j(b)$ $(1 \le j \le n)$ with
complex parameter $b$ by
\begin{equation} \label{eqn:Pjb}
\index{A}{Pjb@$P_j(b)$|main}%
   P_j(b) := \epsilon_j\xzeta_j \square - (2E+n-2b)
              \frac{\partial}{\partial\xzeta_j},
\end{equation}
where we set
\index{A}{2square@\par\indexspace$\square$}%
$ \square = \sum_{j=1}^n \epsilon_j
              \frac{\partial^2}{\partial\xzeta_j^2}$.
By definition, we have
\begin{equation*}
  D_j = \frac{\sqrt{-1}}{2} \epsilon_j P_j(-1)
  \qquad (1\le j \le n).
\end{equation*}

In Theorem \ref{thm:Pj},
we shall see that the differential operator $P_j:=P_j(1)$ is tangential to $C$,
 and that the differential action $d\pi(N_j)$ is given
as
\begin{equation} \label{eqn:NjPj}
  d\pi(N_j) = \frac{\sqrt{-1}}{2} \epsilon_j P_j |_C
  \qquad (1 \le j \le n).
\end{equation}

\section{Commuting differential operators on $C$}
\label{subsec:Pjb}

In this section, we investigate basic properties of the differential operators
\begin{equation*}
  P_j(b) 
   = \epsilon_j x_j \square - (2E+n-2b) \frac{\partial}{\partial x_j}
  \qquad (1 \le j \le n)
\end{equation*}
on $\mathbb{R}^n$ introduced in \eqref{eqn:Pjb},
and in particular, explain why and how the differential operators
$P_j(-1)$ $(1\le j \le n)$ (see \eqref{eqn:Pjb}) induce the
differential operators $P_j(1)$ along the isotropic cone $C$.

Again we recall
\begin{equation*}
   n = n_1+n_2 = (p-1)+(q-1).
\end{equation*}
We also recall from Section \ref{subsec:1.1} that
$\mathbb{R}[x,\frac{\partial}{\partial x}] \equiv
 \mathbb{R}[x_1,\dots,x_n,\frac{\partial}{\partial x_1},\dots,
            \frac{\partial}{\partial x_n}]
$
is the 
\index{B}{Weyl algebra}%
Weyl algebra 
and that $\mathbb{R}[x,\frac{\partial}{\partial x}]^C$
is the subalgebra consisting of differential operators tangential to $C$.

The main result of this section is the following theorem concerning
with the 
\index{B}{fundamental differential operator|main}%
\textit{fundamental differential operators} 
on the isotropic cone defined by
\begin{equation} \label{eqn:Pj}
\index{A}{P_j@$P_j$}%
  P_j := P_j(1)
  \quad (1 \le j \le n).
\end{equation}

\begin{theorem} \label{thm:Pj}
\begin{enumerate}
    \renewcommand{\labelenumi}{\upshape \theenumi)}
\item  
$P_i P_j = P_j P_i$ $(1 \le i,j \le n)$.
\item  
$P_j \in \mathbb{R} [x,\frac{\partial}{\partial x}]^C$
$(1 \le j \le n)$.
\item  
The Lie algebra generated by
\begin{equation*}
  [P_i, x_j] = P_i x_j - x_j P_i
  \quad (1 \le i, j \le n)
\end{equation*}
in the Weyl algebra $\mathbb{R}[x,\frac{\partial}{\partial x}]$ is
isomorphic to $\mathfrak{o}(p-1,q-1)+\mathbb{R}$,
the Lie algebra of the differential action of the conformal linear
transformation group 
\index{A}{CO(Q)@$CO(Q)$}%
$CO(Q)$.
\item  
$P_j |_C$ extends to a self-adjoint operator on $L^2(C)$.
\item  
$\bigl( \displaystyle\sum_{j=1}^n \epsilon_j P_j^2 \bigr) \big|_C = 0$.
\end{enumerate}
\end{theorem}

We shall give a proof of Theorem \ref{thm:Pj} in the following order:
(2), (4), (5), (1), and (3).
An important step for the proof is the following:
\begin{proposition}\label{prop:Pj}
For $u \in C_0^\infty(\mathbb{R}^n\setminus \{0 \})$,
we have:
\begin{enumerate}
    \renewcommand{\labelenumi}{\upshape \theenumi)}
\item  
$P_j(-1) (u \delta(Q)) = (P_j(1)u) \delta(Q)$.

Hence, we have $P_j(-1) \circ T = T \circ P_j(1)$.

\item  
$P_j(1) (uQ) = (P_j(-1)u)Q$.

In particular, $P_j(1)(uQ)|_C = 0$.
\end{enumerate}
\end{proposition}

Admitting Proposition \ref{prop:Pj} for a while,
we give a proof of Theorem \ref{thm:Pj} (2) and (4).

\begin{proof}[Proof of Theorem \ref{thm:Pj} (2)]
For a smooth function $\psi$ defined in an open subset $V$ of $C$,
we extend it to a smooth function $\tilde{\psi}$ in an open subset of 
$\mathbb{R}^n\setminus\{0\}$.
Then,
$P_j(1)\tilde{\psi}|_V$
is independent of the choice of the extension $\tilde{\psi}$,
and is determined by $\psi=\tilde{\psi}|_V$.
In fact, let $\tilde{\psi}_1$, $\tilde{\psi}_2$ be extensions of
$\psi$.
Since $(\tilde{\psi}_1 - \tilde{\psi}_2) |_C = 0$,
we find a smooth function locally defined in a neighborhood of $V$
such that
\begin{equation*}
   \tilde{\psi}_1 - \tilde{\psi}_2 = uQ.
\end{equation*}
Then
$(P_j\tilde{\psi}_1 - P_j\tilde{\psi}_2) |_C
 = (P_j(-1)u)Q |_C = 0
$
by Proposition \ref{prop:Pj} (2).
Therefore, $P_j\tilde{\psi}_1 |_C = P_j\tilde{\psi}_2 |_C$.
Thus, we have a well-defined map
\begin{equation*}
   C^\infty(V) \to C^\infty(V),
   \quad \psi \mapsto P_j(1)\tilde{\psi}|_V.
\end{equation*}
Since this is a sheaf morphism,
it is given by a differential operator on $C$.
Hence, $P_j = P_j(1)$ is tangential to $C$.
Therefore, Theorem \ref{thm:Pj} (2) is shown.
\end{proof}

By Proposition \ref{prop:Pj},
we also get the self-adjointness of $P_j(1)|_C$ as follows.

\begin{proof}[Proof of Theorem \ref{thm:Pj} (4)]
By \eqref{diag:D},
the differential action $d\pi(N_j)$ is characterized by the relation
\begin{equation*}
   D_j \circ T = T \circ d\pi(N_j).
\end{equation*}
On the other hand,
we have proved in Proposition \ref{prop:Pj} (1)
\begin{equation*}
   P_j(-1) \circ T = T \circ (P_j(1)|_C).
\end{equation*}
Since $D_j = \frac{\sqrt{-1}}{2} \epsilon_j P_j(-1)$
by definition,
we get
\begin{equation*}
   d\pi(N_j) = \frac{\sqrt{-1}}{2} \epsilon_j P_j(1) |_C.
\end{equation*}
As $(\pi,L^2(C))$ is a unitary representation of the Lie group
$G = O(p,q)$,
the differential action 
$\sqrt{-1} d\pi(X)$ on $L^2(C)^\infty$ 
extends to a self-adjoint operator on $L^2(C)$ for
any $X \in \mathfrak{g} = \mathfrak{o}(p,q)$.
Hence, $P_j = P_j(1)$ is self-adjoint.
Thus, Theorem \ref{thm:Pj} (4) is proved.
\end{proof}

Before giving a proof of Proposition \ref{prop:Pj},
we set up the notation of the meromorphic continuation of the
distribution $Q(x)_+^\lambda$.
For a complex parameter $\lambda$ with $\operatorname{Re}\lambda>-1$, 
we define a distribution on $\mathbb{R}^n$ by
\begin{equation*}
\index{A}{Qxlambda@$Q(\xzeta)_+^\lambda$|main}%
   Q(\xzeta)_+^\lambda
   := \begin{cases}
         Q(\xzeta)^\lambda  &\text{if $Q(\xzeta)>0$}, \\
         0                 &\text{if $Q(\xzeta)\le0$}.
      \end{cases}
\end{equation*}
Then, $Q(x)_+^\lambda$
 continues as a distribution depending meromorphically on the parameter 
$\lambda \in \mathbb{C}$ 
(see Appendix \ref{subsec:Riesz}).
In particular, as a distribution on $\mathbb{R}^n\setminus\{0\}$,
$Q(\xzeta)_+^\lambda$ has only simple poles located at
$\lambda=-1,-2,-3,\dotsc$.
Since the gamma function $\Gamma(\lambda+1)$ has simple poles exactly
 at the same places $\lambda=-1,-2,\dotsc$,
we see that 
$\frac{Q(\xzeta)_+^\lambda}{\Gamma(\lambda+1)}$
defines a distribution on $\mathbb{R}^n\setminus\{0\}$ 
depending holomorphically on $\lambda\in\mathbb{C}$.

In view of the residue formula \eqref{eqn:resxlmd}
 in Appendix \ref{subsec:Riesz} for the one variable
case
$\frac{t_+^\lambda}{\Gamma(\lambda+1)}$,
we have the following formula of generalized functions on 
$\mathbb{R}^n \setminus \{0\}$ 
(see \cite[Chapter III, \S2.2]{xGeSh}):
\begin{equation*}
   \delta(Q(\xzeta))
   = \frac{Q(\xzeta)_+^\lambda}{\Gamma(\lambda+1)}
     \,\Big|_{\lambda=-1}.
\end{equation*}
Therefore, the map 
\index{A}{TL2CS@$T: L^2(C) \to \mathcal{S}'(\mathbb{R}^{n})$}%
$T: L^2(C) \to \mathcal{S}'(\mathbb{R}^n)$
defined in \eqref{def:T} has the following expression:
\begin{equation}
\label{eqn:TuQ}
   T(u|_C) 
   = \frac{uQ(\xzeta)_+^\lambda}{\Gamma(\lambda+1)}
     \,\Big|_{\lambda=-1}
\end{equation}
for $u\in C_0^\infty(\mathbb{R}^n\setminus\{0\})$.
The proof of Proposition \ref{prop:Pj} will make use of
\eqref{eqn:TuQ}. 
Along this line, we prepare:
\begin{lemSec}\label{lem:fQ}
For $u\in C_0^\infty(\mathbb{R}^n\setminus\{0\})$,
\begin{equation*}
\index{A}{2square@\par\indexspace$\square$}%
   \square(uQ_+^\lambda)
   =  (\square u)Q_+^\lambda + 4\lambda(Eu)Q_+^{\lambda-1}
       + 2\lambda(2\lambda+n-2)uQ_+^{\lambda-1}.
\end{equation*}
\end{lemSec}
Similar formulas also hold if we replace $Q_+$ by $Q$ for positive
integers $\lambda$.
In particular, letting $\lambda=1$, we have
\begin{equation*}
   \square(uQ) 
   = (\square u)Q + (4E+2n)u.
\end{equation*}

\begin{proof}
By the Leibniz rule,
we have
\begin{alignat}{2}
  & \frac{\partial}{\partial\xzeta_j}(uQ_+^\lambda)
    ={} &&\frac{\partial u}{\partial\xzeta_j}Q_+^\lambda
          + 2\lambda\epsilon_j\xzeta_j uQ_+^{\lambda-1},
\label{eqn:fQ1}
\\
  & \frac{\partial^2}{\partial\xzeta_j}(uQ_+^\lambda)
    ={}&& \frac{\partial^2 u}{\partial\xzeta_j^2}Q_+^\lambda
         + 4\lambda\epsilon_j\xzeta_j \frac{\partial u}{\partial\xzeta_j}
         Q_+^{\lambda-1}
\nonumber
\\
  &    &&{}+ 2\lambda\epsilon_j uQ_+^{\lambda-1}
         + 4\lambda(\lambda-1)\xzeta_j^2 uQ_+^{\lambda-2}.
\label{eqn:fQ2}
\end{alignat}
Summing up \eqref{eqn:fQ2} multiplied by the signature $\epsilon_j$ over $j$,
we get Lemma.
\end{proof}

Here is a key formula for the proof of Proposition \ref{prop:Pj}:
\begin{lemSec}\label{lem:Pjkey}
For $u\in C_0^\infty(\mathbb{R}^n\setminus\{0\})$,
we have the following identity as distributions on
$\mathbb{R}^n\setminus\{0\}$:
\begin{equation} \label{eqn:keyPj}
\index{A}{Pjb@$P_j(b)$}%
  P_j(b)(uQ_+^\lambda)
  = \bigl(P_j(b-2\lambda)u\bigr)Q_+^\lambda
    - 4\lambda(\lambda-b)\epsilon_j\xzeta_j u Q_+^{\lambda-1}
\end{equation}
for any $1\le j\le n$, $b\in\mathbb{C}$, and
$\lambda\in\mathbb{C}\setminus\{-1,-2,\dotsc\}$.
\end{lemSec}

\begin{proof}
By \eqref{eqn:fQ1},
we have
\begin{align*}
   (2E+n-2b) \frac{\partial}{\partial x_j} (uQ_+^\lambda)
 ={}& \Bigl( 
     \bigl( 2E+n-2(b-2\lambda) \bigr) 
     \frac{\partial u}{\partial x_j} 
     \Bigr) Q_+^\lambda
\\
  &+ 2\lambda (4\lambda+n-2-2b)\epsilon_j x_j u Q_+^{\lambda-1}
\\
  &+ 4\lambda \epsilon_j x_j(Eu) Q_+^{\lambda-1}.
\end{align*}
Then, \eqref{eqn:keyPj} follows from Lemma \ref{lem:fQ}.
\end{proof}

Now we are ready to prove Proposition \ref{prop:Pj}.

\begin{proof}[Proof of Proposition \ref{prop:Pj}]
Since $\frac{Q_+^\lambda}{\Gamma(\lambda+1)}$
is a distribution on $\mathbb{R}^n\setminus\{0\}$ with 
parameter holomorphically dependent on $\lambda\in\mathbb{C}$,
we have
\begin{equation*}
   P_j(b) \Bigl( u\frac{Q_+^\lambda}{\Gamma(\lambda+1)} \Bigr)
   = \bigl( P_j(b-2\lambda)u \bigr)
     \frac{Q_+^\lambda}{\Gamma(\lambda+1)}
     - (\lambda-b)\epsilon_j\xzeta_j u \frac{Q_+^{\lambda-1}}{\Gamma(\lambda)}
\end{equation*}
for all $\lambda \in \mathbb{C}$.
By letting $b=\lambda$, we have
\begin{equation*}
   P_j(\lambda) \Bigl( u\frac{Q_+^\lambda}{\Gamma(\lambda+1)} \Bigr)
   = \bigl( P_j(-\lambda)u \bigr)
     \frac{Q_+^\lambda}{\Gamma(\lambda+1)}.
\end{equation*}
Further, by letting $\lambda=-1$,
we get
\begin{equation*}
  P_j(-1) \bigl( u\delta(Q) \bigr)
  = \bigl( P_j(1)u \bigr) \delta(Q).
\end{equation*}
By letting $\lambda=1$,
we get
\begin{equation*}
   P_j(1) (uQ) = (P_j(-1)u) Q.
\end{equation*}
Therefore, we have proved  Proposition \ref{prop:Pj}.
\end{proof}

Thus, the proof of Theorem \ref{thm:Pj} (2) is completed.
Next, let us prove the following:
\begin{proposition}\label{prop:relPj}
$\displaystyle\sum_{j=1}^n \epsilon_jP_j(1)^2 = Q \,
\index{A}{2square@\par\indexspace$\square$}%
\square^2$.
\end{proposition}
Admitting Proposition \ref{prop:relPj} for the time being,
we give a proof of Theorem \ref{thm:Pj} (5).

\begin{proof}[Proof of Theorem \ref{thm:Pj} (5)]
We have already shown that
$P_j(1) \in \mathbb{R} [ x,\frac{\partial}{\partial x} ]^C$
$(1 \le j \le n)$.
Then, Theorem \ref{thm:Pj} (5)
 is an immediate consequence of the following identity 
in the algebra
$\mathbb{R} [ x,\frac{\partial}{\partial x} ]^C$.
\end{proof}

To prove Proposition \ref{prop:relPj},
we list some basic relations of the Lie bracket
$[A,B] = AB-BA$ in the Weyl algebra
$\mathbb{R} [x,\frac{\partial}{\partial x}]$:
\begin{claim}\label{claim:Drel}
\quad
\begin{enumerate}
    \renewcommand{\labelenumi}{\upshape \theenumi)}
\item  
$[E,\xzeta_j] = \xzeta_j$.
\item  
$\displaystyle [E,\frac{\partial}{\partial\xzeta_j}]
 = -\frac{\partial}{\partial\xzeta_j}$.
\item  
$[E,\square] = -2 \,\square$.
\item  
$\displaystyle [\square,\xzeta_j] 
 = 2\epsilon_j\frac{\partial}{\partial\xzeta_j}$.
\end{enumerate}
\end{claim}
\noindent
Here, $\xzeta_j$ denotes the multiplication operator by $\xzeta_j$.

\begin{proof}
Straightforward by the Leibniz rule.
\end{proof}

\begin{proof}[Proof of Proposition \ref{prop:relPj}]
In light of the definition
\begin{equation*}
   P_i(b) = \epsilon_i\xzeta_i \square
            - (2E+n-2b)\frac{\partial}{\partial\xzeta_i},
\end{equation*}
we develop $P_i(b)P_j(b)$ as
\begin{equation*}
   P_i(b) P_j(b) 
   = (\textrm{I})+(\textrm{II})+(\textrm{III})+(\textrm{IV}),
\end{equation*}
where $(\textrm{I})\equiv(\textrm{I})_{ij}(b)$ is given by
\begin{align*}
   (\textrm{I})
   :={}& (\epsilon_i\xzeta_i\square)(\epsilon_j\xzeta_j\square)
 \\
    ={}& \epsilon_i\epsilon_j\xzeta_i\xzeta_j\square^2 + 2\epsilon_i \xzeta_i
         \frac{\partial}{\partial\xzeta_j} \,\square,
 \\
\intertext{and similarly,}
   (\textrm{II})
   :={}& -(\epsilon_i\xzeta_i\square)
         \Bigl( (2E+n-2b)\frac{\partial}{\partial\xzeta_j}
         \Bigr) 
 \\
    ={}& -\epsilon_i(2E+n-2b+2)
         \xzeta_i\frac{\partial}{\partial\xzeta_j} \,\square,
 \\
   (\textrm{III})
   :={}& - \Bigl( (2E+n-2b)\frac{\partial}{\partial\xzeta_i} 
           \Bigr)
           (\epsilon_j\xzeta_j\square)
 \\
    ={}& -\epsilon_j(2E+n-2b)
         \Bigl( \xzeta_j\frac{\partial}{\partial\xzeta_i}+\delta_{ij}
         \Bigr)
         \square,
 \\
   (\textrm{IV})
   :={}& \Bigl( (2E+n-2b)\frac{\partial}{\partial\xzeta_i}
         \Bigr)
         \Bigl( (2E+n-2b)\frac{\partial}{\partial\xzeta_j}
         \Bigr)
 \\
    ={}& 
         \bigl( 4E^2+4(n+1-2b)E+(n-2b)(n+2-2b) \bigr)
         \frac{\partial^2}{\partial\xzeta_i\partial\xzeta_j}.
\end{align*}
Here, $\delta_{ij}$ stands for Kronecker's delta.
Now, we take $i=j$ and $b=1$,
and sum up these terms over $j$:
\begin{alignat*}{2}
\allowdisplaybreaks
   & \sum_{j=1}^n \epsilon_j(\textrm{I})_{jj}(1)
   &&{} = Q \,\square^2 + 2E \,\square,
 \\
   & \sum_{j=1}^n \epsilon_j(\textrm{II})_{jj}(1)
   &&{} = -(2E+n)E \,\square,
 \\
   & \sum_{j=1}^n \epsilon_j(\textrm{III})_{jj}(1)
   &&{} = -(2E+n-2)(E+n) \,\square,
 \\
   & \sum_{j=1}^n \epsilon_j(\textrm{IV})_{jj}(1)
   &&{} = \bigl( 4E^2+4(n-1)E+(n-2)n \bigr) \,\square.
\end{alignat*}
Summing up these four equations, we get Proposition \ref{prop:relPj}.
\end{proof}

\medskip

It is easy to see that
the formulas for \textrm{(I)}, \textrm{(II)}, \textrm{(III)}, 
and \textrm{(IV)} used in the previous proof 
lead us also to the following:
\begin{proposition}\label{prop:Pij}
As differential operators on\/ $\mathbb{R}^n$, we have
\begin{equation*}
\index{A}{Pjb@$P_j(b)$}%
   [P_i(b), P_j(b)] = 0
\end{equation*}
for any $1\le i,j \le n$ and $b\in \mathbb{C}$.
\end{proposition}
As a special case of Proposition \ref{prop:Pij},
we have:
\begin{proof}[Proof of Theorem \ref{thm:Pj} (1)]
This follows from Proposition \ref{prop:Pij} by letting $b=1$.
\end{proof}

\medskip

Finally, let us prove Theorem \ref{thm:Pj} (3).
\begin{proof}[Proof of Theorem \ref{thm:Pj} (3)]
We continue the notation \eqref{eqn:ej1n} for $\epsilon_j = \pm1$
$(1 \le j \le n = p+q-2)$. 

Let $E_{ij}$ $(1 \le i, j \le n)$ be the matrix unit.
Then
\[
   \epsilon_i \epsilon_j E_{ij} - E_{ji}
\quad
   (1 \le i < j \le n)
\]
forms a basis of the Lie algebra $\mathfrak{o}(p-1,q-1)$ of the Lie group
 $\Mmax_+ \simeq O(p-1,q-1)$. 
Correspondingly, the natural differential action on
$\mathbb{R}^n = \mathbb{R}^{(p-1)+(q-1)}$ gives a vector field
\begin{equation*}
   X_{ij} := \epsilon_i \epsilon_j 
             x_i \frac{\partial}{\partial x_j} -
             x_j \frac{\partial}{\partial x_i}
  \qquad (1 \le i < j \le n).
\end{equation*}
Likewise, the differential of the dilation is given by the 
\index{B}{Euler vector field}%
Euler vector field
\begin{equation*}
\index{A}{E@$E$}%
   E = \sum_{i=1}^n x_i \frac{\partial}{\partial x_i}.
\end{equation*}
Hence, Theorem \ref{thm:Pj} (3) is an immediate consequence of the
following lemma.
\end{proof}
\begin{lemSec} \label{lem:Pixj}
As differential operators on $\mathbb{R}^n$,
we have
\begin{enumerate}
    \renewcommand{\labelenumi}{\upshape \theenumi)}
\item  
$[P_i,x_j] = X_{ij}$ \quad $(1 \le i < j \le n)$.
\item  
$[P_i,x_i] = -2(E+n-2)$ \quad $(1 \le i \le n)$.
\end{enumerate}
\end{lemSec}

\begin{proof}
By the definition \eqref{eqn:Pjb} of $P_j = P_j(1)$, we have
\begin{align*}
   [P_i,x_j]
   ={}& \epsilon_i x_i [ \, \square,x_j]
    - (2E+n-2) \Bigl[ \frac{\partial}{\partial x_i}, x_j \Bigr]
\\
   & - [2E+n-2, x_j] \frac{\partial}{\partial x_i}.
\\
\intertext{By Claim \ref{claim:Drel}, we have}
   ={}& 2\epsilon_i \epsilon_j x_i \frac{\partial}{\partial x_j}
    - \delta_{ij} (2E+n-2)
    - 2x_j \frac{\partial}{\partial x_i}
\\
   ={}& 2X_{ij} - \delta_{ij} (2E+n-2).
\end{align*}
Hence, Lemma is proved.
\end{proof}

Hence, the proof of Theorem \ref{thm:Pj} is completed.

\begin{remNonumber} \label{rem:Pjout}
Let $R$ be the subalgebra of\/ 
\index{A}{RxC@$\mathbb{R}[x,\frac{\partial}{\partial x}]^C$}%
$\mathbb{R}[x,\frac{\partial}{\partial x}]^C$ 
generated by $x_k$, $E$, and $X_{ij}$ 
$(1 \le k \le n, \, 1 \le i<j \le n)$.
Then, 
\index{A}{P_j@$P_j$}%
$P_j \in \mathbb{R}[x,\frac{\partial}{\partial x}]^C$ but
$P_j \notin R$.

To see this, we say an element
$P$ of the Weyl algebra $\mathbb{R} [x,\frac{\partial}{\partial x}]$
is of degree $k$ if
\begin{equation*}
   [E,P] = kP.
\end{equation*}
If $P,Q$ are of degree $k,l$, respectively,
then $PQ$ is of degree $k+l$ because
\begin{equation*}
   [E,PQ] = [E,P]Q + P[E,Q].
\end{equation*}
Since $x_k$ is of degree $1$,
and $E$ and $X_{ij}$ is of degree $0$,
any element of $R$ is expressed as a linear combination of operators of
non-negative degrees.
Since $P_j$ is of degree $-1$,
we conclude $P_j \notin R$.
\end{remNonumber}

\begin{remNonumber} \label{rem:Pn2}
Our concern here is with the case $n>2$ $($i.e.\ $p+q>4$$)$.
Let us examine the degenerate case where
 $n=2$ and\/ 
$n_1=n_2=1$ $($i.e.\ $p=q=2$$)$. 
Then, the polynomial $Q(x)=x_1^2-x_2^2$ is not irreducible,
and the differential operators $P_1$ and $P_2$ have the following
formulas: 
\begin{align*}
   &P_1 + P_2 = -(x_1+x_2)
                \Bigl(\frac{\partial}{\partial x_1} +
                 \frac{\partial}{\partial x_2}\Bigr)^2,
\\
   &P_1 - P_2 = -(x_1-x_2)
                \Bigl(\frac{\partial}{\partial x_1} -
                 \frac{\partial}{\partial x_2}\Bigr)^2.
\end{align*}
By the change of variables
\begin{equation*}
   y_1 = x_1+x_2, \quad y_2 = x_1-x_2,
\end{equation*}
the isotropic cone $C$ is given by
\begin{equation*}
   \{ (y_1,y_2) \in \mathbb{R}^2 \backslash \{0\}:
       \ \text{$y_1 = 0$ or\/ $y_2 = 0$} \},
\end{equation*}
and we have
\begin{align*}
   &P_1 + P_2 = -4y_1 \frac{\partial^2}{\partial y_1^2},
\\
   &P_1 - P_2 = -4y_2 \frac{\partial^2}{\partial y_2^2}.
\end{align*}
\end{remNonumber}

\section{The unitary inversion operator
  $\mathcal{F}_C=\pi(w_0)$} 
\label{subsec:w}

In this section, we list some important features of the element 
$$
\index{A}{w0@$w_0$}%
w_0 =\begin{pmatrix}  I_p  &  0  \\  0    & -I_q      \end{pmatrix}
\in G,
$$
and investigate key properties of the corresponding unitary operator
\begin{equation*}
\index{A}{FC@$\mathcal{F}_C$}%
   \mathcal{F}_C := \pi(w_0)
\end{equation*}
on $L^2(C)$.

\begin{enumerate}[{\bf {I}1}]
\addtocounter{enumi}{-1}
\index{A}{1zpropertiesI0-I5@\textbf{I0}--\textbf{I5}|(}%
\item
\label{item:I0}
({\it Order two})\enspace
Obviously, the element $w_0$ is of order two.
Therefore, 
$\mathcal{F}_C^2 = \operatorname{id}$
on $L^2(C)$.

\item
\label{item:I01}
(\textit{Cartan involution})\enspace
The 
\index{B}{Cartan involution}%
Cartan involution 
\index{A}{1theta@$\theta$}%
$\theta(g) = {}^t \! g^{-1}$ is given as the
conjugation by $w_0$:
\begin{equation} \label{eqn:Cartanw}
   \theta(g) = w_0  g  w_0^{-1}
\end{equation}
because ${}^t\! g w_0 g = w_0$ for $g \in G$.

\item
\label{item:I4}
({\it Center of $K$})\enspace
$w_0$ lies in the center of $K$.
This is obvious from \eqref{eqn:Cartanw}.
It also follows directly from the definition of $w_0$ in the matrix form.

\item
\label{item:I1}
({\it Bruhat decomposition})\enspace
Retain the notation as in Section \ref{subsec:act}.
Then, 
\begin{equation}\label{eqn:wa}
\index{A}{H@$H$}%
\operatorname{Ad}(w_0)H=-H,
\end{equation}
and therefore
$
 \operatorname{Ad}(w_0) \vert_{\mathfrak a}=-\operatorname{id}
$. 
We also see readily from
\eqref{eqn:NNbar} and \eqref{eqn:Cartanw} that
\begin{equation}\label{eqn:Adwj}
 \operatorname{Ad}(w_0)\overline{N}_j=
 \epsilon_j N_j \quad (1 \le j \le p+q-2),
\end{equation}
and therefore
\begin{equation}\label{eqn:Adw}
 \operatorname{Ad}(w_0)
\index{A}{nsmallmax@$\mathfrak {n}^{\max }$}%
\nmax=
\index{A}{nsmallmaxbar@$\overline {\mathfrak {n}^{\max }}$}%
\overline{\nmax}.
\end{equation}
The Gelfand--Naimark decomposition
\begin{equation*}
   \mathfrak{g} = (\mmax + \mathfrak{a} + \overline{\nmax}) + \nmax
                = \overline{\pmax} + \operatorname{Ad}(w_0)
                  \overline{\nmax}
\end{equation*}
gives, in turn,
the following Bruhat decomposition of $G$ on the group level:
\begin{equation}
 G=\overline{\Pmax} w_0 \overline{\Pmax}\amalg \overline{\Pmax}.
\end{equation}

\item
\label{item:I3}
({\it Restricted root system})
Let $\mathfrak{g} = \mathfrak{k} + \mathfrak{p}$
be the Cartan decomposition, and
we take a maximal 
abelian subalgebra $\mathfrak{b}$ of $\mathfrak{p}$.
Since $\operatorname{Ad}(w_0)$ acts on $\mathfrak{p}$ as 
$-\operatorname{id}$,
$w_0$ acts on the restricted root system
$\Sigma (\mathfrak g, \mathfrak b)$ as 
$-\operatorname{id}$. 
We note that the longest element 
in the Weyl group $W(\Sigma(\mathfrak g, \mathfrak b))$ 
is equal to $-\operatorname{id}$ if $p \ne q$ or if $p = q \in 2\mathbb{Z}$.
In the case $p=q \in 2\mathbb{Z} + 1$, 
the restricted root system
$\Sigma (\mathfrak{g},\mathfrak{b})$
is of type $D_q$ and
$$
 -\operatorname{id}\notin 
 W(\Sigma(\mathfrak g, \mathfrak b))
 \simeq \mathfrak S_q \ltimes (\mathbb Z/2\mathbb Z)^{q-1},
$$
and therefore 
$\operatorname{Ad}(w_0)\vert_{\mathfrak{b}} = -\operatorname{id}$
gives an outer automorphism on $\Sigma(\mathfrak{g},\mathfrak{b})$.

\item
\label{item:I2}
\index{B}{Jordan algebra}%
({\it Jordan algebras})\enspace
We regard $\mathbb R^{p-1, q-1}\simeq \mathbb R^{p+q-2}$
as the semisimple Jordan algebra with 
indefinite quadratic form $Q(\xzeta)$
(see \eqref{def:Q}).
This Jordan algebra is euclidean if $\operatorname{min}(p,q)=2$,
and non-euclidean if $p,q >2$.
The conformal group 
\index{B}{Kantor--Koecher--Tits group}%
({\it Kantor--Koecher--Tits group}) 
of $\mathbb R^{p-1,q-1}$ is equal to the group $G=O(p,q)$, 
and the action of the element $w_0$ 
on $\mathbb R^{p-1,q-1}$ corresponds to 
the 
\index{B}{conformal inversion|main}%
{\it conformal inversion}
$\xzeta \mapsto -\xzeta^{-1}$ (see
 \cite{xfk,xSp}).
Thus, 
we call $w_0$ 
the {\it conformal inversion element}. 
The structure group 
\begin{equation}\label{def:Lp}
 L_+:= \Mmax_+ A \simeq 
\index{A}{CO(Q)@$CO(Q)$}%
CO(Q) = O(p-1,q-1) \times \mathbb{R}_{>0}
\end{equation}
acts
on $\mathbb R^{p-1,q-1}$ by
$x \mapsto e^t mx$ for
$(m,e^{tH}) \in \Mmax_+ \times A$,
and on its dual space by
$
 \xzeta \mapsto e^{-t}\,{}^t m^{-1} \xzeta
$.
\index{A}{1zpropertiesI0-I5@\textbf{I0}--\textbf{I5}|)}%

\item
\label{item:I5}
({\it The action on the minimal representation})\enspace
In the conformal model $(\varpi^{p,q},\overline{V^{p,q}})$
(see Section \ref{subsec:Kpic}),
the whole group $G$ acts on the geometry, and therefore, 
the unitary operator $\varpi^{p,q}(w_0)$ is easy to describe:
\begin{equation*}
   (\varpi^{p,q}(w_0)h)(v_0,\dots,v_{p+q-1})
   = h(v_0,\dots,v_{p-1},-v_p,\dots,-v_{p+q-1}),
\end{equation*}
(see the definition \eqref{def:vpi}).
Then, by using the intertwining operator
$\mathcal{T}: L^2(C) \to \overline{V^{p,q}}$
(see \eqref{def:B}),
we get the formula of $\mathcal{F}_C = \pi(w_0)$
in the Schr\"{o}dinger model
$(\pi,L^2(C))$ as
\begin{equation}\label{eqn:FCTvarpi}
  \mathcal{F}_C = \mathcal{T}^{-1} \circ \varpi^{p,q}(w_0) \circ
  \mathcal{T}.
\end{equation}
However, we do not know how to find
the formulas \eqref{def:Psi0}--\eqref{def:Psi} of $\mathcal{F}_C$
directly from \eqref{eqn:FCTvarpi} and the definition of
$\mathcal{T}$. 
Thus, in order to find the
unitary inversion formulas of $\mathcal{F}_C = \pi(w_0)$, 
we shall take a roundabout course 
(by using the $K$-type decomposition in part) 
in later chapters.

\end{enumerate}

Now, let us study basic properties of the unitary inversion operator
\index{A}{FC@$\mathcal{F}_C$}%
$\mathcal{F}_C=\pi(w_0)$ on $L^2(C)$.
First, we extend $\mathcal{F}_C$ to distribution vectors.

Let $L^2(C)^{-\infty}$ be the continuous dual of the Fr\'{e}chet space
$L^2(C)^\infty$ (Definition \ref{def:smoothvec}), namely,
the space of continuous $\mathbb{C}$-linear maps 
$L^\infty(C) \to \mathbb{C}$.
Elements in 
\index{A}{L2Cminusinfty@$L^2(C)^{-\infty}$|main}%
$L^2(C)^{-\infty}$ 
are called 
\index{B}{distribution vector|main}%
\textit{distribution vectors}.

We introduce a conjugate linear map 
\begin{equation*}
   \iota: L^2(C) \to L^2(C)^{-\infty}
\end{equation*}
characterized by
\begin{equation} \label{eqn:iGel}
   \langle v, \iota(u) \rangle
   := \int_C v(x) \overline{u(x)} d\mu(x)
   \quad\text{for any $v \in L^\infty(C)$}.
\end{equation}
The inclusion
\begin{equation*}
   L^2(C)^\infty \subset L^2(C) \subset L^2(C)^{-\infty}
\end{equation*}
is sometimes referred to as the 
\index{B}{Gelfand triple}%
\textit{Gelfand triple}.

{}From the standard theory of Sobolev spaces and 
the definition of smooth vectors, 
we have the following inclusive relation:
\begin{equation*}
   C_0^\infty(C) \subset 
\index{A}{L2Cinfty@$L^2(C)^\infty$}%
L^2(C)^\infty \subset C^\infty(C).
\end{equation*}
(To see $L^2(C)^\infty \subset C^\infty(C)$,
it is enough to use the $\overline{\Pmax}$-action on $L^2(C)$.)
Then, taking their duals,
we get
\begin{equation} \label{eqn:5C}
   \mathcal{E}'(C) \subset L^2(C)^{-\infty} \subset \mathcal{D}'(C),
\end{equation}
where $\mathcal{D}'(C)$ stands for the space of distributions on $C$,
and $\mathcal{E}'(C)$ for that of compactly supported distributions on
$C$. 

For $g \in G$,
we extend the unitary operator $S = \pi(g)$ on $L^2(C)$ to a
continuous operator $\widetilde{S}$ on $L^2(C)^{-\infty}$ by
\begin{equation} \label{eqn:defStild}
   \langle v,\widetilde{S}u \rangle
   := \langle S^{-1}v,u\rangle
   \quad\text{for $u\in L^2(C)^{-\infty}$ and $v\in L^2(C)^\infty$}.
\end{equation}
Here, we have used the fact that the unitary operator $S^{-1}$ on
$L^2(C)$ induces a continuous map
(we use the same letter $S^{-1}$):  
$$
S^{-1}: L^2(C)^\infty \to L^2(C)^\infty.
$$
Then, the extension $\widetilde{S}$ satisfies
\begin{equation*}
   \widetilde{S} \circ \iota
   = \iota \circ S
   \quad\text{on $L^2(C)$},
\end{equation*}
because for $u \in L^2(C)$ and $v \in L^\infty(C)$ we have
\begin{align*}
   \langle v,\widetilde{S}\iota(u) \rangle
   & = \langle S^{-1}v, \iota(u) \rangle
\\
   & = (S^{-1}v, u)_{L^2(C)}
\\ 
   & = (v, Su)_{L^2(C)}
\\
   & = \langle v, \iota(Su) \rangle.
\end{align*}
Likewise, the Lie algebra $\mathfrak{g}$ acts on $L^2(C)^{-\infty}$ by 
\begin{equation} \label{eqn:difftild}
   \langle v, d\tilde{\pi}(X)u \rangle
   := - \langle d\pi(X)v, u \rangle
  \quad\text{for $u \in L^2(C)^{-\infty}$ and
                 $v \in L^2(C)^{\infty}$}.
\end{equation}
Then, it readily follows that
\begin{equation} \label{eqn:Adtild}
   d\tilde{\pi} (\operatorname{Ad}(g)X)
   = \tilde{\pi}(g) \, d\tilde{\pi}(X) \, \tilde{\pi}(g^{-1})
\end{equation}
for $g\in G$ and $X\in\mathfrak{g}$.

We shall write simply $S$ for $\widetilde{S}$,
and $d\pi(X)$ for $d\tilde{\pi}(X)$ if there is no
confusion. 

For $\xi \in C$,
we denote by 
\index{A}{1delta_\xi@$\delta_\xi$}%
$\delta_\xi$ the 
\index{B}{Dirac delta function}%
Dirac delta function at $\xi \in C$,
namely, 
\begin{equation*}
   \langle v, \delta_\xi \rangle := v(\xi)
\end{equation*}
for a continuous function $v$ on $C$.
Since
\begin{equation*}
   \delta_\xi : L^2(C)^\infty \to \mathbb{C}
\end{equation*}
is a continuous map,
we may regard $\delta_\xi \in L^2(C)^{-\infty}$.
Then, we have
\begin{equation} \label{eqn:defKxi}
   \mathcal{F}_C^{-1} \delta_\xi
   \in L^2(C)^{-\infty}.
\end{equation}
Applying the definition \eqref{eqn:defStild} to $S=\mathcal{F}_C^{-1}$, 
we have for any $v \in L^\infty(C)$,
\begin{align}
   \langle v,\mathcal{F}_C^{-1} \delta_\xi \rangle
   & = \langle \mathcal{F}_C v, \delta_\xi \rangle
\nonumber
\\
   & = (\mathcal{F}_C v) (\xi).
\label{eqn:FCKer}
\end{align}

\begin{remNonumber} \label{rem:defKxi}
In contrast to the (Euclidean) Fourier transform
\index{A}{FRn@$\mathcal{F}_{\mathbb{R}^n}$}%
$\mathcal{F}_{\mathbb{R}^n}$, 
$\mathcal{F}_C^{-1} \delta_\xi$
is a real valued distribution.
In other words, the kernel\/ 
$k(x,\xi) = \frac{1}{(2\pi)^{\frac{n}{2}}} e^{\sqrt{-1} \langle x,\xi\rangle}$
for\/ $\mathcal{F}_{\mathbb{R}^n}$ is not real valued, 
while the kernel $K(x,\xi)$ for\/ $\mathcal{F}_C$ below is real valued.
This reflects the fact that both kernels $k(x,\xi)$ and $K(x,\xi)$ are
characterized by the system of differential equations up to scalar
(see Introduction \ref{subsec:1.2}):
\begin{alignat*}{4}
   & p_j k(\cdot,\xi) 
   &&  = \xi_j k(\cdot,\xi)
   &&\quad \text{in $\mathbb{R}^n$}
   && \quad(1 \le j \le n),
\\
   & P_j K(\cdot,\xi)
   &&  = 4\xi_j K(\cdot,\xi)
   &&\quad \text{in $C$}
   && \quad(1 \le j \le n),
\end{alignat*}
where $p_j = -\sqrt{-1} \frac{\partial}{\partial x_j}$.
A distinguishing feature here is that the self-adjoint differential
operator $p_j$ of first order satisfies
\begin{equation*}
   \overline{p_j} = -p_j,
\end{equation*}
while the self-adjoint differential operator 
\index{A}{P_j@$P_j$}%
$P_j$ of second order
satisfies
\begin{equation*}
   \overline{P_j} = P_j.
\end{equation*}
\end{remNonumber}

For a fixed $\xi\in C$,
$\mathcal{F}_C^{-1}\delta_\xi$ is a distribution on $C$ by
\eqref{eqn:defKxi} and \eqref{eqn:5C}.
We shall see that this distribution is locally integrable on the open
dense subset
\begin{equation*}
   \{ x\in C: \langle x,\xi \rangle \ne 0 \}.
\end{equation*}
Thus, it is convenient to write the distribution
$\mathcal{F}_C^{-1}\delta_\xi$ as a 
\index{B}{generalized function}%
generalized function in the sense
of Gelfand--Shilov \cite{xGeSh} by using the canonical measure
$d\mu(x) = \delta(Q(x))$ on $C$,
that is, 
\begin{equation*}
   (\mathcal{F}_C^{-1}\delta_\xi)(x) = K(x,\xi)d\mu(x)
\end{equation*}
for some generalized function $K(\cdot,\xi)$ on $C$.
Then, the formula \eqref{eqn:FCKer} can be expressed as
\begin{equation} \label{eqn:FCKx}
   (\mathcal{F}_C v)(\xi)
   = \int_C K(x,\xi)v(x)d\mu(x).
\end{equation}

The argument so far assures that $K(\cdot,\xi)$ is a generalized function on
$C$ for each fixed $\xi\in C$.
On the other hand, since the sesqui-linear map
\begin{equation*}
   L^2(C) \times L^2(C) \to \mathbb{C},
   \quad  (v,u) \mapsto (\mathcal{F}_Cv,u)_{L^2(C)}
\end{equation*}
is continuous,
we may regard $K(x,\xi)$ is a generalized function
(or equivalently, $K(x,\xi)d\mu(x)d\mu(\xi)$ is a distribution)
on the direct product manifold $C \times C$
such that
\begin{equation} \label{eqn:FCKxy}
   (\mathcal{F}_Cv,u)_{L^2(C)}
   = \int_{C\times C} v(x) \overline{u(\xi)} K(x,\xi) d\mu(x)
      d\mu(\xi)
\end{equation}
for any $u,v \in L^2(C)^\infty$.

Now, we are ready to state basic properties of the 
\index{B}{unitary inversion operator}%
`unitary inversion operator' 
\index{A}{FC@$\mathcal{F}_C$}%
$\mathcal{F}_C = \pi(w_0)$ as the counterpart of the
properties \textbf{I\ref{item:I0}}, 
\textbf{I\ref{item:I1}} of the conformal inversion
$w_0 \in G$:
\begin{theorem} \label{thm:FC}
\quad
\begin{enumerate}
    \renewcommand{\labelenumi}{\upshape \theenumi)}
\item  
$\mathcal{F}_C$ is an involutive unitary operator on $L^2(C)$.
In particular,
we have
\begin{equation} \label{eqn:Kselfadj}
   K(x,\xi) = \overline{K(\xi,x)}
\end{equation}
as distributions on $C \times C$.

\item  
If $v\in L^2(C)^\infty$,
then $\mathcal{F}_C v \in L^2(C)^\infty$.
In particular,
\begin{equation*}
   \mathcal{F}_C v \in C^\infty(C) \cap L^2(C)
\end{equation*}
for any compactly supported 
$C^\infty$ function $v$ on $C$.
\item  
For\/ $1 \le j \le n$,
the following identities of operators on $L^2(C)^\infty$ hold:
\begin{alignat*}{2}
   &\mathcal{F}_C \circ 4x_j
   &&= P_j \circ \mathcal{F}_C,
\\
   &\mathcal{F}_C \circ P_j
   &&= 4x_j \circ \mathcal{F}_C,
\\
   &\mathcal{F}_C \circ E
   &&= -(E+n-2) \circ \mathcal{F}_C.
\end{alignat*}
These identities hold also on $L^2(C)^{-\infty}$. 
\item  
For any $v \in L^2(C)^\infty$,
\begin{equation*}
   (\mathcal{F}_C v)(\xi)
   = \int_C v(x) K(x,\xi) d\mu(x).
\end{equation*}
\item  
For each fixed
$\xi = (\xi_1,\dots,\xi_n) \in C$,
the generalized function $K(\xi,\cdot)$ solves the system of differential
equations on $C$: 
\begin{equation*}
   P_j K(\xi,\cdot)
   = 4\xi_j K(\xi,\cdot)
   \qquad (1 \le j \le n).
\end{equation*}
In turn,
$K(\cdot,\xi)$ satisfies
\begin{equation*}
   P_j K(\cdot,\xi) = 4\xi_j K(\cdot,\xi)
   \qquad (1 \le j \le n).
\end{equation*}
\end{enumerate}
\end{theorem}

\begin{remNonumber} \label{rem:FC}
We shall see in Theorem \ref{thm:A} that
\begin{equation*}
   K(x,\xi) = K(\xi,x)
\end{equation*}
as distributions on $C\times C$.
\end{remNonumber}

\begin{proof}[Proof of Theorem \ref{thm:FC}]
1)\enspace
The first statement is straightforward
 from definition $\mathcal{F}_C = \pi(w_0)$ because
$(\pi,L^2(C))$ is a unitary representation of $G$ (Fact \ref{fact:B})
and $w_0^2 = 1$.
The second statement follows from \eqref{eqn:FCKxy} and 
\begin{equation*}
   (\mathcal{F}_C v,u)
   = (\mathcal{F}_C v, \mathcal{F}_C^2 u)
   = (v,\mathcal{F}_C u).
\end{equation*}

2)\enspace
The group $G$ preserves $L^2(C)^\infty$,
and in particular, we have
\begin{equation*}
   \mathcal{F}_C(L^2(C)^\infty) = L^2(C)^\infty.
\end{equation*}
Since 
$C_0^\infty(C) \subset L^2(C)^\infty \subset C^\infty(C) \cap L^2(C)$,
we get
$\mathcal{F}_C(C_0^\infty(C)) \subset C^\infty(C) \cap L^2(C)$.

3)\enspace
By \eqref{eqn:Adwj} in \textbf{I\ref{item:I1}},
we have the following identities on $L^2(C)^{\infty}$,
and also on $L^2(C)^{-\infty}$ by \eqref{eqn:Adtild}:
\begin{alignat*}{2}
   &\mathcal{F}_C \circ d\pi(\overline{N}_j) 
   &&= \epsilon_j d\pi(N_j) \circ \mathcal{F}_C,
\\
   &\mathcal{F}_C \circ d\pi(N_j)
   &&= \epsilon_j d\pi(\overline{N}_j) \circ \mathcal{F}_C,
\\
   &\mathcal{F}_C \circ d\pi(H)
   &&= -d\pi(H) \circ \mathcal{F}_C.
\end{alignat*}
Now we recall the formulas of the differential action:
\begin{alignat*}{2}
   &d\pi(\overline{N}_j)
   = 2\sqrt{-1} x_j
   &&\qquad\text{(see \eqref{eqn:barNj})},
\\
   &d\pi(N_j)
   = \frac{\sqrt{-1}}{2} \epsilon_j P_j
   &&\qquad\text{(see \eqref{eqn:NjPj})},
\\
   &d\pi(H)
   = -(E+\frac{n-2}{2})
   &&\qquad\text{(see \eqref{eqn:dpiH})}.
\end{alignat*}
These operators are defined on
$L^2(C)^\infty \subset C^\infty(C)$,
and in turn,
they are defined on
$L^2(C)^{-\infty} \subset \mathcal{D}'(C)$
(see \eqref{eqn:5C}).
Therefore, the third statement follows.

4)\enspace
Since $w_0^{-1}=w_0$,
we have $\mathcal{F}_C^{-1} = \mathcal{F}_C$,
and therefore the statement is equivalent to what we have seen
in \eqref{eqn:FCKx}.

5)\enspace
It follows from the third and fourth statements that we have
(by switching $x$ and $\xi$)
\begin{align*}
\int_C 4\xi_j K(\xi,x) v(\xi) d\mu(\xi)
   ={}& (P_j)_x \int_C K(\xi,x) v(\xi) d\mu(\xi)
\\
   ={}& \int_C (P_j)_x K(\xi,x) v(\xi) d\mu(\xi)
\end{align*}
for any $v \in L^2(C)^\infty$.
Hence, we have shown
\begin{equation*}
   4\xi_j K(\xi,x) = (P_j)_x K(\xi,x).
\end{equation*}
The second assertion follows from $\overline{P_j} = P_j$ as we saw in
Remark \ref{rem:defKxi}.
Thus, Theorem \ref{thm:FC} has been proved.
\end{proof}

If $A$ is a continuous operator on $L^2(C)$,
then $A$ induces a linear map (we use the same letter)
$A: L^2(C)^\infty \to L^2(C)^\infty$ on the space of smooth vectors. 

The next result shows that the intertwining relation in Theorem
\ref{thm:FC} (3) characterizes the operator $\mathcal{F}_C$ up to
scalar: 
\begin{thmSec} \label{thm:FCunique}
Suppose $A$ is a continuous operator on $L^2(C)$ such that
\begin{alignat*}{2}
   & A \circ 4x_j &&= P_j \circ A,
\\
   & A \circ P_j &&= 4x_j \circ A,
\end{alignat*}
on $L^2(C)^\infty$ for $1\le j \le n$.
Then, $A$ is a scalar multiple of\/ $\mathcal{F}_C$.
In particular, $A^2$ is a scalar multiple of the identity operator.
\end{thmSec}

\begin{proof}
We set $B:=\mathcal{F}_C^{-1}\circ A$.
Then, $B$ satisfies
\begin{alignat*}{2}
   & B \circ 4x_j &&= 4x_j \circ B,
\\
   & B \circ P_j &&= P_j \circ B,
\end{alignat*}
on $L^2(C)^\infty$, and consequently,
\begin{alignat*}{2}
   & B \circ d\pi(\overline{N}_j) &&= d\pi(\overline{N}_j) \circ B,
\\
   & B \circ d\pi(N_j) &&= d\pi(N_j) \circ B,
\end{alignat*}
for $1\le j \le n$.
Since the vector space
 $\nmax + \overline{\nmax} = \sum_{j=1}^n \mathbb{R} N_j +
       \sum_{j=1}^n \mathbb{R} \overline{N}_j$
generates the whole Lie algebra $\mathfrak{g}$,
we have
\begin{equation} \label{eqn:Binter}
   B \circ d\pi(X) = d\pi(X) \circ B
\end{equation}
for any $X \in \mathfrak{g}$.
In particular, 
such an operator $B$ preserves the space of $K$-finite vectors,
namely, 
\begin{equation*}
   B(L^2(C)_K) \subset L^2(C)_K.
\end{equation*}
Therefore, the identities \eqref{eqn:Binter} hold also on $L^2(C)_K$. 

On the other hand,
it follows from Fact \ref{fact:A} and Fact \ref{fact:B} that $G$ acts
irreducibly on $L^2(C)$.
Then, 
 $L^2(C)_K$ is an irreducible $(\mathfrak{g},K)$-module. 
Therefore, 
$B$ must be a scalar multiple of the identity by Schur's lemma.
Hence, $A$ is a scalar multiple of $\mathcal{F}_C$.
The last assertion is an immediate consequence of the fact that
$\mathcal{F}_C^2 = \operatorname{id}$
(see \textbf{I\ref{item:I0}}).
\end{proof}

\begin{thmSec} \label{thm:Pjunique}
For a fixed 
$\xi = (\xi_1,\dots,\xi_n) \in \mathbb{R}^n \setminus \{0\}$, 
we consider the system of differential equations:
\begin{equation} \label{eqn:PjPDE}
   P_j \psi = 4 \xi_j \psi
   \quad (1 \le j \le n),
\end{equation}
as distributions on $C$.
\begin{enumerate}
    \renewcommand{\labelenumi}{\upshape \theenumi)}
\item  
If $Q(\xi) \ne 0$, then $\psi = 0$.
\item  
If $Q(\xi) = 0$, then the solution space in $L^2(C)^{-\infty}$ is one-dimensional. 
\end{enumerate}
\end{thmSec}

\begin{proof}
1)\enspace
It follows from \eqref{eqn:PjPDE} that
\begin{equation*}
   ( \sum_{j=1}^n \epsilon_j P_j^2) \psi
   = 4 Q(\xi) \psi.
\end{equation*}
Then, by Theorem \ref{thm:Pj} (5),
the left-hand side vanishes.
Hence, $\psi = 0$ unless $Q(\xi) = 0$.

2)\enspace
Suppose $\xi\in C$.
Taking the Fourier transform $\mathcal{F}_C$ of the differential
equation \eqref{eqn:PjPDE},
we have the following equation
\begin{equation*}
   4x_j(\mathcal{F}_C\psi) = 4\xi_j (\mathcal{F}_C \psi)
   \quad\text{in $L^2(C)^{-\infty}$}
\end{equation*}
from Theorem \ref{thm:FC} (3).
Hence,
\begin{equation} \label{eqn:FPjPDE}
   (x_j - \xi_j) (\mathcal{F}_C \psi) = 0
   \quad (1 \le j \le n).
\end{equation}
Since 
$\mathcal{F}_C \psi \in L^2(C)^{-\infty} \subset \mathcal{D}'(C)$, 
the equations \eqref{eqn:FPjPDE} hold as distributions on $C$.
Then, considering the equations \eqref{eqn:FPjPDE} in a local chart, 
we see that
$\mathcal{F}_C\psi$ is a scalar multiple of the delta function
$\delta_\xi$.
Conversely,
$\delta_\xi$ solves \eqref{eqn:FPjPDE} and 
$\delta_\xi \in L^2(C)^{-\infty}$.
Therefore, $\mathcal{F}_C^{-1} \delta_\xi \in L^2(C)^{-\infty}$ solves
\eqref{eqn:PjPDE}.
Now, Theorem is proved.
\end{proof}

\chapter{$K$-finite eigenvectors in the Schr\"{o}dinger model $L^2(C)$}
\label{sec:K}

In the conformal model (see Section \ref{subsec:Kpic}),
we can find readily explicit $K$-finite vectors.
However, it is far from being obvious to find explicit forms of
$K$-finite vectors in the $L^2$-model for the minimal representation. 
In this chapter, generalizing the idea of \cite[Theorem 5.8]{xkors3}
for the minimal $K$-type,
we find explicit vectors in $L^2(C)$ for \textit{every} $K$-type,
by carrying out the computation involving the integral operator
$\B :
L^2(C) \rarrowsim \overline{V^{p,q}}$ (see Fact \ref{fact:B}).

\section{Result of this chapter}
\label{subsec:K}

Throughout this chapter, we assume $p\ge q \ge 2$ and $p+q \ge 6$. 
For $(l,k) \in \mathbb N^2$, we consider the following two
(non-exclusive) cases:
\begin{equation}\label{def:Cases}
\fbox{\parbox{12em}{$
  \text{\ Case 1}:~
  \frac{p-q}2+l-k\ge 0, \\[1ex]
  \text{\ Case 2}:~
  \frac{p-q}2+l-k\le 0.
$}}
\end{equation}
The case $\frac{p-q}2+l-k=0$ belongs to both 
Cases 1 and 2. 
This overlap will be convenient later because
all the formulas below are the same for both 
Cases 1 and 2 if $(l,k)$ satisfies
$\frac{p-q}{2} + l -k = 0$.

For $(l,k) \in \mathbb N^2$, we define real analytic functions 
$f_{l,k}$ on $\mathbb R_+$ by 
\begin{align}
\label{def:f}
\index{A}{f0lkr@$f_{l,k}(r)$|main}%
  f_{l,k}(r):={}&
  \begin{cases}
     r^{-\frac{q-3}{2}+l}K_{\frac{q-3}{2}+k}(2r)
      &\text{Case 1},
     \\
     r^{-\frac{p-3}{2}+k}K_{\frac{p-3}{2}+l}(2r)
      &\text{Case 2},
  \end{cases}
\\
  ={}& r^{l+k} \times
  \begin{cases}
     \widetilde{K}_{\frac{q-3}{2}+k}(2r)
      &\text{Case 1},
     \\
     \widetilde{K}_{\frac{p-3}{2}+l}(2r)
      &\text{Case 2}.
  \end{cases}
\label{eqn:flkn}   
\end{align}
Here, 
\index{A}{Knuz@$K_\nu(z)$}%
$K_\nu(z)$ is the $K$-Bessel function, i.e.,
 the modified Bessel function of the second kind
(see Appendix \ref{subsec:B}) and
\index{A}{Knuztilde@$\widetilde {K}_\nu (z)$}%
 $\widetilde{K}_\nu(z)=(\frac{z}{2})^{-\nu}K_\nu(z)$ 
is the normalized $K$-Bessel function (see \eqref{eqn:Ktilde}).

By using the polar coordinate \eqref{def:pol}, we
define a linear subspace $H_{l,k}$ of $C^\infty(C)$ consisting of
linear combinations of the following functions:
\begin{equation}
\label{def:H}
  f_{l,k}(r) \phi(\omega) \psi(\eta)
 \quad           (\phi \in \Har{l}{p-1}, 
 \ 
 \psi \in \Har{k}{q-1}).
\end{equation}
Here,
\index{A}{HjRm@$\mathcal{H}^j(\mathbb{R}^m)$}%
$\mathcal{H}^j(\mathbb{R}^m)$
denotes the space of spherical harmonics of degree $j$
(see Appendix \ref{subsec:H}).

We recall from Section~\ref{subsec:act} that there are two key
compact subgroups for the analysis on the minimal representation
$L^2(C)$: 
\begin{alignat}{2}
\index{A}{K1Opq@$K \simeq O(p)\times O(q)$}%
 & K &\  &\simeq O(p) \times O(q), \notag\\
\index{A}{K1Opqdash@$K' \simeq O(p-1)\times O(q-1)$}%
 & K' = K \cap \Mmax_+ &\  &\simeq O(p-1) \times O(q-1).\notag
\end{alignat}
We note that the $K'$-action on $L^2(C)$ is just the pull-back of the
$K'$-action on $C$ (see \eqref{eqn:rM}),
but the $K$-action on $L^2(C)$ is more complicated because $K$ cannot
act on $C$.
Then, here is our main result of this chapter:

\begin{theorem}
\label{thm:K}
  For each pair  $(l,k) \in \mathbb N^2$, we have 
  
{\rm 1)} 
(asymptotic behavior) 
\index{A}{Hlk@$H_{l,k}$|main}%
$H_{l,k} \subset L^2(C)$ for any $l,k \in \mathbb N$.

{\rm 2)}
($K$-type and $K'$-type) 
$H_{l,k} \simeq \Har{l}{p-1} \otimes \Har{k}{q-1}$ as a
$K'$-module. Furthermore, 
$H_{l,k}$ belongs to the $K$-type
$\Har{a(l,k)}{p} \otimes  \Har{a(l,k)+\frac{p-q}2}{q}$ of $L^2(C)$.
Here, we define a non-negative integer $a(l,k)$ by
\begin{equation}\label{def:alk}
 a(l,k):=\max(l,k-\frac{p-q}{2})=
    \begin{cases} l &\text{Case 1}, \\ k-\frac{p-q}2
   &\text{Case 2}.
    \end{cases}  
\end{equation}

{\rm 3)} 
(eigenspace of $\pi(w_0)$)
$\pi(w_0)$ acts on $H_{l,k}$ by the scalar $(-1)^{a(l,k)+\frac{p-q}2}$.  

{\rm 4)} (intertwining operator) Fix any 
$\phi \in \Har{l}{p-1}$, 
\thinspace
$\psi \in \Har{k}{q-1}$, 
and set
\begin{equation}\label{eqn:ulk}
u_{l,k}(r\omega, r\eta):=f_{l,k}(r) \phi(\omega) \psi(\eta)
\in H_{l,k}.
\end{equation}
Then 
\index{A}{TL2C@$\B: L^2(C) \to \overline{V^{p,q}}$}%
$\B :L^2(C) \to \overline \Vpq$ 
has the following form on the subspace $H_{l,k}$:  
  \begin{equation}
  \label{eqn:Tu}
  \B u_{l,k}
  = c_{l,k} I_{l \to a(l,k)}^p(\phi) I_{k\to a(l,k)+\frac{p-q}2}^q (\psi),
  \end{equation} 
where 
\index{A}{Iijm@$I_{i\to j}^m:\mathcal{H}^i(\mathbb{R}^{m-1})
      \to \mathcal{H}^j(\mathbb{R}^m)$}%
$I_{i \to j}^m:\Har{i}{m-1} \to \Har{j}{m}$ 
$(0\le i\le j)$ is an $O(m-1)$-homomorphism  
defined in  Fact \ref{fact:H}, and 
the constant $c_{l,k}$ is given by
\begin{equation}\label{def:clk}
 c_{l,k}:=
 \frac{\sqrt{-1}^{l+k}}{\sqrt{\pi}}
 \times 
 \begin{cases}
\displaystyle
 \frac{\Gamma(\frac{p-q}2+l-k+1) }
      {2^{\frac{p-q}{2}+l-k+3}  \Gamma(\frac{p-2}{2}+l)}
  &\text{Case 1},   \\[2.5ex]
\displaystyle
 \frac{\Gamma(-\frac{p-q}2-l+k+1)}
      {2^{-(\frac{p-q}{2}+l-k)+3}  \Gamma(\frac{q-2}{2}+k)}

 &\text{Case 2}.
 \end{cases}
\end{equation}  

{\rm 5)}
($L^2$-norm)
For any $(l,k) \in \mathbb{N}^2$,
\begin{equation}\label{eqn:TulkL2}
\| \mathcal{T}u_{l,k}\|_{L^2(M)}^2
= \frac{1}{2a(l,k)+p-2}
  \|u_{l,k}\|_{L^2(C)}^2.
\end{equation}
\end{theorem}

A special case (the $k=0$ case) of Theorem \ref{thm:K}
(2) provides us an explicit $K$-finite vector for every $K$-type. 
\begin{corollary}\label{cor:Kfinite}
For each $a\in\mathbb{N}$,
the function
\begin{equation} \label{eqn:Kfa}
   r^{a-\frac{q-2}{2}} K_{\frac{q-2}{2}} (2r) \phi(\omega)
\end{equation}
is a $K$-finite vector for any 
$\phi\in\mathcal{H}^a(\mathbb{R}^{p-1})$.
More precisely, it belongs to the $K$-type
\begin{equation*}
   \mathcal{H}^a(\mathbb{R}^p) \otimes
   \mathcal{H}^{a+\frac{p-q}{2}} (\mathbb{R}^q).
\end{equation*}
\end{corollary}

\begin{proof}
Under our assumption $p\ge q$, $a(l,0)=l$ by the definition
\eqref{def:alk} because $(l,0)$ belongs to Case 1.
Hence, Corollary follows immediately from Theorem \ref{thm:K} (2).
\end{proof}

By using the unitary inner product $(\ ,\ )_M$
(see \eqref{eqn:innerM})
for the conformal model
$(\varpi^{p,q},\overline{V^{p,q}})$ of the minimal representation
$G=O(p,q)$,
the formula \eqref{eqn:TulkL2} in Theorem \ref{thm:K} (5)
can be restated as
\begin{equation} \label{eqn:Tnorm}
  (\mathcal{T}u_{l,k},\mathcal{T}u_{l,k})_M
  = \frac{1}{2} \|u_{l,k}\|_{L^2(C)}^2.
\end{equation}
\begin{remNonumber}
Theorem \ref{thm:K} (2) for $l=k=0$ (or Corollary \ref{cor:Kfinite}
for $a=0$) was
proved in \cite[Theorem 5.8]{xkors3}.
Since $p \ge q$, we are dealing with Case~1 if $l = k = 0$
and $a(0,0) = 0$. In this 
particular case, Theorem~\ref{thm:K} (2) asserts that
$f_{0,0}(r)=r^{-\frac{q-3}2} K_{\frac{q-3}2}(2r)$
belongs to the minimal $K$-type $\Har{0}{p} \otimes \Har{\frac{p-q}2}{q}$
of $(\pi, L^2(C))$.
\end{remNonumber}

\begin{remNonumber} \label{rem:3.1.3}
For $q=2$, 
$\pi$ is essentially a 
\index{B}{highest weight module}%
highest weight module 
(see Section \ref{subsec:1.6} for the $n_2=1$ case).
We note that for $q=2$,
$\Har{k}{q-1}$ is non-zero only if $k=0$ or $1$
(see Appendix \ref{subsec:H}). 
Thus our assumption $p+q\ge 6$ combined with $q=2$ and $k=0,1$ 
implies 
$\frac{p-q}2+l-k=\frac{p+q}2-2+l-k\ge 0$.
Hence, 
the pairs $(l,k)$ $(k=0,1)$ belong automatically to Case 1.
It then follows from the definition \eqref{def:f} that
\begin{align*}
  & f_{l,0}(r) = r^{\frac{1}{2}+l} K_{-\frac{1}{2}} (2r),
\\
  & f_{l,1}(r) = r^{\frac{1}{2}+l} K_{\frac{1}{2}} (2r).
\end{align*}
By using the formula
\[
  K_{\frac{1}{2}}(z) = K_{-\frac{1}{2}}(z)
  = \bigl( \frac{\pi}{2z} \bigr)^{\frac{1}{2}} e^{-z}
\]
(see \eqref{eqn:K12}), we get
$$
f_{l,0}(r)= f_{l,1}(r)
=\frac{\sqrt{\pi}}2 r^l e^{-2r}.
$$ 
We note that $f_{l,0}(r)$ $(= f_{l,1}(r))$ 
coincides with the function 
`$f_{l,l}(r)$' in \cite[(3.2.4)]{xkmano2} up to 
a constant multiple.
(In \cite{xkmano2}, we investigated the Schr\"{o}dinger model for
$O(p,2)$, and the function `$f_{l,l}(r)$' was defined by using the
Laguerre polynomial. 
The definition of `$f_{l,l}(r)$' is different from \eqref{def:f}.)
\end{remNonumber}

Our method to show Theorem \ref{thm:K} is based on the technique of 
\cite[\S 5.6, \S 5.7]{xkors3}.
The key lemma is Lemma \ref{lem:Han}, which 
gives a formula of 
the Hankel transform of the $K$-Bessel function 
with trigonometric parameters
by means of the Gegenbauer polynomial.

The subspace $\bigoplus\limits_{l,k\in\mathbb{N}} H_{l,k}$
is not dense in $L^2(C)$,
but is large enough 
(see Section \ref{subsec:Hlk}) that we can make use of
Theorem \ref{thm:K} for the proof 
of Theorem \ref{thm:C} in Chapter \ref{sec:C}
(see also Section \ref{subsec:C} for its idea).

\section{$K \cap \Mmax$-invariant subspaces $H_{l,k}$} \label{subsec:Hlk}
The subspace $\bigoplus\limits_{l,k\in\mathbb{N}} H_{l,k}$ is not
dense in $L^2(C)$,
but serves as a `skeleton'.
In this section, we try to clarify its meaning.

We begin with
 the branching law
$G \downarrow K$ (see \eqref{eqn:Ktype}) and $K \downarrow K'$
($K'$ denotes $K':=K \cap \Mmax_+\simeq
O(p-1)\times O(q-1)$):
\begin{align}
\index{A}{L2CK@$L^2(C)_K$}%
 L^2(C)_K\simeq{} & \bigoplus_{a=0}^\infty 
 \Har{a}{p} \otimes \Har{a+\frac{p-q}2}{q}  
\label{eqn:3.1.7}
\\
\simeq{} & \bigoplus_{a=0}^\infty
\bigoplus_{l=0}^a \bigoplus_{k=0}^{a+\frac{p-q}2}
\Har{l}{p-1} \otimes \Har{k}{q-1}.
\label{eqn:3.1.8}
\end{align}
The irreducible decomposition \eqref{eqn:3.1.7} shows that
$L^2(C)_K$ is multiplicity-free as a $K$-module.
Hereafter, we identify the $K$-module
$\mathcal{H}^a(\mathbb{R}^p) \otimes 
 \mathcal{H}^{a+\frac{p-q}{2}}(\mathbb{R}^q)$
with the corresponding subspace of $L^2(C)_K$.
Then we observe: 
\begin{enumerate}[{\bf {S}1}]
\index{A}{1zpropertiesS1-S2@\textbf{S1}--\textbf{S2}|(}%
\item \label{S1}({\it $K'$-type}) Fix a pair
$(l,k)\in \mathbb N^2$. 
In light of \eqref{eqn:3.1.7} and \eqref{eqn:3.1.8},
the $K'$-type $\Har{l}{p-1} \otimes \Har{k}{q-1}$
occurs in the $K$-module $\Har{a}{p} \otimes \Har{a+\frac{p-q}{2}}{q}$
if and only if $a \ge a(l,k)$.
Further, 
\index{A}{Hlk@$H_{l,k}$}%
$H_{l,k}$ is characterized as a subspace of $L^2(C)$ satisfying
the following two conditions:
$$
 \begin{cases}W \simeq \Har{l}{p-1}\otimes \Har{k}{q-1} 
              \text{ as $K'$-modules,}\\
              W \subset \Har{a(l,k)}{p}\otimes \Har{a(l,k)+\frac{p-q}{2}}{q}.
 \end{cases}
$$ 

\item \label{S2}({\it $K$-type}) Fix $a \in \mathbb N$.
Then for $(l,k) \in \mathbb{N}^2$,
$$
 \biggl(\bigoplus_{l,k\in \mathbb N}H_{l,k}\biggr)
 \cap \biggl(\Har{a}{p}\otimes
 \Har{a+\frac{p-q}2}{q}\biggr)
$$ 
is non-zero if and only if at least one of $l$ and $k$ attains its maximum
in the set
$\{(l',k')\in\mathbb{N}^2: 0 \le l' \le a, \ 
                         0 \le k' \le a+\frac{p-q}{2}\}$
or equivalently, in the set of
the $K'$-types $(l',k')$ occurring in 
$\Har{a}{p}\otimes \Har{a+\frac{p-q}2}{q}$ 
(see the black dots
\
\setlength{\unitlength}{0.0003in}
\begin{picture}(100,100)(50,-130)
\put(0,0)\blacken\ellipse{150}{150}
\end{picture}
\
 in the figure below).
\begin{figure}[H]
\setlength{\unitlength}{0.0003in}
\begingroup\makeatletter\ifx\SetFigFont\undefined%
\gdef\SetFigFont#1#2#3#4#5{%
  \reset@font\fontsize{#1}{#2pt}%
  \fontfamily{#3}\fontseries{#4}\fontshape{#5}%
  \selectfont}%
\fi\endgroup%
{\renewcommand{\dashlinestretch}{30}
\begin{picture}(7875,8559)(-1800,-10)
\put(2550,6912){\blacken\ellipse{150}{150}}
\put(2550,6912){\ellipse{150}{150}}
\put(3450,6912){\blacken\ellipse{150}{150}}
\put(3450,6912){\ellipse{150}{150}}
\put(4350,6912){\blacken\ellipse{150}{150}}
\put(4350,6912){\ellipse{150}{150}}
\put(5250,6912){\blacken\ellipse{150}{150}}
\put(5250,6912){\ellipse{150}{150}}
\put(6150,6912){\blacken\ellipse{150}{150}}
\put(6150,6912){\ellipse{150}{150}}
\put(6150,6012){\blacken\ellipse{150}{150}}
\put(6150,6012){\ellipse{150}{150}}
\put(6150,5112){\blacken\ellipse{150}{150}}
\put(6150,5112){\ellipse{150}{150}}
\put(6150,4212){\blacken\ellipse{150}{150}}
\put(6150,4212){\ellipse{150}{150}}
\put(6150,3312){\blacken\ellipse{150}{150}}
\put(6150,3312){\ellipse{150}{150}}
\put(6150,2412){\blacken\ellipse{150}{150}}
\put(6150,2412){\ellipse{150}{150}}
\put(6150,1512){\blacken\ellipse{150}{150}}
\put(6150,1512){\ellipse{150}{150}}
\put(6150,612){\blacken\ellipse{150}{150}}
\put(6150,612){\ellipse{150}{150}}
\path(2550,8412)(2550,12)
\path(1875,612)(7650,612)
\path(1950,2712)(7650,8412)
\path(1950,2712)(7650,8412)
\put(6750,4812){\makebox(0,0)[lb]{\smash{{{\SetFigFont{12}{14.4}{\familydefault}{\mddefault}{\updefault}Case 1}}}}}
\put(3450,7512){\makebox(0,0)[lb]{\smash{{{\SetFigFont{12}{14.4}{\familydefault}{\mddefault}{\updefault}Case 2}}}}}
\put(7725,8037){\makebox(0,0)[lb]{\smash{{{\SetFigFont{10}{14.4}{\familydefault}{\mddefault}{\updefault}$k-l=\displaystyle\frac{p-q}{2}$}}}}}
\put(2100,8412){\makebox(0,0)[lb]{\smash{{{\SetFigFont{12}{14.4}{\familydefault}{\mddefault}{\updefault}$k$}}}}}
\put(1200,3237){\makebox(0,0)[lb]{\smash{{{\SetFigFont{10}{14.4}{\familydefault}{\mddefault}{\updefault}${\displaystyle\frac{p-q}{2}}$}}}}}
\put(2100,87){\makebox(0,0)[lb]{\smash{{{\SetFigFont{12}{14.4}{\familydefault}{\mddefault}{\updefault}$0$}}}}}
\put(6075,87){\makebox(0,0)[lb]{\smash{{{\SetFigFont{12}{14.4}{\familydefault}{\mddefault}{\updefault}$a$}}}}}
\put(7875,87){\makebox(0,0)[lb]{\smash{{{\SetFigFont{12}{14.4}{\familydefault}{\mddefault}{\updefault}$l$}}}}}
\put(400,6837){\makebox(0,0)[lb]{\smash{{{\SetFigFont{10}{14.4}{\familydefault}{\mddefault}{\updefault}$a+\displaystyle\frac{p-q}{2}$}}}}}
\end{picture}
}
\caption{}
\label{fig:skt}
\end{figure}

\index{A}{1zpropertiesS1-S2@\textbf{S1}--\textbf{S2}|)}%
\end{enumerate}

\section{Integral formula for the 
 $K \cap \Mmax$-intertwiner}
\label{subsec:intwn}

As a preparation for Theorem \ref{thm:K}, 
we rewrite
the integral formula for 
the $G$-isomorphism
$\B: L^2(C) \to \overline \Vpq (\subset L^2(M))$
(see Fact \ref{fact:B}) applied to spherical harmonics. 
The main result of this section is Lemma \ref{lem:intwn}.

We write $v \in M= S^{p-1} \times S^{q-1}
\subset \mathbb R^{p+q}$ as 
\begin{equation}\label{def:coM}
  v =
  {}^t(v_0, v', v'', v_{p+q-1}), \qquad 
  v_0, v_{p+q-1} \in \mathbb{R}, \
  v' \in \mathbb R^{p-1}, \ v'' \in \mathbb R^{q-1}
\end{equation}
satisfying 
$$
  v_0^2+|v'|^2= |v''|^2+v_{p+q-1}^2=1.
$$

\begin{lemSec}
\label{lem:intwn}
Suppose that  $u \in L^2(C)$ is of the form
$$
u(r\omega, r\eta)=f(r)\phi(\omega)\psi(\eta),
$$
for $\phi \in \Har{l}{p-1}$, $\psi \in \Har{k}{q-1}$, 
and $f\in L^2(\mathbb{R}_+, r^{p+q-5}dr)$
with regard to the polar coordinate \eqref{def:pol}.
Then, $\B u$ is reduced to the following integral transform
of one variable:
for $v\in M$ such that $v_0+v_{p+q-1} > 0$,
\begin{multline}
\label{int:T}
  \bigl(\B u \bigr) (v) =
  e^{-\frac{\sqrt{-1}(l+k)\pi}{2}}
  \frac{|v'|^{-\frac{p-3}2} |v''|^{-\frac{q-3}2}}{v_0+v_{p+q-1}}
  \phi\biggl( \frac{v'}{|v'|} \biggr) \psi \biggl( \frac{v''}{|v''|} \biggr) \\
  \times 
  \int_0^\infty
  f(r) J_{\frac{p-3}2+l}\biggl(\frac{2|v'|r}{v_0+v_{p+q-1}}\biggr)
        J_{\frac{q-3}2+k}\biggl(\frac{2|v''|r}{v_0+v_{p+q-1}}\biggr)
  r^\frac{p+q-4}2 dr,
\end{multline}
where $J_\nu(z)$ denotes the Bessel function.

\end{lemSec}

\begin{proof}
We recall from \eqref{def:B} that
 $\B=(\widetilde{\Psi}^\ast)^{-1}\circ \mathcal F^{-1}
\circ T$. 
Let us compute 
$(\mathcal F^{-1} \circ Tu)(z)$ for
$z \in \mathbb R^{p+q-2}$. 
By the formula \eqref{eqn:Q} and \eqref{def:T} for 
$Tu=u \delta(Q)$, we have
\begin{alignat}{2}
  & 2(\mathcal{F}^{-1}\circ Tu)(z) \notag \\
  ={} & (2\pi)^{-\frac{p+q-2}{2}} \int_0^\infty \int_{S^{p-2}} \int_{S^{q-2}}
     f(r)\phi(\omega)\psi(\eta)
     e^{-\sqrt{-1}(\langle z',r\omega \rangle+\langle z'',r\eta \rangle)}
     r^{p+q-5}drd\omega d\eta  \notag \\
  ={} &\sqrt{-1}^{\,-(l+k)} 
    |z'|^{-\frac{p-3}{2}}| z'' |^{-\frac{q-3}{2}}  
     \phi\biggl(\frac{z'}{|z'|} \biggr)
    \psi\biggl(\frac{z''}{|z''|} \biggr)   \notag  \\
  \label{eqn:intwn}
    &\qquad \times \int_0^\infty f(r)J_{\frac{p-3}{2}+l}(r|z'|)
     J_{\frac{q-3}{2}+k}
    (r|z''|)r^{\frac{p+q-4}{2}}dr.
\end{alignat}
Here, in \eqref{eqn:intwn}, we used the following formula (see 
e.g. \cite[Introduction, Lemma 3.6]{xHe}):
\begin{equation}
  \int_{S^{m-1}} e^{\sqrt{-1} \lambda \langle\eta, \omega\rangle} \phi(\omega) d\omega
= (2\pi)^\frac{m}2 \sqrt{-1}^l \phi(\eta) 
   \frac{J_{l-1+\frac{m}2}(\lambda)}{\lambda^{\frac{m}2-1}}.
\end{equation} 
Then, by the definition \eqref{def:Psiinv} of the pull-back 
$(\Tilde \Psi^\ast)^{-1}$,
we get the desired result \eqref{int:T}.
\end{proof}

\section{$K$-finite vectors $f_{l,k}$ in $L^2(C)$}
\protect\index{A}{f0lkr@$f_{l,k}(r)$}%
\label{subsec:flk}

In this section,
we collect basic results on the real analytic functions $f_{l,k}$ 
defined in \eqref{def:f}.

\begin{lemSec}\label{lem:flkL2}
$f_{l,k} \in L^2(\mathbb{R}_+, \frac{1}{2}r^{p+q-5} dr)$.
\end{lemSec}

\begin{proof}
The $K$-Bessel function decays exponentially at infinity.
The 
\index{B}{asymptotic behavior!$\flk@---, $f_{l,k}$}%
asymptotic formula (see Fact \ref{fact:Bas} (2) in Appendix)
implies
\begin{equation}\label{eqn:flkinfty}
f_{l,k}(r) \sim
  \begin{cases}
     c \ r^{-\frac{q-2}{2}+l} e^{-2r}
     &\text{Case 1},
     \\
     c \ r^{-\frac{p-2}{2}+k} e^{-2r}
     &\text{Case 2},
  \end{cases}
\end{equation}
as $r\to \infty$.
On the other hand, 
since $\widetilde{K}_\nu(r) = O(r^{-2\nu})$ as $r$ tends to $0$ 
(see \eqref{eqn:Ktild0}), 
\begin{equation}\label{eqn:flk0}
 f_{l,k} =
  \begin{cases}
     O(r^{l-k-q+3})  &\text{Case 1}, \\
     O(r^{-l+k-p+3})  &\text{Case 2}, 
  \end{cases}
\end{equation}
as $r$ tends to $0$.
In either case,
$f_{l,k}=O(r^{\frac{-p-q}{2}+3})$ by the definition \eqref{def:Cases}
of Cases 1 and 2.
Hence, we have
$f_{l,k}\in L^2(\mathbb R_+, \frac{1}{2}r^{p+q-5}dr)$ for $p+q > 4$. 
\end{proof}

We shall give a finer estimate on the derivatives of $f_{l,k}$ in
Lemma \ref{lem:flkas}.
The explicit formula of the $L^2$-norm of $f_{l,k}$ is obtained by
the integration formula
\eqref{eqn:KKL2} of $K$-Bessel functions as follows:

\begin{proposition}\label{prop:flkL2}
\begin{align}
& \| f_{l,k}\|^2_{L^2(\mathbb{R}_+, \frac{1}{2}r^{p+q-5}dr)}
\nonumber
\\[1ex]
&=\begin{cases}\displaystyle
    \frac{\Gamma(\frac{p-1}{2}+l)^2
          \Gamma(\frac{p+q-4}{2}+l+k)
          \Gamma(\frac{p-q+2}{2}+l-k)}
         {16\Gamma(p-1+2l)}
    &\text{Case 1},
  \\[2ex]
    \displaystyle
    \frac{\Gamma(\frac{q-1}{2}+k)^2
          \Gamma(\frac{p+q-4}{2}+k+l)
          \Gamma(\frac{q-p+2}{2}+k-l)}
         {16\Gamma(q-1+2k)}
    &\text{Case 2}.
  \end{cases}
\label{eqn:flkL2}
\end{align}
\end{proposition}

Here is another integral formula that we use later:
\begin{lemSec}\label{lem:flkMellin}
\begin{align*}
&\int_0^\infty r^{\frac{p+q-6}{2}+\sqrt{-1}\xzeta} f_{l,k}(r) dr
\\
&=\begin{cases}
    \frac{1}{4}\Gamma(\frac{p+q-4}{4}+\frac{l+k+\sqrt{-1}\xzeta}{2})
               \Gamma(\frac{p-q}{4}+\frac{l-k+1+\sqrt{-1}\xzeta}{2})
    &\text{Case 1},
    \\
    \frac{1}{4}\Gamma(\frac{p+q-4}{4}+\frac{k+l+\sqrt{-1}\xzeta}{2})
               \Gamma(\frac{q-p}{4}+\frac{k-l+1+\sqrt{-1}\xzeta}{2})
    &\text{Case 2}.
  \end{cases}
\end{align*}
\end{lemSec}

\begin{proof}[Proof of Lemma \ref{lem:flkMellin}]
Apply the formula \eqref{eqn:B1} of the Mellin transform of $K$-Bessel
functions.
\end{proof}

In order to compute $\B u_{l,k}$ explicitly by using
 the integral formula \eqref{int:T}, 
we need yet another lemma:
\begin{lemSec}
\label{lem:prK} For a pair $(l,k)\in \mathbb N^2$,
let $f_{l,k}$ be the function on $C$ defined in
\eqref{def:f}. With respect to the coordinate
${}^t(v_0, v', v'', v_{p+q-1}) \in M 
\subset \mathbb R^{p+q}$ (see \eqref{def:coM}), the integral
\begin{equation}\label{eqn:prK}
\int_0^\infty
  f_{l,k}(r) J_{\frac{p-3}2+l}\biggl(\frac{2|v'|r}{v_0+v_{p+q-1}}\biggr)
        J_{\frac{q-3}2+k}\biggl(\frac{2|v''|r}{v_0+v_{p+q-1}}\biggr)
  r^\frac{p+q-4}2 dr
\end{equation}
is equal to:
\begin{align*}
 & \frac{\Gamma(\frac{p-q}2+l-k+1)}
           {2^{\frac{p-q}2+l-k+3}\sqrt{\pi} }
  (v_0+v_{p+q-1}) 
  |v'|^{\frac{p-3}2+l} |v''|^{\frac{q-3}2+k}
  \widetilde{C}_{\frac{p-q}2+l-k}^{\frac{q-2}2+k}(v_{p+q-1}) 
 \quad \text{Case 1,}  \\
 & \frac{\Gamma(-\frac{p-q}2-l+k+1)}
           {2^{-\frac{p-q}2-l+k+3}\sqrt{\pi}}
  (v_0+v_{p+q-1}) 
  |v'|^{\frac{p-3}{2}+l} |v''|^{\frac{q-3}{2}+k}
  \widetilde{C}_{-(\frac{p-q}2+l-k)}^{\frac{p-2}2+l}(v_0)
  \quad \text{Case 2}.
\end{align*}
\end{lemSec}

\begin{proof}[Proof of Lemma \ref{lem:prK}]
We treat Case 1 first.
By the change of variables
 $t:=2r$, 
 the integral \eqref{eqn:prK} amounts to
$$
 \frac{1}{2^{\frac{p+1}2+l}}
 \int_0^\infty t^{\frac{p-1}2+l}
 J_{\frac{p-3}2+l}\biggl(\frac{|v'|t}{v_0+v_{p+q-1}}\biggr)
 J_{\frac{q-3}2+k}\biggl(\frac{|v''|t}{v_0+v_{p+q-1}}\biggr) 
 K_{\frac{q-3}{2}+k}(t) \, dt.
$$
Applying Lemma \ref{lem:Han} with 
$$
\mu:=\frac{p-3}2+l, \quad \nu:=\frac{q-3}2+k, \quad
\cos \theta:=v_0, \quad \cos \phi:=v_{p+q-1},
$$
we get the formula in Case 1.

The proof for Case 2 goes similarly. 
In this case, 
the integral amounts to
\begin{equation*}
\frac{1}{2^{\frac{q+1}{2}+k}}\int_0^\infty t^{\frac{q-1}{2}+k}
J_{\frac{p-3}{2}+l} \Bigl(\frac{|v'|t}{v_0+v_{p+q-1}}\Bigr)
J_{\frac{q-3}{2}+k} \Bigl(\frac{|v''|t}{v_0+v_{p+q-1}}\Bigr)
K_{\frac{p-3}{2}+l} (t)dt
\end{equation*}
by the change of variables $t:=2r$.
Now,
we substitute 
$\mu:=\frac{q-3}2+k$, $\nu:=\frac{p-3}2+l$,
$\cos \theta:=v_{p+q-1}$, and $\cos \phi:=v_0$
into \eqref{eqn:Han}.
\end{proof}

\section{Proof of Theorem \ref{thm:K}}
\label{subsec:prK}

In this section,
we complete the proof of Theorem \ref{thm:K}.
\begin{proof}[Proof of Theorem \ref{thm:K}]
1)
By the isomorphism
$$
L^2(C) \simeq L^2(\mathbb{R}_+, \frac{1}{2} r^{p+q-5} dr)
\hatotimes L^2(S^{p-2}) \hatotimes L^2(S^{q-2})
\quad\text{(see \eqref{eqn:LtwoC})}
$$
in the polar coordinate,
the first statement follows immediately from Lemma \ref{lem:flkL2}.

4)
By \eqref{int:T} and Lemma \ref{lem:prK}, we have
\begin{equation*}
 \B u_{l,k}(v)=
 \begin{cases}
  c_{l,k} \Gamma(\frac{p-2}{2}+l)|v'|^l \phi(\frac{v'}{|v'|})|v''|^k
  \psi(\frac{v''}{|v''|})
  \widetilde{C}_{\frac{p-q}2+l-k}^{\frac{q-2}2+k}(v_{p+q-1})
  &\text{Case 1}, \\
  c_{l,k} \Gamma(\frac{q-2}{2}+k)|v'|^l  \phi(\frac{v'}{|v'|})  
  \widetilde{C}_{-(\frac{p-q}2+l-k)}^{\frac{p-2}2+l}(v_0)|v''|^k
  \psi(\frac{v''}{|v''|}) 
  &\text{Case 2},
 \end{cases}
\end{equation*}
where $c_{l,k}$ is the constant defined in \eqref{def:clk}
and $v={}^t(v_0,v',v'',v_{p+q-1})\in S^{p-1}\times S^{q-1}$.
Now we use the definition \eqref{def:I} that
\begin{equation*}
I_{i\to j}^m (\phi) (x_0,x')
= |x'|^i \phi \Bigl(\frac{x'}{|x'|}\Bigr)
  \tilC_{j-i}^{\frac{m-2}{2}+i} (x_0)
\end{equation*}
for $(x_0,x')\in S^{m-1}$,
and in particular,
for $i=j$,
\begin{equation*}
I_{j\to j}^m (\phi) (x_0,x')
= \Gamma \Bigl(\frac{m-2}{2}+j\Bigr) |x'|^j \phi \Bigl(\frac{x'}{|x'|}\Bigr)
\end{equation*}
by \eqref{eqn:Ijj}.
Thus, the formula \eqref{eqn:Tu} follows.

2)
The first statement is obvious. 
Since
$I_{i\to j}^m$ maps
$\mathcal{H}^i(\mathbb{R}^{m-1})$ into 
$\mathcal{H}^j(\mathbb{R}^m)$
(see Fact \ref{fact:H}),
\eqref{eqn:Tu} implies 
$$
 \B H_{l,k} \subset \Har{a(l,k)}{p}\otimes 
 \Har{a(l,k)+\frac{p-q}2}{q}.
$$

Hence, we have proved the second statement.

3) By \eqref{def:vpi}, the unitary inversion operator
  $\varpi^{p,q}(w_0)$ on $\overline{V^{p,q}}$ is given by
$$
  \Bigl(\varpi^{p,q} (w_0) h\Bigr)(v', v'') = h(v', -v'').
$$
On the other hand, it is easy to see
 $h(-x)=(-1)^jh(x)$ for $h\in \Har{j}{q}$ 
(see Appendix \ref{subsec:H}, \textbf{H\ref{item:H1}}).
Therefore, 
$\varpi^{p,q}(w_0)$ acts on each $K$-type component
$\Har{a}{p} \otimes \Har{a+\frac{p-q}2}{q}$ 
by a scalar multiple $(-1)^{a+\frac{p-q}2}$.
Since
$\mathcal{T}: L^2(C) \to \overline{V^{p,q}}$
intertwines $\pi$ and $\varpi^{p,q}$,
$\pi(w_0)$ acts on $H_{l,k}$ as the scalar
$(-1)^{a+\frac{p-q}{2}}$
because $H_{l,k}$ belongs to the $K$-type
 $\mathcal{H}^{a(l,k)}(\mathbb{R}^p) \otimes 
  \mathcal{H}^{a(l,k)+\frac{p-q}{2}}(\mathbb{R}^q)$.
Thus, the third statement is proved.

5)
We claim that the following formula holds:
\begin{equation}\label{eqn:Tulk2}
\|\mathcal{T}u_{l,k}\|_{L^2(S^{p-1}\times S^{q-1})}^2
= b_{l,k} \Gamma(k+l+\frac{p+q}{2}-2)
  \|\phi\|_{L^2(S^{p-1})}^2 \|\psi\|_{L^2(S^{q-1})}^2,
\end{equation}
where
\begin{equation*}
b_{l,k}
= \begin{cases}
    \displaystyle
    \frac{\pi\Gamma(p-2+2l)\Gamma(\frac{p-q}{2}+l-k+1)}
         {2^{2p+4l} \Gamma(l+\frac{p}{2})^2}
   &\text{Case 1},
  \\[1.5ex]
    \displaystyle
    \frac{\pi\Gamma(q-2+2k)\Gamma(k-l-\frac{p-q}{2}+l)}
         {2^{2q+4k} \Gamma(k+\frac{q}{2})^2}
   &\text{Case 2}.
  \end{cases}
\end{equation*}
To see this,
we begin with \eqref{eqn:Tu}:
\begin{equation*}
\|\mathcal{T}u_{l,k}\|_{L^2(S^{p-1}\times S^{q-1})}^2
= |c_{l,k}|^2 \|I_{l\to a(l,k)}^p (\phi) \|_{L^2(S^{p-1})}^2
  \|I_{k\to a(l,k)+\frac{p-q}{2}}^q (\psi)\|_{L^2(S^{q-1})}^2.
\end{equation*}
By \eqref{eqn:Tu} and the norm formula \eqref{eqn:Iijnorm},
the right-hand side amounts to 
\begin{align*}
& |c_{l,k}|^2 \frac{2^{3-p-2l} \pi \Gamma(p-2+l+a(l,k))}
                   {(a(l,k)-l)! (a(l,k)+\frac{p-2}{2})}
  \|\phi\|_{L^2(S^{p-1})}^2
\\
&\times
  \frac{2^{3-q-2k} \pi \Gamma(q-2+k+a(l,k)+\frac{p-q}{2})}
       {(a(l,k)+\frac{p-q}{2}-k)! (a(l,k)+\frac{p-q}{2}+\frac{q-2}{2})}
  \|\psi\|_{L^2(S^{q-1})}^2.
\end{align*}
Now, substituting \eqref{def:alk} and \eqref{def:clk} into this formula,
we get \eqref{eqn:Tulk2} by elementary computations.

Finally, comparing \eqref{eqn:Tulk2} with Proposition
\ref{prop:flkL2}, 
we obtain \eqref{eqn:TulkL2} by the Gauss duplication formula
\eqref{eqn:dup} of gamma functions.
Hence, Theorem \ref{thm:K} is proved.

\end{proof}

\chapter{Radial part of the inversion}
\label{sec:C}

The goal of this chapter is 
to find the \lq radial part\rq{} $T_{l,k}$ of the unitary inversion operator 
$\mathcal{F}_C : L^2(C) \to L^2(C)$.
The main result of this chapter is Theorem \ref{thm:C}.

\section{Result of this chapter}
\label{subsec:C}

Suppose  $p\ge q\ge 2$ and $p+q\ge 6$.
We recall from Section~\ref{subsec:Sch} that the polar coordinate
of the (generalized) light cone $C$:
$$
 \mathbb{R}_+ \times S^{p-2} \times S^{q-2} 
 \simeq C, \quad (r,\omega,\eta) \mapsto (r\omega, r\eta)
$$
induces the isomorphism of Hilbert spaces:
$$
 L^2(C) \simeq L^2(\mathbb{R}_+, \frac{1}{2}r^{p+q-5}dr)
 \hatotimes L^2(S^{p-2}) \hatotimes L^2(S^{q-2}).
$$
This isomorphism respects the action of the compact group
$$
K':=K\cap \Mmax_+\simeq O(p-1)\times O(q-1).
$$
Since the Hilbert space $L^2(S^{p-2})$ is decomposed into the direct
Hilbert sum of spherical harmonics (see \textbf{H\ref{item:H6}} in Appendix
\ref{subsec:H}):
$$
L^2(S^{p-2}) \simeq \sideset{}{^\oplus}\sum\limits_{j=0}^\infty
\mathcal{H}^j(\mathbb{R}^{p-1}),
$$
and likewise for $L^2(S^{q-2})$,
we have a decomposition of the Hilbert space $L^2(C)$ into the discrete direct sum:
\begin{equation}
\label{eqn:iso}
  L^2(C) \simeq \sideset{}{^\oplus}\sum_{l,k=0}^\infty
  L^2(\mathbb R_+, \frac{1}{2}r^{p+q-5}dr) \otimes \Har{l}{p-1} \otimes
  \Har{k}{q-1}.
\end{equation}
Each summand of \eqref{eqn:iso} is a $K'$-isotypic component.

For each $(l,k)\in\mathbb{N}^2$,
we introduce real analytic function
$K_{l,k}(t)$ by
\begin{equation}
\label{def:Klk}
\index{A}{Klkt@$K_{l,k}(t)$|main}%
  K_{l,k}(t):= 
  \begin{cases}
   4 (-1)^{l+\frac{p-q}2}
   G_{04}^{20}( t^2 | \frac{l+k}2, \frac{-q+3+l-k}2, \frac{-p-q+6-l-k}2, 
                           \frac{-p+3-l+k}2 ) 
   &\text{Case 1},  \\
   4(-1)^k
   G_{04}^{20}( t^2 | \frac{l+k}2, \frac{-p+3-l+k}2, \frac{-p-q+6-l-k}2, 
                           \frac{-q+3+l-k}2 )  
   &\text{Case 2}.
   \end{cases}
\end{equation}
Here, 
$G^{20}_{04}(x|b_1,b_2,b_3,b_4)$ denotes Meijer's $G$-function
(see Appendix \ref{subsec:G} for definition).
For the definition of Cases 1 and 2 with regard to the parameter
 $(l,k)\in\mathbb{N}^2$, 
we recall \eqref{def:Cases}, namely,
\begin{itemize}
\item[]
Case 1: $\frac{p-q}{2}+l-k \ge 0$,
\item[]
Case 2: $\frac{p-q}{2}+l-k \le 0$.
\end{itemize}
The above formulas give the same result in Cases 1 and 2 if
$\frac{p-q}{2}+l-k=0$.
Later, we shall give an integral expression of $K_{l,k}$ 
by means of the Mellin--Barnes type integral
(see Lemma \ref{lem:Klkint}).
The latter expression is more natural in the sense that the formula
 is independent of Cases 1 and 2.

\begin{thmSec}
\label{thm:C}
 {\rm 1)} The unitary inversion operator
$\mathcal{F}_C = \pi(w_0)$ preserves each summand of \eqref{eqn:iso},
 on which $\pi(w_0)$ acts as a form $T_{l,k} \otimes \operatorname{id} 
 \otimes \operatorname{id}$. Here,  
\index{A}{Tlk@$T_{l,k}$|main}%
$T_{l,k}$ 
is a unitary operator
 on the Hilbert space 
$L^2(\mathbb R_+, \frac{1}{2}r^{p+q-5}dr)$.

 {\rm 2)} For each $l, k \in \mathbb N$, 
the unitary operator $T_{l,k}$ 
is given by the integral transform against the real analytic function
$K_{l,k}$ (see \eqref{def:Klk} for definition):
  \begin{equation}
  \label{int:C}
    (T_{l,k}f)(r)= \frac{1}{2}\int_0^\infty K_{l,k}(rr') f(r') r'^{p+q-5}dr'.
  \end{equation}
\end{thmSec}

\begin{remNonumber}[Case $q=2$]
\label{rem:Cq2}
Theorem \ref{thm:C} covers the case $q=2$,
where the representation $\pi$ is essentially a
\index{B}{highest weight module}%
highest weight module 
(see Section \ref{subsec:1.6} for the $n_2=1$ case).  
In this case, 
$\Har{k}{1}$ is non-zero only if $k=0$ or $1$ (see Appendix \ref{subsec:H}
for convention).
As we saw in Remark \ref{rem:3.1.3}, 
the pair $(l,k)$ belongs to Case~1 for any $l\in\mathbb{N}$
because the inequality
$\frac{p-q}{2} + l-k \ge 0$ is implied by
$p \ge 6-q = 4$.
Hence,
by the definition \eqref{def:Klk} of $K_{l,k}(t)$, we have
\begin{align*}  
& K_{l,0}(t)=4(-1)^{l+\frac{p-2}2}
  G_{04}^{20}(t^2 | \frac{l}2, \frac{l+1}2, \frac{-p+4-l}2, \frac{-p+3-l}2),
\\
& K_{l,1}(t)=4(-1)^{l+\frac{p-2}2}
  G_{04}^{20}(t^2 | \frac{l+1}2, \frac{l}{2}, \frac{-p+3-l}2, \frac{-p+4-l}2).
\end{align*} 
In view of the symmetric property of the $G$-function:
$$
  G_{04}^{20}(x| b_1, b_2, b_3, b_4)
=G_{04}^{20}(x| b_2, b_1, b_4, b_3),
$$
the above formulas show
$
 K_{l,0}(t)=K_{l,1}(t)
$.
Applying
 the reduction formula \eqref{eqn:GJ} of
the $G$-function of the form 
$G_{04}^{20}(x| a,a+\frac{1}{2},b,b+\frac{1}{2})$,
we get
\begin{equation}\label{eqn:rdB}
 K_{l,0}(t)=K_{l,1}(t)=
 4(-1)^{l+\frac{p-2}2} t^{-\frac{p-3}2}
  J_{p-3+2l}(4\sqrt{t}). 
\end{equation}
Thus, the integral transform $T_{l,k}$ $(k=0,1)$ is 
 the Hankel transform on $\mathbb R_+$
(after a suitable change of variables).
Therefore, 
Theorem \ref{thm:C} in the case $q=2$
gives the same result with
\cite[Theorem 6.1.1]{xkmano2},
but the proof here is different from that of \cite{xkmano2}.
\end{remNonumber}

\begin{remNonumber}[Comparison with the 
\index{B}{Weil representation}%
Weil representation]
 \label{rem:4.1.3}
We record here the corresponding result for the Schr\"{o}dinger model
$L^2(\mathbb{R}^n)$ of the 
Weil representation $\varpi$ of $\Mp(n, \mathbb R)$
(see also \cite[Remark 6.1.3]{xkmano2}).
We adopt the same normalization of the Weil representation with
\cite{xKaVe}.
Then, the unitary inversion operator $\varpi(w'_0)
= e^{\frac{\sqrt{-1}n\pi}{4}}
\index{A}{FRn@$\mathcal{F}_{\mathbb{R}^n}$}%
\mathcal{F}_{\mathbb{R}^n}$.

{\rm 1)}
According to the $O(n)$-isotypic decomposition of $L^2(\mathbb R^n)$
$$
 L^2(\mathbb R^n)\simeq \sideset{}{^\oplus}
\sum_{l=0}^\infty L^2(\mathbb R_+, r^{n-1}dr)
\otimes 
\index{A}{HjRm@$\mathcal{H}^j(\mathbb{R}^m)$}%
\Har{l}{n},
$$
the unitary inversion operator 
\index{A}{w0'@$w'_0$}%
$\varpi(w'_0)$ decomposes as 
$$
 \varpi(w'_0)\simeq \sideset{}{^\oplus}
\sum_{l=0}^\infty T'_l \otimes \operatorname{id},
$$
by a countable family of unitary operators 
$T'_l$ $(l \in \mathbb N)$ on $L^2(\mathbb R_+, r^{n-1}dr)$.

{\rm 2)} The unitary operator $T'_l$ 
is given by the 
\index{B}{Hankel transform}%
Hankel transform 
$$
 (T'_l f)(r)=\int_0^\infty K'_l(rr')f(r')r'^{n-1}dr',
$$
where the kernel $K'_l(t)$ is defined by
$$
 K'_l(t):= e^{-\frac{\sqrt{-1}(n-1+2l)}{4}\pi}
 \, t^{-\frac{n-2}2}J_{\frac{n-2+2l}2}(t).
$$
\end{remNonumber}

Returning to Theorem \ref{thm:C},
we remark that the group law $w_0^2=1$ in $G$
implies $\pi(w_0)^2= \operatorname{id}$, 
and consequently $T_{l,k}^2= \operatorname{id}$ for every
$l,k\in\mathbb{N}$.
Further, 
 as $\pi(w_0)$ is a unitary operator on $L^2(C)$, 
so is its restriction $T_{l,k}$ on $L^2(\mathbb R_+, \frac 1 2 r^{p+q-5} dr)$
 for every $l,k$.  
Hence, 
Theorem \ref{thm:C} gives
a simple group theoretic proof for 
the Plancherel and reciprocal formulas
on the integral transform involving the 
$G$-functions due to C. Fox \cite{xFo}:

\begin{corollary}[%
\index{B}{Plancherel formula!G-function@---, $G$-function}%
Plancherel formula]
\label{cor:C1}
Let $b_1, b_2, \gamma$ be half-integers such that $b_1 \ge 0$, 
\thinspace
$\gamma \ge 1$, 
\thinspace
$\frac{1-\gamma}2 \le b_2 \le \frac{1}2+b_1$.
Then the integral transform 
$$
  S_{b_1,b_2,\gamma}: f(x) \mapsto \frac{1}\gamma 
\int_0^\infty G_{04}^{20}
((xy)^{\frac{1}\gamma} | b_1, b_2, 1-\gamma-b_1, 1-\gamma-b_2)f(y)dy
$$
is a unitary operator on $L^2(\mathbb R_+)$. 
In particular, we have
\begin{equation} \label{eqn:Sunitary}
  \| S_{b_1,b_2,\gamma} f\|_{L^2(\mathbb R_+)}
=\|f\|_{L^2(\mathbb R_+)}.
\end{equation}
\end{corollary}

\begin{proof}
Set $b_1:=\frac{l+k}2$, 
\thinspace
$b_2:=\frac{-q+3+l-k}2$, 
\thinspace
$\gamma:=\frac{p+q-4}2$,
\thinspace
$x=r^{2\gamma}$, 
\thinspace 
$y=r'^{2\gamma}$.
Then, the assertion is equivalent to the fact that $T_{l,k}$ is
a unitary operator on 
$L^2(\mathbb{R}_+, \frac{1}{2}r^{p+q-5}dr)$.
\end{proof}

\begin{corollary}[%
\index{B}{reciprocal formula!G-function@---, $G$-function}%
Reciprocal formula]
\label{cor:C2}
Retain the notation and the assumption as in Corollary \ref{cor:C1}. 
Then, the unitary operator $S_{b_1,b_2,\gamma}$ is of order two in 
$L^2(\mathbb R_+)$, that is, the following 
reciprocal formula holds for $f \in L^2(\mathbb R_+)$:
\begin{align*}
  (S_{b_1,b_2,\gamma}f)(x)
 &=
 \frac{1}{\gamma} \int_0^\infty G_{04}^{20}
((xy)^{\frac{1}\gamma}|b_1, b_2, 1-\gamma-b_1,1-\gamma-b_2) f(y)dy 
\\
  f(x)
 &=
  \frac{1}{\gamma} \int_0^\infty G_{04}^{20}
((xy)^{\frac{1}{\gamma}}|b_1,b_2,1-\gamma-b_1,1-\gamma-b_2)
     (S_{b_1,b_2,\gamma}f)(y)dy
\end{align*}
\end{corollary}

\begin{remNonumber}[%
\index{B}{Fox's reciprocal formula}%
Fox's reciprocal formula]
\label{rem:4.1.6}
The reciprocal formula 
for the $G$-transform 
was found by C. Fox \cite{xFo} for 
the following generality:
$$
 G_{2n, 2m}^{m,n}\Bigl(x \Big|
  \begin{matrix} a_1, \cdots, a_n, 1-\gamma-a_1, \cdots, 
                1-\gamma-a_n  \\
                b_1, \cdots, b_m, 1-\gamma-b_1, \cdots, 
                1-\gamma-b_m
  \end{matrix}  
 \Bigr).
$$ 
Corollary \ref{cor:C2} is a special case of Fox's formula
corresponding to the case $(m, n)=(2,0)$.
On the other hand,
it follows from Remark \ref{rem:Cq2} that Corollary \ref{cor:C2}
in the special case when $q=2$ yields
the classic reciprocal formula of the Hankel transform
(\cite[\S 14.3]{xWa}, see also \cite{xkmano2} for the connection with
the representation of the conformal group $O(p,2)$).
Our approach gives a new representation theoretic interpretation
(and also a proof) of these
formulas. 
\end{remNonumber}

\begin{remNonumber}[Fourth order differential equation]
Let $\Omega$ be the Casimir element for the Lie algebra $\mathfrak k$. 
Since an element of the Lie algebra $\mathfrak{n}^{\max}$ acts on
smooth vectors of $L^2(C)$ as a differential operator of second order 
(see \eqref{eqn:Dj}),
 the Casimir operator $d\pi(\Omega)$ acts 
 as a differential operator on $C$ of fourth order.
Let us examine what information on $K_{l,k}(t)$ we can obtain from the
Casimir operator $d\pi(\Omega)$.

Since $\operatorname{Ad}(w_0) \Omega =\Omega$,
we have the commutative relation
\begin{equation}\label{eqn:4diffeq}
\pi(w_0) \circ d\pi(\Omega)= d\pi(\Omega) \circ \pi(w_0).
\end{equation}
On the other hand, $\Omega$ commutes with all the elements in $\mathfrak k$, 
in particular with $\mathfrak k'$.
This implies that the differential operator $d\pi(\Omega)$ preserves each 
summand of \eqref{eqn:iso}, and the identity \eqref{eqn:4diffeq} can be
regarded as the identity on each summand of \eqref{eqn:iso}.
Hence,
\eqref{eqn:4diffeq} leads us to the fact that
 the kernel function $K_{l,k}(t)$ 
for $\pi(w_0)$ solves a certain differential equation of order four
for each $(l,k)$.
In view of the formula \eqref{def:Klk} of $K_{l,k}$,
this corresponds to the fact that 
Meijer's $G$-function 
$G_{04}^{20}(x \mid b_1,b_2,b_3,b_4)$
satisfies the fourth order differential equation
(see \eqref{eqn:diffeqG24}):
\begin{equation*}
   \prod_{j=1}^4 (x\frac{d}{dx}-b_j)u = 0.
\end{equation*}

For $q=2$, the situation becomes simpler because
the minimal representation $\pi$ is essentially a 
\index{B}{highest weight module}%
highest weight module.
In this case, the Lie algebra 
$\mathfrak k=\mathfrak{so}(p)\oplus \mathfrak{so}(q)$
contains the center $\mathfrak{so}(q)=\mathfrak{so}(2)$.
Taking a normalized generator $Z$ in $\mathfrak{so}(2)$,
we can show
$$
 d\pi(\Omega)=d\pi(Z)^2+\operatorname{constant}.
$$ 
We note that the differential operator $d\pi(Z)$ is of order two, and 
this in turn corresponds to the fact that the kernel function $K_{l,k}$ 
reduces to the Bessel function (see \eqref{eqn:rdB}) 
solving the Bessel differential equation of order two
(Appendix \ref{subsec:B}).
\end{remNonumber}

The rest of this chapter is devoted to the proof 
of Theorem \ref{thm:C}.
The key properties of the element $w_0$ and the unitary inversion
operator $\pi(w_0)$
that we use in the proof are listed as follows:

1) $w_0$ commutes with $K'$.

2) $\operatorname{Ad}(w_0)H = -H$ (see \eqref{eqn:wa}).

3) $\pi(w_0)|_{H_{l,k}}=\pm \operatorname{id}$
~(Theorem \ref{thm:K} (3)).

The first property (1) gives immediately the proof 
of Theorem \ref{thm:C} (1)
(see Section \ref{subsec:Tlk}).

The second property (2) implies that the kernel function 
$K_{l,k}$ of the unitary operator $T_{l,k}$ is a function essentially of 
one variable $rr'$ (see Section \ref{subsec:F}).

The non-trivial part of Theorem \ref{thm:C} is to prove that 
this one variable function is given by
Meijer's $G$-function $G_{04}^{20}$, namely, by the formula \eqref{def:Klk}.
The property (3) is used for the proof of this part.
Our trick here is
based on  Fact \ref{obsv:skt},
which assures that the desired formula  
\eqref{int:C} holds for any function $f$ once we can prove \eqref{int:C} 
holds for a specific function $f=f_{l,k}$ 
(a $(K, K')$-\textit{skeleton function}).
Actual computations for this are carried out in Section \ref{subsec:prC}.
A technical ingredient of the proof is to show that $K_{l,k}$ is a
tempered distribution 
(see Proposition \ref{prop:Klk}).

Summing up the formulas for $K_{l,k}$ over $l,k$,
we shall give a proof of 
an integral formula of $\pi(w_0)$ for an arbitrary $L^2$-function on $C$
in the next Chapter \ref{sec:A} (see Theorem \ref{thm:A}).

\section{Proof of Theorem \ref{thm:C} (1)}
\label{subsec:Tlk}

As we have explained at the beginning of Section \ref{subsec:C},
the decomposition (\ref{eqn:iso}) of the Hilbert space $L^2(C)$ corresponds to
the branching laws of 
the restriction of the unitary representation $(\pi,L^2(C))$ 
from $G$ to $K'\simeq O(p-1) \times O(q-1)$.
Thanks to the group law
$w_0 m = m w_0$
(see \textbf{I\ref{item:I4}} in Section \ref{subsec:w}), 
the relation $\pi(w_0) \circ \pi(m) = \pi(m) \circ 
\pi(w_0)$ holds for all $m \in K'$.
Therefore, the unitary inversion operator
$\pi(w_0)$ preserves each summand of the decomposition
(\ref{eqn:iso}). 

We now observe that the group $K'$
acts irreducibly on $\Har{l}{p-1} \otimes \Har{k}{q-1}$
(see Appendix \ref{subsec:H}, \textbf{H\ref{item:H2}}).
Therefore, it follows from Schur's lemma that 
$\pi(w_0)$ acts on each summand:
\begin{equation*}
    L^2(\mathbb{R}_+, \frac{1}{2}r^{p+q-5}dr)
    \otimes \mathcal{H}^l(\mathbb{R}^{p-1})
    \otimes \mathcal{H}^k(\mathbb{R}^{q-1})
\end{equation*}
as the form
$T_{l,k} \otimes \operatorname{id}
\otimes \operatorname{id}$ for some operator
$T_{l,k}$ on $L^2(\mathbb R_+, \frac{1}{2}r^{p+q-5}dr)$. 
Such an operator 
$T_{l,k}$ must be unitary because $\pi(w_0)$ is unitary.
Now Theorem \ref{thm:C} (1) has been proved.
\qed

\section{Preliminary results on multiplier operators}
\label{subsec:multop}

We recall the classic theory of multipliers on $\mathbb{R}$,
followed by the observation that the multiplier is determined by the
action on an (appropriate) single
function. 

We write $l(t)$ $(t\in\mathbb{R})$ for the translation operator on
$L^2(\mathbb{R})$, namely,
\begin{equation}\label{eqn:translt}
   (l(t)F)(x) := F(x-t)
   \qquad\text{for $F\in L^2(\mathbb{R})$.}
\end{equation}
An operator $B$ on $L^2(\mathbb{R})$ is called
\textit{translation invariant} if
\begin{equation*}
   B \circ l(t) = l(t) \circ B
   \qquad\text{for any $t\in \mathbb{R}$.}
\end{equation*}

We write 
\index{A}{S1@$\mathcal{S}(\mathbb{R}^n)$}%
$\mathcal{S}(\mathbb{R})$ for the space of rapidly decreasing
$C^\infty$-functions on $\mathbb{R}$
(the Schwartz space endowed with the Fr\'{e}chet
topology), 
and 
\index{A}{S1'@$\mathcal{S}'(\mathbb{R}^n)$}%
$\mathcal{S}'(\mathbb{R})$ for its dual space consisting of
\index{B}{tempered distribution}%
tempered distributions. 
Then, the Fourier transform
$\mathcal{F}: \mathcal{S}(\mathbb{R}) \to \mathcal{S}(\mathbb{R})$,
$g \mapsto \frac{1}{\sqrt{2\pi}} \int_{\mathbb{R}} g(x)
 e^{\sqrt{-1}x\xi} dx$
(see \eqref{eqn:Fourier}) extends to $\mathcal{S}'(\mathbb{R})$ by
\begin{equation} \label{eqn:FourierSp}
   \langle U,\overline{g} \rangle
   = \langle U, \overline{\mathcal{F}g} \rangle
   \qquad\text{for $U \in \mathcal{S}'(\mathbb{R})$,
              $g\in\mathcal{S}(\mathbb{R})$}. 
\end{equation}
For $U \in \mathcal{S}'(\mathbb{R})$ and $f\in\mathcal{S}(\mathbb{R})$, 
we have 
\begin{equation*}
   l(t) (U\ast f) = U\ast (l(t)f)
   \quad\text{for any $t\in\mathbb{R}$}.
\end{equation*}
Hence,
$f \mapsto U_* f$ is translation invariant.
Furthermore,
if $\mathcal{F}U$ is a bounded measurable function,
then $f \mapsto U\ast f$ extends to a bounded operator on
$L^2(\mathbb{R})$ and its operator norm is given by
$\sqrt{2\pi}\|\mathcal{F}U\|_{L^\infty(\mathbb{R})}$
because 
\begin{equation*}
   \|U\ast f\|_{L^2(\mathbb{R})}
   = \|\mathcal{F}(U\ast f)\|_{L^2(\mathbb{R})}
   = \sqrt{2\pi} \|(\mathcal{F}U)(\mathcal{F}f)\|_{L^2(\mathbb{R})}.
\end{equation*}
Conversely, the following theorem also holds:
\begin{fact}\label{obsv:skt}
Let $B$ be a bounded, translation invariant operator 
on $L^2(\mathbb R)$.
Then we have

{\rm 1)} There exists a unique tempered distribution
$U$ whose Fourier transform $\mathcal F U$ is a bounded measurable
function 
such that $B f=U\ast f$
 for any $f \in L^2(\mathbb R)$.

{\rm 2)} 
If, moreover, $B f_0=f_1$ for 
some function $f_0 \in L^2(\mathbb{R})$
such that 
$\mathcal F f_0(\xi) \ne 0$ for a.e.\@ $\xi\in\mathbb{R}$, then
$U=\frac{1}{\sqrt{2\pi}}\mathcal F^{-1}
\bigl( \frac{\mathcal F f_1}{\mathcal F f_0}\bigr)$.
\end{fact}

\begin{proof}
1)\enspace
See Stein--Weiss \cite[Chapter I, Theorem 3.18]{xStWe}, for example.

2)\enspace
It follows from $Bf_0=f_1$ that
$\mathcal{F}(U\ast f_0) = \mathcal{F}f_1$,
and therefore we have
$$
\sqrt{2\pi}(\mathcal{F}U)(\mathcal{F}f_0) = \mathcal{F}f_1.
$$

Hence, the bounded measurable function $\mathcal{F}U$ is determined by
the formula
$$
\mathcal{F}U(\xi) =
\frac{1}{\sqrt{2\pi}}\frac{\mathcal{F}f_1(\xi)}{\mathcal{F}f_0(\xi)}.
$$
\end{proof}

Next, we introduce two linear maps 
\index{A}{1sigma_+@$\sigma_+$}%
$\sigma_+$
and 
\index{A}{1sigma_-@$\sigma_-$}%
$\sigma_-$ by
\begin{equation}\label{def:sig}
 \begin{array}{l}
  \sigma_+:  L^2(\mathbb R_+, \frac{1}{2}r^{p+q-5}dr) \to L^2(\mathbb R),
                    \quad f(r) \mapsto \frac{1}{\sqrt{2}}e^{\frac{p+q-4}{2}x}f(e^x),
  \\
  \sigma_-:  L^2(\mathbb R_+, \frac{1}{2}r^{p+q-5}dr) \to L^2(\mathbb R),
                    \quad f(r) \mapsto \frac{1}{\sqrt{2}}e^{-\frac{p+q-4}{2}x}f(e^{-x}).
 \end{array}
\end{equation}
Then, both $\sigma_+$ and $\sigma_-$ are unitary operators.
Further, clearly we have
\begin{equation}\label{eqn:sigmapm}
(\sigma_- f)(x) = (\sigma_+ f)(-x).
\end{equation}
The inverse map $\sigma_-^{-1}$ is given by
$$
(\sigma_-^{-1}F)(r) = \sqrt{2}r^{-\frac{p+q-4}{2}} F(-\log r)
\quad\text{for $F\in L^2(\mathbb{R})$}.
$$

We define a subspace 
\index{A}{S0@$\mathcal{S}=\sigma_-^{-1}(\mathcal{S}(\mathbb{R}))$}%
$\mathcal S$ of $L^2(\mathbb R_+, \frac{1}{2}r^{p+q-5}dr)$ by
\begin{equation}\label{def:S}
  \mathcal S:= \sigma_-^{-1} \bigl( \mathcal S (\mathbb R) \bigr)
  = \sigma_+^{-1}(\mathcal{S}(\mathbb{R})),
\end{equation}
and endow $\mathcal{S}$ with the topology induced from 
that of the Schwartz space $\mathcal S(\mathbb R)$.
Now let 
\index{A}{S'@$\mathcal{S}'$}%
$\mathcal S'$ be the dual space of $\mathcal S$.  

Since $\sigma_-$ is unitary,
we get the following identity
for 
$F \in L^2(\mathbb{R})$
and 
$u\in L^2(\mathbb{R}_+, \frac{1}{2}r^{p+q-5} dr)$:
\begin{align*}
\langle\sigma_- u,F\rangle_{L^2(\mathbb{R})}
&=
\int_{\mathbb{R}}(\sigma_- u)(x) F(x) dx
\\
&=
\langle u,\sigma_-^{-1} F\rangle_{L^2(\mathbb{R}_+,\frac{1}{2}r^{p+q-5}dr)}.
\end{align*}

Then $\sigma_-$ extends naturally to an isomorphism 
from the dual space $\mathcal S'$ onto $\mathcal S'(\mathbb R)$ 
by the formula
\begin{equation} \label{eqn:Sdual}
  \langle \sigma_- u, F \rangle := \langle u, \sigma_-^{-1} F \rangle
\quad \text{for
  $F \in \mathcal S(\mathbb R),~ u \in \mathcal S'$}.
\end{equation}

Recall from \eqref{def:f} that we have defined a family of real analytic functions
 $f_{l,k}$ on $\mathbb{R}_+$ parametrized by $(l,k)\in\mathbb{N}^2$.
As we saw in Lemma \ref{lem:flkL2}, 
$f_{l,k}\in L^2(\mathbb R_+, \frac{1}{2}r^{p+q-5}dr)$.

For a continuous operator $A$ on $L^2(\mathbb R_+, \frac{1}{2}r^{p+q-5}dr)$, 
we define a continuous operator 
$\Tilde A$ on $L^2(\mathbb R)$ by
\begin{equation}\label{eqn:Atilde}
\Tilde A:= \sigma_- \circ A \circ \sigma_+^{-1}.
\end{equation}
Thus, the following diagram commutes:
\begin{alignat}{3}
  &L^2(\mathbb R_+, \frac{1}{2}r^{p+q-5}dr)  &
  \ \stackrel{A}{\longrightarrow}\ &  L^2(\mathbb R_+, \frac{1}{2}r^{p+q-5}dr)
  \notag
  \\
  \label{diag:A}
  &\qquad \quad \sigma_+ \downarrow && \qquad \quad \sigma_- \downarrow
  \\
  &\qquad \quad L^2(\mathbb R) 
  &\ \stackrel{\Tilde A}{\longrightarrow}\ & \qquad \quad L^2(\mathbb R).
  \notag
\end{alignat}
Since $\sigma_\pm$ are unitary operators,
$A$ is unitary if and only if $\Tilde{A}$ is unitary.

For $\kappa \in \mathcal S'$,
we define an operator $A_\kappa$ by
\begin{equation}\label{eqn:Akappa} 
\index{A}{Akappa@$A_\kappa$}%
  A_\kappa :  \mathcal S \to \mathcal S', \quad 
  f(r) \mapsto \frac{1}{2}\int_0^\infty \kappa(rr') f(r') r'^{p+q-5} dr'.
\end{equation}
It follows from the definition \eqref{def:sig} of $\sigma_+$ and
$\sigma_-$ that
\begin{equation*}
(\sigma_-\kappa \ast \sigma_+ f)(x) 
= \frac{1}{2}e^{-\frac{p+q-4}{2}x} \int_{-\infty}^\infty
  e^{(p+q-4)y} \kappa(e^{-x}e^y) f(e^y) dy.
\end{equation*}
Then, by the change of variables, we have
\begin{equation}
\label{eqn:conv}
  A_\kappa f = \frac{1}{\sqrt{2}}\sigma_-^{-1} (\sigma_- \kappa \ast \sigma_+ f)
 \quad \text{for $f \in \mathcal S$}.
\end{equation}

The following lemma characterizes operators of the form
$A_\kappa$:
\begin{lemSec}
\label{lem:F}
Let $\set{\rho(t)}{t\in \mathbb R}$ 
be a one parameter family of unitary operators on 
$L^2(\mathbb R_+, \frac{1}{2}r^{p+q-5}dr)$ defined by
\begin{equation}
\label{def:rho}
  (\rho(t)f)(r):= e^{-\frac{p+q-4}2 t}f(e^{-t}r).
\end{equation}
Suppose that a unitary operator $T$ on $L^2(\mathbb R_+, \frac{1}{2}r^{p+q-5}dr)$ 
satisfies the following (anti-)commutative relation:
\begin{equation}
\label{eqn:rho}
  T \circ \rho(t) = \rho(-t) \circ T  \quad \text{for any
 $t \in \mathbb R$}. 
\end{equation}
Then, there exists a unique  
distribution $\kappa \in \mathcal S'$
such that $T=A_\kappa$.
\end{lemSec}

\begin{proof}
For a general operator $A$ on $L^2(\mathbb R_+, \frac{1}{2}r^{p+q-5}dr)$,
we shall use the symbol $\Tilde A$ to denote 
$\sigma_- \circ A \circ \sigma_+^{-1}$ as in (\ref{eqn:Atilde}). 
Then we have
\begin{alignat*}{2}
  \widetilde{T\circ \rho(t)}&= 
  \sigma_- \circ (T \circ \rho(t) ) \circ \sigma_+^{-1}&&=
  \Tilde T \circ (\sigma_+ \circ \rho(t)\circ \sigma_+^{-1}),
 \\
 \widetilde{\rho(-t) \circ T}&= 
 \sigma_- \circ (\rho(-t)\circ T ) \circ \sigma_+^{-1}&&=
 (\sigma_- \circ \rho(-t) \circ \sigma_-^{-1})\circ \Tilde T.
\end{alignat*}
On the other hand, by a simple computation, 
we have the following identities:
$$
\sigma_+ \circ \rho(t)\circ \sigma_+^{-1}
= \sigma_- \circ \rho(-t) \circ \sigma_-^{-1}=l(t).
$$
Here, $l(t)$ denotes the translation operator 
\eqref{eqn:translt}. 
Hence, the relation (\ref{eqn:rho}) is equivalent to
$$
  \Tilde T \circ l(t) = l(t) \circ \Tilde T
  \quad\text{for any $t \in \mathbb{R}$},
$$
that is, $\Tilde T$ is a translation invariant bounded operator on $L^2(\mathbb R)$. 
Therefore, 
the operator $\Tilde T$ must be a convolution operator $U \ast$ for 
some tempered distribution $U \in \mathcal S'(\mathbb R)$  
such that its Fourier transform $\mathcal F U$ is
a bounded function (see Fact \ref{obsv:skt} (1)).

Finally, by setting $\kappa:=\sqrt{2}\sigma_-^{-1} U$, 
we have for any
$f\in L^2(\mathbb{R}_+, \frac{1}{2}r^{p+q-5}dr)$,
\begin{align*}
Tf
&= \sigma_-^{-1}\circ \Tilde{T}\circ \sigma_+f
\\
&= \frac{1}{\sqrt{2}}\sigma_-^{-1}(\sigma_-\kappa \ast \sigma_+ f)
\\
&= A_\kappa f
\end{align*}
by \eqref{eqn:conv}.
Therefore,
 $T=A_\kappa$. 
\end{proof}

\section{Reduction to Fourier analysis}
\label{subsec:F}

In this section,
we apply the results of the previous section to the unitary operator
$T_{l,k}$ on $L^2(\mathbb{R}_+,\frac{1}{2}r^{p+q-5}dr)$.
The main result of this section is 
Proposition \ref{prop:F}.

We begin with a refinement of Lemma \ref{lem:flkL2}.

\begin{lemSec}\label{lem:flkas}
\index{B}{asymptotic behavior!$\sigma_-(f_{l,k})$@---, $\sigma_-(f_{l,k})$}%
 $f_{l,k}$ belongs to $\mathcal S$.
In particular,
\index{A}{f0lkr@$f_{l,k}(r)$}%
$f_{l,k}\in L^2(\mathbb{R}_+, \frac{1}{2}r^{p+q-5}dr)$.
\end{lemSec}

\begin{proof}
By the definition of $\mathcal S$
(see \eqref{def:S}), it is sufficient to 
show $\sigma_- f_{l,k} \in \mathcal S
(\mathbb R)$. 
The proof is divided into two steps.

{\bf Step 1}: {\it For any $m\in\mathbb{N}$,
$x^m(\sigma_-f_{l,k})(x)$ is rapidly 
decreasing.}

By the definition \eqref{def:f} of $f_{l,k}$,
we have
\begin{equation} \label{eqn:sigflk}
(\sigma_- f_{l,k})(x)
= \frac{1}{\sqrt{2}}
     e^{-(l+k+\frac{p+q}{2}-2)x}
     \times
  \begin{cases}
     \widetilde{K}_{\frac{q-3}{2}+k}(2e^{-x})
   &\text{Case 1},
  \\
     \widetilde{K}_{\frac{p-3}{2}+l}(2e^{-x})
   &\text{Case 2}.
  \end{cases}
\end{equation}
Therefore, by the asymptotic behavior 
of $K$-Bessel functions 
(see Fact \ref{fact:Bas} 1), 2), respectively),
we get 
\begin{alignat}{2} \label{eqn:flkas1}
(\sigma_- f_{l,k})(x)\sim{}&
  \frac{e^{-x}}{2\sqrt{2}}
  \times
 \begin{cases}
  \Gamma(\frac{q-3}{2}+k)e^{-(\frac{p-q}2+l-k)x} 
  \\
  \Gamma(\frac{p-3}{2}+l)e^{(\frac{p-q}2+l-k)x}  
 \end{cases}
&&\begin{array}{l}
  \text{Case 1},\\
  \text{Case 2}, \vphantom{e^{(\frac{p-q}2+l-k)x}}
 \end{array}
 \qquad \text{as $x\to +\infty$},  
\\ \label{eqn:flkas2}
(\sigma_- f_{l,k})(x)\sim{}&
\frac{1}{\sqrt{2}}
 \times
 \begin{cases}
  e^{-(\frac{p-2}2+l)x}e^{-e^{-x}} 
 \\
  e^{-(\frac{q-2}2+k)x}e^{-e^{-x}} 
 \end{cases}
&&\begin{array}{l}
  \text{Case 1},\\
  \text{Case 2}, \vphantom{e^{-(\frac{q-2}2+k)x}}
 \end{array}
 \qquad \text{as $x\to -\infty$}.
\end{alignat}
Thus, Step 1 is proved.

{\bf Step 2}: {\it For any $n \in 
\mathbb N$, $\frac{d^n}{dx^n}
(\sigma_- f_{l,k}(x))$ is rapidly decreasing.}

Let us prove the above statement using induction on $n$.
We already know from Step 1 that 
the above statement is true for any $l,k\in\mathbb{N}$ in the case $n=0$.
Now assume that the statement is true for 
$n$. Then 
the statement for $n+1$ immediately 
follows from the next claim. 
Now Step 2 has been proved, and 
so has Lemma \ref{lem:flkas}.
\end{proof}

\begin{claim}\label{claim:4.4.2}
For any $l,k \in \mathbb N$, we have 
the following recurrence formula:
\begin{equation} \label{eqn:sigflk2}
\frac{d}{dx} (\sigma_- f_{l,k})=
-(\frac{p+q-4}2+l+k) \sigma_- f_{l,k}+
2 \sigma_- f_{l+1, k+1}.
\end{equation}
\end{claim}

\begin{proof}
For simplicity,
we rewrite the formula \eqref{eqn:sigflk} as
\begin{equation*}
   (\sigma_-f_{l,k})(x) = \frac{1}{\sqrt{2}}
    e^{-ax}\widetilde{K}_\nu(2e^{-x}),
\end{equation*}
where
\begin{equation*}
   a = \frac{p+q-4}{2}+l+k,
   \quad 
   \nu = \begin{cases}
            \frac{q-3}{2}+k  &\text{Case 1}, \\
            \frac{p-3}{2}+l  &\text{Case 2}. \\
         \end{cases}
\end{equation*}
By the recurrence formula \eqref{eqn:diffK2} of the $K$-Bessel function, 
we have
\begin{align*}
   \frac{d}{dx}(\sigma_-f_{l,k})(x)
     &= \frac{1}{\sqrt{2}}
        (-ae^{-ax}\widetilde{K}_\nu(2e^{-x})+2e^{-(a+2)x}
        \widetilde{K}_{\nu+1}(2e^{-x}))
\\
     &= -a(\sigma_-f_{l,k})(x)+2(\sigma_-f_{l+1,k+1})(x).
\end{align*}
Here, we have used the fact that $(k,l)$ belongs to Case 1
(i.e.\ $\frac{p-q}{2}+l-k\ge0$, see \eqref{def:Cases})
if and only if $(k+1,l+1)$ belongs to Case 1,
and likewise for Case 2.
Therefore, we have proved Claim \ref{claim:4.4.2}.
\end{proof}

\begin{remNonumber}
The above proof for $f_{l,k}\in \mathcal{S}$ was straightforward from
the asymptotic behavior of $f_{l,k}$ and its derivatives.
We shall prove in Proposition \ref{prop:Klk}
that $K_{l,k} \in \mathcal S'$, namely, 
$\sigma_- K_{l,k}\in \mathcal S'(\mathbb R)$.
Unlike the above proof for $f_{l,k} \in \mathcal{S}$,
this is not obvious from the asymptotic behavior of $K_{l,k}$
(see Remark \ref{rem:CG}).
In fact,
it follows from \eqref{eqn:Kinfty} that
\begin{align*}
\limsup_{x\to -\infty} e^{-\epsilon|x|}(\sigma_-K_{l,k})(x)
&=
\limsup_{x\to -\infty} e^{-\epsilon|x|-\frac{x}4} 
  \Bigl(\frac{4}{\sqrt{2\pi}}\cos (4 e^{-\frac{x}2}-\frac{2q-3}{4}\pi)\Bigr)
\\
&=+\infty
\end{align*}
if $\epsilon<\frac{1}{4}$.
Thus, the asymptotic behavior of
 $\sigma_- K_{l,k}$ is worse than that of any 
polynomial as $x$ tends to $-\infty$. 
As a result, our proof for $K_{l,k} \in \mathcal{S}'$ 
is more delicate,
and is based on an explicit computation of the Fourier
integral (see Proposition \ref{prop:Klk}).  
We note that $(\sigma_- K_{l,k})(x)$ decays exponentially as $x$ tends to
$+\infty$.
\end{remNonumber}

\begin{proposition}
\label{prop:F}
Let 
\index{A}{Tlk@$T_{l,k}$}%
$T_{l,k}$ be the unitary operator defined in Theorem \ref{thm:C}
(1). 
Then, there exists uniquely a distribution 
\index{A}{1kappa_{l,k}@$\kappa_{l,k}$|main}%
$\kappa_{l,k} \in \mathcal S'$ such that
$T_{l,k}=A_{\kappa_{l,k}}$ 
(see \eqref{eqn:Akappa} for notation).
Namely, we have:
\begin{equation}\label{eqn:F}
 (T_{l,k}f)(r)=
 \frac{1}{2}\int_0^\infty \kappa_{l,k}(rr')f(r')
 r'^{p+q-5}dr'.
\end{equation}
\end{proposition}

\begin{proof}
[Proof of Proposition \ref{prop:F}]
We recall from \eqref{eqn:rA} that the unitary operator 
$\pi(e^{tH})$  on $L^2(C)$ can be written by means of 
the unitary operator $\rho(t)$ (see \eqref{def:rho}) on $L^2(\mathbb R_+, \frac{1}{2}r^{p+q-5}dr)$
as follows:
$$
  \pi(e^{tH})\bigl(f(r)\phi(\omega) \psi(\eta) \bigr)=
  \bigl(\rho(t)f\bigr)(r) \phi(\omega) \psi(\eta),
$$ 
where $f \in L^2(\mathbb R_+, \frac{1}{2}r^{p+q-5}dr)$, \thinspace
$\phi \in \Har{l}{p-1}$, 
and $\psi \in \Har{k}{q-1}$.

Applying $\pi(w_0)$ to the both sides,
together with the definition of $T_{l,k}$ 
(see Theorem \ref{thm:C} (1)), we obtain 
\begin{equation*}
\pi(w_0)\circ \pi(e^{tH})(f(r)\phi(w)\psi(\eta))
= (T_{l,k}\rho(t)f)(r)\phi(w)\psi(\eta).
\end{equation*}
Similarly, applying $\pi(w_0)$ followed by $\pi(e^{-tH})$, we get
\begin{equation*}
\pi(e^{-tH})\circ \pi(w_0)(f(r)\phi(w)\psi(\eta))
= (\rho(-t)T_{l,k}f)(r)\phi(w)\psi(\eta).
\end{equation*}
On the other hand,
it follows from $\operatorname{Ad}(w_0)H=-H$ 
(see \eqref{eqn:wa}) that
$$
  w_0 e^{tH} = e^{-tH} w_0,
$$
and then we have
\begin{equation*}
  \pi(w_0) \circ \pi(e^{tH})
  = \pi(e^{-tH}) \circ \pi(w_0).
\end{equation*}
Therefore, 
$$
  T_{l,k} \circ \rho(t) = \rho(-t) \circ T_{l,k}.
$$
Now, Proposition \ref{prop:F} follows from Lemma \ref{lem:F}.
\end{proof}

\section{Kernel function $K_{l,k}$}
\label{subsec:Klk}
We defined a family of real analytic functions 
\index{A}{Klkt@$K_{l,k}(t)$}%
$K_{l,k}(t)$ in
 \eqref{def:Klk} by means of Meijer's $G$-function $G_{04}^{20}$. 
This section studies basic properties of $K_{l,k}(t)$ 
as a preparation for the proof of Theorem \ref{thm:C} (2). 
The main result here is 
Proposition \ref{prop:Klk}.

We begin with an asymptotic estimate of $K_{l,k}(t)$.

\begin{lemSec}[%
\index{B}{asymptotic behavior!$k_{l,k}$@---, $K_{l,k}$}%
Asymptotic behavior]\label{rem:CG}
$K_{l,k}(t)$
has the following asymptotics 
as $t$ tends to $0$ or $\infty$:

{\rm 1)} As $t$ tends to $0$,
\begin{alignat*}{2}
 K_{l,k}(t)& = 
 \begin{cases}
  O(t^{-q+3+l-k}) \quad \text{Case 1} \\
  O(t^{-p+3-l+k}) \quad \text{Case 2} 
 \end{cases} & & \quad (q>2), \\
 K_{l,k}(t)& = O(t^{l}) & & \quad (q=2).
\end{alignat*}

{\rm 2)} There are some constants $P_1, \cdots, Q_1, \cdots$ such
that
\begin{align}
 K_{l,k}(t)
 ={}&\frac{4}{\sqrt{2\pi}} \, t^{-\frac{2p+2q-9}4} 
  \cos \bigl(4t^\frac{1}2-\frac{2q-3}{4}\pi\bigr)
  (1+P_1 t^{-1}+P_2 t^{-2}+\cdots)  
\nonumber
\\
 &+t^{-\frac{2p+2q-9}4}\sin\bigl(4t^\frac{1}2-\frac{2q-3}{4}\pi\bigr)
  (Q_1 t^{-\frac{1}2}+Q_2 t^{-\frac{3}2}+\cdots),
\label{eqn:Kinfty}
\end{align}
as $t$ tends to $+\infty$.
\end{lemSec}

\begin{proof}
Directly obtained from the asymptotic formula of Meijer's $G$-function
$G_{04}^{20}(x \mid b_1,b_2,1-\gamma-b_1,1-\gamma-b_2)$ given in Lemma
\ref{lem:Gas} in Appendix.
\end{proof}

Next, we give an integral expression of $K_{l,k}(t)$ $(t>0)$,
where
the integral path $L$ will be taken
independently of $l,k\in\mathbb{N}$.
We note that the integrands in \eqref{int:Klk} and \eqref{eqn:Klk2} 
are the same.
The first expression \eqref{int:Klk} is convenient in Case 1,
and the second expression \eqref{eqn:Klk2} is convenient in Case 2
 (see Remark \ref{rem:Klkint}).

\begin{lemSec}\label{lem:Klkint}
Fix a real number $\gamma>-\frac{p+q-5}{2}$,
and let $L$ be a contour that starts at $\gamma-\sqrt{-1}\infty$
and ends at $\gamma+\sqrt{-1}\infty$ and
passes the real axis at a point $s_0$ satisfying
$s_0 < -\frac{p+q-6}{2}$
(see Figure \ref{fig:Ls}).
(Later, we shall assume also that $-\frac{p+q-4}{2} < s_0$ in
Chapter \ref{sec:A}.)
Then, we have
\begin{align}\label{int:Klk}
  K_{l,k}(t)&= 
     \frac{(-1)^{l+\frac{p-q}2}}{\pi \sqrt{-1}}
     \int_L \frac{\Gamma(\frac{l+k-s}2 ) 
                          \Gamma(\frac{-q+3+l-k-s}2 )}
                         {\Gamma(\frac{p+q-4+l+k+s}2) 
                          \Gamma(\frac{p-1+l-k+s}2)}
     t^s ds
\\
 &=
     \frac{(-1)^k}{\pi \sqrt{-1}}
     \int_L \frac{\Gamma(\frac{l+k-s}2 ) 
                          \Gamma(\frac{-p+3-l+k-s}2 )}
                         {\Gamma(\frac{p+q-4+l+k+s}2) 
                          \Gamma(\frac{q-1-l+k+s}2)}
     t^s ds.
\label{eqn:Klk2}
\end{align}
\end{lemSec}

\begin{figure}[H]
\setlength{\unitlength}{0.00033333in}
\begingroup\makeatletter\ifx\SetFigFont\undefined%
\gdef\SetFigFont#1#2#3#4#5{%
  \reset@font\fontsize{#1}{#2pt}%
  \fontfamily{#3}\fontseries{#4}\fontshape{#5}%
  \selectfont}%
\fi\endgroup%
{\renewcommand{\dashlinestretch}{30}
\begin{picture}(10812,9034)(-1700,-10)
\path(9600,8413)(9600,2113)
\path(4800,4363)(4800,4063)
\path(1200,4363)(1200,4063)
\dashline{60.000}(2400,8413)(2400,13)
\path(900,4213)(10800,4213)
\path(2843,6769)(2790,6573)
\path(2859,6764)(3015,6623)
\path(2546,12)(2546,13)(2546,15)
	(2547,19)(2548,26)(2549,35)
	(2550,47)(2552,63)(2555,82)
	(2557,105)(2561,130)(2564,160)
	(2569,192)(2573,227)(2578,265)
	(2584,305)(2589,347)(2595,390)
	(2602,436)(2608,482)(2615,530)
	(2623,578)(2630,628)(2638,679)
	(2646,730)(2655,783)(2664,836)
	(2674,891)(2684,947)(2695,1004)
	(2706,1063)(2718,1124)(2731,1186)
	(2745,1250)(2759,1315)(2774,1381)
	(2790,1447)(2806,1513)(2825,1589)
	(2845,1663)(2864,1732)(2883,1796)
	(2901,1856)(2918,1912)(2935,1964)
	(2952,2012)(2967,2056)(2983,2098)
	(2998,2137)(3012,2174)(3026,2208)
	(3040,2241)(3054,2271)(3067,2300)
	(3079,2327)(3091,2353)(3103,2376)
	(3114,2398)(3123,2417)(3132,2435)
	(3140,2450)(3147,2463)(3153,2474)
	(3157,2482)(3161,2489)(3163,2493)
	(3165,2496)(3166,2497)(3166,2498)
\path(3166,2498)(3166,2499)(3167,2501)
	(3169,2505)(3172,2511)(3176,2519)
	(3181,2531)(3188,2544)(3196,2561)
	(3204,2580)(3214,2602)(3225,2626)
	(3237,2652)(3249,2680)(3261,2709)
	(3274,2739)(3288,2771)(3301,2804)
	(3315,2839)(3329,2875)(3343,2912)
	(3357,2951)(3371,2993)(3386,3036)
	(3401,3082)(3416,3130)(3431,3181)
	(3446,3235)(3461,3290)(3476,3347)
	(3491,3408)(3504,3468)(3517,3525)
	(3528,3579)(3537,3631)(3546,3679)
	(3553,3725)(3559,3769)(3565,3811)
	(3570,3850)(3574,3889)(3578,3925)
	(3581,3960)(3584,3994)(3586,4026)
	(3588,4057)(3590,4085)(3592,4111)
	(3593,4134)(3594,4154)(3594,4171)
	(3595,4185)(3595,4196)(3596,4203)
	(3596,4208)(3596,4211)(3596,4212)
\path(2550,8413)(2550,8412)(2550,8410)
	(2551,8406)(2552,8399)(2553,8390)
	(2554,8378)(2556,8362)(2559,8343)
	(2561,8320)(2565,8295)(2568,8265)
	(2573,8233)(2577,8198)(2582,8160)
	(2588,8120)(2593,8078)(2599,8035)
	(2606,7989)(2612,7943)(2619,7895)
	(2627,7847)(2634,7797)(2642,7746)
	(2650,7695)(2659,7642)(2668,7589)
	(2678,7534)(2688,7478)(2699,7421)
	(2710,7362)(2722,7301)(2735,7239)
	(2749,7175)(2763,7110)(2778,7044)
	(2794,6978)(2810,6912)(2829,6836)
	(2849,6762)(2868,6693)(2887,6629)
	(2905,6569)(2922,6513)(2939,6461)
	(2956,6413)(2971,6369)(2987,6327)
	(3002,6288)(3016,6251)(3030,6217)
	(3044,6184)(3058,6154)(3071,6125)
	(3083,6098)(3095,6072)(3107,6049)
	(3118,6027)(3127,6008)(3136,5990)
	(3144,5975)(3151,5962)(3157,5951)
	(3161,5943)(3165,5936)(3167,5932)
	(3169,5929)(3170,5928)(3170,5927)
\path(3170,5927)(3170,5926)(3171,5924)
	(3173,5920)(3176,5914)(3180,5906)
	(3185,5894)(3192,5881)(3200,5864)
	(3208,5845)(3218,5823)(3229,5799)
	(3241,5773)(3253,5745)(3265,5716)
	(3278,5686)(3292,5654)(3305,5621)
	(3319,5586)(3333,5550)(3347,5513)
	(3361,5474)(3375,5432)(3390,5389)
	(3405,5343)(3420,5295)(3435,5244)
	(3450,5190)(3465,5135)(3480,5078)
	(3495,5017)(3508,4957)(3521,4900)
	(3532,4846)(3541,4794)(3550,4746)
	(3557,4700)(3563,4656)(3569,4614)
	(3574,4575)(3578,4536)(3582,4500)
	(3585,4465)(3588,4431)(3590,4399)
	(3592,4368)(3594,4340)(3596,4314)
	(3597,4291)(3598,4271)(3598,4254)
	(3599,4240)(3599,4229)(3600,4222)
	(3600,4217)(3600,4214)(3600,4213)
\put(2700,8263){\makebox(0,0)[lb]{\smash{{{\SetFigFont{10}{14.4}{\familydefault}{\mddefault}{\updefault}$\gamma+\sqrt{-1}\infty$}}}}}
\put(0,3463){\makebox(0,0)[lb]{\smash{{{\SetFigFont{10}{14.4}{\familydefault}{\mddefault}{\updefault}$-\frac{p+q-5}{2}$}}}}}
\put(4050,3463){\makebox(0,0)[lb]{\smash{{{\SetFigFont{10}{14.4}{\familydefault}{\mddefault}{\updefault}$-\frac{p+q-6}{2}$}}}}}
\put(9225,3763){\makebox(0,0)[lb]{\smash{{{\SetFigFont{10}{14.4}{\familydefault}{\mddefault}{\updefault}$0$}}}}}
\put(2700,88){\makebox(0,0)[lb]{\smash{{{\SetFigFont{10}{14.4}{\familydefault}{\mddefault}{\updefault}$\gamma-\sqrt{-1}\infty$}}}}}
\put(4725,8863){\makebox(0,0)[lb]{\smash{{{\SetFigFont{10}{14.4}{\familydefault}{\mddefault}{\updefault}$s$-plane}}}}}
\put(3300,6163){\makebox(0,0)[lb]{\smash{{{\SetFigFont{10}{14.4}{\familydefault}{\mddefault}{\updefault}$L$}}}}}
\put(1965,3838){\makebox(0,0)[lb]{\smash{{{\SetFigFont{10}{14.4}{\familydefault}{\mddefault}{\updefault}$\gamma$}}}}}
\put(3170,3853){\makebox(0,0)[lb]{\smash{{{\SetFigFont{10}{14.4}{\familydefault}{\mddefault}{\updefault}$s_0$}}}}}
\end{picture}
}
\caption{}
\label{fig:Ls}
\end{figure}

\begin{proof}
The equality \eqref{int:Klk} $=$ \eqref{eqn:Klk2} is an immediate
consequence of the following formula:
\begin{equation}\label{eqn:Gammalk}
\frac{\Gamma(\frac{-q+3+l-k-s}{2})}{\Gamma(\frac{p-1+l-k+s}{2})}
\cdot
\frac{\Gamma(\frac{q-1-l+k+s}{2})}{\Gamma(\frac{-p+3-l+k-s}{2})}
=
(-1)^{\frac{p-q}{2}+l-k},
\end{equation}
which is derived from
$$
\Gamma(z) \Gamma(1-z) = \frac{\pi}{\sin\pi z}.
$$

Let us show \eqref{int:Klk} in Case 1,
and \eqref{eqn:Klk2} in Case 2
(see \eqref{def:Cases} for the definition of Cases 1 and 2).
As meromorphic functions of the variable $s$,
the poles of the numerators in the integrands \eqref{int:Klk} and
\eqref{eqn:Klk2} are given by
\begin{align*}
W_1
 &:= \set{l+k+2a, -q+3+l-k+2a}{a \in \mathbb N},
\\
W_2
 &:= \set{l+k+2a, -p+3-l+k+2a}{a \in \mathbb N},
\end{align*}
respectively.
Then, 
\begin{align*}
&\inf W_1 \ge -\frac{p+q-6}{2}
 \quad\text{in Case 1 \ (i.e. $\frac{p-q}{2}+l-k\ge0$)},
\\
&\inf W_2 \ge -\frac{p+q-6}{2}
 \quad\text{in Case 2 \ (i.e. $\frac{p-q}{2}+l-k\le0$)}.
\end{align*}
Therefore, in either case,
the contour $L$ leaves all these sets $W_1$ and $W_2$ on the right
 because our $L$ passes
the real axis at some point $s_0 < -\frac{p+q-6}2$.

By the definition of Meijer's $G$-function (see \eqref{def:G} in
Appendix, see also Example \ref{ex:G24}), 
we get \eqref{int:Klk} in Case 1 and \eqref{eqn:Klk2} in Case 2 by the
change of variables $s:=2\lambda$.
Hence, Lemma is proved.
\end{proof}

\begin{remNonumber}\label{rem:Klkint}
We shall use the expression
\eqref{int:Klk} in Case 1 and \eqref{eqn:Klk2} in Case 2
later. 
The point here is that there is no
 cancellation
of the poles of the numerator and the denominator of the integrand.
For example, 
the poles of the denominator 
 of the integrand \eqref{int:Klk} 
are given by
$$
V_1 := \{ -p-q+4-l-k-2b, \, -p+1-l+k-2b: b\in\mathbb{N} \}.
$$
Then,
\begin{equation*}
\sup V_1 < \inf W_1 \quad\text{in Case 1},
\end{equation*}
and therefore $V_1 \cap W_1 = \emptyset$.
Similarly, there is no cancellation of the poles of the numerator and
the poles of the denominator of the integrand \eqref{eqn:Klk2} 
in Case 2. 
\end{remNonumber}
           
As $K_{l,k}$ is a real analytic function on $\mathbb{R}_+$,
so is
$(\sigma_- K_{l,k})(x)$ on $\mathbb R$
 (see (\ref{def:sig}) for the definition of $\sigma_-$), 
which in turn is a distribution on $\mathbb R$. 
More strongly, we shall see in Proposition \ref{prop:Klk} that 
$(\sigma_- K_{l,k})(x)$ is a tempered distribution. 

For this, we define a meromorphic function $\psi(\zeta)$ on $\mathbb C$ by
\begin{align}
\psi(\zeta):= {}
&(-1)^{l+\frac{p-q}{2}}
\frac{\Gamma\bigl(\frac{p+q-4}4+\frac{l+k-\sqrt{-1} \zeta}2\bigr)
   \Gamma\bigl(\frac{p-q}4+\frac{l-k+1-\sqrt{-1} \zeta}2\bigr)}
 {\Gamma\bigl(\frac{p+q-4}4+\frac{l+k+\sqrt{-1} \zeta}2\bigr)
 \Gamma\bigl(\frac{p-q}4+\frac{l-k+1+\sqrt{-1} \zeta}2\bigr)}
\label{def:psi}
\\
={}
&(-1)^k
\frac{\Gamma\bigl(\frac{p+q-4}4+\frac{l+k-\sqrt{-1} \zeta}2\bigr)
   \Gamma\bigl(-\frac{p-q}4+\frac{-l+k+1-\sqrt{-1} \zeta}2\bigr)}
 {\Gamma\bigl(\frac{p+q-4}4+\frac{l+k+\sqrt{-1} \zeta}2\bigr)
 \Gamma\bigl(-\frac{p-q}4+\frac{-l+k+1+\sqrt{-1} \zeta}2\bigr)}.
\label{def:psi2}
\end{align}
We shall use \eqref{def:psi} in Case 1,
and \eqref{def:psi2} in Case 2.
The proof of the equality \eqref{def:psi} $=$ \eqref{def:psi2} is the same as
the proof of the equality \eqref{int:Klk} $=$ \eqref{eqn:Klk2}.

\begin{lemSec}\label{lem:espsi}
{\rm 1)} $|\psi(\zeta)|=1$  for $\zeta \in \mathbb R$. 
In particular, the inverse Fourier transform 
$\mathcal F^{-1} \psi$ is defined to be a tempered distribution.

{\rm 2)} $\psi(\zeta)$ is a meromorphic function on $\mathbb{C}$,
and the set of its poles is contained in
$$
\{ -\sqrt{-1} m : m=1,2,3,\ldots \}.
$$

{\rm 3)} For $\eta_1\le\eta\le\eta_2$,
\begin{equation}
\label{eqn:espsi}
|\psi(\xi-\sqrt{-1}\eta)| \sim \left|\frac{\xi}{2}\right|^{-2\eta}
\quad\text{as $|\xi|\to\infty$}.
\end{equation}
\end{lemSec}
\begin{proof}
1) Since $\Gamma(\overline z )=\overline{\Gamma(z)}$ for 
$z \in \mathbb C$, 
we have $|\psi(\zeta)| = 1$.
Therefore $\psi \in \mathcal{S}'(\mathbb{R})$ and thus
$\mathcal{F}^{-1}\psi \in \mathcal{S}'(\mathbb{R})$.

2) The proof is straightforward from the definitions \eqref{def:psi} and
   \eqref{def:psi2} in each case.

3) We recall Stirling's 
\index{B}{asymptotic behavior!gamma function@---, gamma function}%
asymptotic expansion 
of the gamma function
(see \cite[Corollary 1.4.4]{xaar} for example):
\begin{equation}\label{eqn:Stirling}
\left|\Gamma(a+\sqrt{-1}b)\right| = \sqrt{2\pi} |b|^{a-\frac{1}{2}}
e^{-\frac{\pi|b|}{2}} (1+O(|b|^{-1}))
\end{equation}
when $a_1\le a \le a_2$ and $|b| \to \infty$.
Then, for $\alpha \in \mathbb{R}$ and 
$z = x+\sqrt{-1}y$ $(y_1 \le y \le y_2)$, 
\begin{equation}\label{eqn:rationgamma}
\left|\frac{\Gamma(\alpha-\sqrt{-1} z)}{\Gamma(\alpha+\sqrt{-1} z)}\right|
= |x|^{2y} (1+O(|x|^{-1}))
\quad\text{as $|x|\to\infty$},
\end{equation}
where the constant implied by the Bachmann--Landau symbol
$O$ depends only on $\alpha$, $y_1$ and
$y_2$. 
Now, applying \eqref{eqn:rationgamma} to
$z = \frac{1}{2}(\xi-\sqrt{-1}\eta)$ twice,
we get \eqref{eqn:espsi}.
\end{proof}
By the change of variable $s=\sqrt{-1}\zeta-\frac{p+q-4}2$,
the integral formula of $K_{l,k}$
(Lemma \ref{lem:Klkint}) can be restated as follows:
\begin{lemSec}\label{lem:Kpsi}
Let $\gamma>-\frac{p+q-5}{2}$ and $L'$ be an integral path starting 
from $-(\gamma + \frac{p+q-4}{2})\sqrt{-1}-\infty$ and ending 
at 
$-(\gamma + \frac{p+q-4}{2})\sqrt{-1}+\infty$
passing the imaginary axis at some point 
in the open interval $(-\frac{\sqrt{-1}}{2}, 
-\sqrt{-1})$ (see Figure \ref{fig:zeta}).
Then, we have
\begin{equation}
K_{l,k}(t)=\frac{1}{\pi}
\int_{L'} \psi(\zeta)t^{-\frac{p+q-4}2+\sqrt{-1}\zeta}d\zeta,
\end{equation}
or equivalently (see \eqref{def:sig} for the definition 
of $\sigma_-$),
\begin{equation} \label{eqn:4.5.11}
  (\sigma_- K_{l,k})(x)=
  \frac{1}{\sqrt{2}\pi} \int_{L'}  \psi(\zeta) e^{-\sqrt{-1}x\zeta} d\zeta.
\end{equation}
\end{lemSec}

\begin{figure}[H]
\setlength{\unitlength}{0.00033333in}
\begingroup\makeatletter\ifx\SetFigFont\undefined%
\gdef\SetFigFont#1#2#3#4#5{%
  \reset@font\fontsize{#1}{#2pt}%
  \fontfamily{#3}\fontseries{#4}\fontshape{#5}%
  \selectfont}%
\fi\endgroup%
{\renewcommand{\dashlinestretch}{30}
\begin{picture}(12282,10008)(600,-10)
\path(9346,3569)(9162,3608)
\path(9344,3562)(9239,3409)
\path(8220,2112)(7920,2112)
\dashline{60.000}(12270,4512)(3870,4512)
\path(8070,9162)(8070,12)
\path(3870,8412)(12270,8412)
\path(8220,5262)(7920,5262)
\path(3869,4366)(3870,4366)(3872,4366)
	(3876,4365)(3883,4364)(3892,4363)
	(3904,4362)(3920,4360)(3939,4357)
	(3962,4355)(3987,4351)(4017,4348)
	(4049,4343)(4084,4339)(4122,4334)
	(4162,4328)(4204,4323)(4247,4317)
	(4293,4310)(4339,4304)(4387,4297)
	(4435,4289)(4485,4282)(4536,4274)
	(4587,4266)(4640,4257)(4693,4248)
	(4748,4238)(4804,4228)(4861,4217)
	(4920,4206)(4981,4194)(5043,4181)
	(5107,4167)(5172,4153)(5238,4138)
	(5304,4122)(5370,4106)(5446,4087)
	(5520,4067)(5589,4048)(5653,4029)
	(5713,4011)(5769,3994)(5821,3977)
	(5869,3960)(5913,3945)(5955,3929)
	(5994,3914)(6031,3900)(6065,3886)
	(6098,3872)(6128,3858)(6157,3845)
	(6184,3833)(6210,3821)(6233,3809)
	(6255,3798)(6274,3789)(6292,3780)
	(6307,3772)(6320,3765)(6331,3759)
	(6339,3755)(6346,3751)(6350,3749)
	(6353,3747)(6354,3746)(6355,3746)
\path(6355,3746)(6356,3746)(6358,3745)
	(6362,3743)(6368,3740)(6376,3736)
	(6388,3731)(6401,3724)(6418,3716)
	(6437,3708)(6459,3698)(6483,3687)
	(6509,3675)(6537,3663)(6566,3651)
	(6596,3638)(6628,3624)(6661,3611)
	(6696,3597)(6732,3583)(6769,3569)
	(6808,3555)(6850,3541)(6893,3526)
	(6939,3511)(6987,3496)(7038,3481)
	(7092,3466)(7147,3451)(7204,3436)
	(7265,3421)(7325,3408)(7382,3395)
	(7436,3384)(7488,3375)(7536,3366)
	(7582,3359)(7626,3353)(7668,3347)
	(7707,3342)(7746,3338)(7782,3334)
	(7817,3331)(7851,3328)(7883,3326)
	(7914,3324)(7942,3322)(7968,3320)
	(7991,3319)(8011,3318)(8028,3318)
	(8042,3317)(8053,3317)(8060,3316)
	(8065,3316)(8068,3316)(8069,3316)
\path(12270,4362)(12269,4362)(12267,4362)
	(12263,4361)(12256,4360)(12247,4359)
	(12235,4358)(12219,4356)(12200,4353)
	(12177,4351)(12152,4347)(12122,4344)
	(12090,4339)(12055,4335)(12017,4330)
	(11977,4324)(11935,4319)(11892,4313)
	(11846,4306)(11800,4300)(11752,4293)
	(11704,4285)(11654,4278)(11603,4270)
	(11552,4262)(11499,4253)(11446,4244)
	(11391,4234)(11335,4224)(11278,4213)
	(11219,4202)(11158,4190)(11096,4177)
	(11032,4163)(10967,4149)(10901,4134)
	(10835,4118)(10769,4102)(10693,4083)
	(10619,4063)(10550,4044)(10486,4025)
	(10426,4007)(10370,3990)(10318,3973)
	(10270,3956)(10226,3941)(10184,3925)
	(10145,3910)(10108,3896)(10074,3882)
	(10041,3868)(10011,3854)(9982,3841)
	(9955,3829)(9929,3817)(9906,3805)
	(9884,3794)(9865,3785)(9847,3776)
	(9832,3768)(9819,3761)(9808,3755)
	(9800,3751)(9793,3747)(9789,3745)
	(9786,3743)(9785,3742)(9784,3742)
\path(9784,3742)(9783,3742)(9781,3741)
	(9777,3739)(9771,3736)(9763,3732)
	(9751,3727)(9738,3720)(9721,3712)
	(9702,3704)(9680,3694)(9656,3683)
	(9630,3671)(9602,3659)(9573,3647)
	(9543,3634)(9511,3620)(9478,3607)
	(9443,3593)(9407,3579)(9370,3565)
	(9331,3551)(9289,3537)(9246,3522)
	(9200,3507)(9152,3492)(9101,3477)
	(9047,3462)(8992,3447)(8935,3432)
	(8874,3417)(8814,3404)(8757,3391)
	(8703,3380)(8651,3371)(8603,3362)
	(8557,3355)(8513,3349)(8471,3343)
	(8432,3338)(8393,3334)(8357,3330)
	(8322,3327)(8288,3324)(8256,3322)
	(8225,3320)(8197,3318)(8171,3316)
	(8148,3315)(8128,3314)(8111,3314)
	(8097,3313)(8086,3313)(8079,3312)
	(8074,3312)(8071,3312)(8070,3312)
\put(7470,9837){\makebox(0,0)[lb]{\smash{{{\SetFigFont{10}{14.4}{\familydefault}{\mddefault}{\updefault}$\zeta$-plane}}}}}
\put(7620,7887){\makebox(0,0)[lb]{\smash{{{\SetFigFont{10}{14.4}{\familydefault}{\mddefault}{\updefault}$0$}}}}}
\put(0,3717){\makebox(0,0)[lb]{\smash{{{\SetFigFont{10}{14.4}{\familydefault}{\mddefault}{\updefault}$-(\gamma+\frac{p+q-4}{2})\sqrt{-1}-\infty$}}}}}
\put(11370,3792){\makebox(0,0)[lb]{\smash{{{\SetFigFont{10}{14.4}{\familydefault}{\mddefault}{\updefault}$-(\gamma+\frac{p+q-4}{2})\sqrt{-1}+\infty$}}}}}
\put(8955,8637){\makebox(0,0)[lb]{\smash{{{\SetFigFont{10}{14.4}{\familydefault}{\mddefault}{\updefault}$\zeta=-\sqrt{-1}(s+\frac{p+q-4}{2})$}}}}}
\put(8970,9267){\makebox(0,0)[lb]{\smash{{{\SetFigFont{10}{14.4}{\familydefault}{\mddefault}{\updefault}$s=\sqrt{-1}\zeta-\frac{p+q-4}{2}$}}}}}
\put(6309,2006){\makebox(0,0)[lb]{\smash{{{\SetFigFont{10}{14.4}{\familydefault}{\mddefault}{\updefault}$-\sqrt{-1}$}}}}}
\put(6345,5187){\makebox(0,0)[lb]{\smash{{{\SetFigFont{10}{14.4}{\familydefault}{\mddefault}{\updefault}$-\frac{\sqrt{-1}}{2}$}}}}}
\put(9435,3132){\makebox(0,0)[lb]{\smash{{{\SetFigFont{10}{14.4}{\familydefault}{\mddefault}{\updefault}$L'$}}}}}
\end{picture}
}
\caption{}
\label{fig:zeta}
\end{figure}

Now, 
we recall from Section \ref{subsec:multop} 
that $\mathcal S'$ is the dual space of
$\mathcal{S}=\sigma_-^{-1}(\mathcal{S}(\mathbb{R}))$
via the following diagram:
\begin{figure}[H]
\begin{equation*}
   \begin{matrix}
\index{A}{S'@$\mathcal{S}'$}%
       \mathcal{S}'   & \rarrowsim   &\mathcal{S}'(\mathbb{R})
   \\
        \cup  && \cup
   \\
       \sigma_-: L^2(\mathbb{R}_+,\frac{1}{2}r^{p+q-5}dr)
           & \rarrowsim   & L^2(\mathbb{R})
   \\
        \cup  && \cup
   \\
\index{A}{S0@$\mathcal{S}=\sigma_-^{-1}(\mathcal{S}(\mathbb{R}))$}%
       \mathcal{S}   & \rarrowsim   &\mathcal{S}(\mathbb{R}).
   \end{matrix}
\end{equation*}
\renewcommand{\figurename}{Diagram}
\caption{}
\label{diag:sigma}
\end{figure}
We are ready to prove the main result of this section:
\begin{proposition}
\label{prop:Klk}
{\rm 1)} $K_{l,k}$ belongs to $\mathcal S'$.
That is, $\sigma_- K_{l,k} \in \mathcal S'(\mathbb R)$.
 
{\rm 2)} The Fourier transform 
$\mathcal F (\sigma_- K_{l,k})(\zeta)$ of $\sigma_- K_{l,k}$ 
is equal to $\frac{1}{\sqrt{\pi}}\psi(\zeta)$ (see \eqref{def:psi} for definition). 
In particular, 
$|\mathcal F (\sigma_- K_{l,k})(\zeta)|=\frac{1}{\sqrt{\pi}}$
 for 
$\zeta \in \mathbb R$.
\end{proposition}

\begin{proof}
It follows from Lemma \ref{lem:espsi} (1) that
$\psi$ is a tempered distribution,
and therefore, its inverse Fourier transform
 $\mathcal{F}^{-1}\psi \in \mathcal{S}'(\mathbb{R})$.
We also know that $\sigma_- K_{l,k} \in C^\infty(\mathbb{R})$ by
definition.
Let $\mathcal{D}'(\mathbb{R})$ be the space of distributions,
namely, the dual space of $C_0^\infty(\mathbb{R})$.
In light of the inclusion
\begin{equation*}
\mathcal{S}'(\mathbb{R}) \subset \mathcal{D}'(\mathbb{R})
\supset C^\infty(\mathbb{R}),
\end{equation*}
all the statements of Proposition \ref{prop:Klk} will be proved if we
show 
\begin{equation}\label{eqn:Klk1}
\sqrt{\pi}\sigma_- K_{l,k} = \mathcal{F}^{-1} \psi
\quad\text{in $\mathcal{D}'(\mathbb{R})$}, 
\end{equation}
that is,
\begin{equation}\label{eqn:Klk}
\sqrt{\pi}\int_{-\infty}^{\infty}(\sigma_- K_{l,k})(x)
 \overline{\phi(x)}dx=
\int_{-\infty}^{\infty}(\mathcal F^{-1}\psi)(x)
 \overline{\phi (x)}dx
\end{equation}
holds for any test function $\phi \in C_0^\infty(\mathbb R)$.
In fact, \eqref{eqn:Klk1}
 will imply that $\sigma_-K_{l,k} \in \mathcal{S}'(\mathbb{R})$ 
and
$\sqrt{\pi}\mathcal{F}(\sigma_-K_{l,k}) = \psi$
as a tempered distribution.

The key to the proof of \eqref{eqn:Klk} is the integral expression of
$K_{l,k}(t)$ stated in Lemma \ref{lem:Kpsi}. 
By \eqref{eqn:4.5.11}, the left-hand side of \eqref{eqn:Klk}
amounts to
\begin{align*}
& \frac{1}{\sqrt{2\pi}}
\int_{-\infty}^\infty 
\biggl(\int_{L'}\psi(\zeta)e^{-\sqrt{-1}x\zeta}d\zeta\biggr)
\overline{\phi(x)}dx \\
={}& \frac{1}{\sqrt{2\pi}}
\int_{L'}\psi(\zeta)\overline{\biggl(\int_{-\infty}^\infty
\phi(x)e^{\sqrt{-1}x\bar\zeta}dx\biggr)} d\zeta  \\
={}&
\int_{L'}\psi(\zeta)\overline{(\mathcal F \phi)(\bar\zeta)}d\zeta \\
={}&
\int_{-\infty}^\infty 
\psi(\zeta)\overline{(\mathcal F \phi)(\zeta)}d\zeta \\
={}& \text{right-hand side of \eqref{eqn:Klk}}.
\end{align*}
In what follows, 
we explain the above equalities in details:

{\bf First equality} is by 
Fubini's theorem.
For this justification, 
take $a>0$ such that $\operatorname{Supp}\phi\subset[-a,a]$.
Then,
$$
|\overline{\phi(x)} e^{-\sqrt{-1}x\zeta}| \le \|\phi\|_\infty
 e^{a\eta}
$$
for $\zeta = \xi - \sqrt{-1}\eta$ with $\eta>0$.
Here, $\|\phi\|_\infty$ denotes the $L^\infty$ norm.

Since $\gamma > -\frac{p+q-5}{2}$,
we may assume $\zeta=\xi-\sqrt{-1}\eta\in L'$ satisfies
$$
   \eta_1 \le \eta \le \eta_2
$$
for some constants $\eta_1$ and $\eta_2$ such that
$\eta_1>\frac{1}{2}$ if $|\xi|$ is sufficiently large.
Then, there exists a constant $C > 0$ such that
$|\psi(\zeta)| \le C|\xi|^{-2\eta}$
as $|\xi| \to \infty$,
by Lemma \ref{lem:espsi} (3).
Thus, if $|\xi|$ is sufficiently large, we have
$$
 |\psi(\zeta)\overline{\phi(x)}e^{-\sqrt{-1}x\zeta}|
\le C\|\phi\|_\infty |\xi|^{-2\eta_1}.
$$
Hence, $\psi(\zeta)\overline{\phi(x)} e^{-\sqrt{-1}x\zeta}$
is absolutely integrable on $L'\times[-a,a]$ (and therefore, on 
 $L'\times(-\infty,\infty)$).
Thus, we can apply Fubini's theorem.

{\bf Second equality} follows immediately from the 
definition \eqref{eqn:Fourier} of the Fourier transform.

{\bf Third equality} follows from Cauchy's integral formula.
First, we observe that
$\psi(\zeta) \overline{(\mathcal F\phi)(\bar\zeta)}$ is
holomorphic
in the domain between the two integral paths $(-\infty, \infty)$ and
$L'$ since its  poles lie only on
$\set{-\sqrt{-1}m}{m=1,2,\cdots}$ (see Lemma \ref{lem:espsi}(2)). 

Next, let us show
\begin{equation}\label{eqn:Cauchyerr}
\lim_{|\xi|\to\infty} \int_0^{\eta_0}
|\psi(\xi-\sqrt{-1}\eta)
\overline{(\mathcal{F}\phi)(\xi+\sqrt{-1}\eta)}| d\eta = 0
\end{equation}
for a fixed $\eta_0$ $(\ge \gamma+\frac{p+q-4}{2})$
(see Figure \ref{fig:zeta}).
To see \eqref{eqn:Cauchyerr}, 
take $a>0$ as before such that 
 $\operatorname{Supp}\phi \subset [-a, a]$.
Then, by the Paley--Wiener theorem, there exists a constant $C$
such that 
$
  |\overline{\mathcal F \phi(\xi+\sqrt{-1}\eta)}| \le
  Ce^{a\eta}
$.
Now combining this with Lemma \ref{lem:espsi} (3), we get
\begin{equation*}
  |\psi(\zeta)\overline{\mathcal F \phi(\bar\zeta)}|
  \le  C'|\xi|^{-2\eta}
  \quad\text{as $|\xi|\to\infty$}
\end{equation*} 
for $\zeta=\xi-\sqrt{-1}\eta$ and bounded $\eta$.
Hence, \eqref{eqn:Cauchyerr} is proved.
By Cauchy's integral formula, 
we get the third equality.

{\bf Last equality} is by the definition of 
the Fourier transform for tempered distributions:
$$
 (f, g)=(\mathcal F^{-1}f, \mathcal F^{-1}g)
 \quad f \in \mathcal{S}'(\mathbb{R}), \ g\in \mathcal S(\mathbb R).
$$

Hence, we have proved \eqref{eqn:Klk1}.
Now, the proof of Proposition \ref{prop:Klk} is completed.
\end{proof}

\section{Proof of Theorem \ref{thm:C} (2)}
\label{subsec:prC}

In this section, we complete the proof 
of Theorem \ref{thm:C} (2).
For this, it is sufficient to show the following proposition:
\begin{proposition}\label{prop:prC}
$\kappa_{l,k} = K_{l,k}$.
\end{proposition}

Here, we recall that the
 kernel distribution $\kappa_{l,k}$ of 
$T_{l,k}$ is given in Proposition \ref{prop:F}
and that $K_{l,k}$ is defined in \eqref{def:Klk}.

\begin{proof}
The proof makes use of the following:

\begin{lemSec}\label{lem:prC}
Let $\kappa_1, \kappa_2 \in \mathcal S'$. 
If there exists $\phi \in \mathcal S$ such that 
\begin{align}\label{eqn:Fnonzero}
&\mathcal F(\sigma_+ \phi)(\zeta)\ne0 
\quad\text{for any $\zeta\in\mathbb{R}$}, 
\\
&A_{\kappa_1}\phi=A_{\kappa_2}\phi,
\label{eqn:Aphi12}
\end{align}
then $\kappa_1=\kappa_2$.
Here, 
we recall from \eqref{eqn:Akappa} for the definition of $A_\kappa$.
\end{lemSec}

\begin{proof}[Proof of Lemma \ref{lem:prC}]
The identity \eqref{eqn:Aphi12} implies
$$\sigma_- \kappa_1 \ast \sigma_+ \phi
=\sigma_- \kappa_2 \ast \sigma_+ \phi$$
by the formula \eqref{eqn:conv} of $A_\kappa$.
Therefore, we have an identity
$$
 \mathcal F(\sigma_- \kappa_1)(\zeta)\cdot
 \mathcal F (\sigma_+ \phi)(\zeta)
=\mathcal F(\sigma_- \kappa_2)(\zeta) \cdot
 \mathcal F (\sigma_+ \phi)(\zeta)
$$ 
in $\mathcal{S}'(\mathbb{R})$ by taking their Fourier transforms.

On the other hand,
it follows from the assumption that
$\sigma_+\phi\in\mathcal{S}(\mathbb{R})$ 
and its Fourier transform
 $\mathcal F(\sigma_+ \phi) $ 
 does not vanish on $\mathbb{R}$,
we can divide the above identity by $\mathcal F(\sigma_+ \phi)(\zeta)$,
and obtain the following identity in $\mathcal{S}'(\mathbb{R})$:
$$
 \mathcal F(\sigma_- \kappa_1)(\zeta) 
=\mathcal F(\sigma_- \kappa_2)(\zeta).
$$ 
Hence, $\sigma_-\kappa_1 = \sigma_-\kappa_2$,
and in turn, $\kappa_1 = \kappa_2$.
\end{proof}
We want to apply Lemma \ref{lem:prC} with 
$\kappa_1:=\kappa_{l,k}$, \thinspace $\kappa_2:=K_{l,k}$, \thinspace 
$\phi:=f_{l,k}$
(see \eqref{def:f} for the definition).
The conditions in the lemma 
will be verified by the following steps.

{\bf Step 1.} $\kappa_{l,k}, K_{l,k}\in \mathcal S'$. 
These statements have been already proved in
 Propositions \ref{prop:F} and \ref{prop:Klk}.

{\bf Step 2.} $f_{l,k} \in \mathcal S$.
This has been proved in Lemma \ref{lem:flkas}.

{\bf Step 3.} 
{\it $\mathcal F(\sigma_+f_{l,k})(\zeta)$ has no
zero points on $\mathbb R$.}
This statement will follow readily from Claim \ref{clm:4.6.3}.
We note that we have assumed
$p\ge q \ge 2$ and
 $p+q\ge 6$. 
\begin{claim}\label{clm:4.6.3}
We recall from \eqref{def:Cases} the definitions of Cases 1 and 2.
Then,
\begin{align}\label{eqn:GG}
 \mathcal F(\sigma_+ f_{l,k})(\zeta)
={}
& \frac{1}{8\sqrt{\pi}}
  \Gamma(\frac{p+q-4}4+\frac{l+k+\sqrt{-1}\zeta}2)
\nonumber
\\
&
 \times
 \begin{cases}
  \Gamma(\frac{p-q}4+\frac{l-k+1+\sqrt{-1}\zeta}2)
  &\text{Case 1}, \\
  \Gamma(-\frac{p-q}4+\frac{-l+k+1+\sqrt{-1}\zeta}2)
  &\text{Case 2}. 
 \end{cases}
\end{align}
\end{claim} 
\begin{proof}
By the definition \eqref{def:f} of 
$f_{l,k}$ and the definition \eqref{def:sig} of $\sigma_+$,
we have
\begin{equation*}
   (\sigma_+f_{l,k})(x)
   = \frac{1}{\sqrt{2}} e^{(l+k+\frac{p+q}{2}-2)x}
     \times
     \begin{cases}
        \widetilde{K}_{\frac{q-3}{2}+k}(2e^x)
        &\text{Case 1},
     \\
        \widetilde{K}_{\frac{p-3}{2}+l}(2e^x)
        &\text{Case 2}.
     \end{cases}
\end{equation*}
In Case 1, we have
\begin{align*}
   \mathcal{F}(\sigma_+f_{l,k})(\zeta)
   &= \frac{1}{\sqrt{2\pi}\sqrt{2}}
      \int_{-\infty}^\infty
      e^{(l+k+\frac{p+q}{2}-2)x}
      \widetilde{K}_{\frac{q-3}{2}+k} (2e^x)
      e^{\sqrt{-1}x\zeta} dx
\\
   &= \frac{1}{2\sqrt{\pi}}
      \int_0^\infty
      r^{(l+k+\frac{p+q}{2}-3+\sqrt{-1}\zeta)}
      \widetilde{K}_{\frac{q-3}{2}+k}(2r)dr.
\end{align*}
Applying the formula \eqref{eqn:B1} of the 
Mellin transform for the $K$-Bessel function,
we obtain the right-hand side of \eqref{eqn:GG}.
Likewise, in Case 2,
$\mathcal{F}(\sigma_+f_{l,k})(\zeta)$ is equal to
\begin{align*}
&  \frac{1}{\sqrt{2\pi}\sqrt{2}}\int_{-\infty}^\infty 
   e^{(l+k+\frac{p+q}{2}-2)x}
   \widetilde{K}_{\frac{p-3}{2}+l} (2e^x) e^{\sqrt{-1}x\zeta} dx.
\end{align*}
Switching the role of $(p,l)$ and $(q,k)$,
we see \eqref{eqn:GG} holds also in Case 2.
\end{proof}

{\bf Step 4.} $A_{\kappa_{l,k}}f_{l,k}=A_{K_{l,k}}f_{l,k}$.

To see this, we prepare the following explicit formulas. 
As we shall see below,
the proof of (1) is algebraic by using the fact that $\pi(w_0)$ acts on
each $K$-type as $\pm\operatorname{id}$.
On the other hand, the proof of (2) is based on 
an explicit integral computation.

\begin{claim}\label{clm:prC}
Let $a(l,k)$ be as in  \eqref{def:alk}.

{\rm 1)} $A_{\kappa_{l,k}}f_{l,k}=(-1)^{a(l,k)+\frac{p-q}2}f_{l,k}$.

{\rm 2)} $A_{K_{l,k}}f_{l,k}=(-1)^{a(l,k)+\frac{p-q}2}f_{l,k}$.
\end{claim}

\begin{proof}[Proof of Claim \ref{clm:prC}]

1) 
The function 
\index{A}{f0lkr@$f_{l,k}(r)$}%
$f_{l,k}$ belongs to the $K'$-invariant subspace
\index{A}{Hlk@$H_{l,k}$}%
$H_{l,k}$ (see \eqref{def:H}),
and therefore, by Theorem \ref{thm:K} (3), 
we have 
$$
\pi(w_0)(f_{l,k}\phi \psi)=
(-1)^{a(l,k)+\frac{p-q}2}f_{l,k}\phi \psi
$$
for $\phi \in \Har{l}{p-1}$ and $\psi \in \Har{k}{q-1}$.
In light of the definition of $T_{l,k}$ (see Theorem \ref{thm:C} (1)),
this implies
$$
T_{l,k}f_{l,k}=(-1)^{a(l,k)+\frac{p-q}2}f_{l,k}.
$$
By the definition of $A_{\kappa_{l,k}}$ (see Proposition 
\ref{prop:F}),
the first statement follows.

2) First, we treat Case 1, namely, the case
$\frac{p-q}{2}+l-k \ge 0$.
Then, 
\begin{align}\notag
&(A_{K_{l,k}}f_{l,k})(r) 
=\frac{1}{2\pi}\int_0^\infty
 \biggl(\int_{L'}\psi(\zeta)
(rr')^{-\frac{p+q-4}2+\sqrt{-1}\zeta}d\zeta
\biggr)f_{l,k}(r')r'^{p+q-5}dr'  \\  \notag
={}&\frac{1}{2\pi}\int_{L'}
 r^{-\frac{p+q-4}2+\sqrt{-1}\zeta}
 \biggl(\int_0^\infty r'^{\frac{p+q-6}2+\sqrt{-1}\zeta}
 f_{l,k}(r')dr'\biggr)
 \psi(\zeta)d\zeta  \\ \notag
={}&\frac{(-1)^{l+\frac{p-q}2}}{8\pi}
 \int_{L'}
\Gamma\Bigl(\frac{p+q-4}4+\frac{l+k-\sqrt{-1}\zeta}2\Bigr)
\nonumber
\\
&\kern\the\textwidth\llap{$\displaystyle
\times\Gamma\Bigl(\frac{p-q}4+\frac{l-k+1-\sqrt{-1}\zeta}2\Bigr)
r^{-\frac{p+q-4}2+\sqrt{-1}\zeta} d\zeta 
\phantom{==}$}
\nonumber
\\
\label{int:GG}
 ={}&\frac{(-1)^{l+\frac{p-q}2}}{8\pi\sqrt{-1}}
 \int_L \Gamma\Bigl(\frac{l+k-s}2 \Bigr) 
                          \Gamma\Bigl(\frac{-q+3+l-k-s}2 \Bigr)
  r^s ds    \\
 ={}& \frac12(-1)^{l+\frac{p-q}2} G^{20}_{02}
 \Bigl(r^2 \Big| \frac{l+k}2, \frac{-q+3+l-k}2 \Bigr)
\notag \\ 
 ={}&(-1)^{l+\frac{p-q}2} r^{-\frac{q-3}{2}+l}K_{\frac{q-3}{2}+k}(2r) 
\nonumber\\
 ={}& (-1)^{a(l,k)+\frac{p-q}2}f_{l,k}(r).  \notag
\end{align}

Let us explain the above equalities in more details.

{\bf First equality}.
This is by the integral
 expression of $K_{l,k}$ (see Lemma \ref{lem:Kpsi}) 
and the definition \eqref{eqn:Akappa} of $A_{K_{l,k}}$.

{\bf Second equality}.
We recall the upper estimate of $|\psi(\zeta)|$
given in Lemma \ref{lem:espsi} (3)
and the asymptotic behavior of 
$f_{l,k}(r')$
(see \eqref{eqn:flkinfty}
and \eqref{eqn:flk0}).
Then, in light of
\begin{equation*}
 | \psi(\zeta)
r'^{\frac{p+q-6}2+\sqrt{-1}\zeta}f_{l,k}(r')| 
\le{} |\psi(\zeta)| r'^{\frac{p+q-6}2+\eta}
 | f_{l,k}(r')|
\quad\text{for $\zeta=\xi-\sqrt{-1}\eta$},
\end{equation*}
the second equality follows from Fubini's theorem.

{\bf Third equality} is by Lemma \ref{lem:flkMellin}.

{\bf Fourth equality} is from the change of the variable
as before: $s=\sqrt{-1}\zeta-\frac{p+q-4}2$.

{\bf Fifth equality}.  The poles of the integrand 
in \eqref{int:GG} are of the form $l+k+2a$ $(a\in\mathbb{N})$ or 
$-q+3+l-k+2a$ $(a\in\mathbb{N})$.
These lie on the right of the contour $L$ because of
the assumption 
$\frac{p-q}2 +l-k \ge 0$. Hence, the fifth equality
holds by the integral expression of Meijer's $G$-function
(see \eqref{def:G} in Appendix).

{\bf Sixth equality} follows from the reduction formula of 
the $G$-function (see (\ref{eqn:GK})).

{\bf Seventh equality} is by the definition \eqref{def:f} of $f_{l,k}$
and the definition \eqref{def:alk} of $a(l,k)$.

Case 2 can be treated in the same manner.
In this case, the integral 
$$
  \int_L \Gamma( \frac{l+k-s}2 ) \Gamma(\frac{-p+3-l+k-s}2) r^s ds
$$
arises instead of \eqref{int:GG}. But again, by the assumption
$\frac{p-q}2+l-k \le 0$, this defines the $G$-function which reduces to 
$f_{l,k}$ by the same reduction formula.
\end{proof}

{\bf Step 5.} $\kappa_{l,k} = K_{l,k}$.

Steps 3 and 4 imply $\kappa_{l,k} = K_{l,k}$ by Lemma \ref{lem:prC}.
Thus, Proposition \ref{prop:Klk} is proved.
\end{proof}

Now the proof of Theorem \ref{thm:C} finishes.

\chapter{Main theorem}
\label{sec:A}

This chapter is a highlight of this book.
We find an explicit formula for the `Fourier transform'
$\mathcal{F}_C$ on the isotropic cone,
in other words,
 we find an integral kernel 
for the unitary inversion operator $\pi(w_0)$ on the Schr\"{o}dinger
model $L^2(C)$ of the minimal representation.
The main result is Theorem \ref{thm:A}.

\section{Result of this chapter}
\label{subsec:A}

Let $C$ be the conical variety 
$\{\xzeta\in\mathbb{R}^{p+q-2} \setminus \{0\}: Q(\xzeta) = 0\}$
where
$Q(\xzeta) = \xzeta_1^2 + \dots + \xzeta_{p-1}^2 
 - \xzeta_p^2 - \dots - \xzeta_{p+q-2}^2$.
We recall from the Introduction that
the generalized function $K(\xzeta, \xzeta')$ on $C \times C$ 
is defined by the following formula:
\begin{equation}\label{def:K}
\index{A}{Kxx'@$K(x,x')$}%
  K(\xzeta, \xzeta')\equiv K(p,q; \xzeta, \xzeta'):=c_{p,q} 
\Phi_{p,q}(\langle \xzeta, \xzeta' \rangle),
\end{equation}
where $\langle \cdot, \cdot \rangle$ stands for  
the standard (positive definite) inner product on 
$\mathbb R^{p+q-2}$.
The constant $c_{p,q}$ and the 
\index{B}{Bessel distribution}%
\textit{Bessel distribution} $\Phi_{p,q}(t)$ 
are given as follows 
(see \eqref{def:Psi0}, \eqref{def:Psi+}, and \eqref{def:Psi}):
\begin{alignat}{2} \label{def:c}
  c_{p,q}&:= \frac{2(-1)^\frac{(p-1)(p+2)}2}{\pi^\frac{p+q-4}{2}}, \\
\label{def:Phi}
\index{A}{1Phi_{p,q}(t)@$\Phi_{p,q}(t)$}%
  \Phi_{p,q}(t)&:=\begin{cases}
  \Phi^+_{\frac{p+q-6}2}(t) \qquad 
  \text{if $\operatorname{min}(p,q)=2$,} 
 \\ \Psi_{\frac{p+q-6}2}^+(t) \qquad \text{if $p, q >2$ are both even,}
  \\
  \Psi_{\frac{p+q-6}2}(t) \qquad \text{if $p, q>2$ are both odd.}   
   \end{cases}
\end{alignat}

We are ready to state the explicit formula for the 
\index{B}{unitary inversion operator}%
unitary inversion operator:
\begin{thmSec}[Integral formula for the unitary inversion operator]
\label{thm:A}
\quad\newline
Let $(\pi, L^2(C))$ be the Schr\"{o}dinger model of the 
minimal representation of $G=O(p,q)$ for
$p, q\ge 2$ and $p+q\ge 6$ even,
and $w_0 = \begin{pmatrix} I_p & 0 \\ 0 & -I_q\end{pmatrix}$.
Then the unitary operator $\pi(w_0):L^2(C) \to 
L^2(C)$ is given by the following integro-differential operator:
\begin{equation*}
\pi(w_0)u(\xzeta)=\int_C K(\xzeta, \xzeta')u(\xzeta')d\mu(\xzeta'),
\quad  u \in L^2(C).
\end{equation*}
\end{thmSec}
\noindent
The right-hand side of \eqref{int:A} involves a singular integral.
It factors through the Radon transform (see
Section \ref{subsec:Radon}),
and we shall see that the right-hand side of \eqref{int:A}
 is well-defined for 
any compactly supported smooth function $u$ on $C$,
and extends as a unitary operator on $L^2(C)$.

As for the Bessel distribution $\Phi_{p,q}$,
we shall give a Mellin--Barnes type integral formula for $\Phi_{p,q}$ in 
Section \ref{subsec:intPsi}, and the differential
equation satisfied by $\Phi_{p,q}$ in 
Section \ref{subsec:ode}.

Since the action of the maximal parabolic subgroup $\overline{\Pmax}$
on $L^2(C)$ is of a simple form 
(see \eqref{eqn:rM}--\eqref{eqn:rN}), 
Theorem \ref{thm:A} gives an explicit action of the whole group $G$ on
$L^2(C)$ because $G = \overline{\Pmax} \coprod \overline\Pmax w_0 
\overline\Pmax$.

Theorem \ref{thm:A} immediately yields two 
corollaries about the Plancherel formula 
and the reciprocal formula of our integral transform.

\begin{corollary}[%
\index{B}{Plancherel formula!FC@---, $\protect\mathcal{F}_C$}%
Plancherel formula]\label{cor:A1}
Let 
\index{A}{FC@$\mathcal{F}_C$}%
$\mathcal{F}_C: L^2(C) \to L^2(C)$ be an integral transform 
against the kernel $K(\xzeta, \xzeta')$ (see \eqref{def:K}),
namely, 
\begin{equation}\label{int:A}
   (\mathcal{F}_C u)(x)
   :=
   \int_C K(\xzeta, \xzeta')u(\xzeta')d\mu(\xzeta').
\end{equation}
Then $\mathcal{F}_C$ is unitary:
$$
  \|\mathcal{F}_C u\|_{L^2(C)} = \|u\|_{L^2(C)}.
$$
\end{corollary}

Since the group law $w_0^2=1$ in $O(p,q)$ implies $\pi(w_0)^2=
\operatorname{id}$ on $L^2(C)$, we immediately obtain the following 
corollary to Theorem  \ref{thm:A}, which 
can also be viewed as giving the inversion formula 
$\mathcal{F}_C^{-1}=\mathcal{F}_C$.

\begin{corollary}[%
\index{B}{reciprocal formula!FC@---, $\protect\mathcal{F}_C$}%
Reciprocal formula]
\label{cor:A2}
  Retain the notation as in Corollary \ref{cor:A1}. 
  The unitary operator $\mathcal{F}_C$ is of order two
  in $L^2(C)$. Namely, we have
  the following reciprocal relation:
\begin{align*}
 (\mathcal{F}_C u)(\xzeta) 
&   = 
  \int_C K(\xzeta, \xzeta') u(\xzeta') d\mu(\xzeta'),
\\
  u(\xzeta)
&   =
   \int_C K(\xzeta, \xzeta') (\mathcal{F}_C u)
   (\xzeta') d\mu(\xzeta').
\end{align*}
\end{corollary}

\begin{remNonumber}[Comparison with the 
\index{B}{Schr\"{o}dinger}%
Schr\"{o}dinger model of the
\index{B}{Weil representation}%
Weil representation, see Section \ref{subsec:1.6a}]
In the case of the Schr\"{o}dinger model of the Segal--Shale--Weil
representation $\varpi$ of the metaplectic group $\Mp(n,\mathbb{R})$,
the corresponding `inversion' element 
\index{A}{w0'@$w'_0$}%
$w'_0$ acts on
$L^2(\mathbb{R}^n)$ as $e^{\frac{\sqrt{-1}n\pi}{4}}
\index{A}{FRn@$\mathcal{F}_{\mathbb{R}^n}$}%
\mathcal{F}_{\mathbb{R}^n}$,
where $\mathcal{F}_{\mathbb{R}^n}$ denotes
the Fourier transform.
We note that  $(w'_0)^4$ gives the unique non-trivial element $\xi_0$
in the
kernel of the metaplectic covering
$Mp(n,\mathbb{R})\to Sp(n,\mathbb{R})$,
and $\varpi(\xi_0)=-\operatorname{id}$. 
This fact
reflects the identity $\mathcal F_{\mathbb{R}^n}^4= \operatorname{id}$ on $L^2(\mathbb R^n)$.
Thus, the above two
 corollaries can be interpreted as the counterparts to 
the Plancherel formula and the equality
$\mathcal F_{\mathbb{R}^n}^4=\operatorname{id}$ of the Fourier transform 
$\mathcal F_{\mathbb{R}^n}$ on 
$\mathbb R^n$.
\end{remNonumber}

\begin{remNonumber}
In \cite[Corollary 6.3.1]{xkmano2},
we gave a different proof of the same Plancherel
and reciprocal formulas in the case $q = 2$ 
based on analytic continuation
of 
\index{B}{holomorphic semigroup}%
holomorphic semigroup of operators.
\end{remNonumber}

This chapter is organized as follows.
In Section \ref{subsec:Radon}, we analyze the integral transform
\eqref{int:A} by means of the (singular) Radon transform. 
In particular, we prove that the integral transform \eqref{int:A} is
well-defined for $u \in C_0^\infty(C)$ in the sense of distributions. 
The second step of the proof of Theorem \ref{thm:A} 
is to use the restriction from $G$ to
$K'=K\cap \Mmax_+$ (see Section \ref{subsec:act}) and is
to show the $(l, k)$-th spectrum 
of the unitary inversion operator $\mathcal{F}_C$ coincides with 
the radial part $T_{l,k}$ of 
$\mathcal{F}_C$ when restricted to each 
$K'$-isotypic component $\Har{l}{p-1} \otimes \Har{k}{q-1}$
(see Lemma \ref{lem:res}).
The latter operator $T_{l,k}$ was studied in details in the previous
chapter (see Theorem \ref{thm:C}).
Section \ref{subsec:spec} explains a general formula giving 
spectra of a $K'$-intertwining integral operator.
In Section \ref{subsec:prA}, 
we reduce the proof of Theorem \ref{thm:A} to the identity of spectra for specific
$K'$-intertwining integral operators.
This identity is verified in Sections \ref{subsec:5.5} and
\ref{subsec:5.6}.
Some technical parts (e.g.\@ Mellin--Barnes type integral formulas of
Bessel distributions) will be postponed until Chapters \ref{sec:Diffeq}
and \ref{sec:app}.

\section{Radon transform for the isotropic cone $C$}
\label{subsec:Radon}
This section studies the kernel $K(\xzeta,\xzeta')$.
The subtle point in defining $K(\xzeta,\xzeta')$ consists of the
following two facts:
\begin{enumerate}[1)]
\item  
The distribution $\Phi_{p,q}(t)$ is not locally integrable near $t=0$. 
\item  
The level set $\langle \xzeta,\xzeta'\rangle = t$ is not a
regular submanifold in $C \times C$ if $t=0$.
\end{enumerate}
In fact, 
the distribution $\Phi_{p,q}(t)$ involves a linear combination of
distributions $\delta^{(k-1)}(t)$ and $t^{-k}$
($k=1,2,\dots,\frac{p+q-6}{2}$) as we shall see in Section
\ref{subsec:intPsi} on the one hand, and the differential form
\begin{equation*}
   dQ(\xzeta) \wedge dQ(\xzeta') \wedge d(\langle \xzeta,\xzeta' \rangle -
   t)
\end{equation*}
of $\xzeta,\xzeta'$ vanishes if $(\xzeta,\xzeta')$ belongs to the
submanifold
\begin{equation*}
   Y := \{(\xzeta,\xzeta') \in C \times C : \mathbb{R} \xzeta
          = \mathbb{R}(I_{p,q} \xzeta') \}
\end{equation*}
on the other hand. 
Here, we note that
$Y \subset \{(\xzeta,\xzeta'): \langle \xzeta,\xzeta' \rangle = 0 \}$. 

Our idea to give a rigorous definition of $K(\xzeta,\xzeta')$ is
to factorize the transform
\eqref{int:A} by using the singular Radon transform and then to find a
Paley--Winer type theorem of the singular Radon transform.

Let $\delta$ denote the Dirac delta function of one variable.
The 
\index{B}{Radon transform}%
{\it Radon transform} 
of a function $\varphi$ on 
$\mathbb R^{\n}$ is defined by
the formula (see for example \cite[Chapter I]{xGeGrVi}):
\begin{equation}\label{def:Radon}
(R\varphi)(\xzeta, t):= \int_{\mathbb R^{\n}} \varphi(\xzeta')
\delta(t- \langle \xzeta, \xzeta' \rangle) d\xzeta',
\end{equation}
for $\xzeta \in \mathbb R^{\n} \setminus \{0\}$, $t \in \mathbb R$.

The Radon transform $R\varphi$ is well-defined,
for example, for a compactly supported continuous function $\varphi$. 
More generally, $R\varphi$ makes sense if $\varphi$ is a compactly
supported distribution such that the multiplication of two
distributions $\varphi(\xzeta')$ and
$\delta(t-\langle\xzeta,\xzeta'\rangle)$ makes sense.

Now we recall the injective map (see \eqref{def:T})
$$
T: L^2(C) \to \mathcal S'(\mathbb R^{p+q-2}),
\quad u \mapsto u \delta(Q)
$$  
yields a compactly supported distribution 
$Tu$ if $u \in C_0^\infty(C)$. 
In this context, what we need here is the following result:

Let $C_0^k(\mathbb{R})$ denote the space of compactly supported
functions on $\mathbb{R}$ with continuous derivatives up to $k$.
\begin{lemSec}
\label{lem:reg}
Suppose 
\index{A}{C0infty@$C_0^\infty(C)$}%
$u\in C_0^\infty(C)$.

{\rm 0)} The Radon transform
$\index{A}{R@$\mathcal{R}$}%
\mathcal{R}u(x,t):=R(Tu)(\xzeta,t)$ is well-defined and continuous as a function of 
$(\xzeta,t)\in C \times (\mathbb R \setminus \{0\})$.
Moreover, there exists $A>0$ such that
\begin{equation*}
   \operatorname{Supp} \mathcal{R}u \subset
   \{(\xzeta,t)\in C \times (\mathbb{R}\setminus\{0\}):
     t \le A |\xzeta| \},
\end{equation*}
where $|\xzeta| := (\xzeta_1^2+\dots+\xzeta_{p+q-2}^2)^{\frac{1}{2}}$.
In particular, 
$\mathcal{R}u(\xzeta,t)$ vanishes 
if\/ $\vert t\vert$ is sufficiently large for a fixed $\xzeta\in C$.

{\rm 1)} If $p,q>2$ and $p+q\ge 8$, then
$\mathcal{R}u(\xzeta,t)$ extends continuously to $t=0$ and
$\mathcal{R}u(\xzeta,\cdot) \in C_0^k(\mathbb{R})$
where $k := \frac{p+q-8}{2}$.

{\rm 2)} If $\operatorname{min}(p,q)=2$, then
$\mathcal{R}u(\xzeta,t)$ is bounded on
$C \times (\mathbb R\setminus \{0\})$. 

{\rm 3)} If $p,q>2$ and $p+q=6$ (namely, $(p,q)=(3,3)$), then
there exists $C \equiv C(\xzeta)>0$ such that
\begin{equation*}
   | R(Tu)(\xzeta,t)| \le C \left|\log|t| \right|
\end{equation*}
if $t$ is sufficiently small.

\end{lemSec}

\begin{proof}
See \cite{xkmano5}.
\end{proof}

We note that 
$$
\mathcal{R}u(\xzeta,t)
=\int_Cf(\xzeta')\delta(\langle\xzeta,\xzeta'\rangle-t)d\mu(x')
$$
 is well-defined for 
$(\xzeta,t) \in (\mathbb{R}^{p+q}\setminus \{0\}) \times
 (\mathbb{R} \setminus \{0\})$,
but we need here only the case where $\xzeta\in C$.

We rewrite the right-hand side of \eqref{int:A} 
for $u \in C_0^\infty(C)$ as follows:
\begin{align}
(\mathcal{F}_C u)(\xzeta)
:={}& \int_C K(\xzeta,\xzeta')u(\xzeta')d\mu(\xzeta')
\nonumber
\\
={}&
  c_{p,q} \int_{\mathbb R^{p+q-2}}  
  \Phi_{p,q}(\langle \xzeta, \xzeta' \rangle) Tu(\xzeta') d\xzeta'  
\qquad\text{by \eqref{def:K}}
\nonumber
\\
={}&
  c_{p,q} \int_{\mathbb{R}} \int_{\langle\xzeta,\xzeta'\rangle=t}
  \Phi_{p,q} (\langle\xzeta,\xzeta'\rangle) Tu(\xzeta')
  d\omega(\xzeta')dt
\nonumber
\\
={}&
c_{p,q} \int_\mathbb R \Phi_{p,q}(t) \mathcal{R}u(\xzeta, t) dt.   
\label{eqn:Radon}
\end{align}

\begin{lemSec}
\label{lem:Radon}
The right-hand side of 
\eqref{eqn:Radon} is well-defined for $u\in C_0^\infty(C)$.
\end{lemSec}

\begin{proof}
It follows from Theorem \ref{prop:Psiint} which we shall prove later
and from the definition \eqref{def:Phi} of the distribution
$\Phi_{p,q}(t)$ that $\Phi_{p,q}(t)$ has the following decomposition: 
\begin{equation*}
   \Phi_{p,q}(t) = 
\index{A}{1Phi_{p,q}^{\mathrm{reg}}(t)@$\Phi_{p,q}^{\mathrm{reg}}(t)$}%
\Phi_{p,q}^{\mathrm{reg}}(t) 
                 + 
\index{A}{1Phi_{p,q}^{\mathrm{sing}}(t)@$\Phi_{p,q}^{\mathrm{sing}}(t)$}%
\Phi_{p,q}^{\mathrm{sing}}(t),
\end{equation*}
where $\Phi_{p,q}^{\mathrm{reg}}(t)$ and
$\Phi_{p,q}^{\mathrm{sing}}(t)$ are distributions on $\mathbb{R}$
such that
\begin{enumerate}
    \renewcommand{\labelenumi}{\theenumi)}
\item  
$\Phi_{p,q}^{\mathrm{reg}}(t) |t|^{-\epsilon}$
is a locally integrable function on $\mathbb{R}$ for any sufficiently
small $\epsilon\ge 0$,
\item  
\begin{equation*}
\addtolength{\abovedisplayshortskip}{-3ex}
\addtolength{\abovedisplayskip}{-3ex}
   \Phi_{p,q}^{\mathrm{sing}}(t)
  = \begin{cases}
       0
       &\text{if \ $\min(p,q)=2$},
    \\[1ex]
    \displaystyle
       -\sum_{k=1}^{\frac{p+q-6}{2}}
        \frac{(-1)^{k-1}}{2^k(m-k)!} \, \delta^{(k-1)} (t)
       &\text{if \ $p,q > 2$ both even},
    \\[3ex]
    \displaystyle
       -\frac{1}{\pi} \sum_{k=1}^{\frac{p+q-6}{2}}
        \frac{(k-1)!}{2^k(m-k)!} \, t^{-k}
       &\text{if \ $p,q > 2$ both odd}.
    \end{cases}
\end{equation*}
\end{enumerate}
We note that $\Phi_{p,q}^{\mathrm{sing}}(t)\ne0$ only if
$p,q>2$ and $p+q\ge8$.
Combining with Lemma \ref{lem:reg},
we see that the right-hand side of \eqref{eqn:Radon} is well-defined
in all the cases.
\end{proof}

Lemma \ref{lem:Radon} defines a linear map
\begin{equation*}
\index{A}{FC@$\mathcal{F}_C$}%
   \mathcal{F}_C: C_0^\infty(C) \to C^\infty(C),
\end{equation*}
and defines $K(\xzeta,\xzeta')$ as a distribution on the direct product
manifold $C \times C$.
In Section \ref{subsec:prA}, 
we shall see that the image 
$\mathcal{F}_C(C_0^\infty(C))$ 
is contained in $L^2(C)$,
and $\mathcal{F}_C$ extends to a unitary operator on $L^2(C)$,
which in turn equals the unitary inversion operator $\pi(w_0)$.

\begin{remNonumber}
The plane wave expansion 
gives a decomposition of the 
\index{B}{Euclidean Fourier transform}%
Euclidean Fourier transform
\index{A}{FRn@$\mathcal{F}_{\mathbb{R}^n}$}%
$\mathcal{F}_{\mathbb{R}^n}$ on $L^2(\mathbb{R}^n)$
into 
the one-dimensional integral transform 
\index{B}{Mellin transform}%
(Mellin transform)
and the 
\index{B}{Radon transform}%
Radon transform, namely:
$$
  (\mathcal F_{\mathbb{R}^n} u)(\xzeta)
= c_n \langle \Psi, (Ru)(\xzeta, \cdot) \rangle,
$$
where $c_n:= \bigl( \frac{1}{2\pi} \bigr)^\frac{n}2$ and $\Psi(t):=
e^{\sqrt{-1}t}$ (e.g. \cite[Chapter I, \S 1.2]{xGeGrVi}).
In this sense, the formula \eqref{eqn:Radon} can be regarded as 
an analog of the plane wave expansion 
for the unitary operator $\pi(w_0)$ on $L^2(C)$.
\end{remNonumber}

\section{Spectra of $K'$-invariant operators on $S^{p-2}\times
 S^{q-2}$}
\label{subsec:spec}
The expansion into spherical harmonics
$$
L^2(S^{n-1}) \simeq
\sideset{}{^\oplus}\sum_{l=0}^\infty 
\index{A}{HjRm@$\mathcal{H}^j(\mathbb{R}^m)$}%
\mathcal{H}^l(\mathbb{R}^n)
$$
gives a 
\index{B}{multiplicity-free decomposition}%
multiplicity-free decomposition of $O(n)$ into its irreducible
representations (see Appendix \ref{subsec:H}),
and consequently, any $O(n)$-intertwining operator on
$L^2(S^{n-1})$ acts on $\mathcal{H}^l(\mathbb{R}^n)$ as a
scalar multiplication owing to Schur's lemma.
The scalar is given by the Funk--Hecke formula
(see \cite[\S9.7]{xaar}, see also \cite[Lemma 5.5.1]{xkmano2}):
for an integrable function $h$ on the interval $[-1,1]$ and for 
$\phi \in \mathcal{H}^l(\mathbb{R}^n)$,
\begin{equation*}
   \int_{S^{n-1}} h(\langle \omega,\omega'\rangle) 
   \phi(\omega') d\omega'
   = c_{l,n}(h) \phi(\omega),
\end{equation*}
where the eigenvalue $c_{l,n}(h)$ is given by
\begin{equation} \label{eqn:FunkHecke}
   c_{l,n}(h)
   = \frac{2^{n-2} \pi^{\frac{n-2}{2}}l!}{\Gamma(n-2+l)}
     \int_{-1}^1 h(x) \widetilde{C}_l^{\frac{n-2}{2}}(x)
     (1-x^2)^{\frac{n-3}{2}} dx.
\end{equation}
Here, 
\index{A}{Clmuxtilde@$\tilC_l^\mu(x)$}%
$\tilC_l^\mu(x)$ stands for the normalized 
\index{B}{Gegenbauer polynomial}%
Gegenbauer polynomial (see Section \ref{subsec:Ge}).

Likewise, any $K'$-intertwining operator on
$L^2(S^{p-2} \times S^{q-2})$
acts on the subspace $\mathcal{H}^l(\mathbb{R}^{p-1}) \otimes
\mathcal{H}^k(\mathbb{R}^{q-1})$ as a scalar multiplication for each
$k,l \in \mathbb{N}$ 
(we recall $K' \simeq O(p-1) \times O(q-1)$).
In this section, we determine this scalar for specific intertwining
(integral) operators.
In particular, the scalar in Example \ref{ex:Ri} will
be used in the proof of our main theorem
(Theorem \ref{thm:A}).

We begin with a general setup for a $K'$-intertwining operator on
$L^2(S^{p-2} \times S^{q-2})$.
Let $h$ be
 an integrable function of two variables on
$[-1,1]\times [-1,1]$.
We consider the following integral transform:
\begin{multline}
 B_h: C(S^{p-2}\times S^{q-2}) \to 
 C(S^{p-2}\times S^{q-2}), \\
 \varphi(\omega, \eta)\mapsto \int_{S^{p-2}\times S^{q-2}}
 h(\langle \omega, \omega' \rangle, \langle \eta, \eta' \rangle)
 \varphi(\omega', \eta') d\omega' d\eta'.
\end{multline}

\begin{lemSec}\label{lem:spec}
$B_h$ acts on each $K'$-type $\Har{l}{p-1}\otimes \Har{k}{q-1}$ by a 
scalar multiplication of $\alpha_{l,k}(h)\in \mathbb C$. 
The spectrum $\alpha_{l,k}(h)$ is given by the following formulas.

{\rm 1)} If $\operatorname{min}(p,q)=2$, say $q=2$, then
for $k=0,1$,
\begin{equation}\label{eqn:a2}
\alpha_{l,k}(h)=\frac{2^{p-3}\pi^\frac{p-3}2 l!}{\Gamma(p-3+l)}
\int_{-1}^{1}(U_k h)(x)\tilC_l^\frac{p-3}2 (x)(1-x^2)^\frac{p-4}2 dx,
\end{equation}
where we set  
\begin{equation}\label{def:h2}
(U_k h)(x):=h(x,1)+(-1)^k h(x,-1).
\end{equation}
For $k\ge2$,
$\alpha_{l,k}(h)=0$.

{\rm 2)} If $p,q>2$, then
\begin{multline} \label{eqn:a}
 \alpha_{l,k}(h)= \frac{2^{p+q-6}\pi^{\frac{p+q-6}2}l! \, k!}
 {\Gamma(p-3+l)\Gamma(q-3+k)} \\
 \times \int_{-1}^1 \int_{-1}^{1} h(x,y)\tilC_l^{\frac{p-3}2}(x)
 \tilC_k^{\frac{q-3}2}(y)(1-x^2)^\frac{p-4}2
 (1-y^2)^\frac{q-4}2 dx dy.
\end{multline}

\end{lemSec}
\begin{proof}
1)
If $q=2$, then $S^{p-1}\times S^{q-1} = S^{p-1} \coprod S^{p-1}$ 
(disjoint union),
and $\mathcal{H}^k(\mathbb{R}^{q-1})=0$ if $k\ge2$
(see Section \ref{subsec:H}).
Then, the formula \eqref{eqn:a2} is essentially the 
\index{B}{Funk--Hecke formula}%
Funk--Hecke formula 
\eqref{eqn:FunkHecke} for $S^{p-1}$.

2) 
Applying \eqref{eqn:FunkHecke} to each factor, we get \eqref{eqn:a}.
\end{proof}

Let us give some examples of 
the explicit computation of  
spectra $\alpha_{l,k}(h)$.

\begin{example}[%
\index{B}{Riesz potential}%
Riesz potential] \label{ex:Ri}
Consider the following Riesz potential for 
$\operatorname{Re}\lambda > -1$:
\begin{align}\label{def:Ri}
 h_\lambda^\pm (x,y)
 &:=\frac{(x+y)_\pm^\lambda}
 {\Gamma(\lambda+1)},
\nonumber
\\
 &= \begin{cases}
    \displaystyle
       \frac{(x+y) ^\lambda}{\Gamma(\lambda+1)}
       &\text{if \ $\epsilon(x+y)>0$},
     \\[1ex]
       0
       &\text{if \ $\epsilon(x+y)\le0$},
    \end{cases}
\end{align}
where $\epsilon=\pm1$.
Then, the spectrum $\alpha_{l,k}(h_\lambda^\pm)$ for the
$K'$-intertwining operator $B_{h_\lambda^\pm}$
amounts to
\begin{align*}
 \alpha_{l,k}(h_\lambda^\pm)&=
 \frac{2^{1-\lambda}\pi^\frac{p+q-2}2(\pm 1)^{l+k}\Gamma(\lambda+
 \frac{p+q-4}2)}
 {\Gamma(\frac{\lambda+p+q-4+l+k}2)
  \Gamma(\frac{\lambda+p-1+l-k}2)
  \Gamma(\frac{\lambda+q-1-l+k}2)
  \Gamma(\frac{\lambda-l-k+2}2)}\\
&= \frac{(\pm1)^{l+k}}{\pi}
   \Gamma\Bigl(\lambda+\tfrac{p+q-4}2\Bigr)
   \sin \Bigl(\tfrac{\lambda-l-k+2}2 \pi\Bigr)
   \sin \Bigl(\tfrac{\lambda+q-1-l+k}2 \pi\Bigr)
   \gamma_{l,k}(\lambda),
\end{align*}
where we set
\begin{equation}\label{eqn:gammalk}
\gamma_{l,k}(\lambda) :=
2^{1-\lambda} \pi^{\frac{p+q-4}{2}}
\frac{\Gamma(\frac{l+k-\lambda}{2})\Gamma(\frac{-q+3+l-k-\lambda}{2})}
     {\Gamma(\frac{\lambda+p+q+l+k-4}{2})\Gamma(\frac{\lambda+p-1+l-k}{2})}.
\end{equation}
\end{example}
\begin{proof}[Proof of Example \ref{ex:Ri}]
Use \eqref{eqn:a}.
We postpone the actual computation of the integral
(the first equation of $\alpha_{l,k}(h_\lambda^\pm)$)
until Appendix
(see Lemma \ref{lem:frac}
with $\mu=\frac{p-3}2, \nu=\frac{q-3}2$).
In the second equation of $\alpha_{l,k}(h_\lambda^\pm)$, 
we have used the functional equation $\Gamma(z)\Gamma(1-z)=
\frac{\pi}{\sin (z\pi)}$.
\end{proof}

We define a kernel function $h_\lambda(x,y) \equiv h^{p,q}_\lambda(x,y)$ 
with parameter $\lambda$ as follows: 
\begin{equation}\label{def:hpm}
h_\lambda(x,y):=
\frac{\Gamma(-\lambda)}{\Gamma(\lambda+\frac{p+q-4}2)}
\times
\begin{cases}
  (x+y)^\lambda_+ &\text{if $p,q>2$ both even,}\\[1ex]
 \left(\dfrac{(x+y)^\lambda_+}{\tan \lambda \pi}
 +\dfrac{(x+y)^\lambda_-}{\sin \lambda \pi}
 \right) &\text{if $p,q>2$ both odd.}
\end{cases}
\end{equation}

\begin{proposition}
\label{prop:R}
Let $\operatorname{Re}\lambda>-1$
.
For a kernel function $h_\lambda$ (see \eqref{def:hpm}), 
the spectrum $\alpha_{l,k}(h_\lambda)$ given in
Lemma \ref{lem:spec} amounts to
\begin{align}\label{eqn:alkh1}
 \alpha_{l,k}(h_\lambda)
&=
 \frac{(-1)^{l+[\frac{q-3}2]}\pi^\frac{p+q-4}2}
 {2^{\lambda}}
 \frac{\Gamma(\frac{l+k-\lambda}2)
       \Gamma(\frac{-q+3+l-k-\lambda}2)}
      {\Gamma(\frac{\lambda+p+q+l+k-4}2)
       \Gamma(\frac{\lambda+p-1+l-k}2)}
\\
&=
 \frac{(-1)^{k+[\frac{p-3}{2}]} \pi^{\frac{p+q-4}{2}}}{2^\lambda}
 \frac{\Gamma(\frac{l+k-\lambda}{2}) \Gamma(\frac{-p+3-l+k-\lambda}{2})}
      {\Gamma(\frac{p+q-4+l+k+\lambda}{2}) \Gamma(\frac{q-1-l+k+\lambda}{2})}.
\label{eqn:alkh2}      
\end{align}
\end{proposition}

\begin{proof}
The second equation \eqref{eqn:alkh2} follows from the identity
\eqref{eqn:Gammalk} of gamma functions.
Let us show the first equation \eqref{eqn:alkh1}.
In terms of $h_\lambda^\pm$ defined in \eqref{def:Ri},
we rewrite
$h_\lambda$ (see \eqref{def:hpm}) as 
\begin{equation}\label{eqn:hhpm}
 h_\lambda= 
 \frac{\pi}{\Gamma(\lambda+\frac{p+q-4}2)\sin(-\lambda \pi)}
 \times
\begin{cases}
   h_\lambda^+
  &\text{if $p,q$ both even}, \\[1ex]
   \frac{h_\lambda^+}{\tan (\lambda\pi)}+
   \frac{h_\lambda^-}{\sin (\lambda \pi)}
  &\text{if $p,q$ both odd}.
\end{cases}
\end{equation}

Since $\alpha_{l,k}$ is linear, i.e., 
$\alpha_{l,k}(ah+bg)=a\alpha_{l,k}(h)+b\alpha_{l,k}(g)$, \thinspace
$a, b \in \mathbb C$, by \eqref{eqn:hhpm}, we have
\begin{equation*}
 \alpha_{l,k}(h_\lambda)=
  C_{l,k}(\lambda) \gamma_{l,k}(\lambda),
\end{equation*}
where $\gamma_{l,k}(\lambda)$ is the meromorphic function given by
 \eqref{eqn:gammalk}, and
\begin{equation*}
C_{l,k}(\lambda) :=
\frac{\sin \frac{\lambda-l-k+2}2 \pi
 \sin \frac{\lambda+q-1-l+k}2 \pi}{\sin(-\lambda \pi)}
 \times
 \begin{cases}
   1 &\text{if $p,q$ both even},
 \\
   \frac{1}{\tan (\lambda\pi)}+\frac{(-1)^{l+k}}{\sin(\lambda\pi)}
     &\text{if $p,q$ both odd}.
 \end{cases}
\end{equation*}
Hence, the proof of Proposition will be completed by the following claim:
\end{proof}

\begin{claim} \label{clm:5.3.4}
$$
 C_{l,k}(\lambda)=\frac{(-1)^{l+[\frac{q-1}2]}}2.
$$
\end{claim}
\begin{proof}
Let us first consider the case where both $p$ and $q$ are even.
Then, the two integers 
$-l-k+2$ and $q-1-l+k$ have different parities.
Hence,
$$
 C_{l,k}(\lambda)=(-1)^{l+\frac{q-2}2}
 \frac{\sin \frac{\lambda}2 \pi \cos \frac{\lambda}2 \pi}
      {\sin \lambda \pi}=
 \frac{(-1)^{l+\frac{q-2}2}}2.
$$

Next, suppose both $p$ and $q$ are odd.
Then,
\begin{align*}
& \frac{1}{\tan (\lambda\pi)}+\frac{(-1)^{l+k}}{\sin(\lambda\pi)} =
  \begin{cases}
    \frac{1}{\tan\frac{\lambda}2 \pi},  \\
    -\tan \frac{\lambda}2 \pi, 
  \end{cases}
\\
& \frac{\sin \frac{\lambda-l-k+2}2 \pi \sin \frac{\lambda+q-1-l+k}2 \pi}
           {\sin (-\lambda)\pi}
  =\begin{cases}
   (-1)^{-l+\frac{q-1}2} \frac{\sin^2 \frac{\lambda}2 \pi}{\sin \lambda\pi}
  = \frac{(-1)^{-l+\frac{q-1}{2}}}{2}\tan\frac{\lambda}{2}\pi,  \\
   (-1)^{-l+\frac{q+1}2} \frac{\cos^2 \frac{\lambda}2 \pi}{\sin \lambda\pi}
  = \frac{(-1)^{-l+\frac{q+1}{2}}}{2}\frac{1}{\tan\frac{\lambda}{2}\pi},
  \end{cases}
\end{align*}
according as $l+k$ is even (upper row) and odd (lower row).
Thus we have
$$
C_{l,k}(\lambda)
  = \frac{(-1)^{\frac{q-1}2-l}}2
$$
in either case.
Hence, Claim \ref{clm:5.3.4} is verified.
\end{proof}

Let $T$ be the triangular domain in $\mathbb R^2$ given by
$$T:=\set{(x,y) \in \mathbb R^2}
{\ x<1, \  y<1, \  0<x+y},$$
and define a function $g_\lambda(x,y)$ with 
parameter $\lambda \in \mathbb C$ by 
\begin{equation}
 g_\lambda(x,y):=
 \begin{cases}
  (x+y)^\lambda (1-x^2)^\frac{p-4}2 (1-y^2)^\frac{q-4}2 
  &(x,y) \in T, \\
  0  &(x,y) \notin T.
 \end{cases}
\end{equation}

\begin{lemSec}\label{lem:mero}
{\rm 1)} For $\operatorname{Re}\lambda>-1$, 
$g_\lambda$ is a distribution of compact support,
and with holomorphic parameter $\lambda$.
That is, $\langle g_\lambda,\varphi\rangle$ is holomorphic in
$\{\lambda\in\mathbb{C}: \operatorname{Re}\lambda>-1\}$
 for any
$\varphi\in C^\infty(\mathbb{R}^2)$.

{\rm 2)} $g_\lambda$  extends as a distribution with a meromorphic parameter
$\lambda \in \mathbb C$.
That is, $\langle g_\lambda,\varphi\rangle$ is a meromorphic function
with respect to $\lambda\in\mathbb{C}$ for any $\varphi\in
C^\infty(\mathbb{R}^2)$.
\end{lemSec}  

\begin{proof}
The first statement is clear because
$g_\lambda \in L^1(\mathbb{R}^2)$ if $\operatorname{Re}\lambda>-1$. 
For the second statement,
we rewrite $g_\lambda$ as
$$
g_\lambda(x,y)
= (x+y)_+^\lambda (1-x)_+^{\frac{p-4}{2}} (1+x)_+^{\frac{p-4}{2}}
  (1-y)_+^{\frac{q-4}{2}} (1+y)_+^{\frac{q-4}{2}}.
$$
Then, Lemma follows from Bernstein's theorem \cite{xBG}.
\end{proof}

\section{Proof of Theorem \ref{thm:A}}
\label{subsec:prA}

We recall from \eqref{int:A} that $\mathcal{F}_C$ is the linear map
 defined by
$$
u(\xzeta) \mapsto 
\int_C K(\xzeta, \xzeta')u(\xzeta')d\mu(\xzeta').
$$
Since $\mathcal{F}_C$ commutes with the 
$K'$-action ($
\index{A}{K1Opqdash@$K' \simeq O(p-1)\times O(q-1)$}%
K'\simeq O(p-1)\times O(q-1)$), 
$\mathcal{F}_C$ preserves each $K'$-isotypic component of 
$L^2(C)$ given in the decomposition (see \eqref{eqn:iso}):
\begin{equation*}
   L^2(C) \simeq
   \sideset{}{^\oplus}\sum_{l,k=0}^\infty
   L^2(\mathbb{R}_+,\frac{1}{2}r^{p+q-5}dr)
   \otimes \mathcal{H}^l(\mathbb{R}^{p-1})
   \otimes \mathcal{H}^k(\mathbb{R}^{q-1}).
\end{equation*}
On the other hand, 
we have seen in Theorem \ref{thm:C} that
$\pi(w_0)$ also preserves 
each $K'$-isotypic component, and accordingly
has a decomposition:
$$
 \pi(w_0)=\sideset{}{^\oplus}\sum_{l,k=0}^\infty
 T_{l,k} \otimes \operatorname{id} \otimes \operatorname{id},
$$
where 
\index{A}{Tlk@$T_{l,k}$}%
$T_{l,k}$ is a unitary operator 
on $L^2(\mathbb R_+, \frac{1}{2}r^{p+q-5}dr)$ whose 
kernel 
\index{A}{Klkt@$K_{l,k}(t)$}%
$K_{l,k}(t)$ is explicitly given in 
\eqref{def:Klk}.

We shall show the equality $\mathcal{F}_C=\pi(w_0)$ by restricting 
to each $(l,k)$ component, namely,  
\begin{lemSec}\label{lem:res}
For each $l,k \in \mathbb N$, we have
\begin{equation} \label{eqn:res}
 \mathcal{F}_C\vert_{L^2(\mathbb R_+, r^{p+q-5}dr)\otimes 
 \Har{l}{p-1}\otimes \Har{k}{q-1}} =
 T_{l,k}\otimes \operatorname{id}\otimes 
 \operatorname{id}.
\end{equation}
\end{lemSec}
Instead of proving Lemma \ref{lem:res}, 
we shall prove
Lemma \ref{lem:a}
on the spectra $\alpha_{l,k}$ and the
kernel functions $K_{l,k}$, 
which turns out to be equivalent to Lemma \ref{lem:res}. 
For that purpose, we set
\begin{equation}\label{def:h}
 h_{r,r'}(x,y):=c_{p,q}\Phi_{p,q}(rr'(x+y)),
\end{equation}
where
$c_{p,q}$ and $\Phi_{p,q}$ are defined in \eqref{def:c} and \eqref{def:Phi}.
Then, by the definition
\eqref{def:K} of $K(\xzeta,\xzeta')$, we have
\begin{align*}
   K(\begin{pmatrix}r\omega \\ r\eta \end{pmatrix},
   \begin{pmatrix}r'\omega' \\ r' \eta'\end{pmatrix}) 
 & = c_{p,q}\Phi_{p,q}\big(rr'(\langle w,w'\rangle + 
                               \langle \eta, \eta' \rangle)\big)\\
 & = h_{r,r'} (\langle \omega, \omega' \rangle, 
               \langle \eta, \eta' \rangle).
\end{align*}
Suppose 
$f(r)u(\omega,\eta) \in L^2(\mathbb {R}_+, r^{p+q-5}dr)
 \otimes \mathcal{H}^l(\mathbb{R}^{p-1})
 \otimes \mathcal{H}^k(\mathbb{R}^{q-1})$.

\begin{align*}
& \mathcal{F}_C(fu)(r\omega,r\eta) \\
&\  = \int_C K \bigg(
\begin{pmatrix}
  r\omega \\ r\eta
\end{pmatrix}, \xzeta'\bigg)
  (fu)(\xzeta')d\mu(\xzeta') \\
&\ = \frac{1}{2} \int^\infty_0 \int_{S^{p-2}} \int_{S^{q-2}}
     h_{r,r'} (\langle \omega, \omega' \rangle, \langle\eta,\eta'\rangle) 
     f(r')u(\omega',\eta')r^{\prime p+q-5}
        dr'd\omega' d\eta' \\
&\ = \frac{1}{2} \int^\infty_0 (B_{h_{r,r'}}u)
        f(r')r^{\prime p+q-5} dr' .
\end{align*}
Since $B_{h_{r,r'}}u = \alpha_{l,k}(h_{r,r'})u$ 
by Lemma \ref{lem:spec}, 
we have 
$$
 \mathcal{F}_C(fu)(r\omega,r\eta) 
= \frac{1}{2} \int^\infty_0 \alpha_{l,k}(h_{r,r'})f(r')r^{\prime p+q-5} dr'
  u(\omega,\eta).
$$

On the other hand,
it follows from Theorem \ref{thm:C} (2) that
\begin{align*}
&  \big((T_{l,k} \otimes \operatorname{id} \otimes \operatorname{id})(fu)\big)
  (r\omega, r\eta) \\
&\   = (T_{l,k}f)(r) u(\omega,\eta) \\
&\   = \frac{1}{2}\int^\infty_0 K_{l,k}(rr')f(r')r^{\prime p+q-5}
        dr'u(\omega,\eta).
\end{align*}
Since 
$L^2(\mathbb{R}_+,r^{p+q-5}dr)\otimes
 \mathcal{H}^l(\mathbb{R}^{p-1})\otimes
 \mathcal{H}^k(\mathbb{R}^{q-1})
$
is spanned by a linear combination of the function of the form
$f(r)u(\omega,\eta)$, 
Lemma \ref{lem:res} is equivalent to
the following formula between kernel functions:
\begin{lemSec}\label{lem:a}
For each $l,k\in\mathbb{N}$, we have
$$
 \alpha_{l,k}(h_{r,r'})=K_{l,k}(rr').
$$
\end{lemSec}

The proof of Lemma \ref{lem:a} will be given in
the following two Sections,
which will then complete the proof of Theorem \ref{thm:A}.

\section{Proof of Lemma \ref{lem:a} 
(Hermitian case $q=2$)}
\label{subsec:5.5}

This section gives a proof of Lemma \ref{lem:a} in the case 
$\min(p,q)=2$. 
Without loss of generality,
we may and do assume $q=2$.
By the definition \eqref{def:Phi} of $\Phi_{p,2}(t)$,
the definition \eqref{lem:a} of $h_{r,r'}$ amounts to:
$$
h_{r,r'}(x,y)=c_{p,2}\Phi_\frac{p-4}2^+(rr'(x+y)).
$$

Since $\Phi_\frac{p-4}2^+(t)$ is a locally integrable function
supported on $t\ge0$ (see  Theorem \ref{eqn:MB2}),
we have from the definition 
\eqref{def:h2} of the operator $U_k$ $(k=0,1)$:
for $-1\le x \le 1$, 
\begin{align*}
 (U_k h_{r,r'})(x)={}&c_{p,2}\Phi_\frac{p-4}2^+(rr'(x+1))
  +(-1)^k c_{p,2}\Phi_\frac{p-4}2^+ (rr'(x-1))\\
 ={}& c_{p,2}\Phi_\frac{p-4}2^+ (rr'(x+1)).\\
\intertext{Then, by \eqref{def:Psi0}, we have:} 
 ={}&\frac{2^{-\frac{p-8}{4}}(-1)^{\frac{(p-1)(p+2)}{2}}}{\pi^{\frac{p-2}{2}}}(rr')^{-\frac{p-4}4}
  (x+1)^{-\frac{p-4}4}J_\frac{p-4}2 (2\sqrt{2rr'(x+1)}).
\end{align*}
Applying the 
\index{B}{Funk--Hecke formula}%
Funk--Hecke formula \eqref{eqn:a2}, 
the spectrum $\alpha_{l,k}(h_{r,r'})$ amounts to
\begin{align*}
\alpha_{l,k}(h_{r,r'})
&=
 \frac{2^{\frac{3p-4}{4}}(-1)^{\frac{(p-1)(p+2)}{2}}l!}
      {\sqrt{\pi}\Gamma(p-3+l)}
 (rr')^{-\frac{p-4}4}
  \times{}
\\ 
 &\quad \int_{-1}^1 J_{\frac{p-4}2}(2\sqrt{2rr'(x+1)})
 \tilC_l^\frac{p-3}2(x)(1+x)^{\frac{p-4}{4}}(1-x)^\frac{p-4}2 dx \\
 &=4(-1)^{\frac{(p-1)(p+2)}{2}+l}(rr')^{-\frac{p-3}2}J_{p-3+2l}(4\sqrt{rr'}) \\
 &=(-1)^{\frac{p^2}{2}} K_{l,k}(rr').
\end{align*}
Here, the second equality follows from the integral formula
\eqref{eqn:Ge3} of the Gegenbauer polynomials in Appendix
with $\alpha=2\sqrt{2rr'}$ and $\nu=\frac{p-4}{2}$, 
and the last equality follows from 
\eqref{eqn:rdB}.
Since $p$ is even in the case $q=2$,
the right-hand side is equal to $K_{l,k}(rr')$.
Hence, Lemma \ref{lem:a} is proved for $q=2$.
\qed

\section{Proof of Lemma \ref{lem:a} ($p,q>2$)}
\label{subsec:5.6}

This section gives a proof of Lemma \ref{lem:a} in the general case
$p,q>2$. 

First, we give an integral formula of Mellin--Barnes type for    
$h_{r,r'}(x,y)$ (see \eqref{def:h}) 
by means of 
$h_\lambda(x,y)$ (see \eqref{def:hpm} for definition):
Suppose $p,q > 2$. 

\begin{claim}\label{claim:hxy}
  Let $\gamma > -1$ and $L$ be a contour that starts
at $\gamma - \sqrt{-1}\infty$ and ends at
   $\gamma + \sqrt{-1}\infty$.
We assume that $L$ passes the real axis
in the interval $(-\frac{p+q-4}{2},-\frac{p+q-6}{2})$.
Then, we have
\begin{equation}
  \label{eqn:}
 h_{r,r'}(x,y) 
 = \frac{c_{p,q}}{2\pi\sqrt{-1}}\int_L 
   (2rr')^\lambda h_\lambda(x,y) d\lambda.  
\end{equation}
\end{claim}
\begin{proof}
By the definition \eqref{def:hpm}
of $h_\lambda(x,y) \equiv h^{p,q}_\lambda(x,y)$ and the integral formulas
\eqref{eqn:MBPsi+} and \eqref{eqn:MBPsi} of
$\Psi^+_{\frac{p+q-6}{2}}(t)$ and
$\Psi_{\frac{p+q-6}{2}}(t)$ respectively, we have
$$
 \frac{1}{2\pi\sqrt{-1}} \int_L s^\lambda h_\lambda (x,y)d\lambda
 = 
 \begin{cases}
   \Psi^+_{\frac{p+q-6}{2}}\bigl(\frac{s(x+y)}{2}\bigr) & p,q \text{ even,}\\
   \Psi_{\frac{p+q-6}{2}}\bigl(\frac{s(x+y)}{2}\bigr)   & p,q \text{ odd,}\\
 \end{cases}
$$
for $s > 0$.
In either case,
it follows from
the definition\eqref{def:Phi} of 
$\Phi_{p,q}(t)$ that
\begin{equation}\label{eqn:hxyint}
  \frac{1}{2\pi\sqrt{-1}} \int_L s^\lambda h_\lambda (x,y)d\lambda
 = \Phi_{p,q}\Bigl(\frac{s(x+y)}{2}\Bigr). 
\end{equation}
Hence, we get Claim \ref{claim:hxy} by the definition
\eqref{def:h} of $h_{r,r'}(x,y)$.
\end{proof}

By the linearity of $\alpha_{l,k}$,
we have
\begin{align*}
 \alpha_{l,k}(h_{r,r'})={}&
 \frac{c_{p,q}}{2\pi\sqrt{-1}}
 \int_L  \alpha_{l,k}(h_\lambda)
  (2rr')^\lambda d\lambda  \\
 ={}&\frac{(-1)^{l+\frac{p-q}2}}{\pi \sqrt{-1}}
 \int_L  
      \frac{\Gamma(\frac{l+k-\lambda}2) \Gamma(\frac{-q+3+l-k-\lambda}2)}
           {\Gamma(\frac{\lambda+p+q+l+k-4}2)\Gamma(\frac{\lambda+p-1+l-k}2)} 
             (rr')^\lambda d\lambda  \\
 ={}& K_{l,k}(rr').
\end{align*}
Here
, in the second equality, we applied Proposition
\ref{prop:R} and then used the equality 
$(-1)^{l+[\frac{q-1}2]}(-1)^{\frac{(p-1)(p+2)}2}=(-1)^{l+\frac{p-q}2}$,
which follows from
the congruence equality:
\begin{equation*}
   \frac{(p-1)(p+2)}{2} + \left[ \frac{q-1}{2} \right]
   \equiv \frac{p-q}{2}
   \bmod 2
\end{equation*}
under the assumption that $p+q$ is even.
The last equality follows from Lemma \ref{lem:Klkint}.
Hence, we have proved Lemma \ref{lem:a} in the general case
$p,q>2$.
\qed

\begin{proof}[Proof of Theorem \ref{thm:A}]
Now, Lemma \ref{lem:a} is proved in all the cases.
Hence, the proof of Theorem \ref{thm:A} is completed.
\end{proof}

\chapter{Bessel distributions}
\label{sec:Diffeq}

We have seen in the previous chapter (see Theorem \ref{thm:A}) that the 
unitary inversion operator 
$\pi(w_0): L^2(C)\to L^2(C)$ is given by 
the distribution kernel $K(\xzeta,\xzeta')$
which is the composition of the restriction of the bilinear map
\begin{equation*}
   C \times C \to \mathbb{R},
   \qquad
   (\xzeta,\xzeta') \mapsto \langle \xzeta,\xzeta' \rangle
\end{equation*}
and Bessel distributions (see \eqref{def:Psi0}--\eqref{def:Psi}) of
one variable.
In this chapter,
we analyze the distribution kernel from three viewpoints:
integral formulas, power series expansion (including distributions
such as $\delta^{(k)}(x)$ and $x^{-k}$), and differential equations.

Section \ref{subsec:ode} 
 gives a heuristic account on why $K(\xzeta,\xzeta')$ is essentially
of one variable, and why
the Bessel distribution arises in 
$K(\xzeta,\xzeta')$.
The results of Section \ref{subsec:ode} is not used for
 other sections.

\section{Meijer's $G$-distributions} \label{subsec:MGdistr}

In this section, we give a definition of 
\index{B}{Meijer's $G$-distribution|main}%
\textit{Meijer's $G$-distributions} 
which have the following two
properties:
\begin{itemize}
\item[1)]
They are distributions on $\mathbb{R}$.
\item[2)]
The restrictions to the positive half line $\{x>0\}$ are (usual)
Meijer's $G$-functions 
(see Appendix \ref{subsec:G}).
\end{itemize}
The main result of this section is Proposition \ref{prop:Gdistr}.

Let $m$, $n$, $p$ and $q$ be integers with 
$0\le m \le q$ and $0 \le n \le p$.
Suppose moreover that the complex numbers 
$a_1,\dots,a_p$ and $b_1,\dots,b_q$ fulfill the condition
\begin{equation*}
   a_j-b_k \ne 1,2,3,\dotsc
   \quad (j=1,\dots,n; \, k=1,\dots,m).
\end{equation*}
This means that
no pole of the gamma function
$\Gamma(b_j-\lambda)$ $(j=1,\dots,m)$ coincides with any
pole of $\Gamma(1-a_k+\lambda)$ $(k=1,\dots,n)$.
We set
\begin{align}
   & c^* := m+n-\frac{p+q}{2},
\label{eqn:cstar}
\\
   & \mu := \sum_{j=1}^q b_j - \sum_{j=1}^p a_j
            + \frac{p-q}{2} + 1.
\label{eqn:muab}
\end{align}
Throughout this chapter, we assume $c^* \ge 0$.
If $c^* = 0$ then we also assume
\begin{equation} \label{eqn:pqmu}
   p-q < 0  \quad\text{or}\quad p-q > \operatorname{Re}\mu.
\end{equation}
It is easy to see that 
the condition \eqref{eqn:pqmu} allows us to find $\gamma\in\mathbb{R}$
such that
\begin{equation}\label{eqn:gamineq}
   \gamma > -1
   \quad\text{and}\quad
   (q-p) \gamma > \operatorname{Re}\mu.
\end{equation}

\begin{remNonumber} \label{rem:Gdistr1}
The conditions \eqref{eqn:pqmu} and $\gamma>-1$ will not be used when we
define (usual) Meijer's $G$-function as an analytic function in $x>0$ 
(see \eqref{def:G}).
They will be used in showing that Meijer's $G$-distribution
$G(x_+)_L$ given by
the Mellin--Barnes type integral
\eqref{eqn:Gdistr} is a locally integrable function on $\mathbb{R}$ 
if we take an
appropriate contour $L$
(see Proposition \ref{prop:Gdistr} (3)).
\end{remNonumber}

We now take a contour $L$ which starts at the point
$\gamma-\sqrt{-1}\infty$ and finishes at $\gamma+\sqrt{-1}\infty$.
Later, we shall impose the following conditions on $L$:
\begin{texteqn} \label{eqn:Lneg}
   $L$ does not go through any negative integer.
\end{texteqn}
\begin{texteqn} \label{eqn:Lab}
   $L$ leaves $b_j$ $(1\le j\le m)$ to the right,
   and $a_j-1$ $(1 \le j \le n)$ to the left.
\end{texteqn}
We note that the condition \eqref{eqn:Lab} implies:
\begin{equation} \begin{minipage}[b]{0.8\linewidth}
   $L$ does not go through any point in 
 \newline
   $\{ b_j+k: 1\le j\le m, 
    \ 
    k\in\mathbb{N} \}
   \cup \{a_j-1-k: 1\le j\le n, \  k\in\mathbb{N} \}$.
\end{minipage} 
\tag*{(\ref{eqn:Lab})$'$}
\end{equation} \ignorespacesafterend
With these parameters, we define a meromorphic function of $\lambda$ by
\begin{equation} \label{eqn:Ggam}
\index{A}{1Gamma@$\Gamma_{p,q}^{m,n} 
     ( \lambda "| \genfrac{}{}{0pt}{}
                              {a_1,\dots,a_p}
                              {b_1,\dots,b_q}
     )$|main}%
    \Gamma_{p,q}^{m,n} 
     \Bigl( \lambda \Bigm| \genfrac{}{}{0pt}{}
                              {a_1,\dots,a_p}
                              {b_1,\dots,b_q}
     \Bigr)
   :=
     \frac{\prod\limits_{j=1}^m \Gamma(b_j-\lambda) 
           \prod\limits_{j=1}^n \Gamma(1-a_j+\lambda)}
          {\prod\limits_{j=m+1}^q \Gamma(1-b_j+\lambda) 
           \prod\limits_{j=n+1}^p \Gamma(a_j-\lambda)}.
\end{equation}
For $\operatorname{Re}\lambda>-1$,
we set
\begin{equation*}
   x_+^\lambda
   := \begin{cases} x^\lambda &(x>0)   \\
                    0         &(x\le0),
      \end{cases}
\qquad
   x_-^\lambda
   := \begin{cases} 0         &(x\ge0) \\
                    |x|^\lambda   &(x<0).
      \end{cases}
\end{equation*}
Then, $x_+^\lambda$ and $x_-^\lambda$ are locally integrable functions
of the variable $x$ in $\mathbb{R}$,
and extend to distributions with meromorphic parameter $\lambda$ in
the entire complex plane
(see Appendix \ref{subsec:Riesz}).

\begin{proposition} \label{prop:Gdistr}
Let $L$ be a contour satisfying \eqref{eqn:Lneg} and\/
{\upshape(\ref{eqn:Lab})}$'$. 

{\upshape 1)}  
The 
\index{B}{Mellin--Barnes type integral}%
Mellin--Barnes type integral:
\begin{align}
\index{A}{GL@$G(x_+)_L$|main}%
   G(x_+)_L
   & \equiv 
\index{A}{Gpq@$G_{p,q}^{m,n}
     ( x_+ "| \genfrac{}{}{0pt}{}
                           {a_1,\dots,a_p}
                           {b_1,\dots,b_q}
     )_L$|main}%
G_{p,q}^{m,n}
     \Bigl( x_+ \Bigm| \genfrac{}{}{0pt}{}
                           {a_1,\dots,a_p}
                           {b_1,\dots,b_q}
     \Bigr)_L
\nonumber
\\  
   & := \frac{1}{2\pi\sqrt{-1}}
     \int_L \Gamma_{p,q}^{m,n}
     \Bigl( \lambda \Bigm| \genfrac{}{}{0pt}{}
                           {a_1,\dots,a_p}
                           {b_1,\dots,b_q}
     \Bigr)
     x_+^\lambda d\lambda
\label{eqn:Gdistr}
\end{align}
is well-defined as a distribution on $\mathbb{R}$.

Its support is given by
\begin{equation*}
  \operatorname{supp} G(x_+)_L = \{ x\in \mathbb{R}: x \ge 0 \}.
\end{equation*}

{\upshape 2)}
If the contour $L$ satisfies \eqref{eqn:Lab},
then the restriction of $G(x_+)_L$ to the positive half line
$\{ x\in \mathbb{R}: x>0\}$ is a real analytic function,
which coincides with the (usual) $G$-function
$G_{p,q}^{m,n}\Bigl( x\Bigm| \genfrac{}{}{0pt}{0}{a_1,\dots,a_p}
 {b_1,\dots,b_q} \Bigr)$
(see \eqref{def:G} for definition).

{\upshape 3)}
If the contour $L$ is contained in the half plane
$\{\lambda\in\mathbb{C}: \operatorname{Re}\lambda>-1 \}$,
then $G(x_+)_L$ is a locally integrable function on $\mathbb{R}$.
More precisely, there exists $\epsilon_0>0$ such that
$G(x_+)_L \, x_+^{-\epsilon}$ is locally integrable for any $\epsilon$
with $0 \le \epsilon < \epsilon_0$.
\end{proposition}

Likewise, we can define the distribution
\begin{align*}
   G(x_-)_L
   \equiv {}&  G_{p,q}^{m,n}
     \Bigl( x_- \Bigm| \genfrac{}{}{0pt}{}
                           {a_1,\dots,a_p}
                           {b_1,\dots,b_q}
     \Bigr)_L
\\
   :={} & \frac{1}{2\pi\sqrt{-1}} \int_L \Gamma_{p,q}^{m,n}
     \Bigl( \lambda \Bigm| \genfrac{}{}{0pt}{}
                           {a_1,\dots,a_p}
                           {b_1,\dots,b_q}
     \Bigr)
     x_-^\lambda d\lambda
\end{align*}
by using the same contour $L$,
and the support of $G(x_-)_L$ is equal to the negative half line
$\{x\in\mathbb{R}: x\le 0\}$.

\begin{remNonumber} \label{rem:Gdistr}
The distribution $G(x_\pm)_L$ depends on the choice of the contour $L$
even when we assume $L$ satisfies the conditions \eqref{eqn:Lneg} and
\eqref{eqn:Lab}. 
In fact, if $L$ and $L'$ are contours satisfying \eqref{eqn:Lneg} and
\eqref{eqn:Lab}, 
then $G(x_\pm)_L$ may differ from $G(x_\pm)_{L'}$ by
 a distribution supported at $0$, namely, a finite sum of
Dirac's delta function and its derivatives.
This is because the distribution $x_\pm^\lambda$ has
 simple poles at $\lambda=-1,-2,\dotsc$, and consequently, its
residues (see \eqref{eqn:resxlmd} and \eqref{eqn:resxlmd2}) may appear
when we move the contour $L$ across negative integers.
In order to define the $G$-distribution in a unique fashion, 
we need to impose an additional constraint on the contour $L$.
We shall work with concrete examples for this in Section
\ref{subsec:intPsi} where we use Cauchy's integral formula for distributions
with meromorphic parameter.
\end{remNonumber}

In order to prove Proposition \ref{prop:Gdistr},
we need an 
\index{B}{asymptotic behavior!---, $\Gamma_{p,q}^{m,n}\Bigl(\lambda\Bigm"|\genfrac{}{}{0pt}{}{a_1,\dots,a_p}{b_1,\dots,b_q}\Bigr)$}%
asymptotic estimate of the $\Gamma$-factors in the
integrand of \eqref{eqn:Gdistr} as follows:

\begin{lemSec} \label{lem:asymgamma}
For any $\epsilon>0$, there exists a constant $C>0$ such that
$$
\left|
  \Gamma_{p,q}^{m,n} 
     \Bigl( \lambda \Bigm| \genfrac{}{}{0pt}{}
                           {a_1,\dots,a_p}
                           {b_1,\dots,b_q}
     \Bigr)
\right|
\le C e^{-\pi c^*\left|\operatorname{Im}\lambda\right|}
    \left|\operatorname{Im} \lambda\right|
         ^{\operatorname{Re}\mu+(p-q)\gamma-1+\epsilon} 
$$
for any $\lambda\in L$ such that
$\left|\operatorname{Im}\lambda\right|$ is sufficiently large.
Here, $c^*$ and $\mu$ are defined as in \eqref{eqn:cstar} and
\eqref{eqn:muab}, 
and $\gamma = \lim\limits_{\substack{\lambda\in L\\
                  \left|\operatorname{Im}\lambda\right|\to\infty}}
     \operatorname{Re}\lambda$.
\end{lemSec}

\begin{proof}
Fix $a\in\mathbb{C}$.
By Stirling's asymptotic formula \eqref{eqn:Stirling} of the gamma
function, we have
\begin{align*}
   |\Gamma(a-\lambda)|
   & = C_a \left|\operatorname{Im}\lambda\right|
         ^{\operatorname{Re}a-\operatorname{Re}\lambda-\frac{1}{2}}
       e^{-\frac{\pi}{2}\left|\operatorname{Im}\lambda\right|}
       (1+O(\left|\operatorname{Im}\lambda\right|^{-1})),
\\
   |\Gamma(1-a+\lambda)|
   & = C_a \left|\operatorname{Im}\lambda\right|
         ^{-\operatorname{Re}a+\operatorname{Re}\lambda+\frac{1}{2}}
       e^{-\frac{\pi}{2}\left|\operatorname{Im}\lambda\right|}
       (1+O(\left|\operatorname{Im}\lambda\right|^{-1})),
\end{align*}
as $\left|\operatorname{Im}\lambda\right|$ tends to infinity with
$\operatorname{Re}\lambda$ bounded.
Here, the constant $C_a$ is given by
\begin{equation*}
   C_a = \sqrt{2\pi}
       \,  e^{-\frac{\pi}{2}\sgn(\operatorname{Im}\lambda)
              \left|\operatorname{Im} a\right|}.
\end{equation*}

By the definition \eqref{eqn:Ggam} of
$\Gamma_{p,q}^{m,n}\Bigl( \lambda \Bigm| \genfrac{}{}{0pt}{0}{a_1,\dots,a_p}
 {b_1,\dots,b_q} \Bigr)$,
we now get the following asymptotic behavior:
\begin{equation*}
   \left|
\Gamma_{p,q}^{m,n}\Bigl( \lambda \Bigm| \genfrac{}{}{0pt}{0}{a_1,\dots,a_p}
      {b_1,\dots,b_q} \Bigr)
   \right|
   = C' \left|\operatorname{Im}\lambda\right|^s
     e^{-\frac{\pi}{2}t\left|\operatorname{Im}\lambda\right|}
     (1+O(\left|\operatorname{Im}\lambda\right|^{-1})),
\end{equation*}
as $\left|\operatorname{Im}\lambda\right|$ tends to infinity,
where $C'$ is a constant depending on $\operatorname{Im}a_j$ and 
$\operatorname{Im}b_j$, and 
\begin{align*}
   s
   &= \sum_{j=1}^m \operatorname{Re}(b_j-\lambda-\frac{1}{2})
    + \sum_{j=1}^n \operatorname{Re}(-a_j+\lambda+\frac{1}{2})
\\
   &\quad 
    - \sum_{j=m+1}^q \operatorname{Re}(\frac{1}{2}-b_j+\lambda)
    - \sum_{j=n+1}^p \operatorname{Re}(a_j-\lambda-\frac{1}{2})
\\
   &= \operatorname{Re} \mu + (p-q) \operatorname{Re}\lambda - 1,
\\[1ex]
   t
   &= m+n - (q-m) - (p-n)
\\
   &= 2c^*.
\end{align*}
As $\operatorname{Re}\lambda$ converges to $\gamma$ when
$\lambda\in L$ goes to infinity,
we get Lemma \ref{lem:asymgamma}.
\end{proof}

We are ready to give a proof of Proposition \ref{prop:Gdistr}.

\begin{proof}[Proof of Proposition \ref{prop:Gdistr}]
3) 
We begin with the proof of the third statement.
Suppose $L$ is contained in the half plane
$\{ \lambda \in \mathbb{C}:
 \operatorname{Re}\lambda > -1 \}$.
We need to show the integral \eqref{eqn:Gdistr} 
makes sense and gives rise to a
 locally integrable function of $x$.
For the convergence of the integral, 
we shall use Lemma \ref{lem:asymgamma} for the estimate as a function
of $\lambda$. 
The non-trivial part is a uniform estimate in the neighborhood of $x=0$.
Let us consider the interval $0<x\le1$.

Since the contour $L$ has the property: 
$$
\gamma 
= \lim\limits_{\substack{\lambda\in L\\
                         \left|\operatorname{Im}\lambda\right|\to\infty}}
  \operatorname{Re}\lambda > -1,
$$
the assumption 
$L\subset\{\lambda\in\mathbb{C}: \operatorname{Re}\lambda>-1\}$ 
implies $\delta>-1$,
where we set
\begin{equation*}
   \delta := \inf_{\lambda\in L}
   \operatorname{Re}\lambda.
\end{equation*}
Hence, we get
\begin{equation*}
   |x_+^\lambda| \le x^\delta
   \quad\text{for $0 < x \le 1$}.
\end{equation*}
On the other hand, it follows from Lemma \ref{lem:asymgamma} that
\begin{equation*}
   \left| \Gamma_{p,q}^{m,n}
          \Bigl( \lambda\Bigm| \genfrac{}{}{0pt}{0}{a_1,\dots,a_p}
                 {b_1,\dots,b_q} \Bigr)
   \right|
   \le
   \begin{cases}
       C e^{-\pi c^*\left|\operatorname{Im}\lambda\right|}
       &\text{if $c^*>0$},
     \\
       C \left|\operatorname{Im}\lambda\right|^{-1+\epsilon}
             &\text{if $c^*=0$},
   \end{cases}
\end{equation*}
when $\left|\operatorname{Im}\lambda\right|$ is sufficiently large.
Here, we used the inequality
$\operatorname{Re}\mu+(p-q)\gamma < 0$ 
(see \eqref{eqn:gamineq}) in the second case.
Hence, 
$\Gamma_{p,q}^{m,n}\Bigl( \lambda\Bigm| \genfrac{}{}{0pt}{0}{a_1,\dots,a_p}
                 {b_1,\dots,b_q} \Bigr)$
 is absolutely integrable on $L$ in either case.
Therefore, the integration \eqref{eqn:Gdistr} 
converges, giving rise to a function of $x$ which is 
bounded by a scalar
multiple of $x^\delta$ on the interval $0<x\le1$, whence a locally
integrable function of $x$.
Thus, $G(x_+)_L$ is locally integrable.
Similarly, if we set
\begin{equation*}
   \epsilon_0 := 1+\delta
   \quad (>0),
\end{equation*}
then for any $0 \le \epsilon < \epsilon_0$,
$x^{-\epsilon+\delta}$ is locally integrable,
and consequently $G(x_+)_L \, x_+^{-\epsilon}$ is locally integrable.
Hence, 
the third statement of Proposition is proved.

1)
We divide the integral \eqref{eqn:Gdistr} into the sum of the
following two integrals
\begin{equation*}
   \int_L = \int_{L'} + \int_C,
\end{equation*}
where $L'$ is a contour contained in the right half plane
$\{\lambda\in\mathbb{C}: \operatorname{Re}\lambda > -1 \}$,
and $C$ is the closed oriented curve given by $L-L'$
(see Figure \ref{fig:LLC}).

\begin{figure}[H]
\resizebox{12cm}{!}{\input{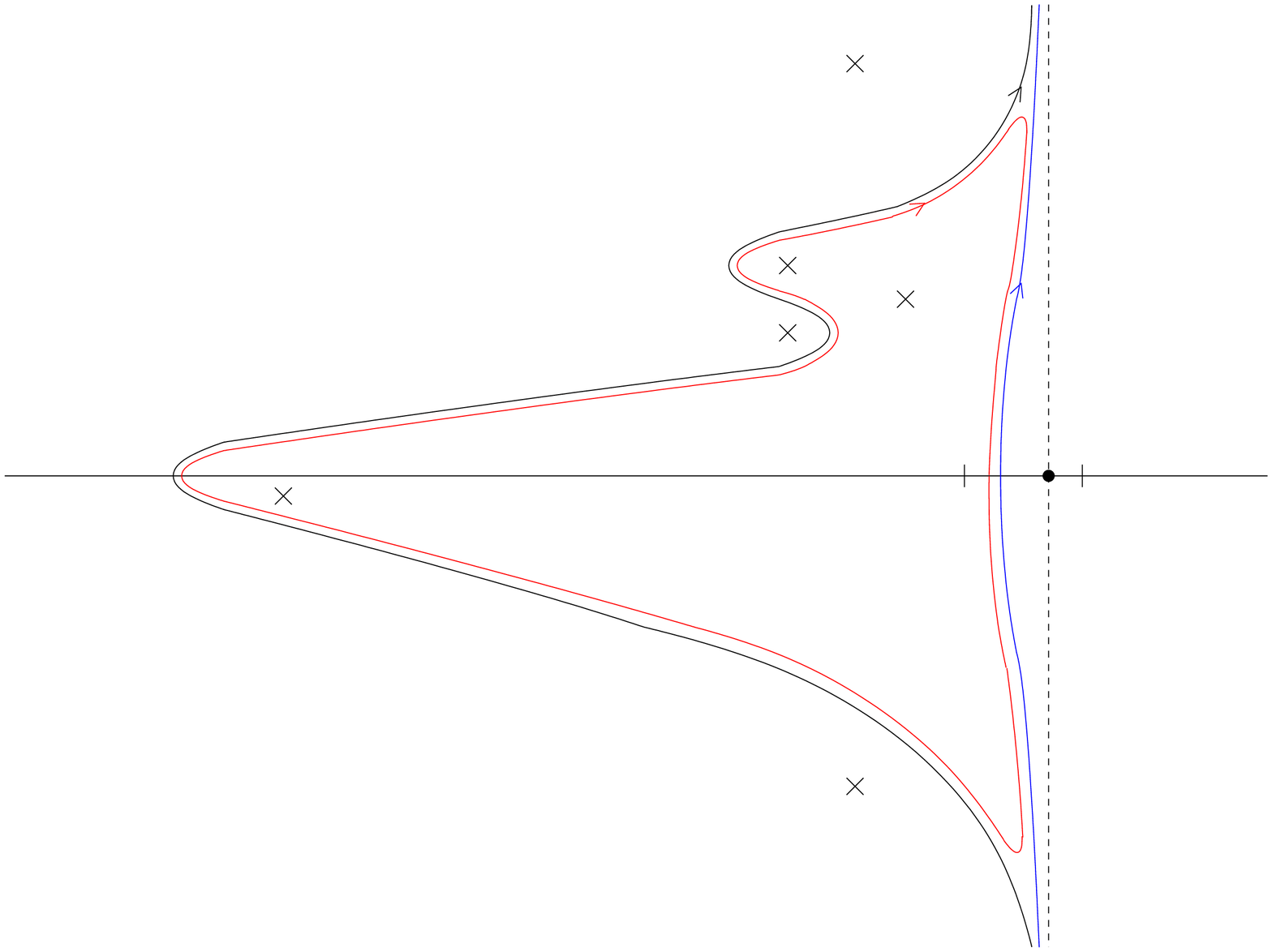}}
\caption{}
\label{fig:LLC}
\end{figure}

Then, 
we have already seen that 
the second term gives a locally integrable function of $x$ 
(the third statement of this proposition).
On the other hand,
the third term is well-defined as a distribution because $C$ is
compact and the integrand is a distribution of $x$ that depends
continuously on $\lambda$ as far as $\lambda$ lies in $C$.
Hence, the first statement is also proved.

2)
This statement is well-known.
See Appendix \ref{subsec:G} for details.
\end{proof}

\section{Integral expression of 
Bessel distributions}
\protect\index{B}{Bessel distribution!integral expression@---, integral expression}%
\label{subsec:intPsi}

In this section, we apply general results on Meijer's
$G$-distributions developed in Section \ref{subsec:MGdistr} to
special cases, and obtain the Mellin--Barnes type integral expression for 
 the distribution kernel of the unitary
inversion operator $\pi(w_0)$.

Let $m$ be a non-negative integer.
We take a contour $L$ such that
\begin{itemize}
\item[1)]
$L$ starts at $\gamma-\sqrt{-1}\infty$,
passes the real axis at some point $s$,
and ends at $\gamma+\sqrt{-1}\infty$.
\item[2)]
$-1<\gamma$ and $-m-1<s<-m$.
\end{itemize}
Likewise, we take a contour $L_0$ (with analogous notation) such that
\begin{itemize}
\item[3)]
$-1<\gamma_0$ and $-1<s_0<0$.
\end{itemize}
For later purpose, we may and do take $\gamma=\gamma_0$.
See Figure \ref{fig:LL}.
\begin{figure}[H]
\setlength{\unitlength}{0.00033333in}
\begingroup\makeatletter\ifx\SetFigFont\undefined%
\gdef\SetFigFont#1#2#3#4#5{%
  \reset@font\fontsize{#1}{#2pt}%
  \fontfamily{#3}\fontseries{#4}\fontshape{#5}%
  \selectfont}%
\fi\endgroup%
{\renewcommand{\dashlinestretch}{30}
\begin{picture}(11274,8439)(-200,-10)
\path(5323,5787)(5322,5787)(5321,5787)
	(5317,5786)(5312,5785)(5305,5783)
	(5294,5781)(5281,5779)(5265,5776)
	(5245,5772)(5222,5767)(5195,5762)
	(5164,5756)(5130,5749)(5092,5741)
	(5051,5733)(5006,5724)(4959,5714)
	(4908,5704)(4854,5693)(4798,5681)
	(4740,5669)(4679,5656)(4617,5643)
	(4553,5630)(4488,5615)(4421,5601)
	(4354,5586)(4285,5571)(4216,5555)
	(4146,5539)(4076,5522)(4005,5506)
	(3934,5488)(3862,5471)(3791,5453)
	(3718,5434)(3645,5415)(3572,5396)
	(3499,5376)(3425,5355)(3350,5334)
	(3275,5312)(3200,5290)(3124,5266)
	(3047,5242)(2971,5217)(2894,5192)
	(2817,5165)(2740,5138)(2663,5110)
	(2588,5082)(2514,5053)(2442,5024)
	(2344,4983)(2253,4942)(2168,4902)
	(2091,4863)(2020,4825)(1956,4789)
	(1899,4754)(1847,4721)(1801,4689)
	(1760,4658)(1724,4628)(1692,4599)
	(1663,4571)(1638,4544)(1616,4518)
	(1597,4492)(1580,4467)(1565,4443)
	(1552,4420)(1542,4397)(1532,4375)
	(1525,4354)(1518,4335)(1513,4316)
	(1509,4299)(1506,4283)(1503,4269)
	(1502,4256)(1500,4246)(1499,4236)
	(1499,4229)(1498,4223)(1498,4219)
	(1498,4215)(1498,4214)(1498,4212)
\path(9012,7962)(9012,7961)(9012,7960)
	(9012,7957)(9011,7953)(9011,7946)
	(9010,7938)(9009,7927)(9008,7913)
	(9006,7897)(9004,7878)(9001,7857)
	(8998,7834)(8994,7808)(8989,7779)
	(8984,7749)(8977,7717)(8969,7684)
	(8961,7649)(8951,7613)(8940,7575)
	(8927,7537)(8913,7498)(8897,7458)
	(8879,7418)(8859,7377)(8837,7335)
	(8812,7293)(8785,7250)(8754,7206)
	(8721,7161)(8684,7116)(8643,7069)
	(8597,7022)(8548,6973)(8493,6924)
	(8434,6873)(8369,6822)(8298,6770)
	(8223,6717)(8142,6664)(8056,6612)
	(7984,6571)(7910,6530)(7834,6491)
	(7758,6454)(7682,6417)(7605,6383)
	(7529,6350)(7453,6318)(7377,6288)
	(7302,6260)(7228,6233)(7154,6207)
	(7081,6182)(7008,6159)(6936,6137)
	(6864,6115)(6793,6095)(6722,6075)
	(6651,6056)(6581,6038)(6512,6021)
	(6443,6004)(6374,5988)(6306,5973)
	(6239,5958)(6173,5944)(6108,5930)
	(6044,5917)(5981,5905)(5921,5893)
	(5862,5882)(5805,5871)(5751,5861)
	(5699,5851)(5650,5843)(5605,5834)
	(5562,5827)(5523,5820)(5488,5814)
	(5456,5809)(5429,5804)(5404,5800)
	(5384,5797)(5367,5794)(5353,5792)
	(5342,5790)(5334,5789)(5329,5788)
	(5326,5787)(5324,5787)(5323,5787)
\path(5323,2637)(5322,2637)(5321,2637)
	(5317,2638)(5312,2639)(5305,2641)
	(5294,2643)(5281,2645)(5265,2648)
	(5245,2652)(5222,2657)(5195,2662)
	(5164,2668)(5130,2675)(5092,2683)
	(5051,2691)(5006,2700)(4959,2710)
	(4908,2720)(4854,2731)(4798,2743)
	(4740,2755)(4679,2768)(4617,2781)
	(4553,2794)(4488,2809)(4421,2823)
	(4354,2838)(4285,2853)(4216,2869)
	(4146,2885)(4076,2902)(4005,2918)
	(3934,2936)(3862,2953)(3791,2971)
	(3718,2990)(3645,3009)(3572,3028)
	(3499,3048)(3425,3069)(3350,3090)
	(3275,3112)(3200,3134)(3124,3158)
	(3047,3182)(2971,3207)(2894,3232)
	(2817,3259)(2740,3286)(2663,3314)
	(2588,3342)(2514,3371)(2442,3400)
	(2344,3441)(2253,3482)(2168,3522)
	(2091,3561)(2020,3599)(1956,3635)
	(1899,3670)(1847,3703)(1801,3735)
	(1760,3766)(1724,3796)(1692,3825)
	(1663,3853)(1638,3880)(1616,3906)
	(1597,3932)(1580,3957)(1565,3981)
	(1552,4004)(1542,4027)(1532,4049)
	(1525,4070)(1518,4089)(1513,4108)
	(1509,4125)(1506,4141)(1503,4155)
	(1502,4168)(1500,4178)(1499,4188)
	(1499,4195)(1498,4201)(1498,4205)
	(1498,4209)(1498,4210)(1498,4212)
\path(9012,462)(9012,463)(9012,464)
	(9012,467)(9011,471)(9011,478)
	(9010,486)(9009,497)(9008,511)
	(9006,527)(9004,546)(9001,567)
	(8998,590)(8994,616)(8989,645)
	(8984,675)(8977,707)(8969,740)
	(8961,775)(8951,811)(8940,849)
	(8927,887)(8913,926)(8897,966)
	(8879,1006)(8859,1047)(8837,1089)
	(8812,1131)(8785,1174)(8754,1218)
	(8721,1263)(8684,1308)(8643,1355)
	(8597,1402)(8548,1451)(8493,1500)
	(8434,1551)(8369,1602)(8298,1654)
	(8223,1707)(8142,1760)(8056,1812)
	(7984,1853)(7910,1894)(7834,1933)
	(7758,1970)(7682,2007)(7605,2041)
	(7529,2074)(7453,2106)(7377,2136)
	(7302,2164)(7228,2191)(7154,2217)
	(7081,2242)(7008,2265)(6936,2287)
	(6864,2309)(6793,2329)(6722,2349)
	(6651,2368)(6581,2386)(6512,2403)
	(6443,2420)(6374,2436)(6306,2451)
	(6239,2466)(6173,2480)(6108,2494)
	(6044,2507)(5981,2519)(5921,2531)
	(5862,2542)(5805,2553)(5751,2563)
	(5699,2573)(5650,2581)(5605,2590)
	(5562,2597)(5523,2604)(5488,2610)
	(5456,2615)(5429,2620)(5404,2624)
	(5384,2627)(5367,2630)(5353,2632)
	(5342,2634)(5334,2635)(5329,2636)
	(5326,2637)(5324,2637)(5323,2637)
\path(4999,5728)(4802,5779)
\path(4995,5723)(4854,5567)
\path(9598,5787)(9598,5786)(9599,5783)
	(9600,5778)(9601,5770)(9603,5759)
	(9606,5745)(9609,5728)(9613,5707)
	(9617,5684)(9622,5658)(9626,5631)
	(9632,5601)(9637,5569)(9642,5536)
	(9648,5502)(9653,5465)(9659,5427)
	(9664,5387)(9670,5344)(9676,5298)
	(9682,5250)(9688,5198)(9694,5142)
	(9699,5084)(9705,5024)(9710,4968)
	(9714,4913)(9718,4860)(9722,4810)
	(9725,4762)(9727,4716)(9730,4673)
	(9731,4632)(9733,4593)(9734,4556)
	(9736,4520)(9737,4485)(9737,4452)
	(9738,4420)(9739,4389)(9739,4360)
	(9739,4333)(9740,4309)(9740,4286)
	(9740,4267)(9740,4251)(9740,4238)
	(9740,4228)(9740,4220)(9740,4216)
	(9740,4213)(9740,4212)
\path(9462,7962)(9462,7961)(9462,7959)
	(9462,7954)(9462,7947)(9462,7938)
	(9462,7924)(9462,7908)(9462,7888)
	(9462,7864)(9462,7838)(9463,7808)
	(9463,7774)(9463,7739)(9463,7701)
	(9464,7661)(9464,7619)(9465,7575)
	(9465,7531)(9466,7485)(9466,7438)
	(9467,7391)(9468,7342)(9469,7293)
	(9470,7242)(9472,7191)(9473,7139)
	(9475,7085)(9477,7030)(9479,6973)
	(9481,6916)(9484,6856)(9487,6796)
	(9490,6735)(9493,6673)(9497,6612)
	(9502,6536)(9507,6464)(9512,6398)
	(9518,6336)(9523,6280)(9528,6229)
	(9533,6182)(9538,6140)(9543,6101)
	(9548,6064)(9553,6031)(9557,6000)
	(9562,5971)(9567,5944)(9571,5919)
	(9575,5896)(9579,5875)(9583,5857)
	(9586,5840)(9589,5826)(9592,5814)
	(9594,5805)(9596,5798)(9597,5793)
	(9597,5789)(9598,5788)(9598,5787)
\path(9598,2637)(9598,2638)(9599,2641)
	(9600,2646)(9601,2654)(9603,2665)
	(9606,2679)(9609,2696)(9613,2717)
	(9617,2740)(9622,2766)(9626,2793)
	(9632,2823)(9637,2855)(9642,2888)
	(9648,2922)(9653,2959)(9659,2997)
	(9664,3037)(9670,3080)(9676,3126)
	(9682,3174)(9688,3226)(9694,3282)
	(9699,3340)(9705,3400)(9710,3456)
	(9714,3511)(9718,3564)(9722,3614)
	(9725,3662)(9727,3708)(9730,3751)
	(9731,3792)(9733,3831)(9734,3868)
	(9736,3904)(9737,3939)(9737,3972)
	(9738,4004)(9739,4035)(9739,4064)
	(9739,4091)(9740,4115)(9740,4138)
	(9740,4157)(9740,4173)(9740,4186)
	(9740,4196)(9740,4204)(9740,4208)
	(9740,4211)(9740,4212)
\path(9462,462)(9462,463)(9462,465)
	(9462,470)(9462,477)(9462,486)
	(9462,500)(9462,516)(9462,536)
	(9462,560)(9462,586)(9463,616)
	(9463,650)(9463,685)(9463,723)
	(9464,763)(9464,805)(9465,849)
	(9465,893)(9466,939)(9466,986)
	(9467,1033)(9468,1082)(9469,1131)
	(9470,1182)(9472,1233)(9473,1285)
	(9475,1339)(9477,1394)(9479,1451)
	(9481,1508)(9484,1568)(9487,1628)
	(9490,1689)(9493,1751)(9497,1812)
	(9502,1888)(9507,1960)(9512,2026)
	(9518,2088)(9523,2144)(9528,2195)
	(9533,2242)(9538,2284)(9543,2323)
	(9548,2360)(9553,2393)(9557,2424)
	(9562,2453)(9567,2480)(9571,2505)
	(9575,2528)(9579,2549)(9583,2567)
	(9586,2584)(9589,2598)(9592,2610)
	(9594,2619)(9596,2626)(9597,2631)
	(9597,2635)(9598,2636)(9598,2637)
\path(9587,5809)(9522,5618)
\path(9598,5803)(9736,5650)
\put(9312,4212){\blacken\ellipse{100}{100}}
\put(9312,4212){\ellipse{100}{100}}
\put(9762,4212){\blacken\ellipse{100}{100}}
\put(9762,4212){\ellipse{100}{100}}
\path(912,4312)(912,4112)
\path(2112,4312)(2112,4112)
\path(12,4212)(11262,4212)
\dashline{60.000}(9312,8412)(9312,12)
\path(9012,4312)(9012,4112)
\path(10212,4312)(10212,4112)
\put(1962,3762){\makebox(0,0)[lb]{\smash{{{\SetFigFont{8}{9.6}{\familydefault}{\mddefault}{\updefault}$-m$}}}}}
\put(387,3762){\makebox(0,0)[lb]{\smash{{{\SetFigFont{8}{9.6}{\familydefault}{\mddefault}{\updefault}$-m\!\!-\!\!1$}}}}}
\put(4512,5862){\makebox(0,0)[lb]{\smash{{{\SetFigFont{11}{13.2}{\familydefault}{\mddefault}{\updefault}$L$}}}}}
\put(9012,4437){\makebox(0,0)[lb]{\smash{{{\SetFigFont{10}{12.0}{\familydefault}{\mddefault}{\updefault}$\gamma$}}}}}
\put(9762,4437){\makebox(0,0)[lb]{\smash{{{\SetFigFont{10}{12.0}{\familydefault}{\mddefault}{\updefault}$s_0$}}}}}
\put(9612,5902){\makebox(0,0)[lb]{\smash{{{\SetFigFont{11}{13.2}{\familydefault}{\mddefault}{\updefault}$L_0$}}}}}
\put(8712,3762){\makebox(0,0)[lb]{\smash{{{\SetFigFont{8}{9.6}{\familydefault}{\mddefault}{\updefault}$-1$}}}}}
\put(10147,3762){\makebox(0,0)[lb]{\smash{{{\SetFigFont{8}{9.6}{\familydefault}{\mddefault}{\updefault}$0$}}}}}
\end{picture}
}
\caption{}
\label{fig:LL}
\end{figure}

Then, we consider the following 
\index{B}{Mellin--Barnes type integral}%
Mellin--Barnes type integrals:
\begin{align}
& 
\index{A}{1Phi@$\Phi_m^+(t)$|main}%
\Phi_m^+(t)
  := \frac{1}{2\pi\sqrt{-1}}
     \int_{L_0} \frac{\Gamma(-\lambda)}{\Gamma(\lambda+1+m)}
     (2t)_+^\lambda d\lambda,
\label{eqn:MB2}
\\
& 
\index{A}{1Psi@$\Psi_m^+(t)$|main}%
\Psi_m^+(t)
  := \frac{1}{2\pi\sqrt{-1}}
     \int_L \frac{\Gamma(-\lambda)}{\Gamma(\lambda+1+m)}
     (2t)_+^\lambda d\lambda,
\label{eqn:MBPsi+}
\\
& 
\index{A}{1Phi@$\Phi_m(t)$|main}%
\Phi_m(t)
  := \frac{1}{2\pi\sqrt{-1}}
     \int_{L_0} \frac{\Gamma(-\lambda)}{\Gamma(\lambda+1+m)}
     \left( \frac{(2t)_+^\lambda}{\tan(\pi\lambda)}
          + \frac{(2t)_-^\lambda}{\sin(\pi\lambda)}\right) d\lambda,
\label{eqn:MBPhi}
\\
& 
\index{A}{1Psi@$\Psi_m(t)$|main}%
\Psi_m(t)
  := \frac{1}{2\pi\sqrt{-1}} 
     \int_L \frac{\Gamma(-\lambda)}{\Gamma(\lambda+1+m)}
     \left( \frac{(2\lambda)_+^\lambda}{\tan(\pi\lambda)}
          + \frac{(2t)_-^\lambda}{\sin(\pi\lambda)} \right) d\lambda.
\label{eqn:MBPsi}
\end{align}
We shall see that these integrals are special cases of
\eqref{eqn:Gdistr} and define distributions on $\mathbb{R}$. 
The next theorem is the main result of this section,
which will be derived from Proposition \ref{prop:Gdistr}
by applying the reduction formula of Meijer's $G$-functions.
\begin{theorem} \label{prop:Psiint}
{\upshape 1)} $\Phi_m^+(t)$ and $\Phi_m(t)$ are locally integrable functions
on $\mathbb{R}$.
Furthermore, for a sufficiently small
$\epsilon>0$,
$\Phi_m^+(t)|t|^{-\epsilon}$ and
$\Phi_m(t)|t|^{-\epsilon}$
are also locally integrable.

{\upshape 2)}
\begin{align}
& \Psi_m^+(t)
  = \Phi_m^+(t) - \sum_{k=1}^{m}
    \frac{(-1)^{k-1}}{2^k(m-k)!} \delta^{(k-1)}(t).
\label{eqn:PPmp}
\\
& \Psi_m(t)
  = \Phi_m(t) - \frac{1}{\pi} \sum_{k=1}^m
    \frac{(k-1)!}{2^k(m-k)!} t^{-k}.
\label{eqn:PPm}
\end{align}
See \eqref{eqn:xkminus} in Appendix for the definition of the
distribution $t^{-k}$. 
In particular,
$\Psi_m^+$ and $\Psi_m$ are defined as functionals on the space
\index{A}{C0kR@$C_0^{k}(\mathbb{R})$}%
$C_0^{m-1}(\mathbb{R})$ of compactly supported functions on
$\mathbb{R}$ with continuous derivatives up to $m-1$ if $m\ge1$.
\end{theorem}

The rest of this section is devoted to the proof of
 Theorem \ref{prop:Psiint}. 

Before regarding the integrals \eqref{eqn:MB2}--\eqref{eqn:MBPsi} 
 as those for distributions,
we consider the classic cases, namely,
their restrictions to $\mathbb{R}\setminus\{0\}$,
which are real analytic functions.

Let $L_i$ $(i=1,2,3)$ be contours that start at
$\gamma_i-\sqrt{-1}\infty$ and end at $\gamma_i+\sqrt{-1}\infty$,
and pass the real axis at some point $s_i$.
We assume
\begin{align*}
   &-\frac{m}{2} < \gamma_1,
   && s_1 < 0,
\\
   &-\frac{m}{2} < \gamma_2,
   && s_2 < -m,
\\
   &
   && s_3 < -m.
\end{align*}
Then, we have the following integral expressions of Bessel functions.
Although the results are classical,
we shall give a proof to illustrate the idea of passing from Bessel
functions to \textit{Bessel distributions}.
The proof below is based on the integral expressions
 of Meijer's $G$-functions
(see Appendix \ref{subsec:G}, see also Proposition \ref{prop:Gdistr} (2)):
\begin{lemSec} \label{lem:JYKmint}
{\upshape 1)}
For $t>0$,
\begin{align}
\index{A}{Jnuztilde@$\protect\widetilde{J}_\nu(z)$}%
\widetilde{J}_m(2\sqrt{2t}) 
&= (2t)^{-\frac{m}{2}}
 J_m(2\sqrt{2t})
\nonumber
\\
&= \frac{1}{2\pi\sqrt{-1}} \int_{L_1}
   \frac{\Gamma(-\lambda)}{\Gamma(\lambda+m+1)}
   (2t)_+^\lambda d\lambda.
\label{eqn:Jmint}
\end{align}

{\upshape 2)}
For $t>0$,
\begin{align}
\index{A}{Ynuztilde@$\protect\widetilde{Y}_\nu(z)$}%
\widetilde{Y}_m(2\sqrt{2t}) 
&= (2t)^{-\frac{m}{2}}
 Y_m(2\sqrt{2t})
\nonumber
\\
&= \frac{1}{2\pi\sqrt{-1}} \int_{L_2}
   \frac{\Gamma(-\lambda)}{\Gamma(\lambda+m+1)}
   \frac{(2t)_+^\lambda}{\tan(\pi\lambda)} d\lambda.
\label{eqn:Ymint}
\end{align}

{\upshape 3)}
For $t<0$,
\begin{align}
\index{A}{Knuztilde@$\protect\widetilde{K}_\nu(z)$}%
\widetilde{K}_m(2\sqrt{2|t|}) 
&= (2|t|)^{-\frac{m}{2}}
 K_m(2\sqrt{2|t|})
\nonumber
\\
&= \frac{(-1)^{m+1}}{4\sqrt{-1}} \int_{L_3}
   \frac{\Gamma(-\lambda)}{\Gamma(\lambda+m+1)}
   \frac{(2t)_-^\lambda}{\sin(\pi\lambda)} d\lambda.
\label{eqn:Kmint}
\end{align}
\end{lemSec}

\begin{proof}[Proof of Lemma \ref{lem:JYKmint}]
Each of the first equalities is by the definition of
the normalized Bessel functions
\index{A}{Jnuztilde@$\widetilde {J}_\nu (z)$}%
$\widetilde{J}_m$, 
\index{A}{Ynuztilde@$\widetilde {Y}_\nu (z)$}%
$\widetilde{Y}_m$, and 
\index{A}{Knuztilde@$\widetilde {K}_\nu (z)$}%
$\widetilde{K}_m$ given in
\eqref{eqn:Jtilde}, \eqref{eqn:Ytilde}, and \eqref{eqn:Ktilde},
respectively. 
Let us verify the second equalities (the integral formulas for the
Bessel functions).

1)
By the reduction formula \eqref{eqn:G} of the $G$-function
$G_{02}^{10}$, we have
\begin{equation*}
   (2t)^{-\frac{m}{2}} J_m (2\sqrt{2t})
   = G_{02}^{10} (2t \bigm| 0,-m)
\end{equation*}
for $t>0$.
Then, by the integral expression \eqref{def:G} of the $G$-function
$G_{02}^{10}$, we have
\begin{equation*}
   G_{02}^{10} (2t \bigm| 0,-m)
   = \frac{1}{2\pi\sqrt{-1}} \int_{L_1}
     \frac{\Gamma(-\lambda)}{\Gamma(1+m+\lambda)}
     (2t)^\lambda d\lambda
\end{equation*}
for $t>0$.
Hence, \eqref{eqn:Jmint} is proved.

2)
By the reduction formula \eqref{eqn:GY} of the $G$-function
$G_{13}^{20}$, we have
\begin{equation*}
   (2t)^{-\frac{m}{2}} Y_m(2\sqrt{2t})
   = G_{13}^{20} \left(2t \Bigm|
   \genfrac{}{}{0pt}{}{-m-\frac{1}{2}}{-m,0,-m-\frac{1}{2}} \right),
\end{equation*}
for $t>0$.
Then, by Example \ref{ex:G2013},
we have
\begin{align*}
  & G_{13}^{20}
   \left(2t \Bigm|
   \genfrac{}{}{0pt}{}{-m-\frac{1}{2}}{-m,0,-m-\frac{1}{2}} \right)
\\
  & = \frac{1}{2\pi\sqrt{-1}} \int_{L_2}
      \frac{\Gamma(-m-\lambda)\Gamma(-\lambda)}
           {\Gamma(m+\frac{3}{2}+\lambda)\Gamma(-m-\frac{1}{2}-\lambda)} 
      (2t)^\lambda d\lambda.
\end{align*}
Now, \eqref{eqn:Ymint} is deduced from this formula and the
following identity:
\begin{align*}
& 
\frac{\Gamma(-m-\lambda)\Gamma(\lambda+1+m)}
       {\Gamma(m+\frac{3}{2}+\lambda)\Gamma(-m-\frac{1}{2}-\lambda)}
 = \frac{1}{\tan\pi\lambda}
\quad\text{for any $m\in\mathbb{Z}$}.
\end{align*}
Here, the last identity is an elementary consequence of the formula
$\Gamma(z)\Gamma(1-z) = \frac{\pi}{\sin\pi z}$.

3)
By the reduction formula \eqref{eqn:GK} of the $G$-function
$G_{02}^{20}$, we have
\begin{equation*}
   (2|t|)^{-\frac{m}{2}} K_m (2\sqrt{2|t|})
   = \frac{1}{2} G_{02}^{20} \Bigl(2|t| \Bigm| 0,-m\Bigr).
\end{equation*}
Suppose $t<0$.
Then, again by the integral expression \eqref{def:G} of $G_{02}^{20}$,
the right-hand side amounts to
\begin{equation*}
   \frac{1}{4\pi\sqrt{-1}} \int_{L_3}
   \Gamma(-\lambda)\Gamma(-m-\lambda)(2t)_-^\lambda d\lambda.
\end{equation*}

Then, \eqref{eqn:Kmint} follows from the identity:
\begin{equation*}
   \Gamma(\lambda+1+m)\Gamma(-m-\lambda)
   = \frac{\pi}{\sin(-\pi(\lambda+m))}
   = \frac{(-1)^{m+1}\pi}{\sin\pi\lambda}.
\end{equation*}
Thus, all the statements of Lemma \ref{lem:JYKmint} are proved.
\end{proof}

The integrals in Lemma \ref{lem:JYKmint} do not depend on the choice
of $L_i$ $(i=1,2,3)$ as ordinary functions on
$\mathbb{R}\setminus\{0\}$.
However, as we mentioned in Remark \ref{rem:Gdistr}, they  
depend on the choice of $L_i$ as distributions on $\mathbb{R}$ because
the poles of the distributions $t_\pm^\lambda$ are located at
$\lambda=-1,-2,-3,\dotsc$.

To avoid this effect,
we need to impose more constraints on the contours $L_i$.
Thus, let us assume further $-1<s_1$ and $-m-1<s_i$ $(i=2,3)$.
Moreover, we assume $-1<\gamma_j$ $(j=1,2)$.
That is, we shall assume from now that the integral paths $L_i$ $(i=1,2,3)$
are under the following constraints:
\begin{alignat}{4}
& -1 < \gamma_1,
&&  &&-1 &&< s_1 < 0,
\label{eqn:L1cond}
\\
& -1 < \gamma_2,
\qquad\qquad
&& -m&&-1 &&< s_2 < -m,
\label{eqn:L2cond}
\\
& 
&& -m&&-1 &&< s_3 < -m.
\label{eqn:L3cond}
\end{alignat}
Then, the right-hand sides of
\eqref{eqn:Jmint}--\eqref{eqn:Kmint} define distributions on
$\mathbb{R}$,
which
are independent of the choice of
the integral paths $L_i$ $(i=1,2,3)$
subject to \eqref{eqn:L1cond}--\eqref{eqn:L3cond}.

\begin{proof}[Proof of Theorem \ref{prop:Psiint}]
The first statement is a special case of Proposition \ref{prop:Gdistr}
(3).

Let us show the second statement.
The contour $L$ used in \eqref{eqn:MBPsi+} and \eqref{eqn:MBPsi}
 meets the constraints 
\eqref{eqn:L2cond} and \eqref{eqn:L3cond}, and
can be used
as $L_2$ and $L_3$.
Likewise, the contour $L_0$ used in \eqref{eqn:MB2} and
\eqref{eqn:MBPhi} can be used as
$L_1$. 
Further, we shall assume that the contour $L_0$
 coincides with $L$ when
$\left|\operatorname{Im}\lambda\right|$ is sufficiently large.

The integrand of \eqref{eqn:MB2} has poles at 
$\lambda=-1,-2,\dots,-m$ inside the closed contour $L_0-L$,
and its residue is given by
\begin{equation*}
   \underset{\lambda=-k}{\operatorname{res}}
   \frac{\Gamma(-\lambda)}{\Gamma(\lambda+1+m)}
   (2t)_+^\lambda
   = \frac{(-1)^{k-1}}{2^{k}(m-k)!}
     \,  \delta^{(k-1)}(t)
\end{equation*}
for $k=1,2,\dots,m$ by \eqref{eqn:resxlmd}.
Therefore, by Cauchy's integral formula,
we have 
\begin{align*}
& \frac{1}{2\pi\sqrt{-1}} \left( \int_{L_0} - \int_L \right)
   \frac{\Gamma(-\lambda)}{\Gamma(\lambda+1+m)}
   (2t)_+^\lambda d\lambda
\\
&= \sum_{k=1}^{m} 
   \frac{(-1)^{k-1}}{2^k(m-k)!} \,  \delta^{(k-1)}(t)
\end{align*}
as distributions.
Hence, \eqref{eqn:PPmp} is proved.

Next, let us prove \eqref{eqn:PPm}.
We recall from \eqref{eqn:xpLaurent} that
the Laurent expansion of the distribution $t_\pm^\lambda$ at
$\lambda=-k$ $(k=1,2,\dotsc)$ is given by
\begin{align*}
   &t_+^\lambda
    = \frac{1}{\lambda+k}  \,  \frac{(-1)^{k-1}}{(k-1)!}
      \,  \delta^{(k-1)}(t) + t_+^{-k} + \dotsb ,
\\
   &t_-^\lambda
    = \frac{1}{(\lambda+k)(k-1)!} 
      \,  \delta^{(k-1)}(t) + t_-^{-k} + \dotsb .
\end{align*}
Combining with the Taylor expansions at $\lambda=-k$
$(k=1,2,\dots,m)$:
\begin{align*}
    \sin\pi\lambda
   &  = (-1)^k \pi(\lambda+k) + \dotsb ,
\\
    \tan \pi\lambda
   &  = \pi(\lambda+k) + \dotsb ,
\\
    \frac{\Gamma(-\lambda)2^\lambda}{\Gamma(\lambda+1+m)}
   &  = b_0 + b_1(\lambda+k) + \dotsb ,
\end{align*}
where $b_0 = \frac{(k-1)!}{2^k(m-k)!}$,
we have
\begin{align*}
& \frac{\Gamma(-\lambda)}{\Gamma(\lambda+1+m)}
  \left(
  \frac{(2t)_+^\lambda}{\tan(\pi\lambda)}
  +\frac{(2t)_-^\lambda}{\sin(\pi\lambda)}
  \right) 
\\
= {}
& \frac{b_0 \bigl( (-1)^{k-1}\delta^{(k-1)}(t)+(-1)^k\delta^{(k-1)}(t)
        \bigr)} 
       {\pi(k-1)!}
  \,  \frac{1}{(\lambda+k)^2}
\\
& + \left( \frac{b_0(t_+^{-k}+(-1)^kt_-^{-k})}{\pi}
           + \frac{b_1 \bigl( (-1)^{k-1} \delta^{(k-1)}(t)+(-1)^k
                   \delta^{(k-1)}(t) \bigr)}
                  {\pi(k-1)^2}
    \right)
    \frac{1}{\lambda+k}
  + \dotsb
\\
= {}
& \frac{(k-1)! \, t^{-k}}{2^k(m-k)! \, \pi} \, 
  \frac{1}{\lambda+k}
  + O(1),
\end{align*}
as $\lambda$ tends to $-k$.

Therefore, by \eqref{eqn:MBPhi} and \eqref{eqn:MBPsi},
we have
\begin{align*}
  \Psi_m(t)-\Phi_m(t)
& = \frac{1}{2\pi\sqrt{-1}} \left(\int_L-\int_{L_0}\right)
    \frac{\Gamma(-\lambda)}{\Gamma(\lambda+1+m)}
    \left( \frac{(2t)_+^\lambda}{\tan(\pi\lambda)}
           + \frac{(2t)_-^\lambda}{\sin(\pi\lambda)}
    \right)
    d\lambda
\\
& = -\sum_{k=1}^m
    \frac{(k-1)!}{2^k(m-k)! \, \pi} t^{-k}.
\end{align*}
Hence, \eqref{eqn:PPm} is proved. 
Now, we have completed the proof of Theorem \ref{prop:Psiint}.
\end{proof}

\begin{remNonumber}[%
\index{B}{Bessel distribution}%
Bessel distributions]
 \label{rem:PhiPsi}
We shall use the symbols
\begin{alignat}{2}
& \widetilde{J}_m(2\sqrt{2t_+})
&& = (2t)_+^{-\frac{m}{2}} J_m(2\sqrt{2t_+}),
\label{eqn:Jmdist}
\\
& \widetilde{K}_m(2\sqrt{2t_+})
&& = (2t)_+^{-\frac{m}{2}} K_m(2\sqrt{2t_+}),
\label{eqn:Kmdist}
\\
& \widetilde{Y}_m(2\sqrt{2t_-})
&& = (2t)_-^{-\frac{m}{2}} Y_m(2\sqrt{2t_-}),
\label{eqn:Ymdist}
\end{alignat}
to denote the distributions defined by the right-hand sides of
\eqref{eqn:Jmint}--\eqref{eqn:Kmint} and by the contours $L_i$
$(i=1,2,3)$ satisfying \eqref{eqn:L1cond}--\eqref{eqn:L3cond},
respectively.

It is noteworthy 
 that $\widetilde{J}_m(2\sqrt{2t_+})$ is a locally integrable function
 on $t\in\mathbb{R}$ in view of the Taylor expansion \eqref{eqn:Jtilde},
but $\widetilde{K}_m(2\sqrt{2t_+})$ and
$\widetilde{Y}_m(2\sqrt{2t_-})$ are not
(see \eqref{eqn:Ymzero} and \eqref{eqn:Kmzero}).
Then, by the above proof of Theorem \ref{prop:Psiint},
we have
\begin{align}
\index{A}{1Phi@$\Phi_m^+(t)$}%
\Phi_m^+(t)
&= \widetilde{J}_m(2\sqrt{2t_+})
\nonumber
\\
&= (2t)_+^{-\frac{m}{2}} J_m(2\sqrt{2t_+}),
\label{eqn:Bes1}
\\
\index{A}{1Psi@$\Psi_m^+(t)$}%
\Psi_m^+(t)
&= \widetilde{J}_m(2\sqrt{2t_+}) - \sum_{k=1}^{m}
   \frac{(-1)^{k-1}}{2^k(m-k)!} \,  \delta^{(k-1)}(t)
\nonumber
\\
&= (2t)_+^{-\frac{m}{2}} J_m (2\sqrt{2t_+}) -
   \sum_{k=1}^{m}
   \frac{(-1)^{k-1}}{2^k(m-k)!} \,  \delta^{(k-1)}(t),
\label{eqn:Bes2}
\\
\index{A}{1Psi@$\Psi_m(t)$}%
\Psi_m(t)
&= \widetilde{Y}_m (2\sqrt{2t_+}) + \frac{2(-1)^{m+1}}{\pi}
   \widetilde{K}_m(2\sqrt{2t_-})
\nonumber
\\
&= (2t)_+^{-\frac{m}{2}} Y_m(2\sqrt{2t_+})
   + \frac{2(-1)^{m+1}}{\pi} (2t)_-^{-\frac{m}{2}}
   K_m(2\sqrt{2t_-}).
\label{eqn:Bes3}
\end{align}
An alternative way to define the distributions
$\widetilde{K}_m(2\sqrt{2t_+})$
and
$\widetilde{Y}_m(2\sqrt{2t_-})$
is to use the infinite sum expressions 
\eqref{eqn:Ymzero} and \eqref{eqn:Kmzero}.
Let $t_\pm^{-k}$ be the regularized distributions given by the second
terms in \eqref{eqn:xpLaurent} and \eqref{eqn:7.1.4}.
In light of \eqref{eqn:Ymzero} and \eqref{eqn:Kmzero},
we define as \emph{distributions}.
\begin{align*}
   \widetilde{Y}_m(2\sqrt{2t_+})
   ={}
   &
   -\frac{1}{\pi} \sum_{k=1}^m
   \frac{(k-1)!}{2^k(m-k)!} t_+^{-k}
   +\frac{1}{\pi} \widetilde{J}_m
   (2\sqrt{2t_+}) \log (2t_+)
\\
   &
   -\frac{1}{\pi} \sum_{l=0}^\infty
   (-1)^l (2t_+)^l
   \frac{\psi(m+l+1)+\psi(l+1)}{l!(m+l)!},
\\
   \frac{2(-1)^{m+1}}{\pi} \widetilde{K}_m(2\sqrt{2t_-})
   ={}
   &
   \frac{-1}{\pi} \sum_{k=1}^m
   \frac{(-1)^k (k-1)!}{2^k(m-k)!} t_-^{-k}
   +\frac{1}{\pi} \widetilde{I}_m (2\sqrt{2t_-}) \log(2t_-)
\\
   &
   -\frac{1}{\pi} \sum_{l=0}^\infty (2t_-)^l
   \frac{\psi(m+l+1)+\psi(l+1)}{l!(m+l)!}.
\end{align*}
If we set
\begin{align}\label{eqn:Phim power}
\index{A}{1Phi@$\Phi_m(t)$}%
   \Phi_m(t)
   :={}
   &
   \frac{-1}{\pi} \sum_{l=0}^\infty (-2t)^l
   \frac{\psi(m+l+1)+\psi(l+1)}{l!(m+l)!}
\nonumber
\\
   &
   +\frac{1}{\pi} (\widetilde{J}_m(2\sqrt{2t_+})
   \log(2t_+) + \widetilde{I}_m(2\sqrt{2t_-}) \log(2t-1)),
\end{align}
then $\Phi_m(t)$ is a locally integrable function,
and
\begin{equation*}
   \widetilde{Y}_m(2\sqrt{2t_+}) + \frac{2(-1)^{m+1}}{\pi}
   \widetilde{K}_m(2\sqrt{2t_-})
\\
   = \phi_m(t) - \frac{1}{\pi} \sum_{k=1}^m
     \frac{(k-1)!}{2^k(m-k)!} t^{-k}.
\end{equation*}
This gives another explanation of \eqref{eqn:PPm}.
In this book,
we have adopted the integral expression \eqref{eqn:MBPhi} for the
definition of $\Phi_m(t)$ in place of the expansion \eqref{eqn:Phim power}.
\end{remNonumber}

\section{Differential equations for Bessel distributions}
\label{subsec:ode}

The kernel $K(\xzeta, \xzeta')$ of the unitary inversion operator
$\pi(w_0)$ 
is given by means of the Bessel distribution. 
In this section, 
we shall give a heuristic account on why the Bessel function 
arises in the kernel function. 

We begin with the observation that generic $L_+$-orbits on $C\times C$ are 
of codimension one
(see Lemma \ref{lem:codim}).
This would force that the kernel function of any
$L_+$-intertwining operator 
(see \eqref{def:Lp} for the definition of $L_+ = \Mmax_+ A$)
 should be a function 
of one variable $\langle \xzeta, \xzeta' \rangle$.
The second step is to
 make use of 
the differential equation arising from  
$\operatorname{Ad}(w_0)\nmax=\overline \nmax$
(see \eqref{eqn:Adw}). 
Since $\nmax$ acts on $L^2(C)^\infty$ as differential operators of
second order,
we get a differential equation of second order that the kernel
distribution $K(\xzeta,\xzeta')$ must satisfy
(see Proposition \ref{prop:ode}).
The technical point here is that we have avoided using the Casimir
operator of $K$ because it acts on $L^2(C)^\infty$ as a differential
operator of fourth order.

The argument here was the clue for us to find a coarse form of the
kernel $K(\xzeta,\xzeta')$.
Though we did not use the results of this section
for the actual proof of our main theorem,
 we think that the idea here is still helpful
to get a coarse solution to
 a similar problem (see Problem \ref{prob:Uncertainty} (2)) 
in other settings,
namely, to 
find the integral kernel of 
the unitary inversion operator $\pi(w_0)$ of the minimal 
representation of other groups.

Let $\theta : g \mapsto {}^t \! g^{-1}$ be the Cartan involution of $G$.
Since $g \in G = O(p,q)$ satisfies
${}^tg I_{p,q} g = I_{p,q}$ where
$I_{p,q}
=\begin{pmatrix}  I_p &  0 \\
                                0   &  -I_q
   \end{pmatrix}$, 
we have
${}^t g^{-1} = I_{p,q} g I_{p,q}^{-1}$.
Since $w_0 = I_{p,q}$, we get
\begin{equation}
\label{eqn:theta}
  \theta (g) = w_0 g w^{-1}_0.
\end{equation}

We let $L_+ = \Mmax_+ A$
act on $C$ by
$$
 me^{tH}\cdot \xzeta = e^t m \xzeta
$$
for $m \in \Mmax_+ \simeq O(p-1,q-1)$ and 
$a = e^{tE} \in A$ (see Section~\ref{subsec:act} for notation).

\begin{lemSec}\label{lem:theta}
The kernel function $K(\xzeta, \xzeta')$ 
of the unitary operator $\pi(w_0)$ satisfies the following
functional equation:
\begin{equation}
\label{eqn:thetal}
K(\xzeta, \xzeta') = K(\theta(l)\xzeta, l\xzeta')
\quad \text{for all $l \in L_+$.}
\end{equation}
\end{lemSec}
\begin{proof}
Building on the unitary representation
$(\pi, L^2(C))$,
we define another unitary representation $\pi^\theta$
on $L^2(C)$ by the following twist:
$$
  \pi^\theta (g):= \pi(\theta(g)).
$$
Then, \eqref{eqn:theta} implies
that $\pi(w_0): L^2(C) \to L^2(C)$ 
is an intertwining operator from $(\pi, L^2(C))$ to $(\pi^\theta, L^2(C))$.
In particular, we have,
\begin{equation}
\label{eqn:lw}
  \pi^\theta(l) \circ \pi(w_0) = \pi (w_0) \circ \pi(l)
  \quad \text{for any $l \in L_+$.}
\end{equation}
For $l = m \in \Mmax_+$,
we recall from \eqref{eqn:nbar} that
$$
 (\pi(m)u)(\xzeta) = u({}^t m\xzeta) \quad\text{for 
  $\xzeta \in C$}.
$$
Hence, for any $u \in L^2(C)$, we have
\begin{equation*}
   \int_C K({}^t\theta(m)\xzeta,\xzeta') u(\xzeta') d\mu(\xzeta') 
   = \int_C K(\xzeta,\xzeta'') u({}^t m\xzeta'') d\mu(\xzeta'').
\end{equation*}
Since $d\mu$ is $\Mmax_+$-invariant,
the right-hand side is equal to
\begin{equation*}
\int_C K(\xzeta,{}^tm^{-1}\xzeta') u(\xzeta') d\mu(\xzeta').
\end{equation*}
Since $u$ is arbitrary, the kernel function must coincide:
$$
 K({}^t\theta(m)\xzeta,\xzeta') = K(\xzeta, {}^tm^{-1}\xzeta').
$$
Replacing ${}^t\theta(m)\xzeta$ with $\xzeta$, we have
$$
 K(\xzeta,\xzeta') = K(m\xzeta, {}^tm^{-1}\xzeta')
 \quad \text{for any $m \in \Mmax_+$}.
$$
Thus, \eqref{eqn:thetal} holds for $l \in \Mmax_+$.
For $l = a := e^{tH} \in A$,
we recall from \eqref{eqn:rA} that 
$$
 (\pi(a)u)(\xzeta) = e^{-\frac{p+q-4}{2}t}u(e^{-t}\xzeta)
 \ (\xzeta \in C).
$$
Since $\pi^\theta(a) = \pi(a^{-1})$,
the equation \eqref{eqn:lw} amounts to
$$
 \pi(w_0) = \pi(a) \circ \pi(w_0) \circ \pi(a).
$$
Hence, for any $u \in L^2(C)$, we have
\begin{equation*}
 \int_C K(\xzeta,\xzeta')u(\xzeta')d\mu(\xzeta') 
  = e^{-(p+q-4)t}
   \int_C K(e^{-t}\xzeta,\xzeta'')u(e^{-t}\xzeta'')d\mu(\xzeta'').
\end{equation*}
By the formula \eqref{eqn:Q} of the measure $d\mu$
in the polar coordinate, we have
$$
 d\mu(\xzeta'') = e^{(p+q-4)t}d\mu(\xzeta')\quad 
\text{for $\xzeta' = e^{-t}\xzeta''$}.
$$
Thus, the right-hand side equals
$\int_C K(e^{-t}\xzeta,e^t\xzeta')u(\xzeta')d\mu(\xzeta')$.
Hence, we have
$$
 K(\xzeta,\xzeta') = K(e^{-t}\xzeta, e^t\xzeta')
 \quad \text{for any $t \in \mathbb{R}$}
$$
and therefore
$$
 K(\xzeta,\xzeta') = K(\theta(a)\xzeta, a\xzeta')
 \quad \text{for any $a \in A$}.
$$
Now, Lemma \ref{lem:theta} is proved.
\end{proof}

Now let $\Mmax_+$ act on the direct product manifold
$C \times C$ by the formula:
$$
  \Mmax_+ \times C \times C \to C \times C,
  \quad 
  (l, \xzeta, \xzeta') \mapsto (\theta(l)\xzeta, l \xzeta').
$$
Furthermore, we define the level set of $C \times C$ by
$$
  H_t := \set{(\xzeta, \xzeta') \in C \times C}
                    {\langle \xzeta, \xzeta' \rangle =t}
  \quad\text{for $t \in \mathbb R$}.
$$
Here, $\langle \cdot,\cdot\rangle$ is
 the standard positive definite
inner product 
on $\mathbb{R}^{p+q-2}$. Then we have:

\begin{lemSec}\label{lem:codim}
{\rm 1)} The level set $H_t$ is stable under the $\Mmax_+$-action.

{\rm 2)} Moreover $H_t$ is a single $\Mmax_+$-orbit for any non-zero $t$.
\end{lemSec}

\begin{proof}
1)\ For $\xzeta,\xzeta' \in \mathbb{R}^{p+q-2}\ (\subset \mathbb{R}^{p+q})$
and $l \in \Mmax_+ \simeq O(p-1,q-1)$, we have
\begin{align*}
  \langle \theta(l)\xzeta, l\xzeta'\rangle
 & = \langle w_0 lw^{-1}_0 \xzeta, l \xzeta'\rangle \\
 & = {}^t\xzeta'{}^tlw_0 lw^{-1}_0\xzeta \\
 & = {}^t\xzeta' w_0w^{-1}_0 \xzeta \\
 & = \langle \xzeta, \xzeta'\rangle .
\end{align*}
Hence, $H_t$ is $\Mmax_+$-stable.

2)\ We replace $(p-1,q-1)$ by $(p,q)$, and
consider the $G$-action on
$$
\index{A}{C0@$\widetilde{C}$}%
 \widetilde C := \{(\xzeta_0,\cdots , \xzeta_{p+q-1}):
 \xzeta^2_0 +\dots + \xzeta^2_{p-1} - \xzeta^2_p - \dots -
 \xzeta^2_{p+q-1} = 0 \}
$$
in place of the $\Mmax_+$-action on $C$
(this change allows us to use the notation $\Nmax$
and $\overline\Nmax$ in Section~\ref{subsec:act}).
Then, we recall from \eqref{eqn:GMN} that $G$ acts transitively on
$\widetilde C$ and the isotropy subgroup at 
$e_0 + e_{p+q-1} = {}^t(1,0,\cdots, 0,1)$
is given by $\Mmax_+\Nmax$ (see \eqref{eqn:GMN}).
For $t\ne0$, we take any $(\xzeta,\xzeta')\in H_t$.
Then, we find $g_1\in G$ such that
$g_1\xzeta'={}^t(1,0,\dots,0,1)$.
We write 
$\theta(g_1)\xzeta = {}^t(x_0,x,x_{p+q-1})
 \in \mathbb{R}\oplus\mathbb{R}^{p+q-2}\oplus\mathbb{R}$.
Then,
\begin{equation*}
  t = \langle \xzeta,\xzeta' \rangle
    = \langle \theta(g_1)\xzeta,g_1\xzeta' \rangle
    = x_0 + x_{p+q-1}.
\end{equation*}

Let us consider the orbit of 
$\theta(\Mmax_+\Nmax) = \Mmax_+ \overline{N}^{\max}$
on $\widetilde C$.
In view of \eqref{eqn:nbar}, we have
\begin{equation} \label{eqn:nbaction}
 \overline{n}_b
 \begin{pmatrix}
   x_0 \\ x \\ x_{p+q-1}
 \end{pmatrix}
 = 
 \begin{pmatrix}
   x_0 - {}^txw_0b \\ x \\ x_{p+q-1} + {}^txw_0b 
 \end{pmatrix}
+ \frac{x_0+x_{p+q-1}}{2}
  \begin{pmatrix}
    -Q(b) \\ 2b \\ Q(b)
  \end{pmatrix},
\end{equation}
for $b = {}^t(b_1,\cdots,b_{p+q-2})$ and
    $x = {}^t(x_1,\cdots,x_{p+q-2}) \in \mathbb{R}^{p+q-2}$. 
If $x_0 + x_{p+q-1} \ne 0$
and ${}^t(x_0,\cdots,x_{p+q-1}) \in \widetilde C$, we set
$$
 b := \frac{-x}{x_0+x_{p+q-1}}.
$$
Since
$x_0^2-x_{p+q-1}^2 = -Q(x)$, we have
\begin{equation*}
   b = \frac{x_0-x_{p+q-1}}{Q(x)}x
   \qquad\text{and}\qquad
   Q(b) = \frac{x_0-x_{p+q-1}}{x_0+x_{p+q-1}}.
\end{equation*}
Then, we have
$$
 \overline{n}_b
 \begin{pmatrix}
   x_0 \\ x \\ x_{p+q-1}
 \end{pmatrix}
=
\begin{pmatrix}
   x_0 \\ 0 \\ x_{p+q-1}
\end{pmatrix}
+\frac{Q(x)}{2(x_0+x_{p+q-1})}
\begin{pmatrix}
  1 \\ 0 \\ -1
\end{pmatrix}
= \frac{x_0+x_{p+q-1}}{2}
\begin{pmatrix}
  1 \\ 0 \\ 1
\end{pmatrix}.
$$
Now, we set
$g:=\theta(\overline{n}_b)g_1$.
Then,
\begin{equation*}
   (\theta(g)\xzeta,g\xzeta')
   = (\overline{n}_b\theta(g_1)\xzeta, \theta(\overline{n}_b)g_1\xzeta') 
   = (\frac{t}{2} \begin{pmatrix} 1\\ 0\\ 1 \end{pmatrix},
                  \begin{pmatrix} 1\\ 0\\ 1 \end{pmatrix}).
\end{equation*}
Hence, the second statement is proved.
\end{proof}

Thus, we have the following proposition by Lemma \ref{lem:theta} and 
Lemma \ref{lem:codim}.
\begin{proposition}\label{prop:Konevariable}
$K(\xzeta, \xzeta')\big|_{C\times C \setminus H_0}$ is  of the form
\begin{equation}\label{eqn:KPsi}
  K(\xzeta, \xzeta')= \Psi (\langle \xzeta, \xzeta' \rangle)
\end{equation}
for some function $\Psi(t)$ defined on $\mathbb R
\setminus \{0\}$. 
\end{proposition}

By lifting the inversion relation 
$\operatorname{Ad}(w_0)\nmax=\overline \nmax$ 
(see \eqref{eqn:Adw})
in the Lie algebra $\mathfrak g$ to the 
actions on $L^2(C)_K$, we get the 
differential equation satisfied by $\Psi$. 
More precisely, 
\begin{proposition}\label{prop:ode}
$\Psi(t)$ satisfies the following 
ordinary differential equation 
on $\mathbb R\setminus \{0\}$:
\begin{equation}\label{eqn:ode}
  t \frac{d^2 \Psi}{dt^2}+\frac{p+q-4}2 \frac{d\Psi}{dt} +2\Psi =0.
\end{equation}
\end{proposition}

\begin{proof}
It follows from 
$\operatorname{Ad}(w_0)\overline{N}_j = \epsilon_j N_j$
(see \eqref{eqn:Adwj}) that
\begin{equation}  \label{eqn:Ad} \pi(w_0)\circ d\pi (\overline{N}_j) = 
 \epsilon_j d\pi(N_j) \circ \pi(w_0).
\end{equation}
We recall from Sections \ref{subsec:act} and \ref{subsec:Pjb} that 
\begin{align*}
   &d\pi(\overline{N}_j) = 2\sqrt{-1} \, \xzeta_j
   && \text{(see \eqref{eqn:rN})},
\\
   &d\pi(N_j) = \frac{\sqrt{-1}}{2} \epsilon_j P_j
   && \text{(see \eqref{eqn:NjPj})},
\end{align*}
where 
\index{A}{P_j@$P_j$}%
$P_j$ 
are the 
\index{B}{fundamental differential operator}%
fundamental differential operators on the isotropic cone $C$ such that
\begin{equation} \label{eqn:Pjpsi}
   P_j\psi
   = \Bigl( \epsilon_j \xzeta_j\square-(2E+p+q-4)
     \frac{\partial}{\partial\xzeta_j} \Bigr) \tilde{\psi}|_C
\end{equation}
if $\psi=\tilde{\psi}|_C$ for a function defined in a neighborhood of
$C$.
Then, \eqref{eqn:Ad} leads us to the functional equation
for any test function $u(\xzeta')$:
\begin{equation} \label{eqn:weq}
   4\int_C\Psi(\langle\xzeta,\xzeta'\rangle)\xzeta'_j
   u(\xzeta')d\mu(\xzeta')
   = \int_C
     \bigl( P_j\Psi(\langle\xzeta,\xzeta'\rangle) \bigr)
     u(\xzeta')d\mu(\xzeta').
\end{equation}
Now, in view of \eqref{eqn:Pjpsi},
\begin{align*}
  P_j\Psi(\langle\xzeta,\xzeta'\rangle)
  ={}& \epsilon_j \xzeta_jQ(\xzeta') \Psi''(\langle\xzeta,\xzeta'\rangle)
 \\
  &   -\xzeta'_j (2\langle\xzeta,\xzeta'\rangle\Psi''
                           (\langle\xzeta,\xzeta'\rangle)
      + (p+q-4) \Psi'(\langle\xzeta,\xzeta'\rangle))
 \\
  ={}& -\xzeta'_j
     (2t\Psi''(t)+(p+q-4)\Psi'(t))
     \Big|_{t=\langle\xzeta,\xzeta'\rangle}.
\end{align*}
Hence, the functional equation \eqref{eqn:weq} for any
$u\in L^2(C)$ implies that $\Psi$ satisfies the following differential
equation: 
\begin{equation*}
  4\Psi(t) 
  = -\bigl( 2t\Psi''(t)+(p+q-4)\Psi'(t) \bigr).
\end{equation*}
Thus, Proposition \ref{prop:ode} is proved.
\end{proof}

We write $\theta = t\frac{d}{dt}$.
As $\theta^2 = t^2 \frac{dt^2}{dt} + t \frac{d}{dt}$,
\eqref{eqn:ode} is equivalent to
\begin{equation*}
   (\theta^2 + \frac{p+q-6}{2} \theta + 2t) \Psi(t) = 0.
\end{equation*}
Finally,
let us see directly the differential equations that
 the 
\index{B}{Bessel distribution!differential equation@---, differential equation}%
Bessel distributions
\index{A}{1Phi@$\Phi _m^+(t)$}%
$\Phi_m^+(t)$, 
\index{A}{1Psi@$\Psi _m^+(t)$}%
$\Psi_m^+(t)$, and 
\index{A}{1Psi@$\Psi _m(t)$}%
$\Psi_m(t)$
(see \eqref{eqn:MB2}, \eqref{eqn:MBPsi+}, and \eqref{eqn:MBPsi} for definition)
satisfy.
It is easy to see that
$\Phi_m^+(t)$, $\Psi_m^+(t)$, and $\Psi_m(t)$ solve the following
differential equation 
 (in an ordinary sense)
\begin{equation}\label{eqn:Bdiffeq}
   (\theta^2 + m\theta + 2t)u = 0
   \quad\text{on $\mathbb{R} \setminus \{0\}$}.
\end{equation}
Of course, 
this fits well with what Proposition \ref{prop:ode} asserts.
On the other hand,
the distribution $\Phi_m(t)$
(see \eqref{eqn:MBPhi} for definition)
does not appear in the kernel function $K(\xzeta,\xzeta')$.
We note that, 
as distributions, 
\begin{align*}
   & (\theta^2 + 2m\theta + 2t) \sum_{k=1}^m
     \frac{(-1)^{k-1}}{2^k(m-k)!} \, \delta^{(k-1)} (t) = 0,
\\
   & (\theta^2 + 2m\theta + 2t) \sum_{k=1}^m
     \frac{(k-1)!}{2^k(m-k)!} \, t^{-k} = \frac{1}{(m-1)!}.
\end{align*}
In particular,
\index{A}{1Phi@$\Phi _m(t)$}%
$\Phi_m(t)$ (see \eqref{eqn:MBPhi}) does not solve
\eqref{eqn:Bdiffeq}, but solves the third order
differential equation on $\mathbb{R} \setminus \{0\}$:
\begin{equation} \label{eqn:Phimt}
   \theta(\theta^2 + 2m\theta + 2t) \Phi_m(t) = 0.
\end{equation}

\chapter{Appendix: special functions}
\label{sec:app}

We have seen that various special functions arise naturally in the
analysis on the minimal representations.
Some of their fundamental properties (e.g.\ integral formulas,
differential equations, etc.) have been used in the proof of
the unitary inversion formulas.
Conversely, representation theoretic properties are reflected as
algebraic relations (e.g.\ functional equations) of such special
functions. 
Further, different models of the same representation yield 
functional equations 
connecting special functions arising from each model.

For the convenience of the reader, we collect the formulas 
and the properties of special functions that were used in the previous
chapters.

\section{Riesz distribution $x_+^\lambda$} 
\protect\index{B}{Riesz distribution|main}%
\protect\index{A}{xlambda@$x_+^\lambda$|main}%
\label{subsec:Riesz}

A distribution $f_\lambda$ on $\mathbb{R}$ with parameter
$\lambda\in\mathbb{C}$ is said to be 
\textit{a distribution with meromorphic parameter} $\lambda$ if the
pairing 
$$
\langle f_\lambda, \varphi\rangle
$$
is a meromorphic function of $\lambda$ for any test function
$\varphi\in C_0^\infty(\mathbb{R})$.
We say $f_\lambda$ has a pole at $\lambda=\lambda_0$ if
$\langle f_\lambda,\varphi\rangle$ has a pole at $\lambda=\lambda_0$
for some $\varphi$.
Then, taking a residue at $\lambda=\lambda_0$,
we get a distribution:
\begin{equation*}
   C_0^\infty(\mathbb{R}) \to \mathbb{C},
   \quad
   \varphi \mapsto \underset{\lambda=\lambda_0}{\operatorname{res}}
   \langle f_\lambda,\varphi\rangle,
\end{equation*}
which we denote by
$\underset{\lambda=\lambda_0}{\operatorname{res}}f_\lambda$.

By Cauchy's integral formula,
if $C$ is a contour surrounding $\lambda=\lambda_0$,
then we have
\begin{equation*}
  \underset{\lambda=\lambda_0}{\operatorname{res}}
  \langle f_\lambda,\varphi\rangle
  = \frac{1}{2\pi\sqrt{-1}}
    \int_C \langle f_\lambda,\varphi\rangle d\lambda,
\end{equation*}
and in turn we get an identity of distributions:
\begin{equation*}
   \underset{\lambda=\lambda_0}{\operatorname{res}}
   f_\lambda
   = \frac{1}{2\pi\sqrt{-1}}
    \int_C f_\lambda \, d\lambda.
\end{equation*}
A classic example of distributions with meromorphic parameter is the
Riesz distribution $x_+^\lambda$ defined as a locally integrable
function (and hence a distribution):
\begin{equation*}
   x_+^\lambda
   = \begin{cases} x^\lambda  &(x>0) \\
                   0          &(x\le 0)
     \end{cases}
\end{equation*}
for $\lambda\in\mathbb{C}$ such that $\operatorname{Re}\lambda > -1$.
Then, $x_+^\lambda$ extends meromorphically to the entire complex
plane,
and all the poles are located at $\lambda=-1,-2,\dotsc$.
The residue is given by
\begin{equation} \label{eqn:resxlmd}
   \underset{\lambda=-k}{\operatorname{res}}
   x_+^\lambda
   = \frac{(-1)^{k-1}}{(k-1)!}
   \, \delta^{(k-1)}(x)
\end{equation}
for $k=1,2,3,\dotsc$.
Here,
$\delta(x)$ is the 
\index{B}{Dirac delta function}%
Dirac delta function. 
To see this, we set
\begin{equation*}
   \varphi_N(x)
   := \varphi(x) - \sum_{k=1}^N
      \frac{\varphi^{(k-1)}(0)}{(k-1)!} \, x^{k-1}.
\end{equation*}
Then,
\begin{align*}
  \langle x_+^\lambda,\varphi\rangle
  = {}&\sum_{k=1}^N \frac{\varphi^{(k-1)}(0)}{(k-1)!}
     \int_0^1 x^{\lambda+k-1} dx
    + \int_0^1 x^\lambda \varphi_N(x) dx
\nonumber
\\
  &  + \int_1^\infty x^\lambda \varphi(x)dx
\nonumber
\\
  = {}&\sum_{k=1}^N \frac{1}{\lambda+k}
     \frac{\varphi^{(k-1)}(0)}{(k-1)!}
     + \int_0^1 x^\lambda
     \varphi_N(x) dx
\nonumber
\\
   &  + \int_1^\infty x^\lambda \varphi(x) dx.
\end{align*}
The first two terms have a simple pole at
$\lambda=-k$ with residue
\begin{equation*}
   \frac{\varphi^{(k-1)}(0)}{(k-1)!}
   = \left\langle \frac{(-1)^{k-1}}{(k-1)!}
     \, \delta^{(k-1)}(x),\varphi(x) \right\rangle, 
\end{equation*}
the second term is holomorphic if
$\operatorname{Re}\lambda > -N-1$
because
$\varphi_N(x) = O(x^N)$,
and the last term is an entire function of $\lambda$ because $\varphi$
is compactly supported.
Hence, \eqref{eqn:resxlmd} is proved.
Likewise,
\begin{equation*}
   x_-^\lambda
   := \begin{cases}
         |x|^\lambda &(x<0) \\
         0           &(x\ge0)
      \end{cases}
\end{equation*}
extends a distribution with meromorphic parameter $\lambda$ and all
the poles are located at $\lambda=0,-1,-2,\dotsc$.
They are simple poles with
\begin{equation} \label{eqn:resxlmd2}
   \underset{\lambda=-k}{\operatorname{res}}
   x_-^\lambda = \frac{\delta^{(k-1)}(x)}{(k-1)!}.
\end{equation}
We write the Laurent expansions of $x_+^\lambda$ and $x_-^\lambda$ at
$\lambda=-k$ $(k=1,2,3,\dotsc)$ as follows:
\begin{align}
   & x_+^\lambda
     = \frac{(-1)^{k-1}}{\lambda+k} \delta^{(k-1)}(x)
       + x_+^{-k} + (\lambda+k)x_+^{-k} \log x_+ + \dotsb,
\label{eqn:xpLaurent}
\\
   & x_-^\lambda
     = \frac{1}{\lambda+k} \delta^{(k-1)}(x)
       + x_-^{-k} + (\lambda+k)x_-^{-k} \log x_- + \dotsb.
\label{eqn:7.1.4}
\end{align}
Then, 
\index{A}{xplusminusk@$x_+^{-k}$}%
$x_+^{-k}$ and $x_-^{-k}$ are tempered distributions supported
on the half lines $x\ge0$ and $x\le0$, respectively.
We note that they are not homogeneous as distributions.
Then, the sum $x_+^\lambda+(-1)^\lambda x_-^\lambda$ becomes a
distribution with holomorphic parameter $\lambda$ in the entire
complex plane because
\begin{equation*}
   \underset{\lambda=-k}{\operatorname{res}}
   (x_+^\lambda + (-1)^\lambda x_-^\lambda) = 0
\end{equation*}
for $k=1,2,3,\dotsc$.
We now define a distribution
\begin{equation} \label{eqn:xkminus}
\index{A}{xminusk@$x^{-k}$}%
   x^{-k} := \bigl(x_+^\lambda + (-1)^\lambda x_-^\lambda\bigr) \big|_{\lambda=-k}.
\end{equation}
This distribution is homogeneous,
and coincides with $x_+^{-k} + (-1)^k x_-^k$.

For $k=1$, $x^{-1}$ is the distribution that gives 
\index{B}{Cauchy's principal value|main}%
Cauchy's principal value:
\begin{equation*}
   \langle x^{-1},\varphi \rangle
   = \lim_{\epsilon\downarrow0}
     \left(\int_{-\infty}^{-\epsilon} +
           \int_{\epsilon}^\infty\right) 
     \frac{\varphi(x)}{x} dx.
\end{equation*}
This formula is valid for any $\varphi\in C_0(\mathbb{R})$.
Likewise, $x^{-k}$ extends to a functional on the space
\index{A}{C0kR@$C_0^{k}(\mathbb{R})$}%
$C_0^{k-1}(\mathbb{R})$ of compactly supported functions on
$\mathbb{R}$ with continuous derivatives up to $k-1$.
See the textbook \cite{xGeSh} of Gelfand and Shilov for a nice
introduction to these distributions.

\section{Bessel functions $J_\nu, I_\nu, K_\nu, Y_\nu$}
\label{subsec:B}
For $\operatorname{Re}\nu > 0$, the series 
\begin{equation*}
\index{A}{Jnuz@$J_\nu(z)$|main}%
  J_\nu(z):=\Bigl(\frac{z}{2} \Bigr)^\nu
 \sum_{j=0}^\infty \frac{(-1)^j (\frac{z}{2})^{2j}}{j! \,  \Gamma(j+\nu+1)}
\end{equation*}
converges in the entire complex plane.
Its sum $J_\nu(z)$ is called the 
\index{B}{Bessel function!of the first kind@---, of the first kind}%
\textit{Bessel function of the first kind} of
order $\nu$
 (see \cite[\S 3.54]{xWa}).
$J_\nu(z)$ extends meromorphically to $\nu \in \mathbb{C}$ by the
Poisson integration formula:
\begin{equation*}
   J_\nu(z)
   = \frac{1}{\sqrt{\pi} \, \Gamma(\nu+\frac{1}{2})}
     \Bigl( \frac{z}{2} \Bigr)^\nu
     \int_{-1}^1 e^{\sqrt{-1}zt} (1-t^2)^{\nu-\frac{1}{2}} dt
\end{equation*}
and solves the 
\index{B}{Bessel differential equation}%
Bessel differential equation:
\begin{equation*}
   \Bigl( z^2 \frac{d^2}{dz^2} + z \frac{d}{dz} + (z^2-\nu^2)
   \Bigr) u = 0.
\end{equation*}

We set
\begin{equation*}
\index{A}{Ynuz@$Y_\nu(z)$|main}%
  Y_\nu(z):=\frac{J_\nu(z)\cos \nu\pi - J_{-\nu}(z)}{\sin \nu \pi}.
\end{equation*}
If $\nu$ is an integer, say $\nu=m$, then this definition reads as
$$
  Y_m(z) := \lim_{\nu \to m}
                   \frac{J_\nu(z)\cos \nu\pi - J_{-\nu}(z)}
                           {\sin \nu \pi}.
$$
$Y_\nu$ is known as the 
\index{B}{Bessel function!of the second kind@---, of the second kind}%
\textit{Bessel function of the second kind} 
 or
\textit{Neumann's function}.

Further, we define two more functions by
\begin{align*}
I_\nu(z)
:={}& e^{-\frac{\sqrt{-1}\nu\pi}{2}} J_\nu(e^{\frac{\sqrt{-1}\pi}{2}}z)
\\
={}& \Bigl(\frac{z}{2}\Bigr)^\nu
   \sum_{j=0}^\infty
   \frac{(\frac{z}{2})^{2j}}{j! \, \Gamma(j+\nu+1)},
\\
K_\nu(z)
:={}& \frac{\pi}{2\sin\nu\pi} (I_{-\nu}(z) - I_{\nu}(z)).
\end{align*}
Both of them solve the following differential equation:
\begin{equation*}
z^2 \frac{d^2u}{dz^2} + z \frac{du}{dz} - (z^2+\nu^2) u = 0.
\end{equation*}
\index{A}{Inuz@$I_\nu(z)$|main}%
$I_\nu(z)$ 
is known as the 
\index{B}{modified Bessel function!of the first kind@---, of the first kind}%
\textit{modified Bessel function of the first kind},
and is real when $\nu\in\mathbb{R}$ and $z>0$.

\index{A}{Knuz@$K_\nu(z)$|main}%
$K_\nu(z)$ 
is known as the 
\index{B}{modified Bessel function!of the third kind@---, of the third kind}%
\textit{modified Bessel function of the third kind} 
or \textit{Basset's function}.
Clearly we have
\begin{equation*}
K_{-\nu}(z) = K_{\nu}(z).
\end{equation*}
We call $J_{\nu}$, 
$Y_{\nu}$, $I_{\nu}$, and $K_{\nu}$ simply as 
\index{B}{JBessel@$J$-Bessel function|main}%
\textit{$J$-Bessel},
\index{B}{YBessel@$Y$-Bessel function|main}%
\textit{$Y$-Bessel},
\index{B}{IBessel@$I$-Bessel function|main}%
\textit{$I$-Bessel}, 
and 
\index{B}{K-Bessel function@$K$-Bessel function|main}%
\textit{$K$-Bessel functions}.

The $K$-Bessel function satisfies the following formula
(see \cite[\textbf{II}, \S7.11 (22)]{xerdHigherI}):
\begin{equation*}
   \Bigl(\frac{d}{z\,dz}\Bigr)^m (z^{-\nu} K_\nu(z))
   = (-1)^m z^{-\nu-m} K_{\nu+m}(z).
\end{equation*}
This formula may be stated as
\begin{equation*}
   \Bigl(-\frac{2d}{z\,dz}\Bigr)^m \widetilde{K}_\nu(z)
   = \widetilde{K}_{\nu+m}(z)
\end{equation*}
in terms of the normalized $K$-Bessel function \eqref{eqn:Ktilde}.
By the change of variables $z=2e^{-x}$,
the $m=1$ case amounts to:
\begin{align} \label{eqn:diffK2}
   \frac{d}{dx} (e^{-ax} \widetilde{K}_\nu(2e^{-x}))
   = {}
   & {}-ae^{-ax}\widetilde{K}_\nu(2e^{-x})
\nonumber
\\
   & {}+2e^{-(a+2)x} \widetilde{K}_{\nu+1}(2e^{-x}).
\end{align}

The $K$-Bessel functions $K_{\nu}(z)$ reduce to combinations of
elementary functions if $\nu$ is half of an odd integer.
For $n\in\mathbb{N}$ we have
\begin{align*}
K_{n+\frac{1}{2}}(z)
&= \Bigl(\frac{\pi}{2z}\Bigr)^{\frac{1}{2}}
   e^{-z} \sum_{j=0}^n \frac{(n+j)!}{j!(n-j)!} \frac{1}{(2z)^j}
\\
&= (-1)^n \Bigl(\frac{\pi}{2z}\Bigr)^{\frac{1}{2}} z^{n+1}
   \Bigl(\frac{d}{z \, dz}\Bigr)^n \frac{e^{-z}}{z}.
\end{align*}
For instance, if $n=0$, we have
\begin{equation}
\label{eqn:K12}
K_{\frac{1}{2}}(z) = \left(\frac{\pi}{2z}\right)^{\frac{1}{2}} e^{-z}.
\end{equation}

The following renormalization is sometimes convenient:
\begin{align}
\index{A}{Jnuztilde@$\protect\widetilde{J}_\nu(z)$|main}%
\widetilde{J}_\nu(z)&:= (\frac{z}2)^{-\nu}J_\nu(z)
 =\sum_{j=0}^\infty 
 \frac{(-1)^j (\frac{z}2)^{2j}}{j! \, \Gamma(\nu+j+1)},  
\label{eqn:Jtilde}
\\
\index{A}{Inuztilde@$\protect\widetilde{I}_\nu(z)$|main}%
\widetilde{I}_\nu(z)&:= (\frac{z}2)^{-\nu}I_\nu(z)
 =\sum_{j=0}^\infty 
 \frac{ (\frac{z}2)^{2j}}{j! \, \Gamma(j+\nu+1)},  
\label{eqn:Itilde}
\\
\index{A}{Ynuztilde@$\protect\widetilde{Y}_\nu(z)$|main}%
\widetilde{Y}_\nu(z)&:=(\frac{z}2)^{-\nu}Y_\nu(z),
\label{eqn:Ytilde}
\\
\index{A}{Knuztilde@$\protect\widetilde{K}_\nu(z)$|main}%
\widetilde{K}_\nu(z)&:=(\frac{z}2)^{-\nu}K_\nu(z).
\label{eqn:Ktilde}
\end{align}
By the Taylor expansion as above,
we see that both $\widetilde{J}_\nu(z)$ and $\widetilde{I}_\nu(z)$ are
holomorphic function of $z$ in the entire complex plane.

$\widetilde{J}_\nu(z)$ and 
$\widetilde{Y}_\nu(z)$ are linearly independent of
each other (whether $\nu$ is an integer or not) and form
a basis of the space of solutions to the following 
differential equation:
\begin{equation}\label{eqn:Bde}
  z \frac{d^2 u}{dz^2} + (2\nu+1) \frac{du}{dz} + zu =0,
\end{equation}
or equivalently,
\begin{equation*}
   (\theta^2 + 2\nu\theta + z^2)u = 0,
\end{equation*}
where $\theta:=z\frac{d}{dz}$.
On the other hand,
$\widetilde{I}_\nu(z)$ and $\widetilde{K}_\nu(z)$ solves
\begin{equation*}
   (\theta^2 + 2\nu\theta - z^2)u = 0.
\end{equation*}

In terms of Meijer's $G$-functions (Appendix \ref{subsec:G}) or the
Barnes hypergeometric function ${}_p F_q$ (see \eqref{eqn:BFpq}),
we shall have the following expressions:
(see \eqref{eqn:G}--\eqref{eqn:GY}):
\begin{align*}
\index{A}{Jnuztilde@$\protect\widetilde{J}_\nu(z)$}%
   \widetilde{J}_\nu(z)
   &= G_{02}^{10} \biggl( \frac{z^2}{4} \biggm|
      0, -\nu \biggr)
    = \frac{1}{\Gamma(\nu+1)} \, {}_0F_1
      \Bigl( \nu+1; -\frac{z^2}{4} \Bigr),
\\
\index{A}{Ynuztilde@$\protect\widetilde{Y}_\nu(z)$}%
   \widetilde{Y}_\nu(z)
   &= G_{13}^{20} \biggl( \frac{z^2}{4} \biggm|
      \begin{matrix}
         -\nu-\frac{1}{2}  \\
         -\nu,0,-\nu-\frac{1}{2}
      \end{matrix}
      \biggr),
\\
\index{A}{Knuztilde@$\protect\widetilde{K}_\nu(z)$}%
   \widetilde{K}_\nu(z)
   &= \frac{1}{2} 
      G_{02}^{20} \biggl( \frac{z^2}{4} \biggm|
      0,-\nu \biggr).
\end{align*}
The Mellin--Barnes type integral expression of
$\widetilde{J}_\nu(z)$, $\widetilde{Y}_\nu(z)$, and
$\widetilde{K}_\nu(z)$ is also given in Lemma \ref{lem:JYKmint}.

For $m=1,2,3,\dotsc$, the infinite sum expressions of $Y_m(z)$ and
$K_m(z)$ (or $\widetilde{Y}_m(z)$ and $\widetilde{K}_m(z)$) at $z=0$
are given in \cite[\textbf{II}, \S7.2, (31) and (37)]{xerdHigherI},
which may be stated as follows:
\begin{align}
  \widetilde{Y}_m(z)
   ={}& -\frac{1}{\pi} \sum_{k=1}^m (\frac{z}{2})^{-2k}
      \frac{(k-1)!}{(m-k)!}
\nonumber
\\
  &   + \frac{2}{\pi} \widetilde{J}_m(z) \log (\frac{z}{2})
\nonumber
\\
  & - \frac{1}{\pi} \sum_{l=0}^\infty (-1)^l (\frac{z}{2})^{2l}
   \,  \frac{\psi(m+l+1)+\psi(l+1)}{l!(m+l)!}.
\label{eqn:Ymzero}
\\
  \widetilde{K}_m(z)
   ={}& \frac{1}{2} \sum_{k=1}^m (-1)^{m-k} (\frac{z}{2})^{-2k}
      \frac{(k-1)!}{(m-k)!}
\nonumber
\\
  & + (-1)^{m+1} \widetilde{I}_m(z) \log (\frac{z}{2})
\nonumber
\\
  & + \frac{1}{2} (-1)^m \sum_{l=0}^\infty (\frac{z}{2})^{2l}
   \,  \frac{\psi(m+l+1)+\psi(l+1)}{l!(m+l)!}.
\label{eqn:Kmzero}
\end{align}
Here, the function 
\index{A}{1ypsi(z)@$\psi(z)$}%
$\psi(z)$ is the logarithmic derivative of the
gamma function:
\begin{equation*}
  \psi(z):= \frac{d\log\Gamma(z)}{dz}
          = \frac{\Gamma'(z)}{\Gamma(z)}.
\end{equation*}
The $\psi$ function is meromorphic with simple poles at
$z=0,-1,-2,\dotsc$.

Next, we summarize the 
\index{B}{asymptotic behavior!Bessel function@---, Bessel function}%
asymptotic behaviors of the Bessel functions:
\begin{fact}[{see \cite[Chapter 4]{xaar}, 
\cite[Chapter VII]{xWa}}] \label{fact:Bas}
The asymptotic behaviors of the Bessel functions 
at $z=0, \infty$ are given by

{\rm 1)} As $z$ tends to $0$, 
$J_\nu(z), I_\nu(z) =
O(z^{\nu})$.

For $\nu>0$,
\begin{equation}\label{eqn:Ktild0}
   \widetilde{K}_\nu(z)
   = \frac{\Gamma(\nu)}{2} \Bigl(\frac{z}{2}\Bigr)^{-2\nu}
     + o(z^{-2\nu})
\quad\text{as $z \to 0$}.
\end{equation}

{\rm 2)} As $z$ tends to infinity
\begin{alignat*}{2}
J_\nu(z) \sim {}
&\sqrt{\frac{2}{\pi z}}
 \biggl(\cos \Bigl(z-\frac{\nu\pi}{2}-\frac{\pi}{4}\Bigr)
 \sum_{j=0}^\infty \frac{(-1)^j(\nu,2j)}{(2z)^{2j}}
\\
&\qquad -\sin \Bigl(z-\frac{\nu\pi}{2} - \frac{\pi}{4}\Bigr)
 \sum_{j=0}^\infty \frac{(-1)^j(\nu,2j+1)}{(2z)^{2j+1}}\biggr)
&& \quad(|\arg z| < \pi),
\\
Y_\nu(z) \sim {}
&\sqrt{\frac{2}{\pi z}}
 \biggl(\sin \Bigl(z-\frac{\nu\pi}{2}-\frac{\pi}{4}\Bigr)
 \sum_{j=0}^\infty (-1)^j\frac{(\nu,2j)}{(2z)^{2j}}
\\
&\qquad +\cos \Bigl(z-\frac{\nu\pi}{2} - \frac{\pi}{4}\Bigr)
 \sum_{j=0}^\infty (-1)^j\frac{(\nu,2j+1)}{(2z)^{2j+1}}\biggr)
&& \quad (|\arg z| < \pi),
\\
I_\nu(z) \sim {}
&\frac{e^z}{\sqrt{2\pi z}}
\sum_{j=0}^\infty \frac{(-1)^j(\nu,j)}{(2z)^j}
 +\frac{e^{-z+(\nu+\frac{1}{2})\sqrt{-1}\pi}}{\sqrt{2\pi z}}
 \sum_{j=0}^\infty \frac{(\alpha,j)}{(2z)^j}
&& \quad (-\frac{\pi}{2} < \arg z < \frac{3}{2}\pi),
\\
K_\nu(z) \sim {}
&\sqrt{\frac{\pi}{2z}} \, e^{-z}
 \Bigl(1+\sum_{j=1}^\infty \frac{(\nu,j)}{(2z)^j}\Bigr)
&& \quad (|\arg z| < \frac{3\pi}{2}).
\end{alignat*}
In particular, we have
\begin{equation*}
\widetilde{K}_\nu(2z) = \frac{\sqrt{\pi}}{2} \, e^{-2z}
\, z^{-\nu-\frac{1}{2}} (1+O(\frac{1}{z}))
\quad\text{as $z\to\infty$}.
\end{equation*}
Here, we have used 
\index{B}{Hankel's notation|main}%
Hankel's notation:
\begin{align*}
\index{A}{1alpha@$(\alpha,j)$|main}%
(\alpha,j) := {}
&
(-1)^j \frac{(\frac{1}{2}-\alpha)_j (\frac{1}{2}+\alpha)_j}{j!}
\\
= {}
&\frac{(4\alpha^2-1^2)(4\alpha^2-3^2)\cdots(4\alpha^2-(2j-1)^2)}
      {2^{2j} j!}.
\end{align*}

\end{fact}

Finally, we list some integral
formulas for the Bessel functions:

\begin{enumerate}[{\bf {B}1}]
\item
\label{item:B1} 
\index{A}{1zpropertiesB1-B3@\textbf{B1}--\textbf{B3}|(}%
(the 
\index{B}{Mellin transform!of Bessel function@---, of Bessel function}%
Mellin transform 
of $K$-Bessel functions, see
\cite[p.~684]{xGrRy}).
For
$\operatorname{Re}(\mu+1\pm\nu)>0$ and $\operatorname{Re}a>0$,
\begin{equation*}
  \int_0^\infty t^\mu K_\nu(at)dt=2^{\mu-1}a^{-\mu-1}
  \Gamma(\frac{1+\mu+\nu}{2})\Gamma(\frac{1+\mu-\nu}{2}).
\end{equation*}
Equivalently, we have
\begin{equation}
\label{eqn:B1}
  \int_0^\infty t^s \widetilde{K}_\nu(at)dt=2^{s-1}a^{-s-1}
  \Gamma(\frac{1+s}{2})\Gamma(\frac{1+s}{2}-\nu).
\end{equation}

\item
\label{item:B2}
Formula of the 
\index{B}{Hankel transform}%
Hankel transform due to W. Bailey 
\cite{xbail}\ (see also \cite[\S~19.6 (8)]{xerdIntII}).
\begin{alignat}{2}
     &\int_0^{\infty}
     t^{\lambda-1}J_{\mu}(a t) J_{\nu}(b t) K_{\rho}(c t) d t
\notag
\\
    &\quad = \frac {2^{\lambda-2}
              a^{\mu} b^{\nu}
              \Gamma({\frac 1 2}(\lambda+ \mu + \nu-\rho))
              \Gamma({\frac 1 2}(\lambda+ \mu + \nu+\rho))}
              {c^{\lambda+ \mu + \nu} \Gamma (\mu+1)\Gamma(\nu+1)}
\notag
\\
\label{eqn:B2}
     &\qquad \times
       F_4(\frac 1 2 (\lambda+ \mu + \nu -\rho),
           \frac 1 2 (\lambda+ \mu + \nu +\rho);
            \mu+1, \nu+1;
            -\frac {a^2}{c^2}, -\frac {b^2}{c^2}).  
\end{alignat}
Here, $F_4$ is 
\index{B}{Appell's hypergeometric function}%
Appell's hypergeometric function of two variables
(see \eqref{def:F_4}).

\item
\label{item:B3}
(see \cite[\S7, 14.2 (36)]{xerdIntII})
For $\operatorname{Re}(\alpha+\beta)>0$ and
$\operatorname{Re}(\rho\pm\mu\pm\nu+1)>0$,
\begin{align*}
& 2^{\rho+2} \Gamma(1-\rho) \int_0^\infty K_\mu(\alpha t) K_\nu(\beta t)
  t^{-\rho} dt
\\
&= \alpha^{\rho-\nu-1} \beta^\nu
 \trF
 (\frac{1+\nu+\mu-\rho}{2}, \frac{1+\nu-\mu-\rho}{2};
  1-\rho; 1-\frac{\beta^2}{\alpha^2})
\\
&\times \Gamma (\frac{1+\nu+\mu-\rho}{2})
        \Gamma (\frac{1+\nu-\mu-\rho}{2})
        \Gamma (\frac{1-\nu+\mu-\rho}{2})
        \Gamma (\frac{1-\nu-\mu-\rho}{2}).
\end{align*}

In particular, we have
\begin{equation}
\int_0^\infty K_\mu (2t)^2 t^{2s-1} dt
= \frac{\Gamma(s)^2 \Gamma(s+\mu) \Gamma(s-\mu)}{8\Gamma(2s)}.
\label{eqn:KKL2}
\end{equation}
\index{A}{1zpropertiesB1-B3@\textbf{B1}--\textbf{B3}|)}%
\end{enumerate}

\section{Associated Legendre functions $P_\nu^\mu$}
\label{subsec:L}
The 
\index{B}{associated Legendre function|main}%
associated Legendre functions 
on the interval 
$(-1, 1)$ is defined as the special value of the hypergeometric
function: 
\begin{equation}
\label{def:L}
\index{A}{Pnumux@$P_\nu^\mu(x)$|main}%
  P_\nu^\mu(x)=
  \frac{1}{\Gamma(1-\mu)}
  \Bigl( \frac{1+x}{1-x} \Bigr)^{\frac{\mu}{2}}
  \trF\bigl(-\nu, \nu+1; 1-\mu; \frac{1-x}{2} \bigr).
\end{equation}
The associated Legendre functions satisfy the following functional relation:
\begin{equation}
\label{eqn:L1}
  \frac{d}{dx}\Bigl( (1-x^2)^{-\frac{\mu}{2}}P_\nu^\mu(-x)\Bigr)
=(1-x^2)^{-\frac{\mu+1}{2}}P_\nu^{\mu+1}(-x),
\end{equation}
which is derived from the following recurrence relation
(see \cite[\S 8.733 (1)]{xGrRy}):   
$$
  (1-x^2)\frac{d}{dx}P_\nu^\mu(x)
=-\sqrt{1-x^2}P_\nu^{\mu+1}(x)-\mu xP_\nu^\mu(x).
$$

Integral formulas for the associated Legendre functions:

\begin{enumerate}[{\bf {L}1}]
\item
\label{item:L1}
\index{A}{1zpropertiesL1-L2@\textbf{L1}--\textbf{L2}|(}%
(see \cite[p.~803]{xGrRy})
Formula of the Riemann--Liouville integral:
$\operatorname{Re}\lambda<1$, $\operatorname{Re}\mu>0$, $0<y<1$,
\begin{equation}
\label{eqn:RL}
  \frac{1}{\Gamma(\mu)}
  \int_0^y (y-x)^{\mu-1}\bigl(x(1-x)\bigr)^{-\frac{\lambda}{2}}
  P_\nu^\lambda (1-2x)dx = \bigl(y(1-y)\bigr)^{\frac{\mu}{2}-\frac{\lambda}{2}}
  P_\nu^{\lambda-\mu}(1-2y).
\end{equation}

\item
\label{item:L2}
(see \cite[p.~798]{xGrRy}) \
For $2\operatorname{Re}\lambda>|\operatorname{Re}\mu|$,
\begin{equation}
\label{eqn:LGam}
  \int_{-1}^1 (1-x^2)^{\lambda-1} P_\nu^\mu(x)dx
= \frac{\pi2^\mu \Gamma(\lambda+\frac{\mu}{2})\Gamma(\lambda-\frac{\mu}{2})
            }
            {\Gamma(\lambda+\frac{\nu}{2}+\frac{1}{2})
             \Gamma(\lambda-\frac{\nu}{2})
             \Gamma(\frac{-\mu+\nu+2}2)
             \Gamma(\frac{-\mu-\nu+1}2)
            }.
\end{equation}
\index{A}{1zpropertiesL1-L2@\textbf{L1}--\textbf{L2}|)}%
\end{enumerate}

\section{Gegenbauer polynomials $C_l^\mu$}
\label{subsec:Ge}
Definition of the 
\index{B}{Gegenbauer polynomial|main}
Gegenbauer polynomials:
For $l \in \mathbb{N}$, we define
\begin{equation}
\label{def:Ge}
\index{A}{C1mu@$C_l^\mu(x)$|main}%
   C_l^\mu(x):=
   \frac{(-1)^l}{2^l}
   \frac{\Gamma(2\mu+l)\Gamma(\mu+\frac{1}{2})
            }
            {\Gamma(2\mu)\Gamma(\mu+l+\frac{1}{2})
            }
   \frac{(1-x^2)^{\frac{1}{2}-\mu}}{l!}
   \frac{d^l}{dx^l} \bigl( (1-x^2)^{\mu+l-\frac{1}{2}} \bigr).
\end{equation}

Slightly different from the usual notation in the literature, 
we adopt the following normalization of 
the Gegenbauer polynomial: 
\begin{equation}\label{def:noGe}
\index{A}{Clmuxtilde@$\tilC_l^\mu(x)$}%
 \tilC_l^\mu(x)
:=\Gamma(\mu)C_l^\mu(x).
\end{equation}
By using 
\index{B}{Gauss's duplication formula}%
Gauss's duplication formula
\begin{equation}\label{eqn:dup}
\Gamma(2\mu) = 2^{2\mu-1} \pi^{-\frac{1}{2}} \Gamma(\mu)
\Gamma(\mu+\frac{1}{2}),
\end{equation}
the definition \eqref{def:Ge} may be stated as
\begin{equation}\label{eqn:Getilde}
\tilC_l^\mu(x)
= \frac{(-1)^l
  \Gamma(2\mu+l)\sqrt{\pi}}{2^{2\mu+l-1}l! \, \Gamma(\mu+l+\frac{1}{2})} 
  (1-x^2)^{-\mu+\frac{1}{2}} \frac{d^l}{dx^l}
  \bigl((1-x^2)^{\mu+l-\frac{1}{2}}\bigr).
\end{equation}
The special value at $\mu=0$ is given by the limit formula
(see \cite[\S 3.15.1 (14)]{xerdHigherI}):
\begin{equation}
 \tilC_l^0(\cos \theta)=\lim_{\mu\to 0}
 \Gamma(\mu)C_l^\mu(\cos \theta)=\frac{2\cos (l\theta)}l.
\end{equation}

On the other hand,
the special value at $l=0$ is given by
\begin{equation*}
\tilC_0^\mu(x) = \Gamma(\mu).
\end{equation*}

Connection with Gauss's hypergeometric function
(see \cite[\S 3.15 (3)]{xerdHigherI}):
\begin{align}
\label{eqn:Gehy}
  \tilC_l^\mu(x)
  &=\frac{\Gamma(l+2\mu)\Gamma(\mu)}{\Gamma(l+1)\Gamma(2\mu)}
  \trF(l+2\mu, -l; \mu+\frac{1}2; \frac{1-x}2 )
\nonumber
\\
&= \frac{\Gamma(l+2\mu)\Gamma(\mu)}{\Gamma(l+1)\Gamma(2\mu)}
   \trF(\frac{l+2\mu}{2},-\frac{l}{2}; \mu+\frac{1}{2}; 1-x^2).
\end{align}
Here, the second equation is derived from
 the formula of quadratic transformation for 
hypergeometric function (see \cite[\S 2.11 (2)]{xerdHigherI}):
$$
  \trF(a, b; a+b+\frac{1}2; 4x(1-x))=\trF(2a, 2b; a+b+\frac{1}2; x).
$$

By using Kummer's transformation formula for the hypergeometric functions:
$$
  \trF(\alpha, \beta; \gamma; z)=(1-z)^{\gamma-\alpha-\beta}
  \trF(\gamma-\alpha, \gamma-\beta; \gamma; z),
$$
one can obtain the following relationship 
between the Gegenbauer polynomials and the associated Legendre functions.
\begin{equation}
\label{eqn:GeL}
  \tilC_l^\mu (x) =
  \frac{\sqrt{\pi}\Gamma(2\mu+l)}{2^{\mu-\frac{1}2}\Gamma(l+1)}
  (1-x^2)^{\frac{1}{4}-\frac{\mu}{2}}
  P_{\mu+l-\frac{1}{2}}^{\frac{1}{2}-\mu}(x), \quad -1<x<1.
\end{equation}

Integral formulas for the Gegenbauer polynomials:

\begin{enumerate}[{\bf {Ge}1}] 
\item\label{item:Ge1}
\index{A}{1zpropertiesGe1-Ge4@\textbf{Ge1}--\textbf{Ge4}|(}%
(Orthogonality relations; see \cite[\S3.15.1 (17)]{xerdHigherI})\quad 
For $\operatorname{Re}\mu>-\frac{1}{2}$,
\begin{equation}
\label{eqn:Ge1}
   \int_{-1}^1 \tilC_l^\mu (x) \tilC_m^\mu(x) 
  (1-x^2)^{\mu-\frac{1}{2}}dx
  = \begin{cases}
  0 &\text{if $l\neq m$},  \\
 \frac{2^{1-2\mu}\pi \Gamma(l+2\mu)}{(l+\mu)\Gamma(l+1)}   
    &\text{if $l=m$}.
  \end{cases}
\end{equation}

\item
\label{item:Ge2}
(see \cite[\S 7.321]{xGrRy})\quad
For $\operatorname{Re}\mu>-\frac{1}{2}$,
\begin{equation}
\label{eqn:Ge2}
  \int_{-1}^1 (1-x^2)^{\mu-\frac{1}{2}} 
 e^{\sqrt{-1}ax} \tilC_l^\mu(x)dx
=\frac{\pi 2^{1-\mu} \Gamma(2\mu+l)}{\Gamma(l+1)}
  a^{-\mu}J_{\mu+l}(a).
\end{equation}

\item 
\label{item:Ge3}
(see \cite[Lemma 8.5.2]{xkmano2})\quad
For $\alpha \in \mathbb C$, $\operatorname{Re}\nu > -1$, 
and $l\in \mathbb N$,
\begin{multline}
\label{eqn:Ge3}
 \int_{-1}^1 J_\nu(\alpha\sqrt{x+1}) 
 \tilC_l^{\nu+\frac{1}{2}}(x) (1+x)^{\frac{\nu}{2}}(1-x)^{\nu} dx= \\
 \frac{2^\frac{3}2 (-1)^l \sqrt{\pi}\Gamma(2\nu+l+1)}
 {\alpha^{\nu+1}l!}
 J_{2\nu+2l+1}(\sqrt{2}\alpha).
\end{multline}

\item
\label{item:Ge4}
For $\operatorname{Re}\nu>-\frac{1}{2}$ and
$\operatorname{Re}\lambda>-1$, 
\begin{align}
&\int_{-x}^1(x+y)^\lambda \tilC_k^\nu(y) (1-y^2)^{\nu-\frac{1}{2}} dy
\nonumber
\\
&=
\frac{\sqrt{\pi}\Gamma(2\nu+k)\Gamma(\lambda+1)}{2^{\nu-\frac{1}{2}}k!} 
(1-x^2)^{\frac{\lambda}{2}+\frac{\nu}{2}+\frac{1}{4}}
P_{\nu+k-\frac{1}{2}}^{-\lambda-\nu-\frac{1}{2}} (-x).
\label{eqn:xyC}
\end{align}
\index{A}{1zpropertiesGe1-Ge4@\textbf{Ge1}--\textbf{Ge4}|)}%
\end{enumerate}
This formula \eqref{eqn:xyC} is essentially the integration formula
\eqref{eqn:RL} for the 
\index{B}{associated Legendre function}%
associated Legendre functions.
For the sake of completeness, we give a proof:
\begin{align*}
&\text{The left-hand side of \eqref{eqn:xyC}}
\\
={}& \frac{\sqrt{\pi}\Gamma(2\nu+k)}{2^{\nu-\frac{1}{2}}k!}
   \int_{-x}^1 (1-y^2)^{\frac{\nu}{2}-\frac{1}{4}}
   (x+y)^\lambda P_{\nu+k-\frac{1}{2}}^{\frac{1}{2}-\nu}
   (y) dy
 \qquad\text{by \eqref{eqn:GeL}}
\\
={}& \frac{2^{\lambda+1}\sqrt{\pi}\Gamma(2\nu+k)}{k!}
   \int_0^{\frac{1+x}{2}} \bigl((1-t)t\bigr)^{\frac{\nu}{2}-\frac{1}{4}}
   \Bigl(\frac{x+1}{2}-t\Bigr)^\lambda
   P_{\nu+k-\frac{1}{2}}^{\frac{1}{2}-\nu}
   (1-2t)dt
\\
={}& \frac{2^{\lambda+1}\sqrt{\pi}\Gamma(2\nu+k)\Gamma(\lambda+1)}{k!}
   \Bigl(\frac{1-x^2}{4}\Bigr)^{\frac{\lambda}{2}+\frac{\nu}{2}+\frac{1}{4}}
   P_{\nu+k-\frac{1}{2}}^{-\lambda-\nu-\frac{1}{2}} (-x)
 \qquad\text{by \eqref{eqn:RL}}
\\
={}& \text{the right-hand side of \eqref{eqn:xyC}}.
\end{align*}

\section{Spherical harmonics $\mathcal{H}^j(\mathbb{R}^m)$ 
and branching laws}
\label{subsec:H}

A 
\index{B}{spherical harmonics|main}%
\textit{spherical harmonics} 
$f$ of degree $j=0,1,2,\ldots$ is the
restriction to the unit sphere $S^{m-1}\subset\mathbb{R}^m$ of a
homogeneous harmonic polynomials of degree $j$ in $\mathbb{R}^m$.
Equivalently, $f$ is a smooth function satisfying the differential
equation:
$$
\Delta_{S^{m-1}} f = -j(j+m-2)f.
$$
The space of spherical harmonics of degree $j$ is denoted by
$$  
\index{A}{HjRm@$\mathcal{H}^j(\mathbb{R}^m)$|main}%
\mathcal{H}^j(\mathbb{R}^m):=
       \set{ f \in C^\infty(S^{m-1})}{\Delta_{S^{m-1}}f=-j(j+m-2)f}.
$$
When $m=1$, it is convenient to set:
$$
  \Har{0}{1}:= \mathbb C \mathbf{1}, \quad 
  \Har{1}{1}:=\mathbb C \text{sgn}, \quad
  \Har{j}{1}:=0 \quad (j\ge 2). 
$$

The following facts are well-known (see \cite[Introduction]{xHe},
\cite{xTa}):
\begin{enumerate}[{\bf {H}1}]
\index{A}{1zpropertiesH1-H3@\textbf{H1}--\textbf{H6}|(}%
\item
\label{item:H1}
  For $f \in \Har{j}{m}$, $f(-x)= (-1)^j f(x)$.

\item
\label{item:H2}
  $O(m)$ acts irreducibly on $\Har{j}{m}$.

\item
\label{item:H3}
 $\Har{j}{m}$ is still irreducible as an $SO(m)$-module if $m\geq 3.$

\item
\label{item:H4}
  $\mathcal{H}^j(\mathbb{R}^2)= \mathbb{C}e^{\sqrt{-1}j\theta} 
   \oplus \mathbb{C}e^{-\sqrt{-1}j\theta},~j\geq 1$ as $SO(2)$-modules, 
 where $\theta = \tan^{-1}\frac{y}{x},~ (x, y) \in \mathbb{R}^2.$

\item
\label{item:H5}
  $\Har{j}{m}\big|_{O(m-1)}
     \simeq \bigoplus_{i=0}^j \Har{i}{m-1}$ as $O(m-1)$-modules. 

\item
\label{item:H6}
The Hilbert space $L^2(S^{m-1})$ decomposes into a direct sum of the
space of spherical harmonics:
$$
L^2(S^{m-1}) \simeq \sideset{}{^\oplus}\sum_{j=0}^\infty
   \mathcal{H}^j(\mathbb{R}^m).
$$
Here, $\sideset{}{^\oplus}\sum$ stands for the Hilbert completion of the
algebraic direct sum
$\bigoplus_{j=0}^\infty \mathcal{H}^j(\mathbb{R}^m)$.
\index{A}{1zpropertiesH1-H3@\textbf{H1}--\textbf{H6}|)}%
\end{enumerate}

Let $(x_0, x) \in \mathbb R^m$, 
\thinspace
$x \in \mathbb R^{m-1}$
be a coordinate of $\mathbb R^m$. 
Then, the branching law \textbf{H\ref{item:H5}} 
is explicitly constructed by the $O(m-1)$-intertwining operator
\begin{equation*}
\index{A}{Iijm@$I_{i\to j}^m:\mathcal{H}^i(\mathbb{R}^{m-1})
      \to \mathcal{H}^j(\mathbb{R}^m)$}%
I_{i\to j}^m:
\mathcal{H}^i(\mathbb{R}^{m-1})
\to \mathcal{H}^j(\mathbb{R}^m)
\end{equation*}
as follows (see \cite[Chapter III]{xTa}):
\begin{fact}
\label{fact:H}
For $0 \le i \le j$ and $\phi \in \Har{i}{m-1}$, we define a function 
$I_{i\to j}^m \phi$ on $S^{m-1}$ by
\begin{equation}
\label{def:I}
  (I_{i\to j}^m (\phi))(x_0, x) :=
  |x|^i \phi\biggl( \frac{x}{|x|} \biggr) 
  \widetilde{C}_{j-i}^{\frac{m-2}2+i}(x_0).
\end{equation}
Here, $\widetilde{C}_l^\nu(z)$ is the 
normalized 
\index{B}{Gegenbauer polynomial}%
Gegenbauer polynomial (see \eqref{def:noGe}).
Then, 
\begin{enumerate}
    \renewcommand{\labelenumi}{{\upshape\theenumi)}}
\item 
$I_{i\to j}^m(\phi)\in\Har{j}{m}$.
\item 
$I_{i\to j}^m$ gives
an injective $O(m-1)$-homomorphism from
$\Har{i}{m-1}$ to $\Har{j}{m}$.
\item 
($L^2$-norm)
\begin{equation}\label{eqn:Iijnorm}
  \|I_{i\to j}^m (\phi) \|_{L^2(S^{m-1})}^2
 = \frac{2^{3-m-2i} \pi \Gamma(m-2+i+j)}
        {(j-i)! \, (j+\frac{m-2}2 )}
  \| \phi \|_{L^2(S^{m-2})}^2.
\end{equation}
\end{enumerate}
\end{fact}

\begin{proof}
We use the following coordinate:
\begin{equation}
  \label{eqn:Sm}
 [-1,1] \times S^{m-2} \to S^{m-1},\quad
(r,\eta) \mapsto \omega = (r,\sqrt{1-r^2}\eta).  
\end{equation}
Then, the standard volume form 
$d\omega$ on the unit sphere $S^{m-1}$ is given by
$(1-r^2)^{\frac{m-3}{2}}d\eta dr$.
Therefore,
$$
\|I_{i\to j}^m(\phi)\|_{L^2(S^{m-1})}^2
= \int_{-1}^1 \int_{S^{m-2}} (1-r^2)^i |\phi(\eta)|^2
|\tilC_{j-i}^{\frac{m-2}{2}+i} (r)|^2
(1-r^2)^{\frac{m-3}{2}} d\eta \, dr.
$$

Now, apply \eqref{eqn:Ge1}.
\end{proof}

We illustrate the intertwining operator $I_{ij}$ by the two extremal
 cases, 
$i=0$ and $i=j$:
\begin{example} \label{ex:Iij}
{\upshape 1)}
 The case $i=0$.  
Then,
\begin{equation}\label{eqn:I0j}
(I_{0\to j}^m\mathbf{1})(x_0, x)=\tilC_j^\frac{m-2}2(x_0)
\end{equation}
 is 
the generator of $O(m-1)$-invariant vectors in $\Har{j}{m}$,
where $\mathbf{1}$ is the constant function on $S^{m-1}$.

{\upshape 2)}
The case $i=j$.
Then, we have simply
\begin{equation}\label{eqn:Ijj}
I_{j\to j}^m(\phi) (x_0,x)
= \Gamma(m) |x|^j \phi (\frac{x}{|x|}).
\end{equation}
\end{example}

\section{Meijer's $G$-functions $G_{p,q}^{m,n}\Bigl( x \Bigm| 
\vcenter{\hbox{$a_1, \ \cdots, \ a_p$}
\hbox{$\,b_1, \ \cdots, \ b_q$}}\Bigr)$}
\label{subsec:G}

Let $m,n,p$ and $q$ be integers with $0 \le m \le q$, $0\le n \le p$
and 
$$
c^* := m+n-\frac{p+q}{2} \ge 0.
$$
Suppose further that the complex numbers $a_1,\dots,a_p$ and
$b_1,\dots,b_q$ fulfill the condition:
\begin{equation} \label{eqn:abinteger}
   a_j-b_k \ne 1,2,3,\dotsc
   \quad (j=1,\dots n; \, k=1,\dots,m).
\end{equation}
Then,
\index{B}{Meijer's $G$-function|main}%
{\it Meijer's $G$-function} 
of order $(m,n,p,q)$ is defined by the
\index{B}{Mellin--Barnes type integral}
Mellin--Barnes type integral
(see \cite[\S1.19, \S5.3]{xerdHigherI}, \cite[{\bf I}, \S 1]{xMe}, \cite[\S8.2]{xPBM}):
for $x>0$,
\begin{multline}\label{def:G}
\index{A}{Gpqmn@$G_{p,q}^{m,n} ( x "| 
        \genfrac{}{}{0pt}{}{a_1, \cdots,  a_p}
        {b_1,  \cdots, b_q}
        )$|main}%
  G_{p,q}^{m,n}\biggl( x \biggm| 
  \begin{matrix}
   a_1,& \cdots, & a_p  \\
   b_1, & \cdots, & b_q 
   \end{matrix}
   \biggr)    \\
  :=\frac{1}{2\pi\sqrt{-1}}
   \int_L
 \frac{\prod\limits_{j=1}^m \Gamma(b_j-\lambda)\prod\limits_{j=1}^n \Gamma(1-a_j+\lambda)}
 {\prod\limits_{j=m+1}^q \Gamma(1-b_j+\lambda) \prod\limits_{j=n+1}^p \Gamma(a_j-\lambda)}
   x^\lambda d\lambda,
\end{multline}
where an empty product is interpreted as 
$1$. 

The contour $L$ starts at the point
$\gamma-\sqrt{-1}\infty$
($\gamma$ is a real number satisfying \eqref{eqn:pqgamma} below if
$c^*=0$), leaving all the poles of the integrand of the forms
\begin{equation} \label{eqn:bpoles}
   \lambda = b_j, b_j+1, b_j+2, \dotsc
   \quad (1 \le j \le m)
\end{equation}
to the right, and all the poles of the forms
\begin{equation*}
   \lambda = a_j-1, a_j-2, a_j-3, \dotsc
   \quad (1 \le j \le n)
\end{equation*}
to the left of the contour and finishing at the point
$\gamma+\sqrt{-1}\infty$.

Here, the condition on the real number $\gamma$ is given by
\begin{equation} \label{eqn:pqgamma}
   (q-p)\gamma > \operatorname{Re} \mu,
\end{equation}
where we set
\begin{equation*}
   \mu := \sum_{j=1}^q b_j - \sum_{j=1}^p a_j
          + \frac{p-q}{2} + 1.
\end{equation*}
It follows from the asymptotic behavior of the gamma factors 
(see Lemma \ref{eqn:gamineq}) that 
the integral \eqref{def:G} converges and is independent of $\gamma$ 
if one of the following conditions holds:

1) $\ c^* > 0,   \ \vert\operatorname{arg}x\vert < c^*\pi$;

2) $\ c^* \ge 0, \ \vert\operatorname{arg}x\vert = c^*\pi$, \ 
   $(q-p)\gamma > \operatorname{Re}\mu$.

In particular,
the $G$-function extends holomorphically to the complex domain
$\vert\operatorname{arg}x\vert < c^*\pi$
if $c^*>0$.

The $G$-function is symmetric in the parameters $a_1,\dots,a_n$,
likewise in $a_{n+1},\dots,a_p$, in $b_1,\dots,b_m$, and in
$b_{m+1},\dots,b_q$.

Obvious changes of variables in the integral give
\begin{align*}
&   x^s 
  G_{p,q}^{m,n}\biggl( x \biggm| 
  \begin{matrix}
   a_1, & \cdots, & a_p  \\
   b_1, & \cdots, & b_q 
   \end{matrix}
   \biggr)    
  =
  G_{p,q}^{m,n}\biggl( x \biggm| 
  \begin{matrix}
   a_1+s, & \cdots, & a_p+s  \\
   b_1+s, & \cdots, & b_q+s 
   \end{matrix}
   \biggr),
\\
&
  G_{p,q}^{m,n}\biggl( x^{-1} \biggm| 
  \begin{matrix}
   a_1, & \cdots, & a_p  \\
   b_1, & \cdots, & b_q 
   \end{matrix}
   \biggr)    
   =
  G_{q,p}^{n,m}\biggl( x \biggm| 
  \begin{matrix}
   1-b_1, & \cdots, & 1-b_q  \\
   1-a_1, & \cdots, & 1-a_p 
   \end{matrix}
   \biggr).
\end{align*}

The $G$-function
$
G_{p,q}^{m,n} 
\Bigl(x\Bigm| \begin{matrix}a_1,\dots,a_p\\ b_1,\dots,b_q\end{matrix}\Bigr)
$
satisfies the differential equation
(see \cite[\S5.4 (1)]{xerdHigherI}):
\begin{equation} \label{eqn:Gdiffeq}
\Bigl((-1)^{p-m-n} x \prod_{j=1}^p (x\frac{d}{dx}-a_j+1)
- \prod_{j=1}^q (x\frac{d}{dx}-b_i)\Bigr) u = 0.
\end{equation}
If $p<q$, the only singularities of \eqref{eqn:Gdiffeq} are
$x=0,\infty$;
$x=0$ is a regular singularity,
$x=\infty$ an irregular one.
For example,
$G_{04}^{20} (x \mid b_1,b_2,b_3,b_4)$ satisfies
the fourth order differential equation:
\begin{equation}
\label{eqn:diffeqG24}
\prod_{j=1}^4 (x\frac{d}{dx} - b_j) u = 0.
\end{equation}

The condition \eqref{eqn:abinteger} implies that none of the poles of 
$\Gamma(b_j-\lambda)$ $(j=1,2,\dots,m)$ coincides with any of the
poles of $\Gamma(1-a_k+\lambda)$ $(k=1,\dots,n)$.
Suppose further that
\begin{equation*}
   b_j-b_k \ne 0, \pm1, \pm2, \dotsc
   \quad (1\le j< k\le m).
\end{equation*}
Then the integrand (as an ordinary function for $x>0$) has simple
poles at the points \eqref{eqn:bpoles}.
(We note that as a distribution of $x$, $x^\lambda$ has simple poles
at $\lambda=-1,-2,-3,\dotsc$,
and the analysis involved is more delicate;
see Sections \ref{subsec:MGdistr} and \ref{subsec:intPsi}.)
For $p\le q$, by the residue calculus, we obtain 
(see \cite[\textbf{I}, (7)]{xMe}):
\begin{align}
&  G_{p,q}^{m,n}\biggl( x \biggm| 
  \begin{matrix}
   a_1,& \cdots, & a_p  \\
   b_1, & \cdots, & b_q 
   \end{matrix}
   \biggr)    
\nonumber
\\
&= \sum_{k=1}^m
  \frac{\prod\limits_{\substack{j=1\\ j\ne k}}^m 
        \Gamma(b_j-b_k)
        \prod\limits_{j=1}^n \Gamma(1+b_k-a_j)}
       {\prod\limits_{j=m+1}^q 
        \Gamma(1+b_k-b_j) 
        \prod\limits_{j=n+1}^p \Gamma(a_j-b_k)}
  \,  x^{b_k}
\nonumber
\\
& \times {}_pF_{q-1}
  (1+b_k-a_1,\dots,1+b_k-a_p; 1+b_k-b_1,
   \overset{k}{\hat{\cdots}}\, , 1+b_k-b_q;
   (-1)^{p-m-n}x).
\label{eqn:Ghg}
\end{align}
Here, 
\index{A}{F5pq@${}_pF_q$|main}%
${}_pF_q$ 
denotes the 
\index{B}{Barnes generalized hypergeometric function|main}%
Barnes generalized hypergeometric function: 
\begin{equation} \label{eqn:BFpq}
   {}_pF_q(\alpha_1,\dots,\alpha_p; \beta_1,\dots,\beta_q; x)
   = \sum_{k=0}^\infty
     \frac{x^k \prod\limits_{j=1}^p 
           \alpha_j(\alpha_j+1) \cdots (\alpha_j+k-1)}
          {k!  \prod\limits_{j=1}^q 
           \beta_j(\beta_j+1) \cdots (\beta_j+k-1)}.
\end{equation}
For example,
$\trF (\alpha_1,\alpha_2;\beta_1;x)$ is the Gauss hypergeometric
function, and
\begin{equation*}
   {}_0F_1(\beta;x)
   = \sum_{k=0}^\infty
     \frac{x^k}{k! \, \beta(\beta+1)\cdots(\beta+k-1)}.
\end{equation*}

Similarly, for $q\le p$, 
if $a_j-a_k \ne 0, \pm1, \pm2, \dotsc$ $(1\le j<k \le n)$,
we have
\begin{align}
 &  G_{p,q}^{m,n}\biggl( x \biggm| 
  \begin{matrix}
   a_1, & \cdots, & a_p  \\
   b_1, & \cdots, & b_q 
   \end{matrix}
   \biggr)    
\nonumber
\\
&= \sum_{k=1}^n
  \frac{\prod\limits_{\substack{j=1\\ j\ne k}}^n
        \Gamma(a_k-a_j)
        \prod\limits_{j=1}^n \Gamma(b_j-a_k+1)}
       {\prod\limits_{j=n+1}^p 
        \Gamma(a_j-a_k+1) 
        \prod\limits_{j=m+1}^q \Gamma(a_k-b_j)}
   \, x^{a_k-1}
\nonumber
\\
 & \times {}_qF_{p-1}
   (1+b_1-a_k,\dots,1+b_q-a_k; 1+a_1-a_k, 
    \overset{k}{\hat{\cdots}} \, , 1+a_p-a_k; (-1)^{q-m-n} x^{-1}).
\label{eqn:Ghg2}
\end{align}

For $p\le q$,
it follows from \eqref{eqn:Ghg} that
\begin{equation} \label{eqn:Gzero}
  G_{p,q}^{m,n}\biggl( x \biggm| 
  \begin{matrix}
   a_1, & \cdots, & a_p  \\
   b_1, & \cdots, & b_q 
   \end{matrix}
   \biggr)    
  = O(|x|^{\min(\operatorname{Re}b_1,\dots,\operatorname{Re}b_m)})
\end{equation}
as $x\to0$ (see also \cite[\textbf{I}, \S5.4.1 (8)]{xerdHigherI},
but there is a typographical error:
$\max\operatorname{Re}b_h$ loc.\ cit.\ should be 
$\min\operatorname{Re}b_h$).
On the other hand, the 
\index{B}{asymptotic behavior!G-function@---, $G$-function}
asymptotic expansion 
of
$G_{p,q}^{m,n}(x)$ $(p\le q)$ for large $x>0$ that we need in this book
 is the following case:
\begin{fact}[{\cite[\textbf{VII}, Theorem 17]{xMe}}] \label{fact:Glarge}
Let $m,p$ and $q$ be integers satisfying
$$
0\le p\le q-2 \quad\text{and}\quad  p+1\le m\le q-1.
$$

Then the $G$-function $G_{p,q}^{m,0}(x)$ possesses the following
asymptotic expansion for large $x>0\,${\upshape:}
\begin{equation*}
   G_{p,q}^{m,0}(x)
   \sim A_{\phantom{m,}q}^{m,0}
        H_{p,q}(xe^{(q-m)\pi\sqrt{-1}})
      + \bar{A}_{\phantom{m,}q}^{m,0}
        H_{p,q}(xe^{(m-q)\pi\sqrt{-1}}).
\end{equation*}
\end{fact}

Here, $H_{p,q}(z)$ is a function that possesses the following
expansion (see \cite[\textbf{I}, (25)]{xMe}):
\begin{align*}
   H_{p,q}(z) ={}
   & \exp\left((p-q)z^{\frac{1}{q-p}}\right) z^\theta 
     \left(\frac{(2\pi)^{\frac{q-p-1}{2}}}{\sqrt{q-p}}
           +\frac{M_1}{z^{\frac{1}{q-p}}}
           +\frac{M_2}{z^{\frac{2}{q-p}}}
           + \dotsb
     \right),
\end{align*}
where $M_1, M_2, \dotsc$ are constants,
and $\theta$ is given by
\begin{equation*}
   \theta := \frac{1}{q-p}
   \biggl( \frac{p-q+1}{2} + \sum_{j=1}^q b_j - \sum_{j=1}^p a_j
   \biggr),
   \quad \text{\cite[\textbf{I}, (23)]{xMe}}.
\end{equation*}
The coefficients $A_{\phantom{m,}q}^{m,0}$ and
$\bar{A}_{\phantom{m,}q}^{m,0}$ are given by
\begin{alignat*}{2}
  & 
\index{A}{Aqm0@$A_{\ \ q}^{m,0}$}%
    A_{\phantom{m,}q}^{m,0}
    := (-2\pi\sqrt{-1})^{m-q}
       \, e^{-(b_{m+1}+\dots+b_q)\pi\sqrt{-1}},
  && \quad\text{\cite[\textbf{II}, (45)]{xMe}},
\\
  & 
\index{A}{Aqm0bar@$\bar{A}_{\ \ q}^{m,0}$}%
    \bar{A}_{\phantom{m,}q}^{m,0}
    := (2\pi\sqrt{-1})^{m-q}
       \, e^{(b_{m+1}+\dots+b_q)\pi\sqrt{-1}},
  && \quad\text{\cite[\textbf{II}, (46)]{xMe}}.
\end{alignat*}

\begin{example} \label{ex:G2013}
For $(m,n,p,q)=(2,0,1,3)$, $c^*=0$.
We take $\gamma$ such that
\begin{equation*}
   \gamma > \frac{1}{2} \operatorname{Re}(b_1+b_2+b_3-a_1).
\end{equation*}
Then, we have an integral expression:
\begin{equation*}
 G_{13}^{20} \left( x \biggm| \genfrac{}{}{0pt}{}{a_1}{b_1,b_2,b_3}\right)
= \frac{1}{2\pi\sqrt{-1}}
  \int_L \frac{\Gamma(b_1-\lambda)\Gamma(b_2-\lambda)}
              {\Gamma(1-b_3+\lambda)\Gamma(a_1-\lambda)}
  x^\lambda d\lambda,
\end{equation*}
where the integral path $L$ starts from $\gamma-\sqrt{-1}\infty$,
leaves $b_1, b_2$ to the right and ends at $\gamma+\sqrt{-1}\infty$ 
(see Figure \ref{fig:G24}).
\end{example}
\begin{figure}[H]
\setlength{\unitlength}{0.00033333in}
\begingroup\makeatletter\ifx\SetFigFont\undefined%
\gdef\SetFigFont#1#2#3#4#5{%
  \reset@font\fontsize{#1}{#2pt}%
  \fontfamily{#3}\fontseries{#4}\fontshape{#5}%
  \selectfont}%
\fi\endgroup%
{\renewcommand{\dashlinestretch}{30}
\begin{picture}(4824,8440)(-5000,-10)
\path(1275,3669)(1419,3526)
\path(1263,3657)(1199,3457)
\path(12,4213)(4812,4213)
\dashline{60.000}(2412,8413)(2412,13)
\path(2266,12)(2266,13)(2266,15)
	(2265,19)(2264,26)(2263,35)
	(2262,47)(2260,63)(2257,82)
	(2255,105)(2251,130)(2248,160)
	(2243,192)(2239,227)(2234,265)
	(2228,305)(2223,347)(2217,390)
	(2210,436)(2204,482)(2197,530)
	(2189,578)(2182,628)(2174,679)
	(2166,730)(2157,783)(2148,836)
	(2138,891)(2128,947)(2117,1004)
	(2106,1063)(2094,1124)(2081,1186)
	(2067,1250)(2053,1315)(2038,1381)
	(2022,1447)(2006,1513)(1987,1589)
	(1967,1663)(1948,1732)(1929,1796)
	(1911,1856)(1894,1912)(1877,1964)
	(1860,2012)(1845,2056)(1829,2098)
	(1814,2137)(1800,2174)(1786,2208)
	(1772,2241)(1758,2271)(1745,2300)
	(1733,2327)(1721,2353)(1709,2376)
	(1698,2398)(1689,2417)(1680,2435)
	(1672,2450)(1665,2463)(1659,2474)
	(1655,2482)(1651,2489)(1649,2493)
	(1647,2496)(1646,2497)(1646,2498)
\path(2262,8413)(2262,8412)(2262,8410)
	(2261,8406)(2260,8399)(2259,8390)
	(2258,8378)(2256,8362)(2253,8343)
	(2251,8320)(2247,8295)(2244,8265)
	(2239,8233)(2235,8198)(2230,8160)
	(2224,8120)(2219,8078)(2213,8035)
	(2206,7989)(2200,7943)(2193,7895)
	(2185,7847)(2178,7797)(2170,7746)
	(2162,7695)(2153,7642)(2144,7589)
	(2134,7534)(2124,7478)(2113,7421)
	(2102,7362)(2090,7301)(2077,7239)
	(2063,7175)(2049,7110)(2034,7044)
	(2018,6978)(2002,6912)(1983,6836)
	(1963,6762)(1944,6693)(1925,6629)
	(1907,6569)(1890,6513)(1873,6461)
	(1856,6413)(1841,6369)(1825,6327)
	(1810,6288)(1796,6251)(1782,6217)
	(1768,6184)(1754,6154)(1741,6125)
	(1729,6098)(1717,6072)(1705,6049)
	(1694,6027)(1685,6008)(1676,5990)
	(1668,5975)(1661,5962)(1655,5951)
	(1651,5943)(1647,5936)(1645,5932)
	(1643,5929)(1642,5928)(1642,5927)
\path(1642,5927)(1642,5926)(1641,5924)
	(1639,5920)(1636,5914)(1632,5906)
	(1627,5894)(1620,5881)(1612,5864)
	(1604,5845)(1594,5823)(1583,5799)
	(1571,5773)(1559,5745)(1547,5716)
	(1534,5686)(1520,5654)(1507,5621)
	(1493,5586)(1479,5550)(1465,5513)
	(1451,5474)(1437,5432)(1422,5389)
	(1407,5343)(1392,5295)(1377,5244)
	(1362,5190)(1347,5135)(1332,5078)
	(1317,5017)(1304,4957)(1291,4900)
	(1280,4846)(1271,4794)(1262,4746)
	(1255,4700)(1249,4656)(1243,4614)
	(1238,4575)(1234,4536)(1230,4500)
	(1227,4465)(1224,4431)(1222,4399)
	(1220,4368)(1218,4340)(1216,4314)
	(1215,4291)(1214,4271)(1214,4254)
	(1213,4240)(1213,4229)(1212,4222)
	(1212,4217)(1212,4214)(1212,4213)
\path(1646,2498)(1646,2499)(1645,2501)
	(1643,2505)(1640,2511)(1636,2519)
	(1631,2531)(1624,2544)(1616,2561)
	(1608,2580)(1598,2602)(1587,2626)
	(1575,2652)(1563,2680)(1551,2709)
	(1538,2739)(1524,2771)(1511,2804)
	(1497,2839)(1483,2875)(1469,2912)
	(1455,2951)(1441,2993)(1426,3036)
	(1411,3082)(1396,3130)(1381,3181)
	(1366,3235)(1351,3290)(1336,3347)
	(1321,3408)(1308,3468)(1295,3525)
	(1284,3579)(1275,3631)(1266,3679)
	(1259,3725)(1253,3769)(1247,3811)
	(1242,3850)(1238,3889)(1234,3925)
	(1231,3960)(1228,3994)(1226,4026)
	(1224,4057)(1222,4085)(1220,4111)
	(1219,4134)(1218,4154)(1218,4171)
	(1217,4185)(1217,4196)(1216,4203)
	(1216,4208)(1216,4211)(1216,4212)
\put(1062,3163){\makebox(0,0)[rb]{\smash{{{\SetFigFont{10}{14.4}{\familydefault}{\mddefault}{\updefault}$L$}}}}}
\put(2112,88){\makebox(0,0)[rb]{\smash{{{\SetFigFont{10}{14.4}{\familydefault}{\mddefault}{\updefault}$\gamma-\sqrt{-1}\infty$}}}}}
\put(2112,8263){\makebox(0,0)[rb]{\smash{{{\SetFigFont{10}{14.4}{\familydefault}{\mddefault}{\updefault}$\gamma+\sqrt{-1}\infty$}}}}}
\put(2622,3838){\makebox(0,0)[lb]{\smash{{{\SetFigFont{10}{14.4}{\familydefault}{\mddefault}{\updefault}$\gamma$}}}}}
\put(1611,5504){\makebox(0,0)[lb]{\smash{{{\SetFigFont{10}{14.4}{\familydefault}{\mddefault}{\updefault}$\times$}}}}}
\put(1611,2714){\makebox(0,0)[lb]{\smash{{{\SetFigFont{10}{14.4}{\familydefault}{\mddefault}{\updefault}$\times$}}}}}
\put(1892,5203){\makebox(0,0)[lb]{\smash{{{\SetFigFont{10}{14.4}{\familydefault}{\mddefault}{\updefault}$b_1$}}}}}
\put(1892,2398){\makebox(0,0)[lb]{\smash{{{\SetFigFont{10}{14.4}{\familydefault}{\mddefault}{\updefault}$b_2$}}}}}
\end{picture}
}
\caption{}
\label{fig:G24}
\end{figure}

\begin{example} \label{ex:G24}
If $p=0$, the $G$-function is denoted by
$G_{0,q}^{m,0}(x \mid b_1,\dots,b_q)$.
The $G$-function that we use most frequently in this book is of type 
$G_{04}^{20}$.
Again, we have $c^*=0$.
Then, we have an integral expression:
\begin{equation*}
   G_{04}^{20}(x \mid  b_1,b_2,b_3,b_4)
   = \frac{1}{2\pi\sqrt{-1}} \int_L
     \frac{\Gamma(b_1-\lambda)\Gamma(b_2-\lambda)}
          {\Gamma(1-b_3+\lambda)\Gamma(1-b_4+\lambda)}
     x^\lambda d\lambda,
\end{equation*}
where $L$ starts from $\gamma-\sqrt{-1}\infty$,
leaves $b_1$, $b_2$ to the right, and ends at $\gamma+\sqrt{-1}\infty$
(see Figure \ref{fig:G24}) for $\gamma\in\mathbb{R}$ such that
\begin{equation*}
   \gamma > \frac{1}{4}
   (\operatorname{Re}(b_1+b_2+b_3+b_4)-1).
\end{equation*}
\end{example}

In Section \ref{subsec:Klk}, we need the
following lemma on the asymptotic behavior:
\begin{lemSec}\label{lem:Gas}
The asymptotic behavior of the 
$G$-functions 
$G^{20}_{04}(x \mid b_1, b_2, 1-\gamma-b_1, 1-\gamma-b_2)$
are given as follows:

{\rm 1)} As $x$ tends to $0$,
 $G^{20}_{04}(x \mid b_1, b_2, 1-\gamma-b_1, 1-\gamma-b_2)
 = O(x^{\operatorname{min}(b_1,b_2)})$.

{\rm 2)} As $x$ tends to $\infty$,
\begin{align}
 &G^{20}_{04}(x \mid b_1, b_2, 1-\gamma-b_1, 1-\gamma-b_2)
\nonumber
\\
 ={}&-\frac{1}{\sqrt{2\pi}} \, x^{\frac{1-4\gamma}8}
 \cos \bigl(4x^\frac{1}4-(\gamma+b_1+b_2+\frac{1}4)\pi\bigr)
 (1+P_1 x^{-\frac{1}2}+P_2 x^{-1}+\cdots) 
\nonumber
\\
&+x^{\frac{1-4\gamma}8}
\sin \bigl(4x^\frac{1}4-(\gamma+b_1+b_2+\frac{1}4)\pi\bigr)
 (Q_1 x^{-\frac{1}4}+Q_2 x^{-\frac{3}4}+\cdots).
\label{eqn:Gasy}
\end{align} 
Here, $P_1, \cdots, Q_1, \cdots$ are the constants independent 
of $x$.
\end{lemSec}

\begin{proof}
1) This estimate is a special case of \eqref{eqn:Gzero}.

2) 
We apply Fact \ref{fact:Glarge} to the case
\begin{equation*}
   (m,p,q) = (2,0,4), \quad
   (b_1,b_2,b_3,b_4) = (b_1,b_2,1-\gamma-b_1, 1-\gamma-b_2).
\end{equation*}
Then, the coefficients 
$A_{\phantom{2}4}^{20}$, $\bar{A}_{\phantom{2}4}^{20}$ and the
constant $\theta$ amount to
\begin{align*}
   &A_{\phantom{2}4}^{20} = -\frac{1}{4\pi^2}
                            \, e^{(2\gamma-2+b_1+b_2)\pi\sqrt{-1}},
 \quad
    \bar{A}_{\phantom{2}4}^{20} = -\frac{1}{4\pi^2}
                            \, e^{-(2\gamma-2+b_1+b_2)\pi\sqrt{-1}},
\\
   & \theta = \frac{1-4\gamma}{8}.
\end{align*}
The expansion of $H_{0,4}(xe^{\pm2\pi\sqrt{-1}})$ is given by
\begin{equation*}
   H_{0,4}(xe^{\pm2\pi\sqrt{-1}})
   = e^{\mp(4x^{\frac{1}{4}}-\frac{1-4\gamma}{4}\pi)\sqrt{-1}}
     \, x^{\frac{1-4\gamma}{8}}
     \biggl( \frac{(2\pi)^{\frac{3}{2}}}{2}
            \pm \frac{M_1}{\sqrt{-1}x^{\frac{1}{4}}}
            + \dotsb 
     \biggr).
\end{equation*}
Hence, $G_{04}^{20}(x)$ has the following asymptotic expansion:
\begin{multline*}
   -\frac{1}{4\pi^2}
    \, e^{-(4x^{\frac{1}{4}}-\pi(\gamma+b_1+b_2-\frac{7}{4}))\sqrt{-1}} 
    \, x^{\frac{1-4\gamma}{8}}
    \biggl( \frac{(2\pi)^{\frac{3}{2}}}{2}
           + \frac{M_1}{\sqrt{-1}x^{\frac{1}{4}}} + \dotsb
    \biggr)
\\
   -\frac{1}{4\pi^2}
    \, e^{(4x^{\frac{1}{4}}-\pi(\gamma+b_1+b_2-\frac{7}{4}))\sqrt{-1}} 
    \, x^{\frac{1-4\gamma}{8}}
    \biggl( \frac{(2\pi)^{\frac{3}{2}}}{2}
           - \frac{M_1}{\sqrt{-1}x^{\frac{1}{4}}} + \dotsb
    \biggr),
\end{multline*}
which is expressed as the right-hand side of \eqref{eqn:Gasy} by
virtue of the formulas
$e^{c\pi\sqrt{-1}} + e^{-c\pi\sqrt{-1}}
 = 2\cos(c\pi)$
and
$e^{c\pi\sqrt{-1}} - e^{-c\pi\sqrt{-1}}
 = 2\sqrt{-1}\sin(c\pi)$.
\end{proof}

Finally, we list the reduction formulas of $G$-functions that are
used in this book:
\begin{alignat}{3}\label{eqn:G}
  G_{0 2}^{10} (x \mid a, b)={} & x^{\frac{1}2 (a+b)} J_{a-b}(2x^\frac{1}2)
 &&\quad \text{\cite[\S 5.6 (3)]{xerdHigherI}},   \\
\label{eqn:GK}
  G_{02}^{20} (x \mid a, b) ={}& 2x^{\frac{1}2 (a+b)} K_{a-b}(2x^\frac{1}2)
 && \quad \text{\cite[\S 5.6 (4)]{xerdHigherI}},  \\
\label{eqn:GJ}
  G_{04}^{20}(x \mid a, a+\frac{1}2, b, b+\frac{1}2 ) ={} & 
       x^{\frac{1}2 (a+b)}J_{2(a-b)}( 4x^\frac{1}4)
 &&\quad \text{\cite[\S 5.6 (11)]{xerdHigherI}},  \\
\label{eqn:GY}
  G_{13}^{20} \biggl( x \biggm| 
    \begin{matrix} 
      a-\frac{1}2   \\
     a, b, a-\frac{1}2 
     \end{matrix}
     \biggr)={}& x^{\frac{1}2(a+b)} Y_{b-a}(2x^\frac{1}2)
 && \quad \text{\cite[\S 5.6 (23)]{xerdHigherI}}.
\end{alignat}

\section{Appell's hypergeometric functions $F_1, F_2, F_3, F_4$}
\label{subsec:Ahy}

\index{B}{Appell's hypergeometric function|main}%
Appell's hypergeometric functions 
(in two variables) $F_1, F_2, F_3, F_4$ are
defined by the following double power series:
\begin{alignat}{1}
\label{def:F_1}
\index{A}{F1@$F_1(\alpha, \beta, \beta', \gamma; x, y)$|main}%
  F_1(\alpha, \beta, \beta', \gamma; x, y)&:=
   \sum_{m,n=0}^\infty  \frac{(\alpha)_{m+n} (\beta)_m (\beta')_n}
               {(\gamma)_{m+n} m! \, n!}   x^m y^n,    \\
\index{A}{F2@$F_2(\alpha, \beta, \beta', \gamma, \gamma'; x, y)$|main}%
  F_2(\alpha, \beta, \beta', \gamma, \gamma'; x, y)&:=
   \sum_{m, n=0}^\infty  \frac{(\alpha)_{m+n}(\beta)_m (\beta')_n}          
           {(\gamma)_m (\gamma')_n  m! \, n!}   x^m y^n,   
\label{def:F_2}
\\
\label{def:F_3}
\index{A}{F3@$F_3(\alpha, \alpha', \beta, \beta', \gamma; x, y)$|main}%
  F_3(\alpha, \alpha', \beta, \beta', \gamma; x, y)&:=
    \sum_{m, n=0}^\infty \frac{(\alpha)_m (\alpha')_n (\beta)_m (\beta')_n}
             {(\gamma)_{m+n} m! \, n!}   x^m y^n,     \\
\label{def:F_4}
\index{A}{F4@$F_4(\alpha, \beta, \gamma, \gamma'; x, y)$|main}%
  F_4(\alpha, \beta, \gamma, \gamma'; x, y)&:=
   \sum_{m, n=0}^\infty  \frac{(\alpha)_{m+n} (\beta)_{m+n}}
            {(\gamma)_m (\gamma')_n  m! \, n!}  x^m y^n.
\end{alignat}  

\begin{enumerate}[{\bf {Ap}1}]
  \item  
\label{item:Ap1}
\index{A}{1zpropertiesAp1-Ap5@\indexspace\par\textbf{Ap1}--\textbf{Ap5}|(}%
Reduction from $F_3$ to $F_1$ (\cite[\S 5.11, (11)]{xerdHigherI}):
\begin{alignat}{1}
\label{Appell 1}
  F_3(\alpha, \alpha', \beta, \beta', \alpha+\alpha'; x,y)=
 (1-y)^{-\beta'} F_1(\alpha, \beta, \beta', \alpha+\alpha'; x, \frac{y}{y-1}).
\end{alignat}

 \item
\label{item:Ap2}
 Reduction from $F_3$ to $\trF$ 
(\cite[\S 5.10, (4)]{xerdHigherI}):
\begin{alignat}{1}
\label{eqn:F3hy}
 F_3(\alpha, \gamma-\alpha, \beta, \gamma-\beta; \gamma; x, y)
=(1-y)^{\alpha+\beta-\gamma}  \trF(\alpha, \beta, \gamma; x+y-xy).
\end{alignat}

\item
\label{item:Ap3}
  Reduction formula of $F_4$ (\cite[\S 5.10, (8)]{xerdHigherI}):
\begin{alignat}{1}
&F_4(\alpha, \beta; 1+ \alpha-\beta, \beta;
 \frac{-x}{(1-x)(1-y)},\frac{-y}{(1-x)(1-y)})
\nonumber
\\
={}
&(1-y)^\alpha \trF(\alpha, \beta; 1+\alpha-\beta; \frac{-x(1-y)}{1-x}).
\label{eqn:F4hy}
\end{alignat}

\item
\label{item:Ap4} 
Single 
\index{B}{integral of Euler's type}%
integral of Euler's type 
for $F_1$ (\cite[\S 5.8.2, (5)]{xerdHigherI}):
For 
$ \operatorname{Re} \alpha >0$ and $\operatorname{Re}(\gamma-\alpha)>0$,
\begin{multline}
\label{Appell 3}
\kern-1em
F_1(\alpha, \beta, \beta', \gamma; x, y)=  
\\
 \frac{\Gamma(\gamma)}{\Gamma(\alpha)\Gamma(\gamma-\alpha)}
 \int_0^1 u^{\alpha-1} (1-u)^{\gamma-\alpha-1} (1-ux)^{-\beta}
 (1-uy)^{-\beta'} du .
\end{multline}

\item
\label{item:Ap5}
 Double integral of Euler's type 
for $F_3$ (\cite[\S 5.8.1, (3)]{xerdHigherI}):
\begin{multline}
\label{Appell 4}
 F_3(\alpha, \alpha', \beta, \beta', \gamma; x, y)=
 \frac{\Gamma(\gamma)}{\Gamma(\beta) \Gamma(\beta') 
       \Gamma(\gamma-\beta-\beta')}   \\
 \times \iint_D u^{\beta-1} v^{\beta'-1} (1-u-v)^{\gamma-\beta-\beta'-1}
                          (1-ux)^{-\alpha} (1-vy)^{-\alpha'} du dv  \\
 \operatorname{Re} \beta>0, \  \operatorname{Re} \beta'>0, \  
 \operatorname{Re} (\gamma-\beta-\beta')>0,
\end{multline}
where $D:=\set{(u, v) \in \mathbb R^2}{ u \ge 0, \   v \ge 0,  \  u+v\le 1}$.
\index{A}{1zpropertiesAp1-Ap5@\indexspace\par\textbf{Ap1}--\textbf{Ap5}|)}%
\end{enumerate}

\section{Hankel transform with trigonometric parameters}
\label{subsec:Han}

This section presents an integral formula \eqref{eqn:Han}
on the Hankel transform with two trigonometric parameters.
In the 
\index{B}{conformal model}%
conformal model \cite{xkors1} (i.e.\ the solution space to the
Yamabe equation) of the minimal representation,
$K$-finite vectors can be explicitly expressed in terms of spherical
harmonics (e.g.\ Gegenbauer's polynomials).
On the other hand, in the 
\index{B}{L2model@$L^2$-model}%
$L^2$-model (the 
\index{B}{Schr\"{o}dinger model}%
Schr\"{o}dinger model)
which is obtained by the Euclidean Fourier transform of the conformal model 
(or the 
\index{B}{N-picture@$N$-picture}%
$N$-picture in a terminology of representation theory),
it is not easy to find explicit $K$-finite vectors.
The formula \eqref{eqn:Han}
bridges these two
models and gives an explicit formula of $K$-finite vectors in the
Schr\"{o}dinger model
(see the proof of Lemma \ref{lem:prK}).

Since we have not found this formula in the literature,
we give a proof here for the sake of completeness. 
The method here is a generalization of the argument in 
\cite[\S 5.6, 5.7]{xkors3}.

\begin{lemSec}\label{lem:Han}
The following integral formula on the 
\index{B}{Hankel transform}%
Hankel transform 
holds:
\begin{align}\notag
&\int_0^\infty t^{\mu+1}
J_\mu \biggl(\frac{t\sin \theta}{\cos \theta+\cos \phi}\biggr)
J_\nu \biggl(\frac{t \sin \phi} {\cos \theta+\cos \phi}\biggr) 
K_\nu(t)dt  \\  \label{eqn:Han}
={}&\frac{2^{\nu-1}}{\sqrt{\pi}}
\Gamma(\mu-\nu+1)
(\cos \theta+\cos \phi)\sin^\mu \theta
\sin^\nu \phi \
\widetilde{C}_{\mu-\nu}^{\nu+\frac{1}2}(\cos \phi).
\end{align}
\end{lemSec}

\begin{proof}
By Baily's formula \eqref{eqn:B2} of the Hankel transform,
the left-hand side of \eqref{eqn:Han} equals
\begin{multline}\label{eqn:HanF4}
  \frac{\Gamma(\mu+\nu+1)}{\Gamma(\nu+1)}
   \frac{2^\mu \sin^\mu\theta \sin^\nu \phi}
        {(\cos \theta+\cos \phi)^{\mu+\nu}} \\
 \times F_4 (\mu+1, \mu+\nu+1; \mu+1, \nu+1;
  -(\frac{\sin\theta}{\cos \theta+\cos\phi})^2, 
   - (\frac{\sin\phi}{\cos \theta+\cos\phi})^2).
\end{multline}
Here $F_4$ denotes Appell's hypergeometric function (see \eqref{def:F_4}).
Thus, the proof of Lemma \ref{lem:Han} will be completed if we show
the following:
\begin{claim}\label{clm:F4Ge}
We have
\begin{multline}\label{eqn:F4Ge}
 F_4 (\mu+1, \mu+\nu+1; \mu+1, \nu+1;
  -(\frac{\sin\theta}{\cos \theta+\cos\phi})^2, 
   - (\frac{\sin\phi}{\cos \theta+\cos\phi})^2) \\
=\frac{(\cos \theta+\cos \phi)^{\mu+\nu+1}}{2^{\mu-\nu+1}\sqrt{\pi}}
 \frac{\Gamma(\mu-\nu+1)\Gamma(\nu+1)}
      {\Gamma(\mu+\nu+1)}
  \tilC_{\mu-\nu}^{\nu+\frac{1}2}(\cos \phi).
\end{multline}
\end{claim}
\end{proof}

Claim \ref{clm:F4Ge} is essentially a restatement of \cite[Lemma 5.7]{xkors3}.
For the convenience of the reader, we include its proof here.

\begin{proof}[Proof of Claim \ref{clm:F4Ge}]
We recall a quadratic transformation for hypergeometric functions
(see \cite[\S~2.11 (32)]{xerdHigherI}):
\begin{equation}\label{eqn:quad}
\trF(\alpha, \beta; 1+\alpha-\beta; z)
= (1-z)^{-\alpha} 
\trF(\frac{\alpha}2, \frac{\alpha + 1 - 2 \beta}2;
 1+\alpha -\beta; \frac{-4 z}{(1-z)^2}).
\end{equation}
Combining the reduction formula \eqref{eqn:F4hy} with 
\eqref{eqn:quad},
and using the symmetry of $a$ and $b$; 
$(c,x)$ and $(d, y)$ in $F_4(a,b;c,d; x,y)$,
we have
\begin{multline}
\label{eqn:F4hq}
  F_4(\alpha, \beta; \alpha, 1-\alpha+\beta; 
  \frac{-x}{(1-x)(1-y)}, \frac{-y}{(1-x)(1-y)}) \\
  = 
  \Biggl(\frac{(1-x)(1-y)}{1-xy} \Biggr)^\beta
  \trF (\frac{\beta}2, \frac{1-2\alpha+\beta}2; 1-\alpha+\beta;
  \frac{4y(1-x)(1-y)}{(1-xy)^2}
  ).
\end{multline}
Consider the change of variables from $(x,y)$ to $(\theta,\phi)$ by
the following identities:
$$
\frac{x}{(1-x)(1-y)} =\biggl(\frac{\sin\theta}{\cos \theta+\cos \phi}\biggr)^2,
 \quad \frac{y}{(1-x)(1-y)}=\biggl(\frac{\sin\phi}{\cos \theta+\cos \phi}\biggr)^2
$$
such that $(x,y)=(0,0)$ corresponds to $(\theta,\phi)=(0,0)$.
Then, a simple computation shows
$$
  \frac{1-xy}{(1-x)(1-y)} =
  \frac{2}{\cos \theta+\cos \phi},
  \quad 
  \frac{4y(1-x)(1-y)}{(1-xy)^2}=\sin^2\phi.
$$
Now, we set
\begin{equation*}
\alpha=\mu+1, \quad \beta=\mu+\nu+1,  
\end{equation*}
in \eqref{eqn:F4hq}. 
Then,
the left-hand side of \eqref{eqn:F4Ge} amounts to
\begin{equation}
\label{eqn:hymn}
 \biggl(\frac{2}{\cos \theta+\cos \phi} \biggr)^{-\mu-\nu-1}
  \trF
  \bigl(\frac{\mu+\nu+1}2, -\frac{\mu-\nu}2; \nu+1; \sin^2\phi \bigr).
\end{equation}
By using \eqref{eqn:Gehy}, \eqref{eqn:hymn} is expressed as
$$
\biggl(\frac{\cos\theta+\cos\phi}{2}\biggr)^{\mu+\nu+1}
\frac{\Gamma(\mu-\nu+1)\Gamma(2\nu+1)}
     {\Gamma(\mu+\nu+1)\Gamma(\nu+\frac{1}{2})} 
\tilC_{\mu-\nu}^{\nu+\frac{1}{2}} (\cos\phi).
$$
By using Gauss's duplication formula \eqref{eqn:dup},
we get Claim. 
\end{proof}

\section{Fractional integral of two variables}
\label{subsec:frac}

In Section \ref{subsec:spec},
we find explicit eigenvalues of intertwining operators on
$L^2(S^{p-2} \times S^{q-2})$
based on the 
\index{B}{Funk--Hecke formula}%
Funk--Hecke formula 
(see Example \ref{ex:Ri}).
The following lemma is the substantial part of computation there.
\begin{lemSec}
\label{lem:frac} For\/ 
$\operatorname{Re}\mu, \operatorname{Re}\nu>-\frac{1}{2}$, \thinspace
$\operatorname{Re} \lambda >-1$ and 
$l, k \in \mathbb N$, 
we have the following formula for the fractional integral:
\begin{multline} 
\label{int:frac}
   \int_{-1}^{1}\int_{-1}^1 \frac{(x+y)_\pm^\lambda}{\Gamma(\lambda+1)}   
   \tilC_l^\mu(x) \tilC_k^\nu(y)
   (1-x^2)^{\mu-\frac{1}{2}} (1-y^2)^{\nu-\frac{1}{2}} dx dy   \\    
  = \frac{b 2^{1-\lambda} \Gamma(\lambda+\mu+\nu+1) 
               }
               {\Gamma(\frac{\lambda+2\mu+2\nu+l+k+2}{2})
                \Gamma(\frac{\lambda+2\mu+l-k+2}{2})
                \Gamma(\frac{\lambda+2\nu-l+k+2}{2})
                \Gamma(\frac{\lambda-l-k+2}{2})
                },
\end{multline}   
where 
$$
b:= \frac{(\pm 1)^{l+k}\pi^2}{2^{2\mu+2\nu}}
   \frac{\Gamma(2\mu+l)\Gamma(2\nu+k)}{l! \, k!}
$$
is a constant independent of $\lambda$. 
\end{lemSec}

\begin{proof}
The left-hand side of \eqref{int:frac} amounts to
\begin{align*}
 & \int_{-1}^1\Bigl( \int_{-x}^1
      \frac{(x+y)^\lambda}{\Gamma(\lambda+1)} 
      \tilC_k^\nu(y) (1-y^2)^{\nu-\frac{1}{2}}dy \Bigr)
      \tilC_l^\mu(x)(1-x^2)^{\mu-\frac{1}{2}}dx \\
={}&\frac{\sqrt{\pi}\Gamma(2\nu+k)}
     {2^{\nu-\frac{1}2} k!}   
     \int_{-1}^1\Bigl( (1-x^2)^{\frac{\lambda}{2}+\frac{\nu}{2}+\frac{1}{4}}
     P_{\nu+k-\frac{1}{2}}^{-\lambda-\nu-\frac{1}{2}}(-x)\Bigr)
     \tilC_l^\mu(x) (1-x^2)^{\mu-\frac{1}{2}} dx  \\
={}& \frac{2^{-\nu-2\mu-l+\frac{3}2}\pi}{\Gamma(\mu+l+\frac{1}{2})}
      \frac{\Gamma(2\nu+k)\Gamma(2\mu+l)}{k! \, l!}
     \int_{-1}^1 
     (1-x^2)^{\frac{\lambda}{2}+\mu+\frac{\nu}{2}-\frac{l}{2}-\frac{1}{4}}
     P_{\nu+k-\frac{1}{2}}^{-\lambda-\nu+l-\frac{1}{2}}(-x) dx  \\
={}& \frac{2^{1-\lambda-2\mu-2\nu}\pi \Gamma(2\mu+l)\Gamma(2\nu+k)
\Gamma(\lambda+\mu+\nu+1)}
 {l! \, k! \, \Gamma(\frac{\lambda+2\mu+2\nu+l+k+2}{2})
 \Gamma(\frac{\lambda+2\mu+l-k+2}{2})
 \Gamma(\frac{\lambda+2\nu-l+k+2}{2})
 \Gamma(\frac{\lambda-l-k+2}{2})}.
\end{align*}
Hence, the right-hand side of \eqref{int:frac} follows.
Some remarks on each equality are given in turn:

{\bf First equality} follows from \textbf{Ge\ref{item:Ge4}} in Appendix
\ref{subsec:Ge}. 

{\bf Second equality}. First, we made use of 
the integral by parts because we have (see \eqref{def:Ge})
$$
  \tilC_l^\mu(x) (1-x^2)^{\mu-\frac{1}{2}}
=\frac{(-1)^l}{2^{2\mu+l-1} ~l!}
  \frac{\Gamma(2\mu+l)\sqrt{\pi}}{\Gamma(\mu+l+\frac{1}{2})}
 \frac{d^l}{dx^l} \Bigl( (1-x^2)^{\mu+l-\frac{1}{2}} \Bigr).
$$
Then, we applied the functional relation 
$$
  \frac{d^l}{dx^l}
     \Bigl( (1-x^2)^{\frac{\lambda}{2}+\frac{\nu}{2}+\frac{1}{4}} 
      P_{\nu+k-\frac{1}{2}}^{-\lambda-\nu-\frac{1}{2}}(-x) \Bigr)   
= (1-x^2)^{\frac{\lambda}{2}+\frac{\nu}{2}-\frac{l}{2}+\frac{1}{4}}
   P_{\nu+k-\frac{1}{2}}^{-\lambda-\nu+l-\frac{1}{2}}(-x),
$$
which is obtained by iterating \eqref{eqn:L1}.  

{\bf Third equality}.
We applied the integral formula \eqref{eqn:LGam}
after changing the variable $x \mapsto -x$.

Therefore,
the proof of Lemma \ref{lem:frac} is completed.
\end{proof}

\def\cprime{'}
\providecommand{\href}[2]{#2}

\bigskip
\noindent
Toshiyuki Kobayashi\\
\textit{Home address}: 
Graduate School of Mathematical Sciences,
  the University of Tokyo, 
  3-8-1 Komaba, Meguro, Tokyo, 153-8914, Japan \\
{E-mail address}: {toshi@ms.u-tokyo.ac.jp}\\
\textit{Current address}: 
Department of Mathematics, FAS, 
Harvard University, 
One Oxford Street, Cambridge MA 02138 USA\\
\\
Gen Mano\\
Graduate School of Mathematical Sciences,
  the University of Tokyo, 
  3-8-1 Komaba, Meguro, Tokyo, 153-8914, Japan \\
{E-mail address}: {gmano@ms.u-tokyo.ac.jp}

\printindex{A}{List of Symbols}
\printindex{B}{Index}

\end{document}